\documentclass[12pt,twoside,english]{article}
\usepackage[T1]{fontenc}
\usepackage{units}
\usepackage{amssymb}
\usepackage{graphicx}

\makeatletter
\newcommand{\lyxaddress}[1]{
\par {\raggedright #1
\vspace{1.4em}
\noindent\par}
}

\usepackage{a4wide}
\usepackage{graphicx}
\usepackage{fancyhdr}
\usepackage{amsmath}
\usepackage{amssymb}
\usepackage{subfigure}
\usepackage{setspace}
\usepackage{verbatim} 
\numberwithin{equation}{section}

\newcommand{\vek}[1]{\mathchoice{\displaystyle\boldsymbol#1}
{\textstyle\boldsymbol#1}{\scriptstyle\boldsymbol#1}
{\scriptscriptstyle\boldsymbol#1}}

\makeatother

\usepackage{babel}
\begin{document}

\title{Higher-order meshing of implicit geometries---part I:\\Integration
and interpolation in cut elements}

\author{T.P. Fries, S. Omerovi{\'c}, D. Sch{\"o}llhammer, J. Steidl}

\maketitle

\lyxaddress{\begin{center}
Institute of Structural Analysis\\
Graz University of Technology\\
Lessingstr. 25/II, 8010 Graz, Austria\\
\texttt{www.ifb.tugraz.at}\\
\texttt{fries@tugraz.at}
\end{center}}
\begin{abstract}
An accurate implicit description of geometries is enabled by the level-set
method. Level-set data is given at the nodes of a higher-order background
mesh and the interpolated zero-level sets imply boundaries of the
domain or interfaces within. The higher-order accurate integration
of elements cut by the zero-level sets is described. The proposed
strategy relies on an automatic meshing of the cut elements. Firstly,
the zero-level sets are identified and meshed by higher-order interface
elements. Secondly, the cut elements are decomposed into conforming
sub-elements on the two sides of the zero-level sets. Any quadrature
rule may then be employed within the sub-elements. The approach is
described in two and three dimensions without any requirements on
the background meshes. Special attention is given to the consideration
of corners and edges of the implicit geometries.

\end{abstract}
\newpage{}\tableofcontents{}\newpage{}

\section{Introduction\label{sec:Introduction}}

\begin{figure}
\centering

\subfigure[Conforming mesh]{\includegraphics[width=4cm]{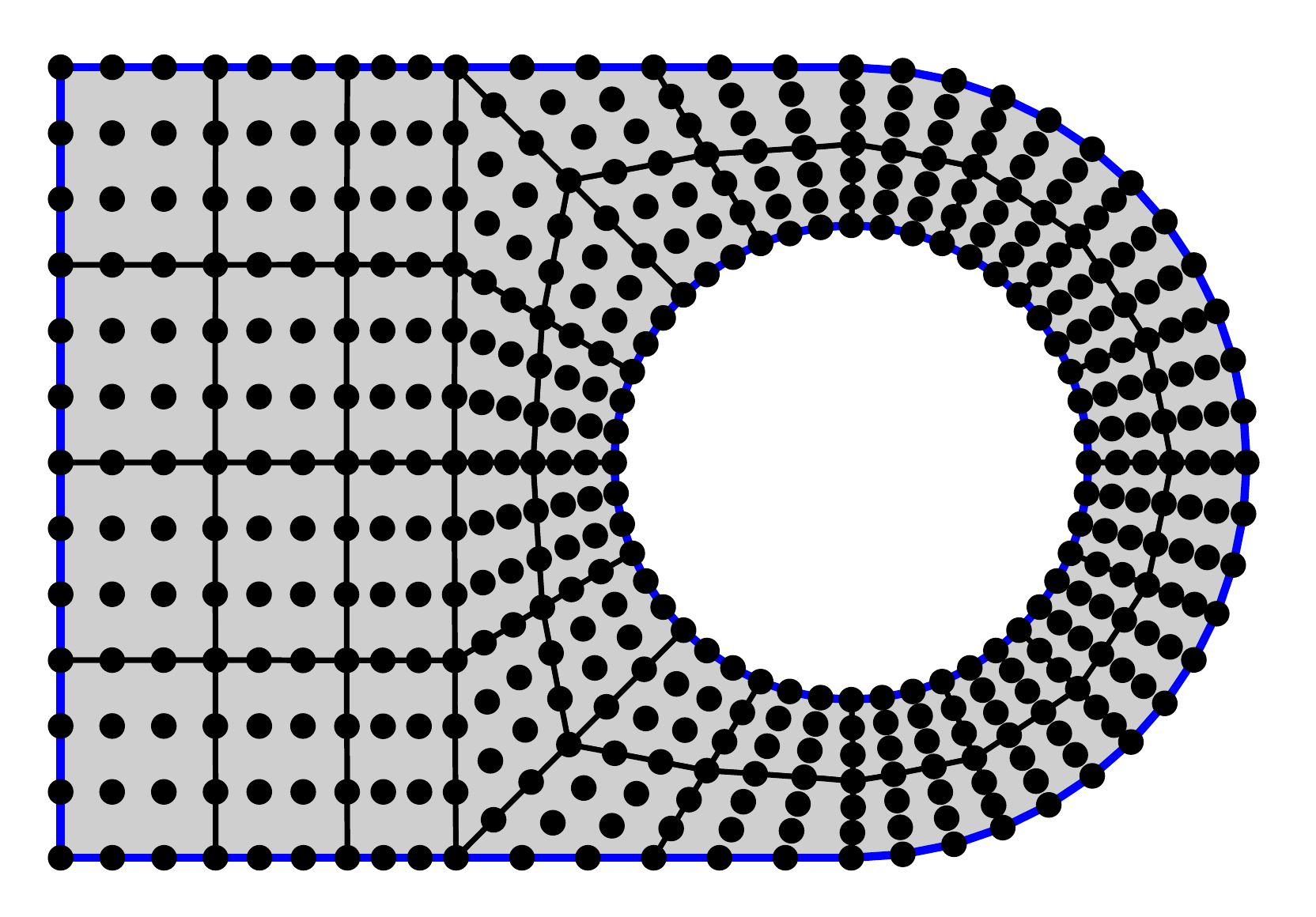}}\;\subfigure[Internal boundary]{\includegraphics[width=4cm]{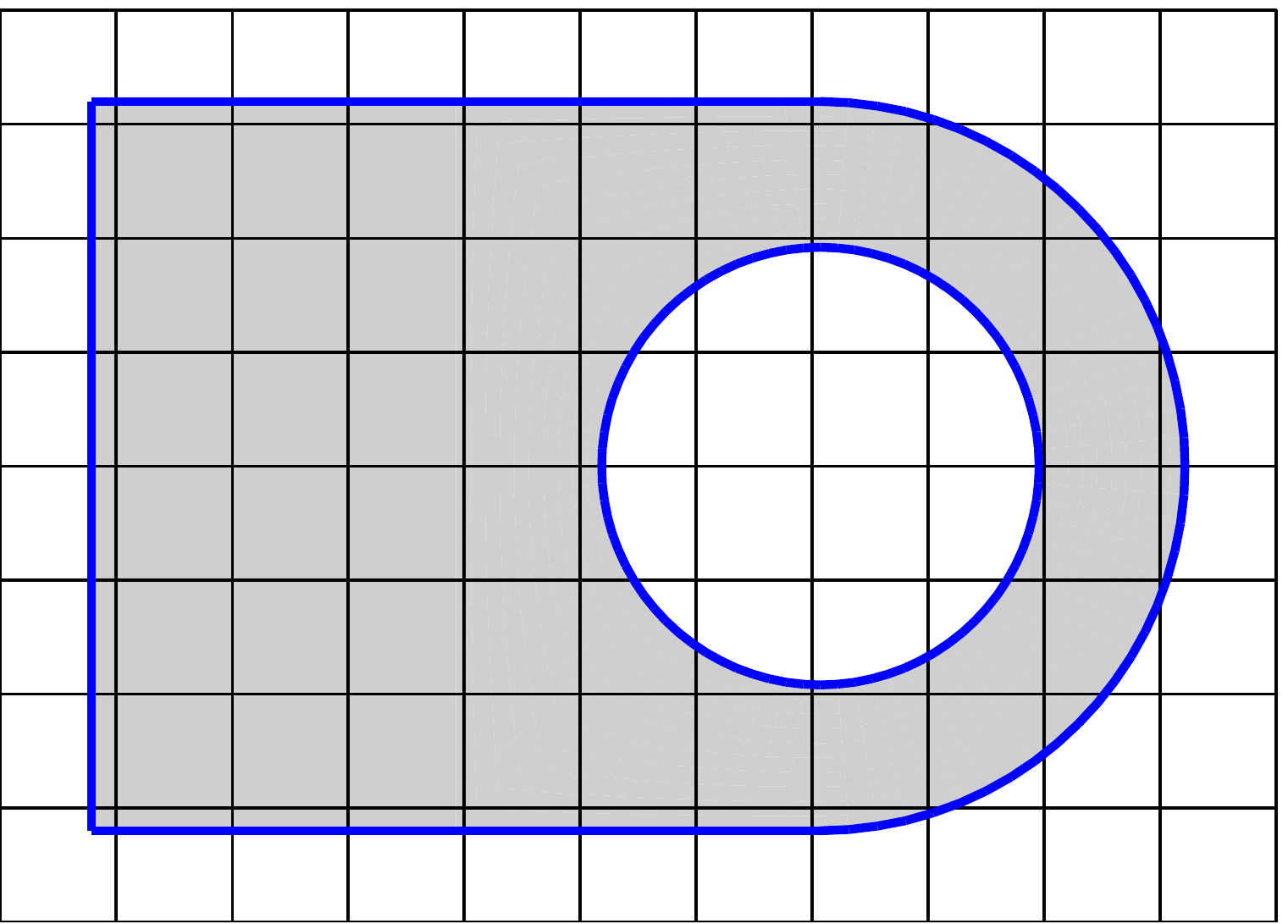}}\;\subfigure[Cracks]{\includegraphics[width=4cm]{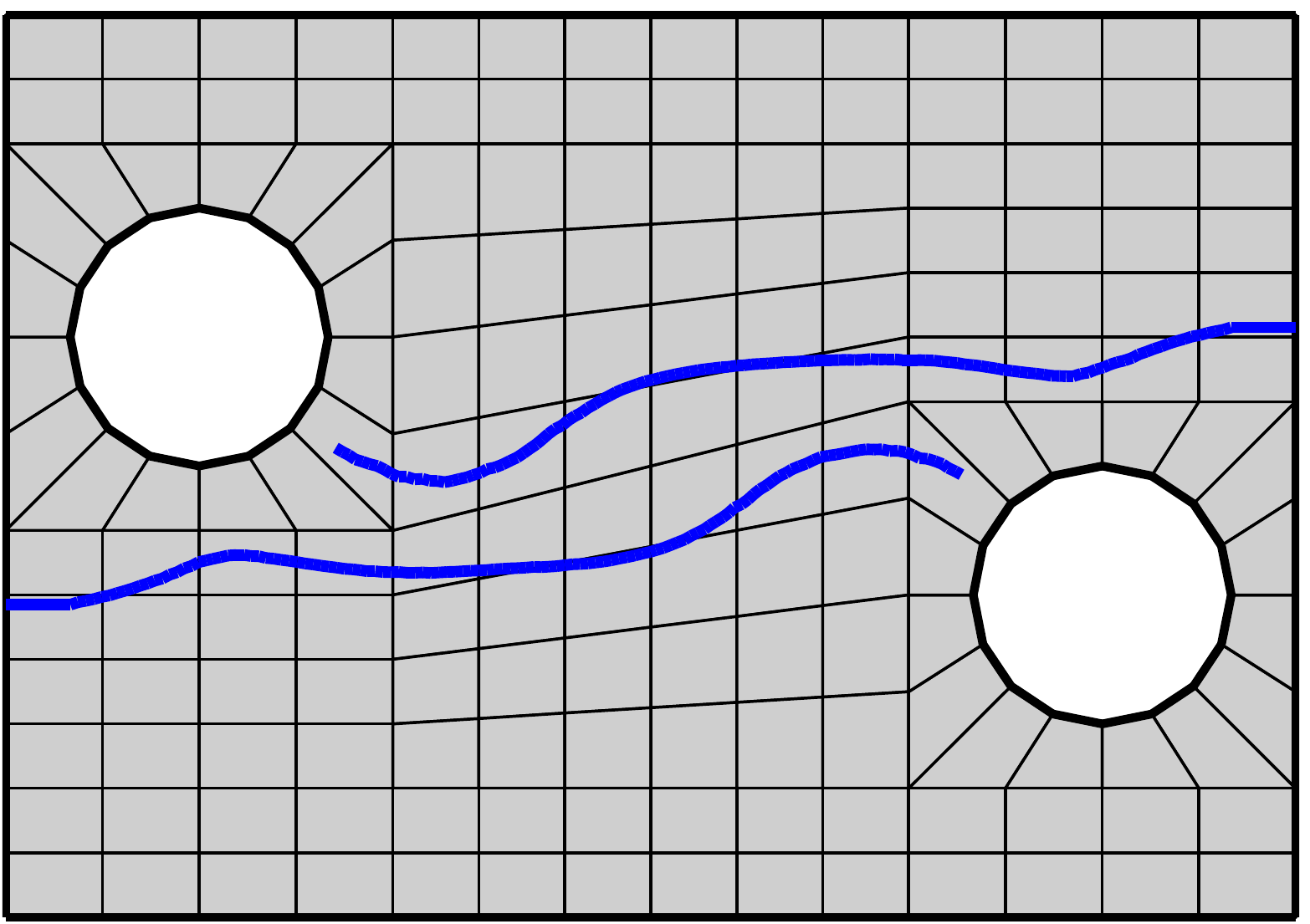}}\;\subfigure[Two-phase flows]{\includegraphics[width=4cm]{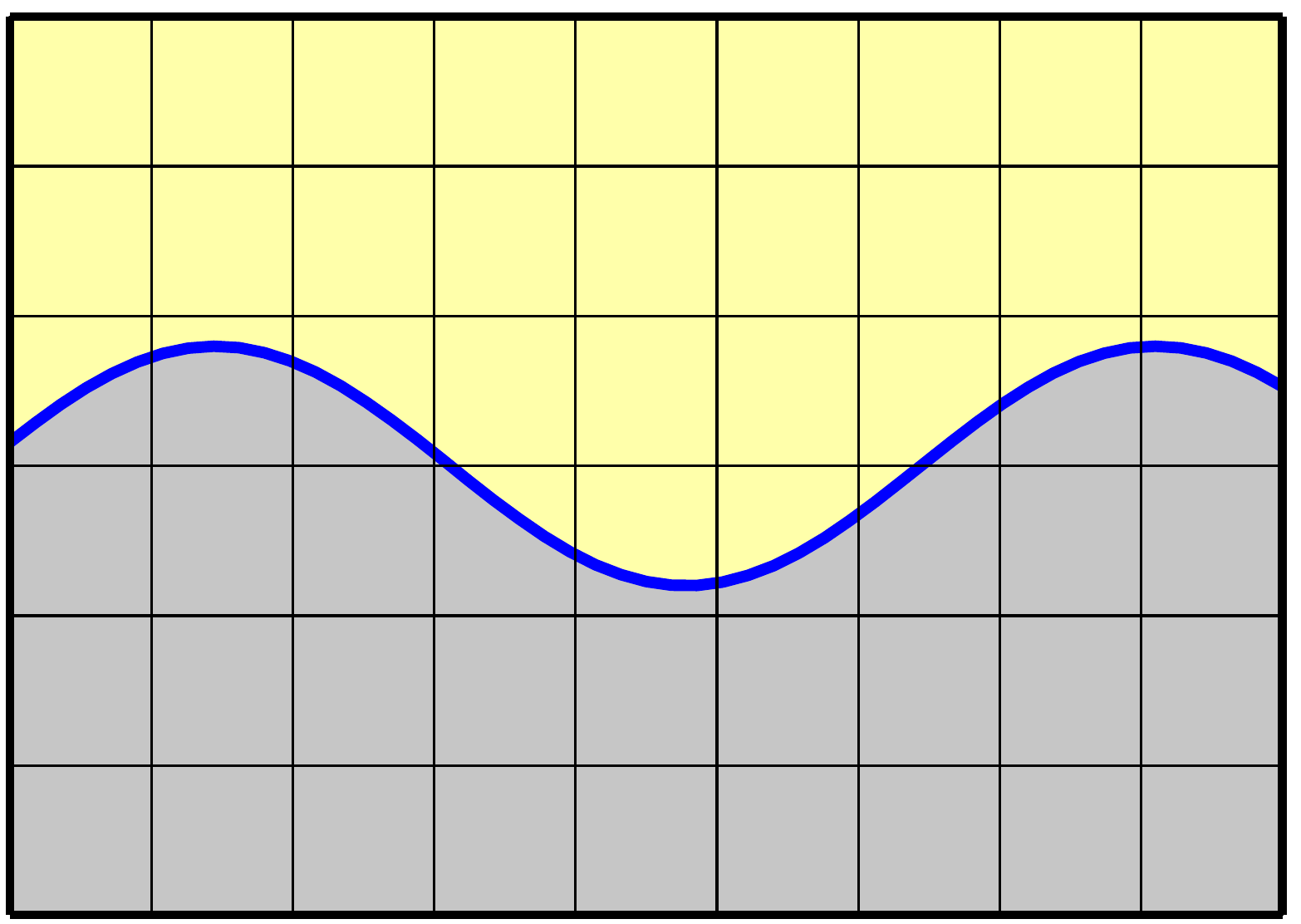}}

\caption{\label{fig:VisIntro}(a) Conforming higher-order mesh used in classical
$p$-FEM, (b) internal boundaries in FDMs, (c) and (d) internal interfaces,
e.g.~in the XFEM.}
\end{figure}

The approximation of boundary value problems (BVPs) based on the finite
element method (FEM) has achieved an enormous importance in engineering,
physics, and related fields. The goal to achieve higher-order accurate
approximations of BVPs can be traced back to the early days of the
FEM, see e.g.~\cite{Szabo_2004a} or the text books \cite{Solin_2003a,Bathe_1996a,Belytschko_2000b,Zienkiewicz_2000d}.
This is typically labeled $p$-FEM and there are numerous references
found. In fact, the isoparametric concept leads to a conceptually
simple approach to achieve optimal results. Therefore, an accurate
geometry representation based on higher-order elements which consider
boundaries and interfaces is necessary, see Fig.~\ref{fig:VisIntro}(a)
for an example. The generation of such meshes is, however, not a simple
task. Moreover, elements may have to be refined during the analysis,
for instance in the context of adaptivity and convergence studies.
The interplay of the FEM software, the meshing tool, and the CAD program
is far from trivial and hardly automated, in particular using higher-order
elements. Things become even worse in case of moving boundaries and
interfaces. 

Hence, it is not surprising that a number of methods have been developed
which do not require conforming meshes. In the case where \emph{boundaries}
are not meshed explicitly, we use the term ``fictitious domain methods''
(FDMs) and there is a large variety of variants. These methods are
characterized by the fact that the physical domain in which a BVP
is formulated, is completely embedded in a fictitious domain, also
called background mesh, see Fig.~\ref{fig:VisIntro}(b). The background
mesh is often structured and no geometry dependent meshing is needed.
The shape functions are provided by the background mesh. Methods falling
into this class include the unfitted or cut finite element method
\cite{Burman_2010a,Burman_2012a,Burman_2014a,Hansbo_2002a},\textbf{
}finite cell method \cite{Abedian_2013b,Duester_2008a,Parvizian_2007a,Schillinger_2014a,Schillinger_2014b},
Cartesian grid method \cite{Uzgoren_2009a,Ye_1999a},\textbf{ }immersed
interface method \cite{Leveque_1994a}, virtual boundary method \cite{Saiki_1996a},
embedded domain method \cite{Lohner_2007a,Neittaanmaki_1995a} etc.

Another class of methods has been developed for problems with \emph{interfaces}
rather than boundaries. Often, simple meshes are chosen to describe
the overall geometry without taking care that internal interfaces
such as material interfaces or crack surfaces are meshed conformingly,
see Fig.~\ref{fig:VisIntro}(c) and (d). Non-standard approximation
spaces are often constructed near the interfaces. In this class, we
mention extended or generalized finite element methods (XFEM/GFEM),
see e.g.~\cite{Belytschko_1999a,Moes_1999a,Fries_2009b} for the
XFEM and \cite{Strouboulis_2000a,Strouboulis_2000b} for the GFEM.
They consider for inner-element jumps and kinks by adding enrichment
functions based on the partition of unity concept \cite{Babuska_1994a,Babuska_1997a,Melenk_1996a}.
The XFEM has also been used in the realm of FDMs e.g.~in \cite{Sukumar_2001c}.
Other mesh independent approximations for non-smooth solutions involve
the weak element method \cite{Goldstein_1986a,Rose_1975a}, the ultra
weak variational method \cite{Cessenat_1998a}, the least-squares
method as proposed in \cite{Monk_1999a}, and the global-local FEM
as described in \cite{Mote_1971a,Noor_1986a}. Finally, we mention
the manifold method \cite{Chen_1998b,Shi_1992a} which may be used
for inner-element discontinuities.

A shared property of FDMs and XFEM-related approaches is the numerical
integration in elements cut by boundaries or interfaces. Approaches
for the integration of elements with internal boundaries and interfaces
may be distinguished based on the fact whether they rely on a decomposition
of the cut elements into sub-elements or not. For the first class,
the standard approach is to recursively decompose a cut element into
polygonal sub-cells until the desired accuracy is obtained \cite{Abedian_2013a,Moumnassi_2011a,Dreau_2010a}.
This typically leads to a very large number of integration points
in the context of higher-order approximations \cite{Stazi_2003a,Zi_2003a,Legrain_2012a,Laborde_2005a,Dreau_2010a}.
It was already noted in \cite{Legay_2005a,Cheng_2009a,Fries_2015a,Sala_2012a}
that a decomposition into sub-elements with curved, higher-order edges
or faces is a strategy which consistently enables the generation of
higher-order accurate integration rules with only a modest number
of integration points. Of course, the effort for generating these
integration points is larger than for the polygonal approaches with
reduced accuracy. The other class of approaches is built by methods
that do \emph{not} decompose the cut elements \cite{Ventura_2006a,Ventura_2009a,Mueller_2013a,Mousavi_2011a}.
Typically, the number of integration points is small, however, the
generation of proper integration weights is often involved and requires
the solution of small systems of equations e.g.~in the context of
moment fitting and Lasserre's technique. The extendability to three
dimensions, higher-order accuracy, general integrands, and to the
presence of corners and edges often poses problems in this class of
approaches.

The higher-order accurate numerical integration in three-dimensional,
cut elements is subject of some very recent contributions and interested
readers are referred to the references given therein: For those approaches
relying on element decompositions, \cite{Fries_2015a} presents a
very general approach in the context of the level-set method which
applies for all element types and has no restriction on the background
meshes. The adaption to the finite cell method and spline based boundary
representations is reported in \cite{Stavrev_2016a}. ``Smart octrees''
are proposed in \cite{Kudela_2016a} and are employed for structured
quadrilateral and hexahedral background meshes. Among the approaches
which avoid decompositions, \cite{Joulaian_2016a} presents an approach
based on moment fitting. Finally, the approach in \cite{Lehrenfeld_2016a}
manipulates the background mesh such that the interfaces are meshed
conformingly. Comparisons of these advanced approaches with the classical
use of (polygonal) recursive spacetrees are reported in \cite{Kudela_2016a,Stavrev_2016a}
and clearly demonstrate the advantages of the recent approaches.

\begin{figure}
\centering

\subfigure[Domain]{\includegraphics[width=4cm]{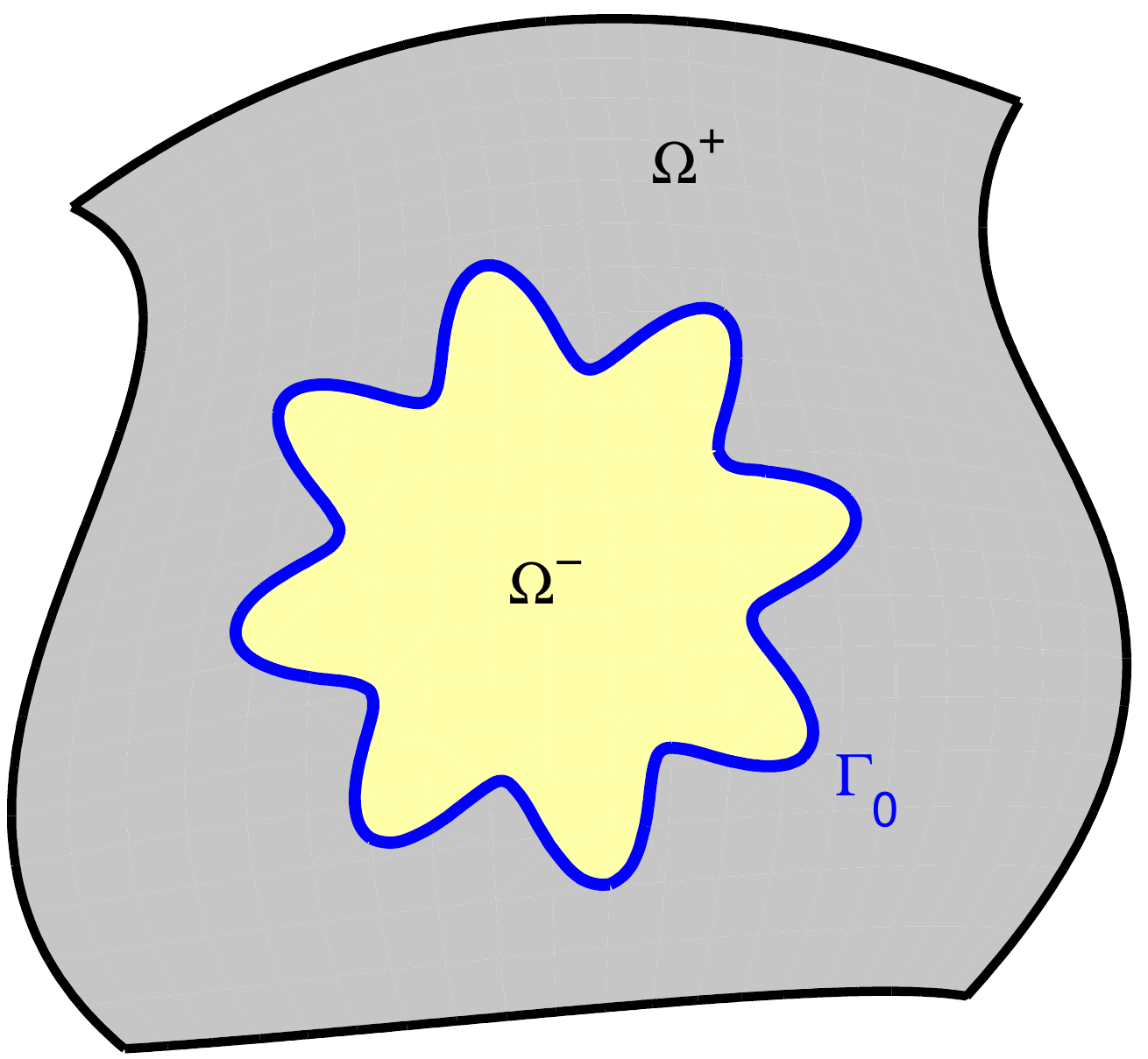}}\qquad\subfigure[Background mesh]{\includegraphics[width=4cm]{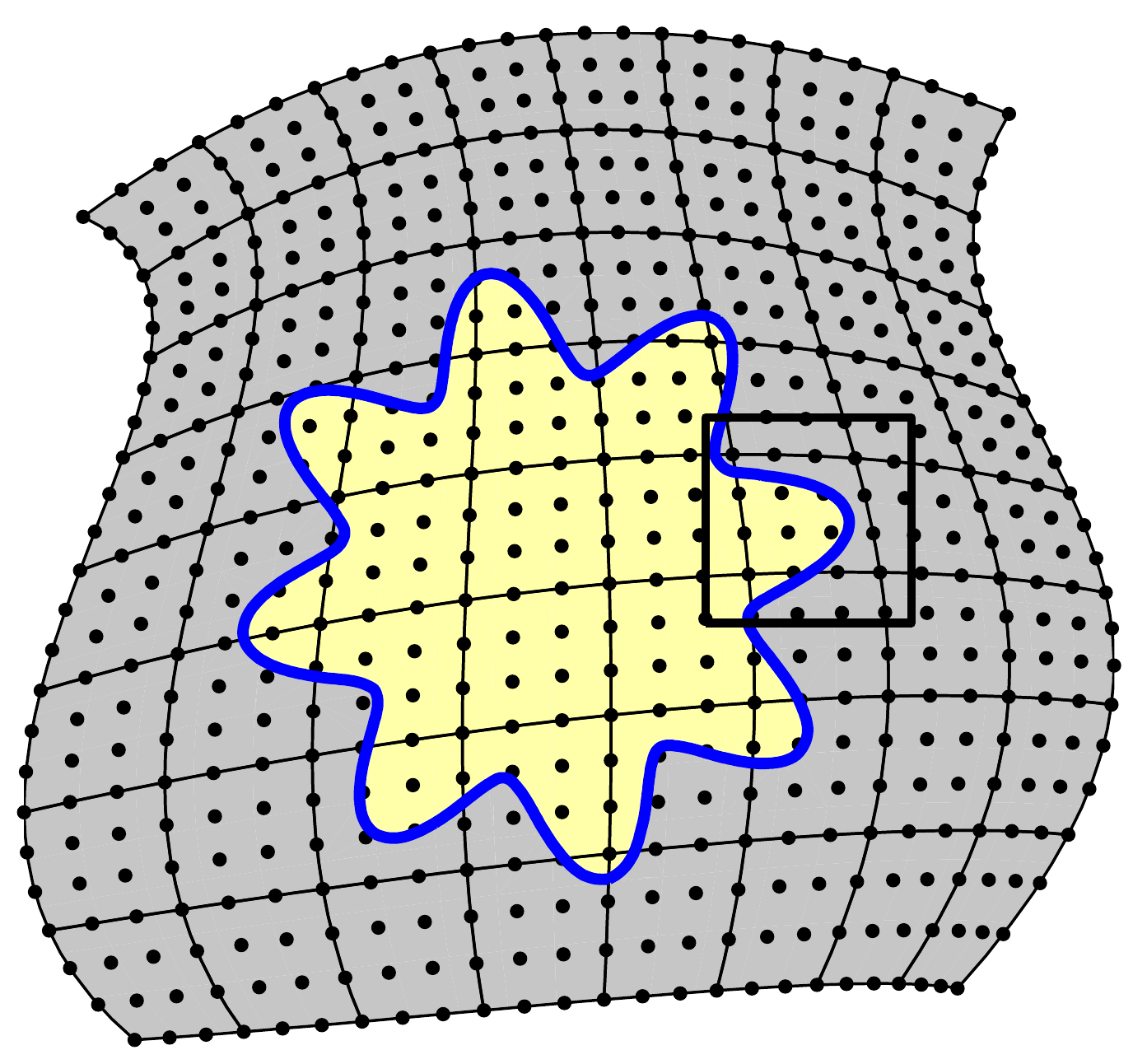}}\qquad\subfigure[Detail]{\includegraphics[width=4cm]{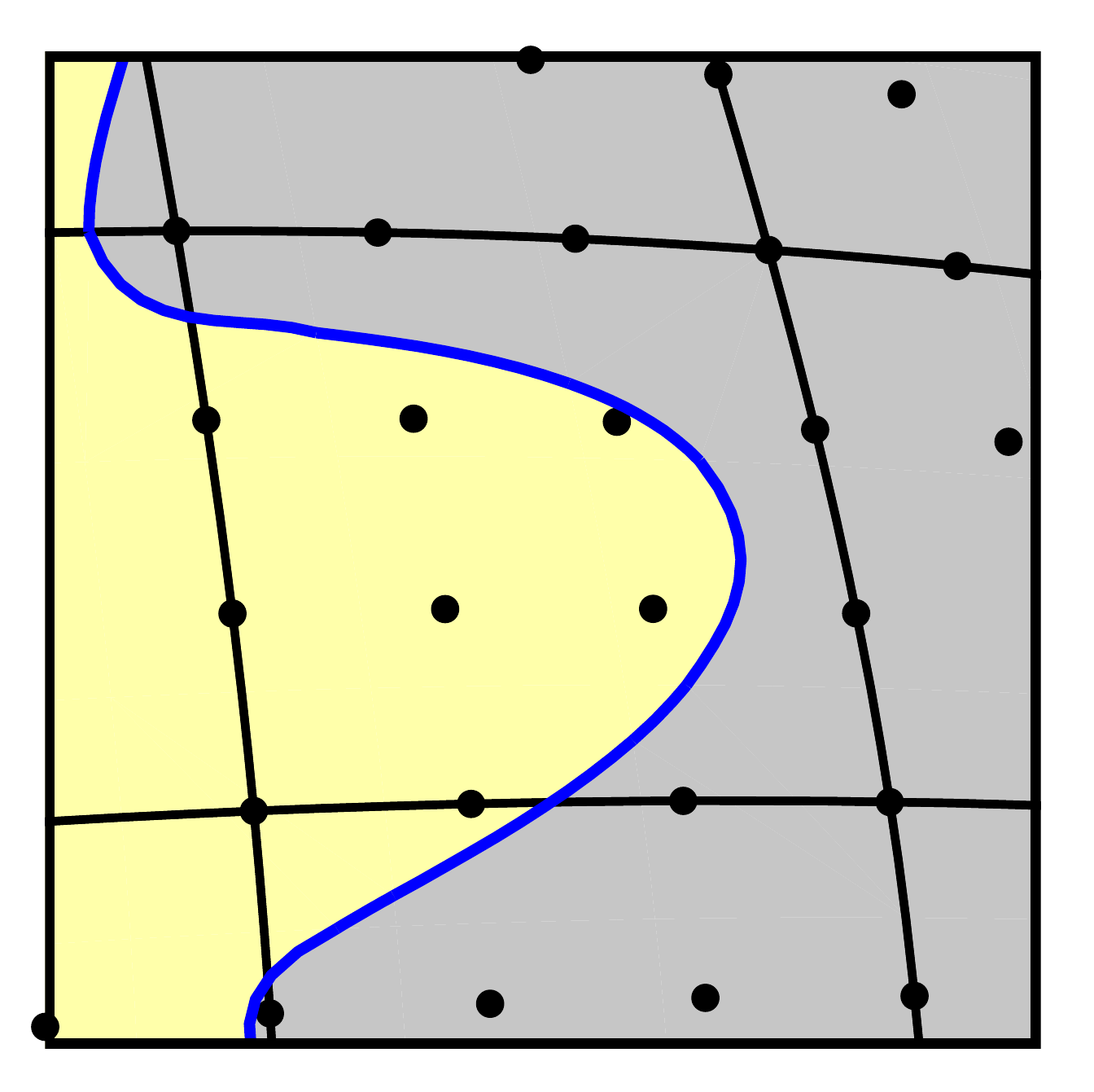}}

\caption{\label{fig:VisStrategyA}(a) Arbitrary domain with an interface described
by a level-set function, (b) higher-order background mesh, (c) detail
showing the higher-order element cut by a curved zero-level set.}
\end{figure}

\begin{figure}
\centering

\subfigure[Decomposed elements]{\includegraphics[width=4cm]{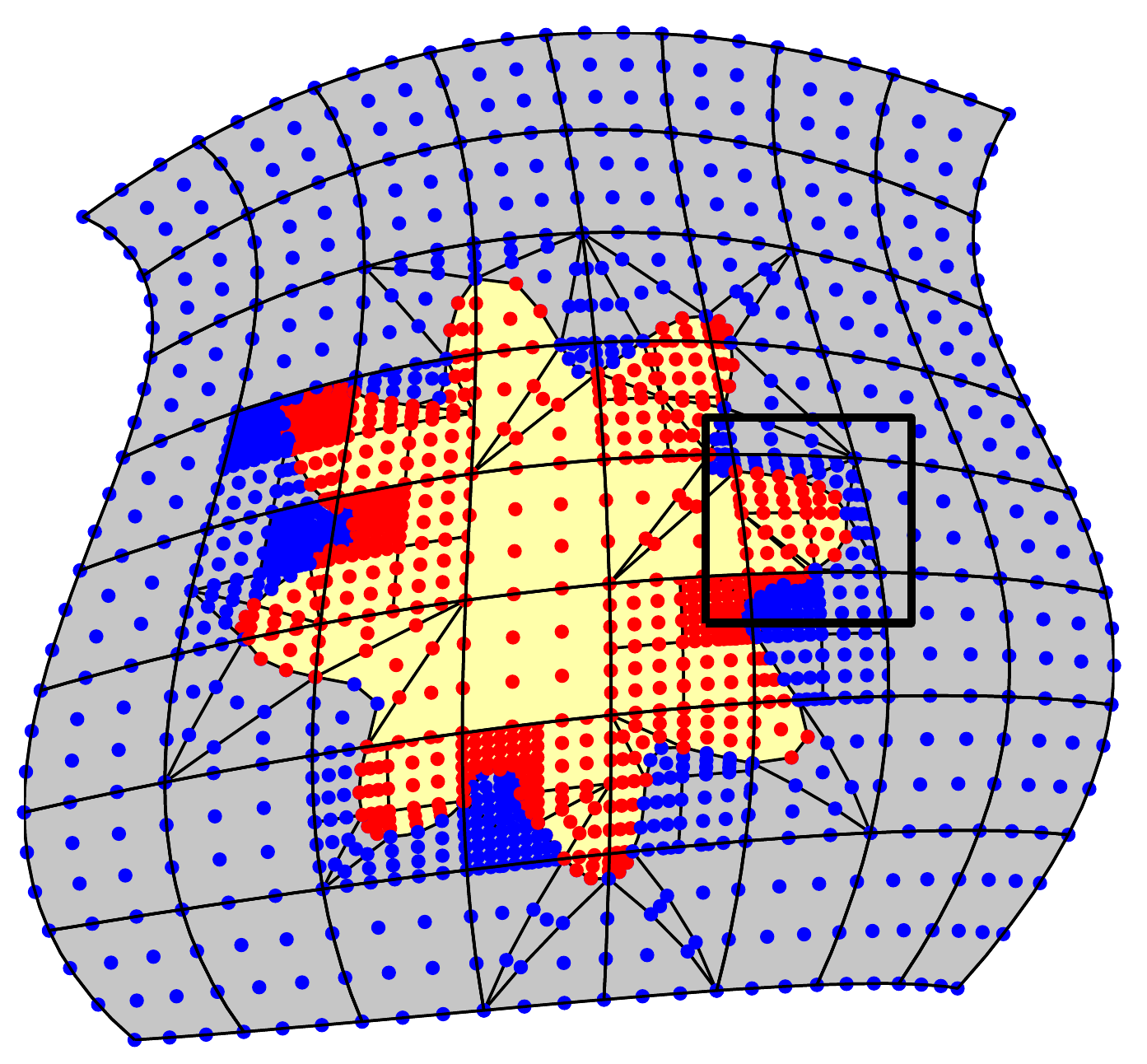}}\qquad\subfigure[Detail: elements]{\includegraphics[width=4cm]{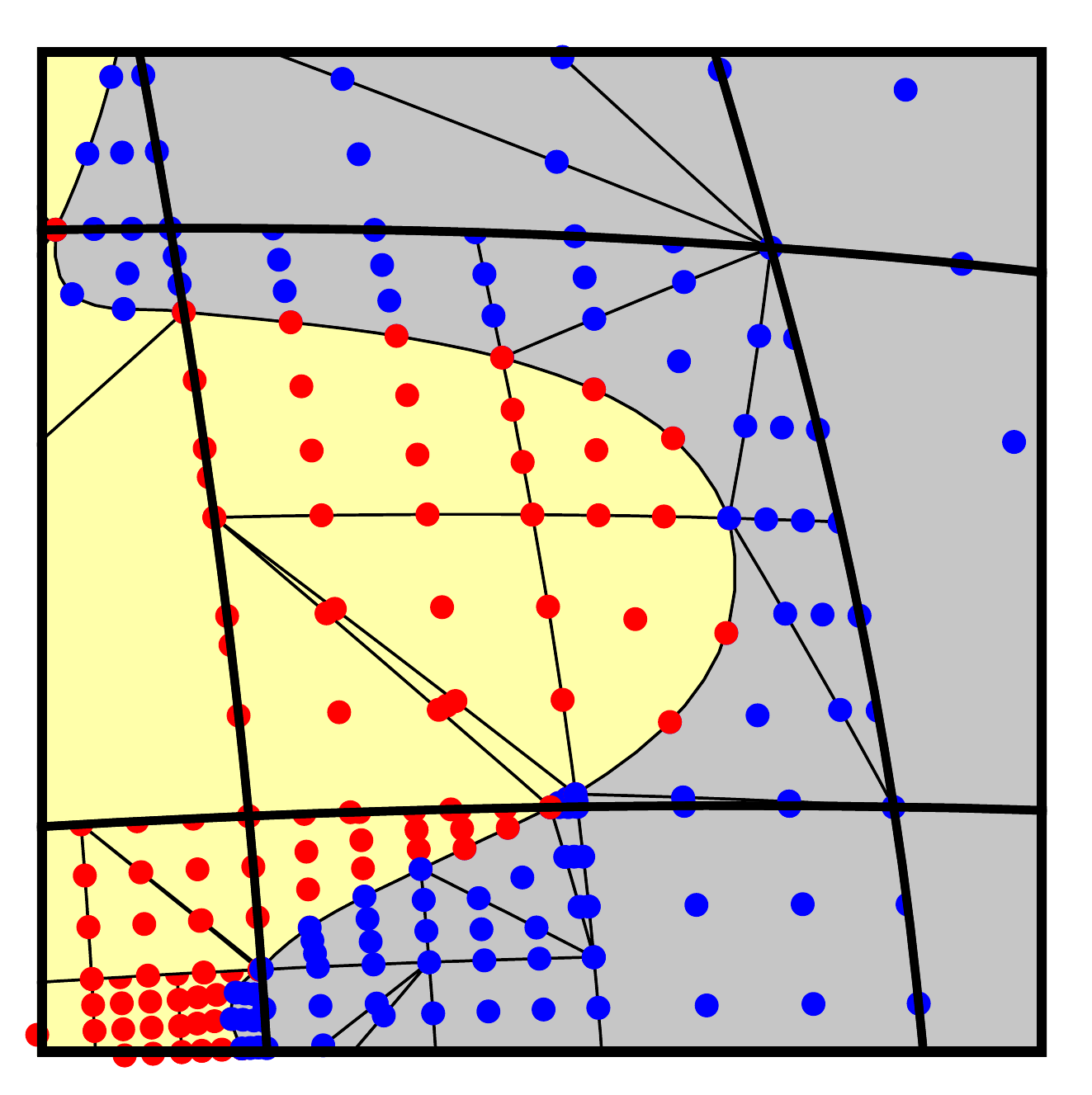}}\qquad\subfigure[Detail: int.~points]{\includegraphics[width=4cm]{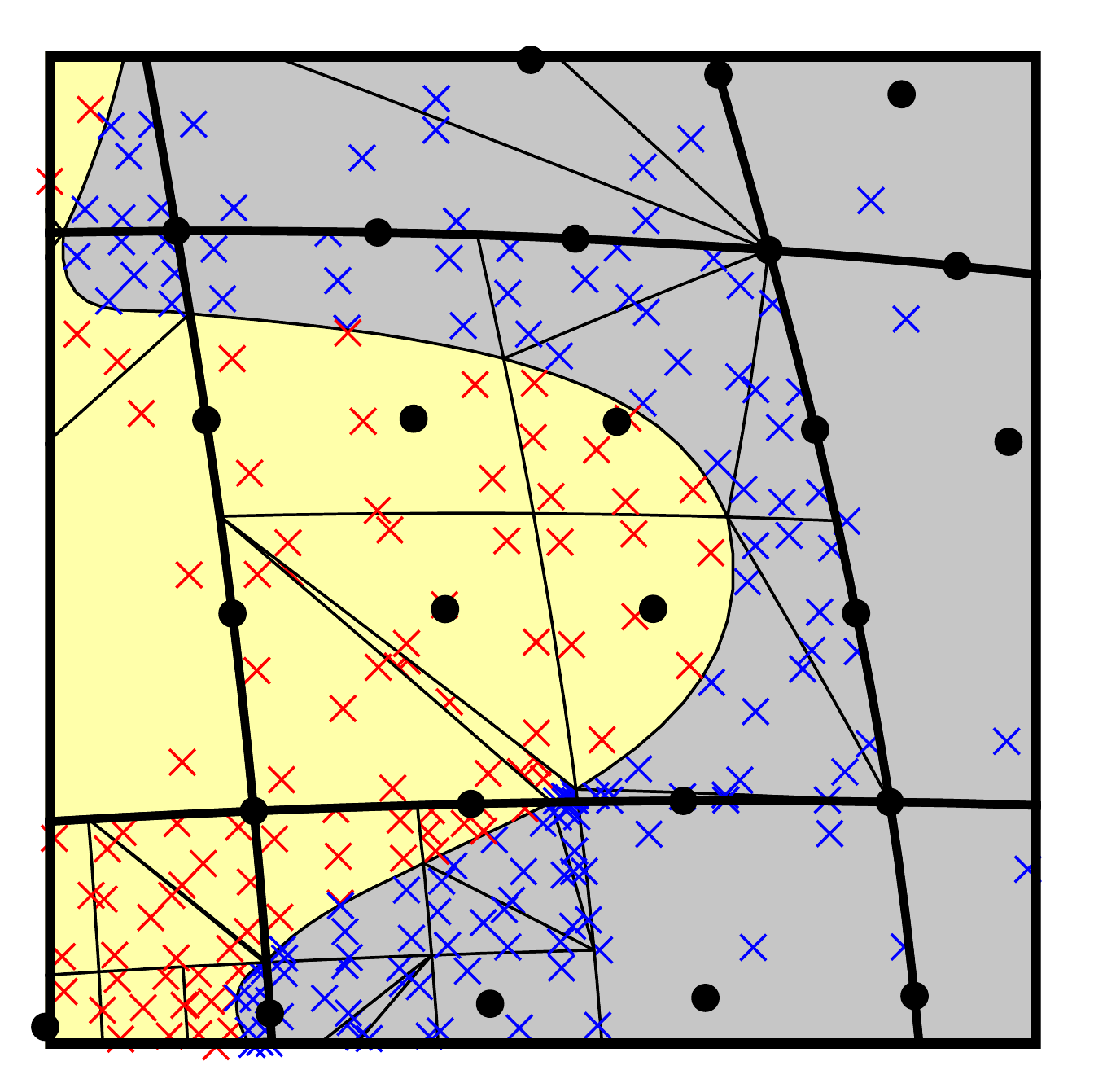}}

\caption{\label{fig:VisStrategyB}(a) Decomposed elements, (b) detail showing
sub-elements in a cut background element (the red and blue dots are
the element nodes), (c) integration points (crosses) mapped to the
sub-elements (the black dots are the element nodes of the background
elements).}
\end{figure}

Herein, we propose a new strategy for the numerical integration of
elements cut by boundaries and interfaces which relies on the decomposition
of cut elements into isoparametric, higher-order sub-elements. Let
us assume the abstract domain represented by Fig.~\ref{fig:VisStrategyA}(a)
featuring an internal interface $\Gamma_{0}$ implicitly defined by
the level-set method. It does not matter here whether $\Gamma_{0}$
is an \emph{interface} between two materials in $\Omega^{+}$ and
$\Omega^{-}$ or a \emph{boundary} if e.g.~$\Omega^{-}$ is considered
to be a void region. A higher-order background mesh is introduced
as shown in Fig.~\ref{fig:VisStrategyA}(b) and the level-set data
is only given at the nodes. The implied zero-level set is, in general,
curved within the elements, see Fig.~\ref{fig:VisStrategyA}(c).
The task is to integrate on one or both sides of the interface individually
with higher-order accuracy. Therefore, the zero-level sets are approximated
by higher-order interface elements (reconstruction) and higher-order
sub-elements result on the two sides of the interface (decomposition).
See Fig.~\ref{fig:VisStrategyB}(a) for the resulting subelements
in the overall mesh and Fig.~\ref{fig:VisStrategyB}(b) for a detail.
Standard integration points, e.g.~Gauss points may now be mapped
into the sub-elements and used in FDMs or XFEM-related methods for
the integration in the corresponding background element, see Fig.~\ref{fig:VisStrategyB}(c).
It is also shown how corners and edges may be considered properly
by the proposed method using several level-set functions and performing
the decomposition successively for each. Note that previous works
in \cite{Fries_2015a,Fries_2016a} by (partly) the same authors do
not offer this simple concept for corners and edges. Furthermore,
we find that the approach herein unifies the treatment of the reconstruction
and decomposition step as they now both rely on element generation,
i.e.~meshing.

This is the first paper in a sequence of contributions where it is
shown that the resulting, automatically generated meshes, herein only
used for the proper placement of integration points in cut elements
and interpolation, may also be used for directly approximating BVPs.
In the second part, this will be shown for BVPs on zero-level sets,
i.e.~on manifolds. Therefore, only the mesh of reconstructed interface
elements will be needed. In the third part, BVPs with inner-element
boundaries and interfaces are considered using meshes consisting of
regular and cut background elements, the latter decomposed into conforming
sub-elements as described herein. Further studies will then outline
how fully automatic, higher-order analyses are enabled for domains
in the context of constructive solid geometries.

The paper is organized as follows: Section \ref{sec:BackgroundMeshAndLevelSetData}
describes some preliminaries and defines the higher-order background
meshes used herein and the level-set method. In Section \ref{sec:Reconstruction},
the identification and meshing of the zero-level set is described
in detail (reconstruction). The decomposition into higher-order, conforming
sub-elements is detailed in Section \ref{sec:Decomposition}. The
proper consideration of corners and edges with several level-set functions
is outlined in Section \ref{sec:CornersAndEdges}. Numerical results
are presented at the end of each section. Finally, the paper concludes
in Section \ref{sec:Conclusions} with a summary and outlook. An appendix
is provided where mappings are detailed which are frequently used
throughout this work.

\section{Background mesh and level-set data\label{sec:BackgroundMeshAndLevelSetData}}

A two or three-dimensional domain of interest, $\Omega\in\mathbb{R}^{d}$,
$d=\{2,3\}$ is defined implicitly based on the level-set method \cite{Osher_2003a,Sethian_1999b}.
That is, a scalar function $\phi\left(\vek x\right)$, $\vek x\in\Omega$
is given, called level-set function, which is positive on one side
of the boundary/interface and negative on the other. The function
is continuous and its zero-level set defines the position of a smooth
boundary or interface. Note that $\phi\left(\vek x\right)$ does not
have to be a signed distance function herein. In Section \ref{sec:CornersAndEdges},
we shall also consider the situation where boundaries and interfaces
are, in fact, \emph{not} smooth but feature edges and corners which
is of high practical relevance. Then, we suggest to use several level-set
functions for the description.

The domain is fully immersed in a background mesh composed by higher-order
Lagrange elements. The background mesh may be unstructured and feature
curved elements. The focus herein is on triangular elements in 2D
and tetrahedral elements in 3D. The level-set data is only given at
the nodes of the background mesh, $\phi_{i}=\phi\left(\vek x_{i}\right)$,
and is interpolated by standard finite element shape functions $N_{i}\left(\vek x\right)$
in-between, hence, $\phi^{h}\left(\vek x\right)=\sum_{i}N_{i}\left(\vek x\right)\cdot\phi_{i}$.
That is, once the nodal level-set data is given, no (further) use
of any analytic knowledge of the level-set functions is made. Obviously,
the background mesh enables a higher-order accurate description of
the domain as long as the boundaries and interfaces are sufficiently
smooth.

\begin{figure}
\centering

\subfigure[Situation in 2D]{\includegraphics[width=15cm]{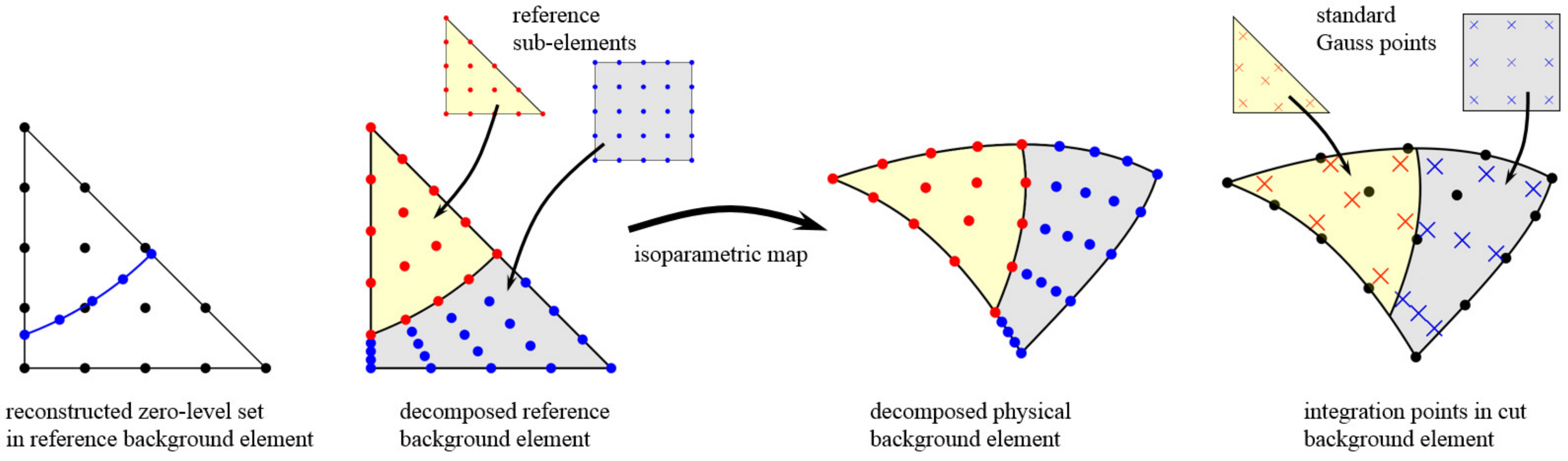}}\\\subfigure[Situation in 3D]{\includegraphics[width=15cm]{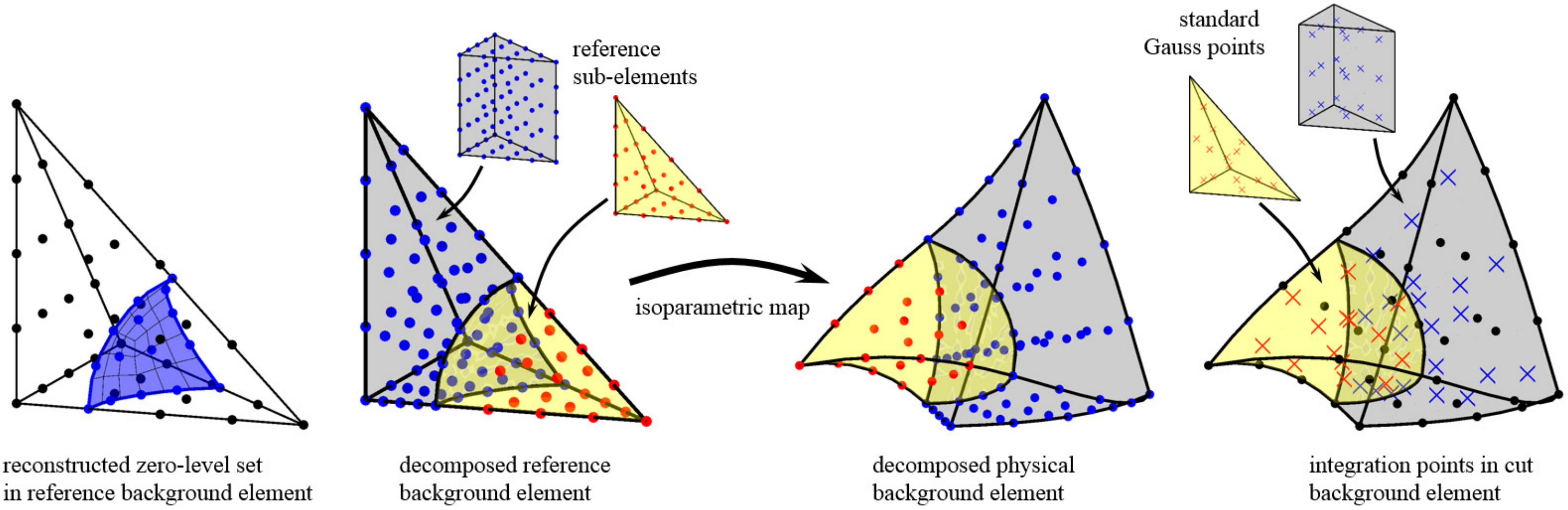}}

\caption{\label{fig:VisOverviewRemesh}The proposed procedure to obtain integration
points in cut elements in (a) 2D and (b) 3D.}
\end{figure}

The task is to obtain integration points which allow a higher-order
accurate integration in the implicit geometry. Therefore, the background
elements are decomposed into sub-elements which conform with the interfaces
and automatically consider also for corners and edges. It is then
simple to use any desired integration rule (e.g.~Gauss quadrature)
within each sub-element. The decomposition is successively applied
in each element for the involved level-set functions. This is a two-step
procedure for each zero-level set that cuts an element: (1) Reconstruction,
which is the meshing of the zero-level set with higher-order interface
elements. (2) Decomposition of cut elements into higher-order sub-elements
on the two sides of the interface elements. These two steps are depicted
in the left half of Fig.~\ref{fig:VisOverviewRemesh}(a) and (b)
for the two- and three-dimensional case, respectively. It is important
to note that both steps are performed first in each cut \emph{reference}
element. This is, in fact, crucial for the ability to handle arbitrarily
curved background elements and perform the decomposition successively
for several level-set functions. Thereafter, (3) the sub-elements
are mapped to the corresponding physical background element and (4)
any desired quadrature points are mapped to these sub-elements, see
the right half of Fig.~\ref{fig:VisOverviewRemesh}. 

Herein, the focus is on triangular and tetrahedral elements, but the
extension to quadrilateral and hexahedral elements is truly straightforward
and described in \cite{Fries_2015a} in a similar context. For example,
one may simply decompose a reference hexahedral element into tetrahedra
which is not a problem because the faces are flat in the reference
domain. Then, the decomposition is realized in each reference tetrahedron
as described below.

\section{Reconstruction\label{sec:Reconstruction}}

Let us focus on a reference element which is cut by the zero-level
set of a level-set function. The nodal values of the level-set function
at the element nodes are taken from the physical element, i.e.~$\phi(\vek r_{i})=\phi(\vek x_{i})$.
The task is to approximate the zero-level set of $\phi^{h}(\vek r)$
by means of a higher-order interface element. That is, roots have
to be identified on the zero-level set which are then the element
nodes of the interface element. In-between, the zero-level set is
only approximated. It is important to note that the location of the
inner element nodes is not unique as shall be seen below.

The description here is partly along the lines of \cite{Fries_2015a,Fries_2016a},
however, the following issues in the context of reconstructions are
new: (i) The generation of start values for the root-finding based
on (cubic) Hermite interpolation in triangular elements, (ii) the
use of tailored mappings for predicting suitable start values in tetrahedral
elements, leading to optimal accuracy in contrast to results reported
in \cite{Fries_2015a,Fries_2016a}, and (iii) the study of integration
\emph{and} interpolation properties of the resulting reconstructed
surface elements.

\subsection{Valid level-set data and recursion\label{sub:ValididyAndRecursion}}

\begin{figure}
\centering

\subfigure[2D, triangle]{\includegraphics[width=4.3cm]{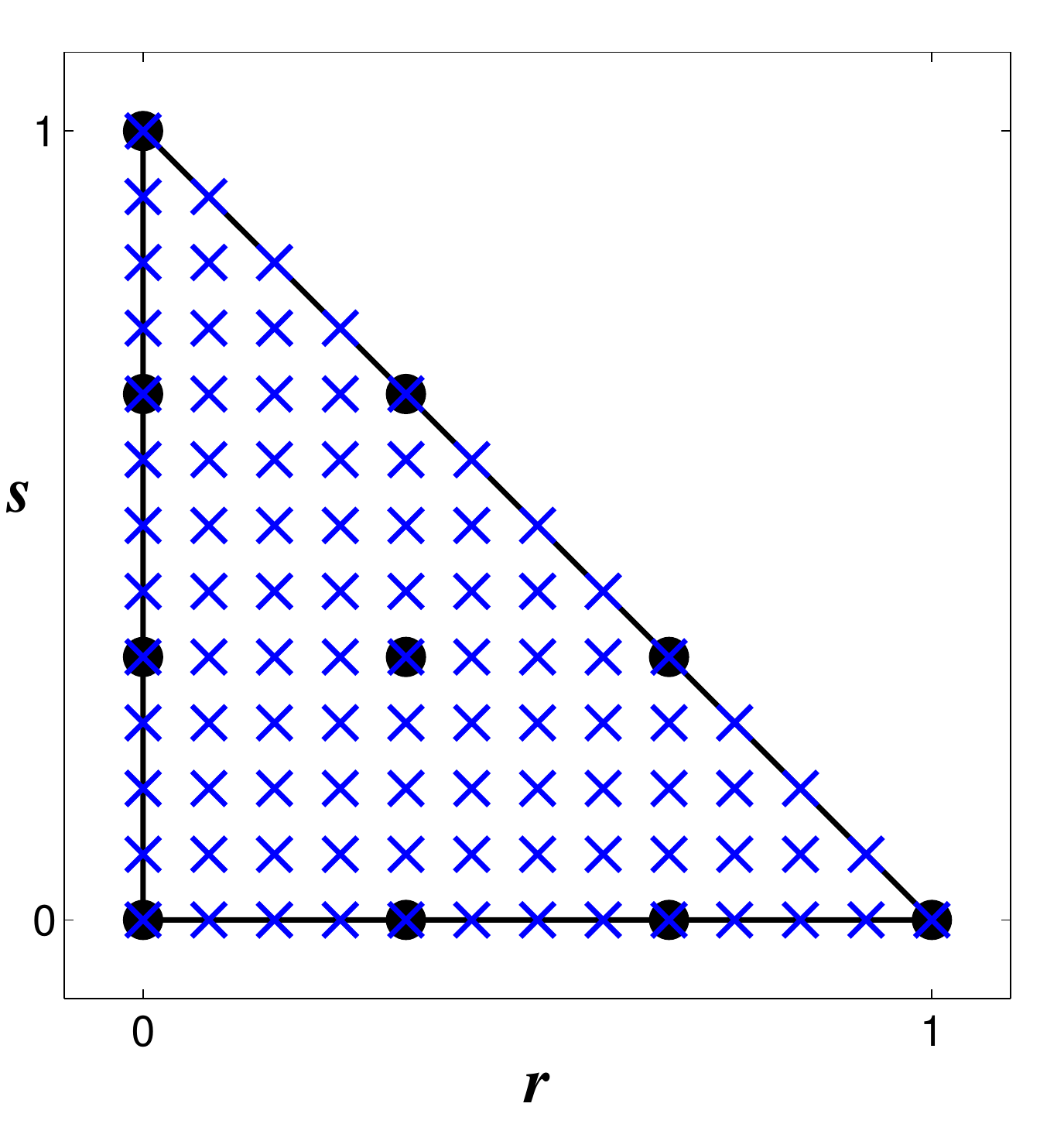}}\qquad\subfigure[3D, tetrahedron]{\includegraphics[width=5cm]{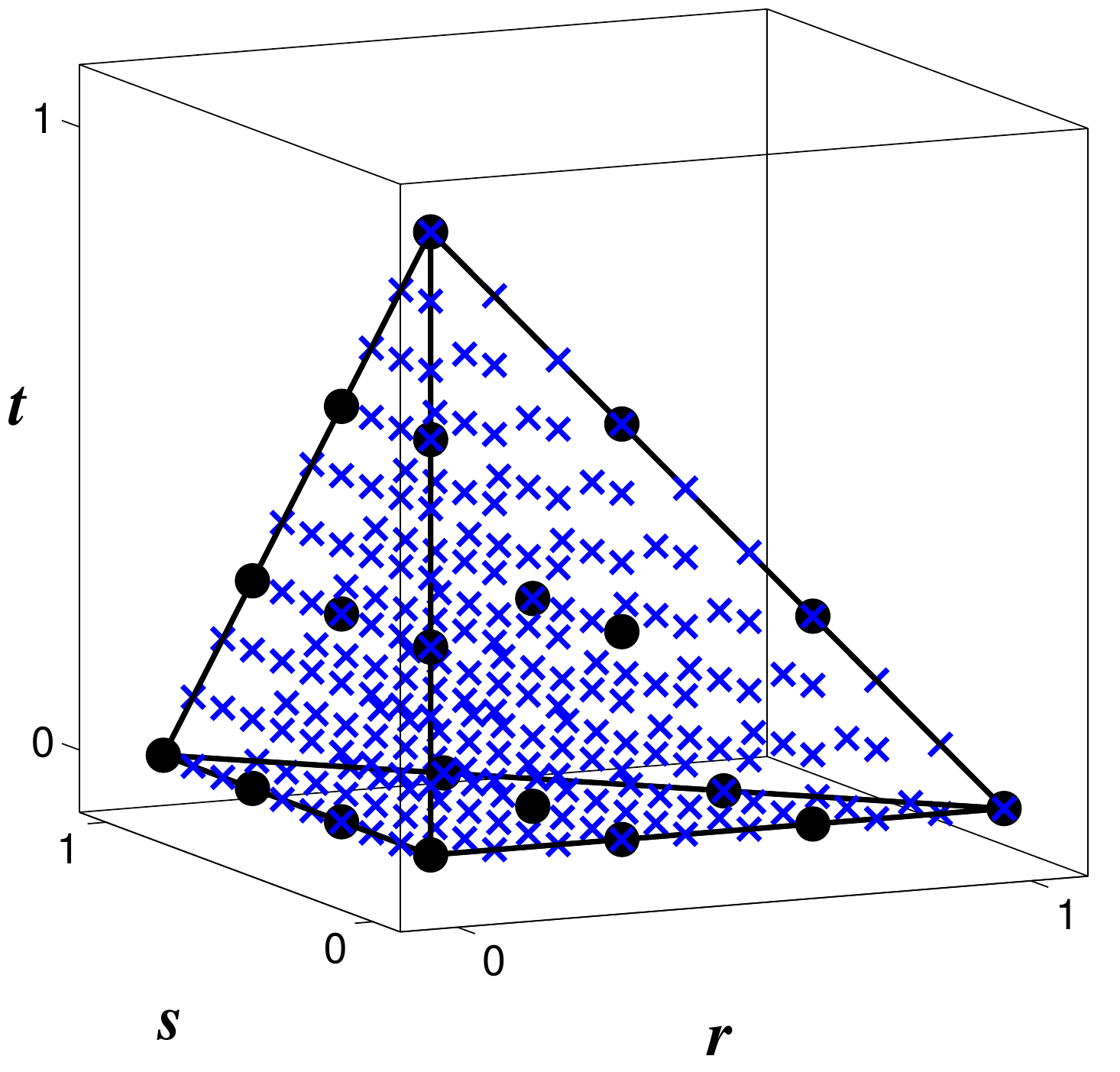}}

\caption{\label{fig:SampleGrids}Sample grids for different reference background
elements of order $3$. The black dots are elements nodes and blue
crosses are sample points.}
\end{figure}

\begin{figure}
\centering

\subfigure[]{\includegraphics[width=4cm]{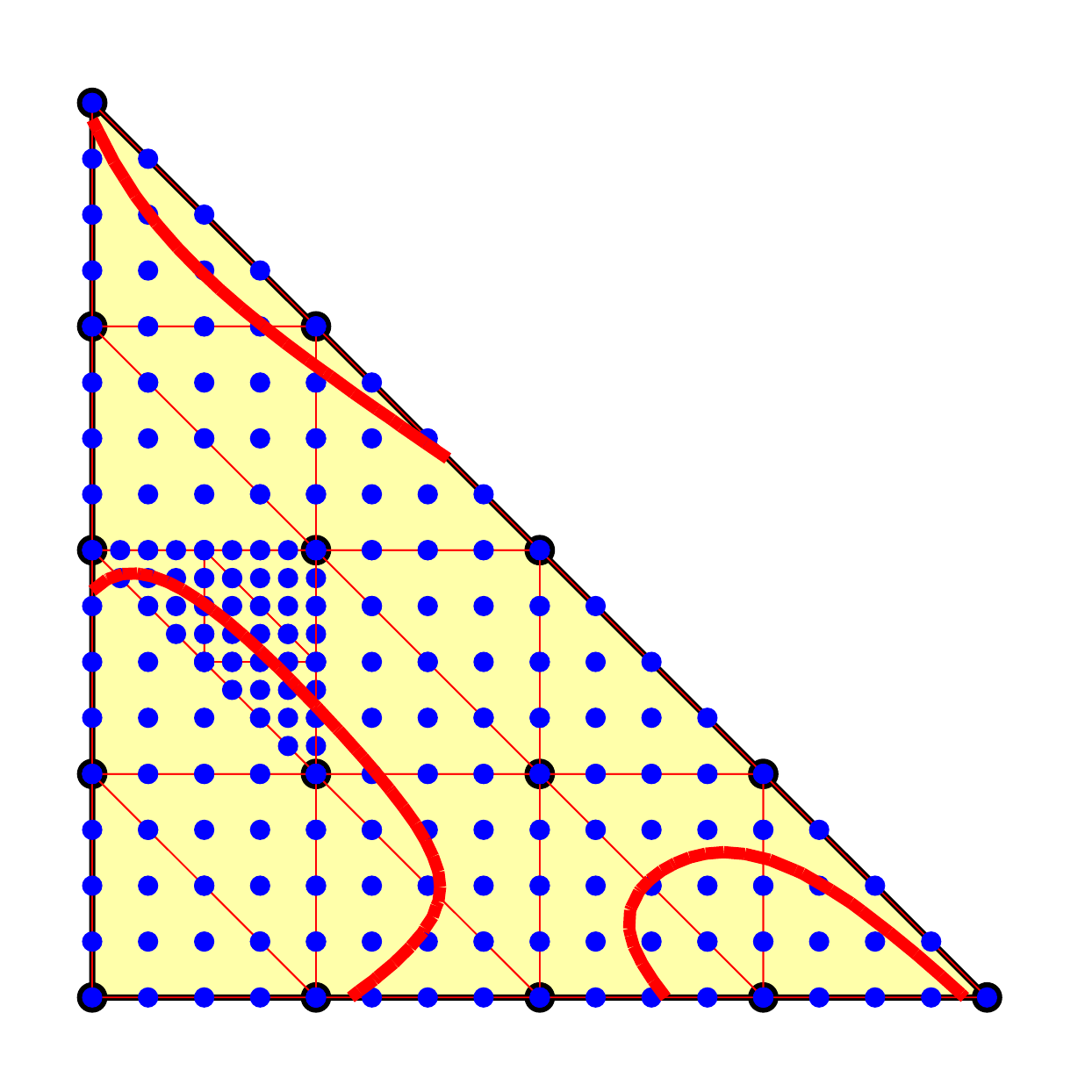}}\subfigure[]{\includegraphics[width=4cm]{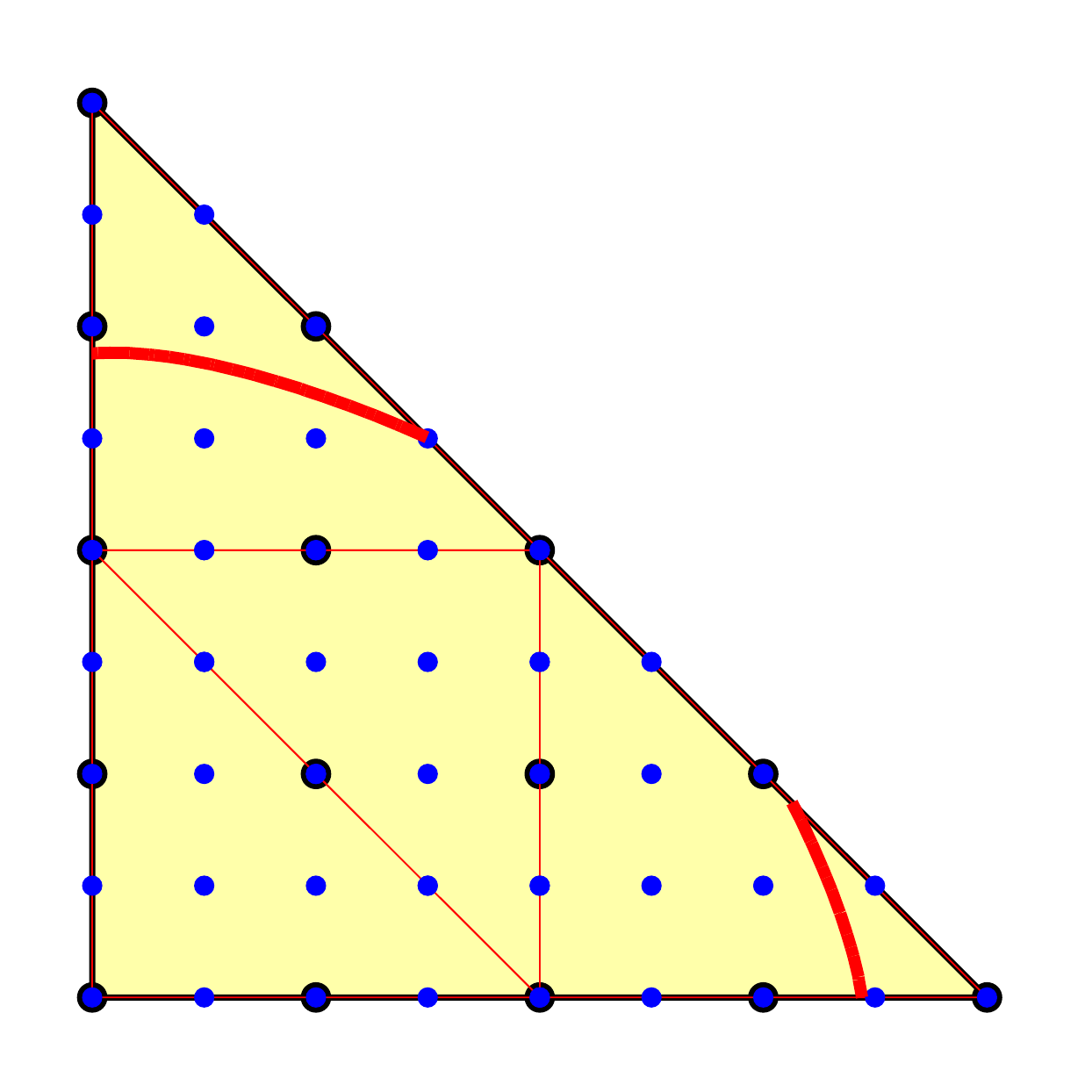}}\subfigure[]{\includegraphics[width=4cm]{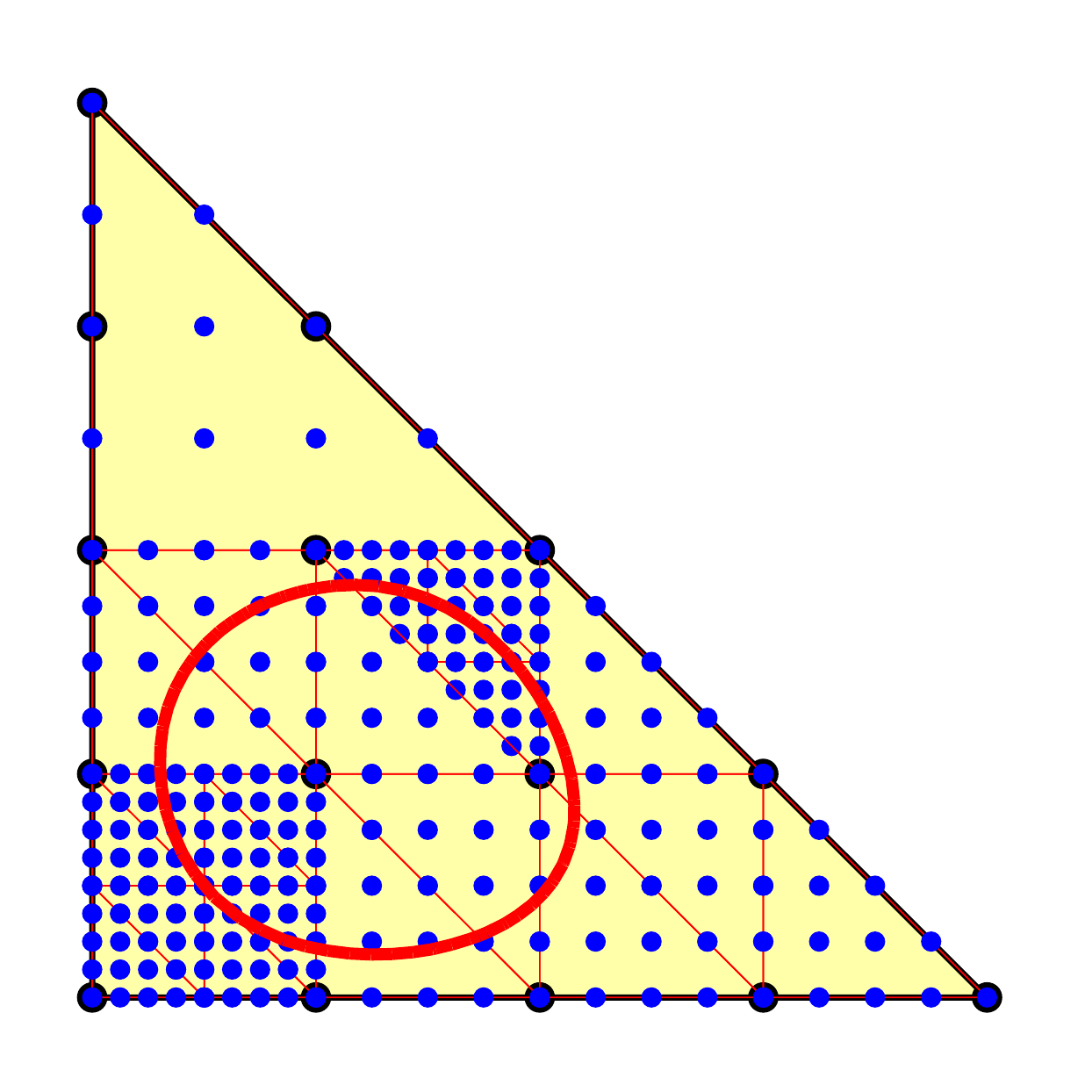}}\subfigure[]{\includegraphics[width=4cm]{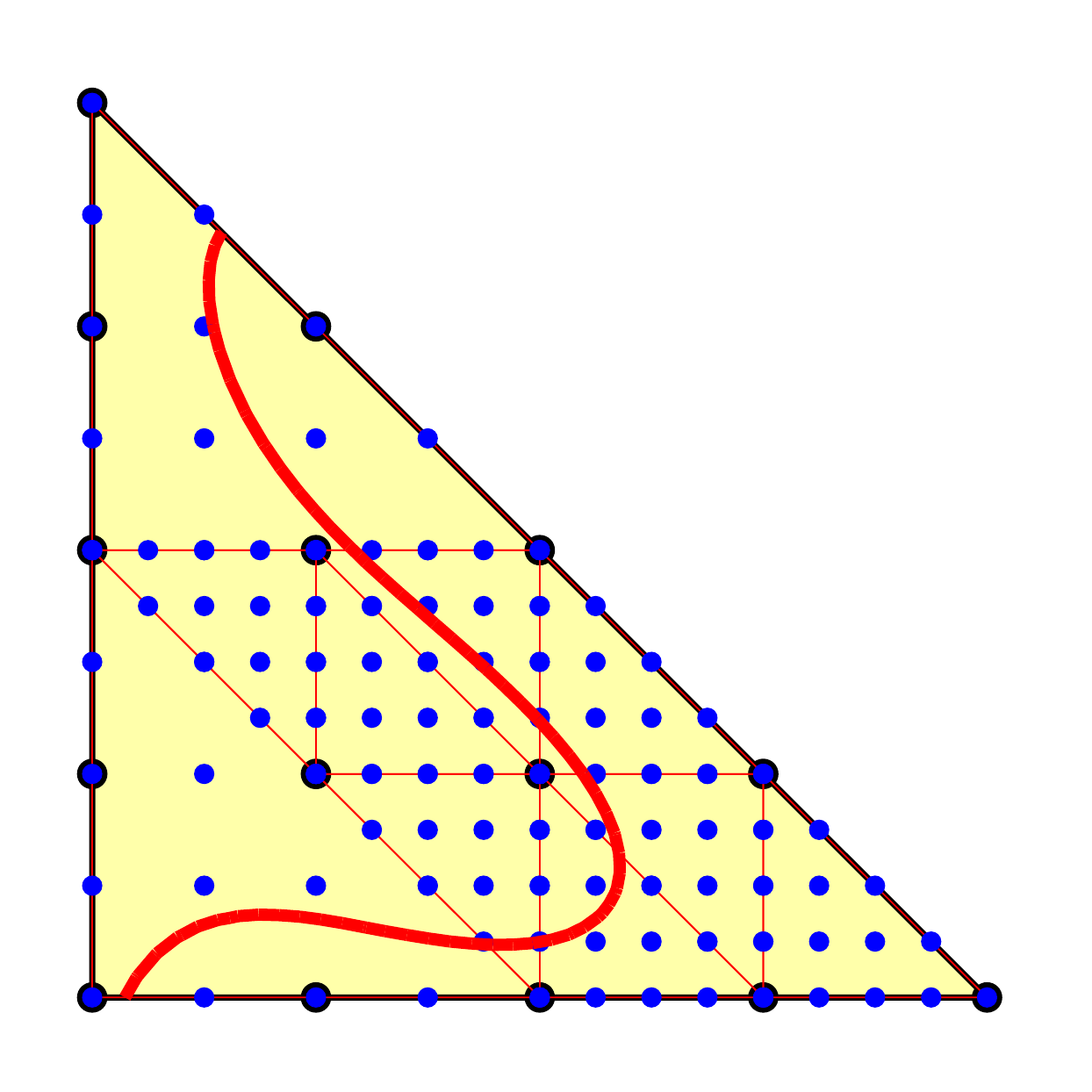}}

\caption{\label{fig:Vis2dIsoLinesInElem}Some complicated zero-level sets in
triangular elements. Note that after recursive refinements, valid
level-set data is obtained in the refined elements.}
\end{figure}

\begin{figure}
\centering

\subfigure[]{\includegraphics[width=6cm]{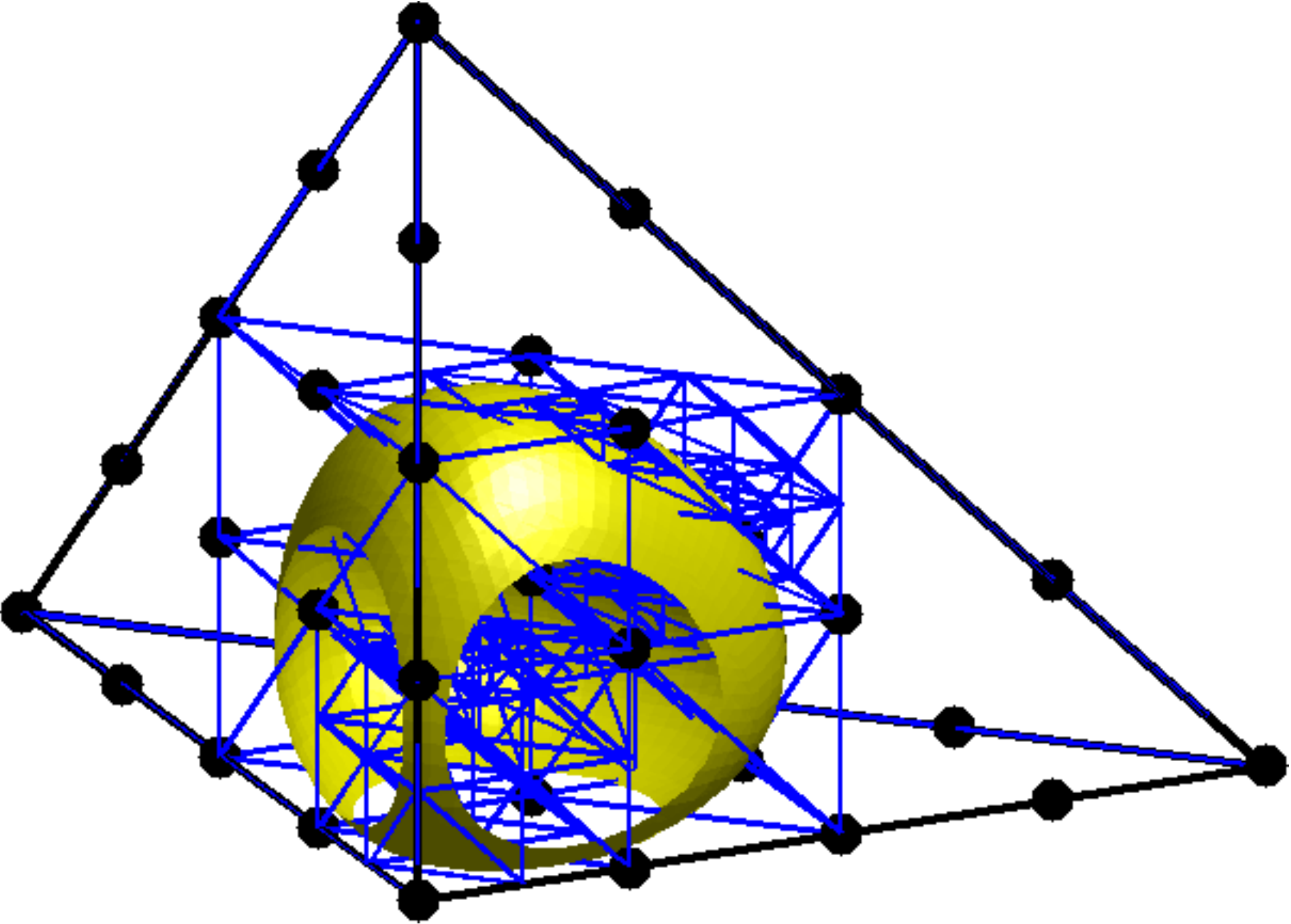}}\subfigure[]{\includegraphics[width=6cm]{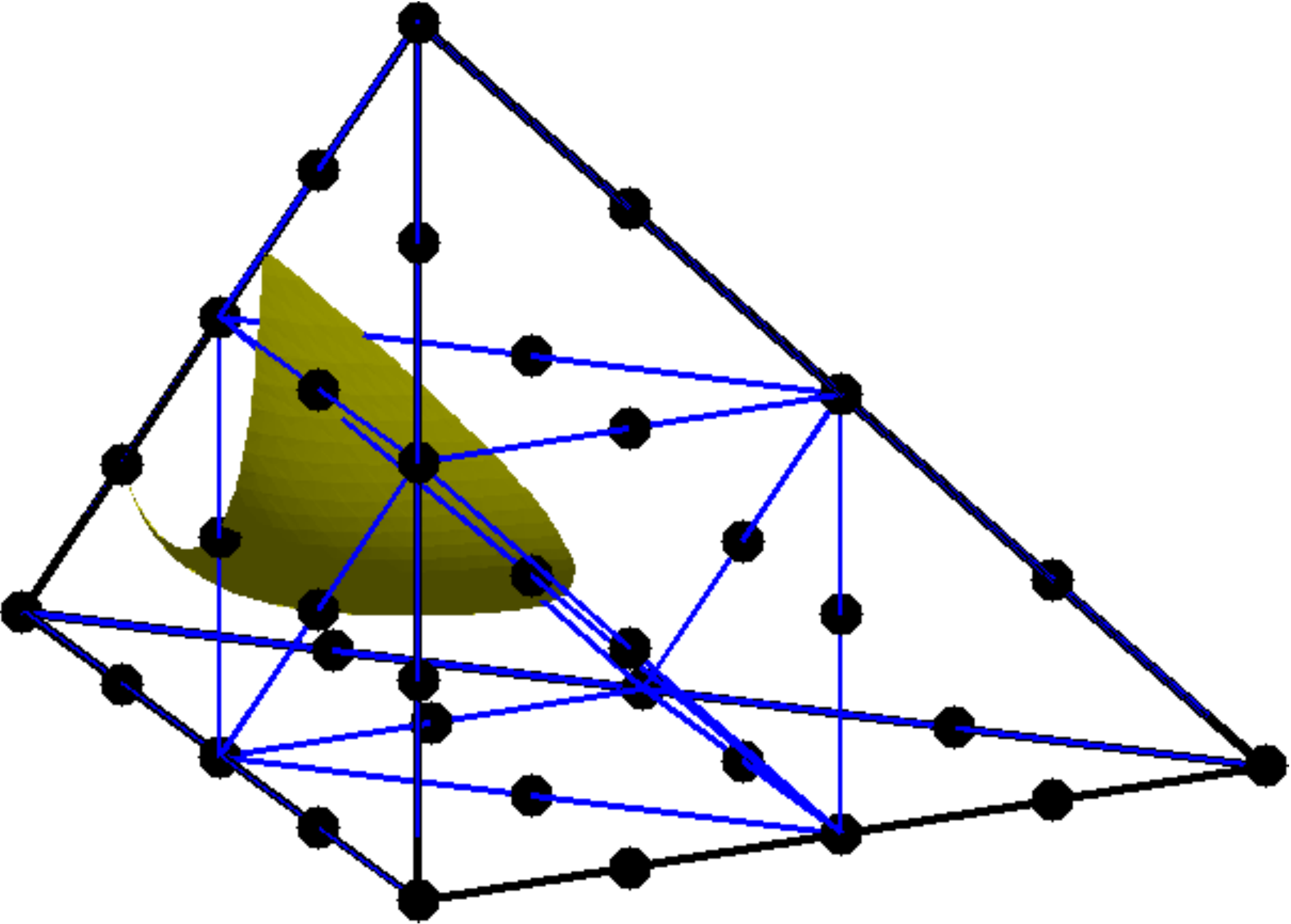}}

\caption{\label{fig:Vis3dIsoSurfInElem}Complicated zero-level sets in tetrahedral
elements and the required recursive refinements to achieve valid level-set
data.}
\end{figure}

It is obvious that higher-order elements may feature very complicated
topologies of zero-level sets. For example, element edges may be cut
several times, or the zero-level set is completely inside an element,
see Fig.~\ref{fig:Vis2dIsoLinesInElem} in 2D and Fig.~\ref{fig:Vis3dIsoSurfInElem}
in 3D. Therefore, a recursive procedure may be required (typically
only in very few elements) until ``valid level-set data'' inside
the element is obtained. In two dimensions, by ``valid'' we refer
to the situation where (i) each element edge is only cut once, (ii)
the overall number of cut edges must be two, and (iii) if no edge
is cut then the element is completely uncut. In three dimensions,
the conditions (i) to (iii) are checked on each element face plus
the additional condition that (iv) if no face is cut then the element
is completely uncut. Whether these conditions are fulfilled is checked
based on a sample grid in the reference element, see Fig.~\ref{fig:SampleGrids}.
Only the sign of the level-set data at the sample points is considered
for these checks.

If the level-set data is not valid the reference element is recursively
refined and the described procedures apply to the refined elements,
see Figs.~\ref{fig:Vis2dIsoLinesInElem} and \ref{fig:Vis3dIsoSurfInElem}.
The distinction between valid and invalid level-set data is done \emph{before
}the reconstruction and decomposition steps. However, the validity
check does not necessarily ensure that the subsequent reconstruction
and decomposition are successful. For example, when a resulting sub-element
after the decomposition features a negative Jacobian, a recursive
refinement is still triggered although the initial level-set data
was possibly valid in the above sense. This may, for example, be the
case for strongly curved zero-level sets, see Fig.~\ref{fig:Vis2dIsoLinesInElem}(d).
A level-set function must not be zero right at the corner nodes of
the background element; if that is the case, the level-set value is
simply perturbed in the range of $10^{-13}$ without any noticeable
change in the results.

\subsection{Topological cases\label{sub:Topologies}}

\begin{figure}
\centering

\subfigure[triangle]{\includegraphics[height=4cm]{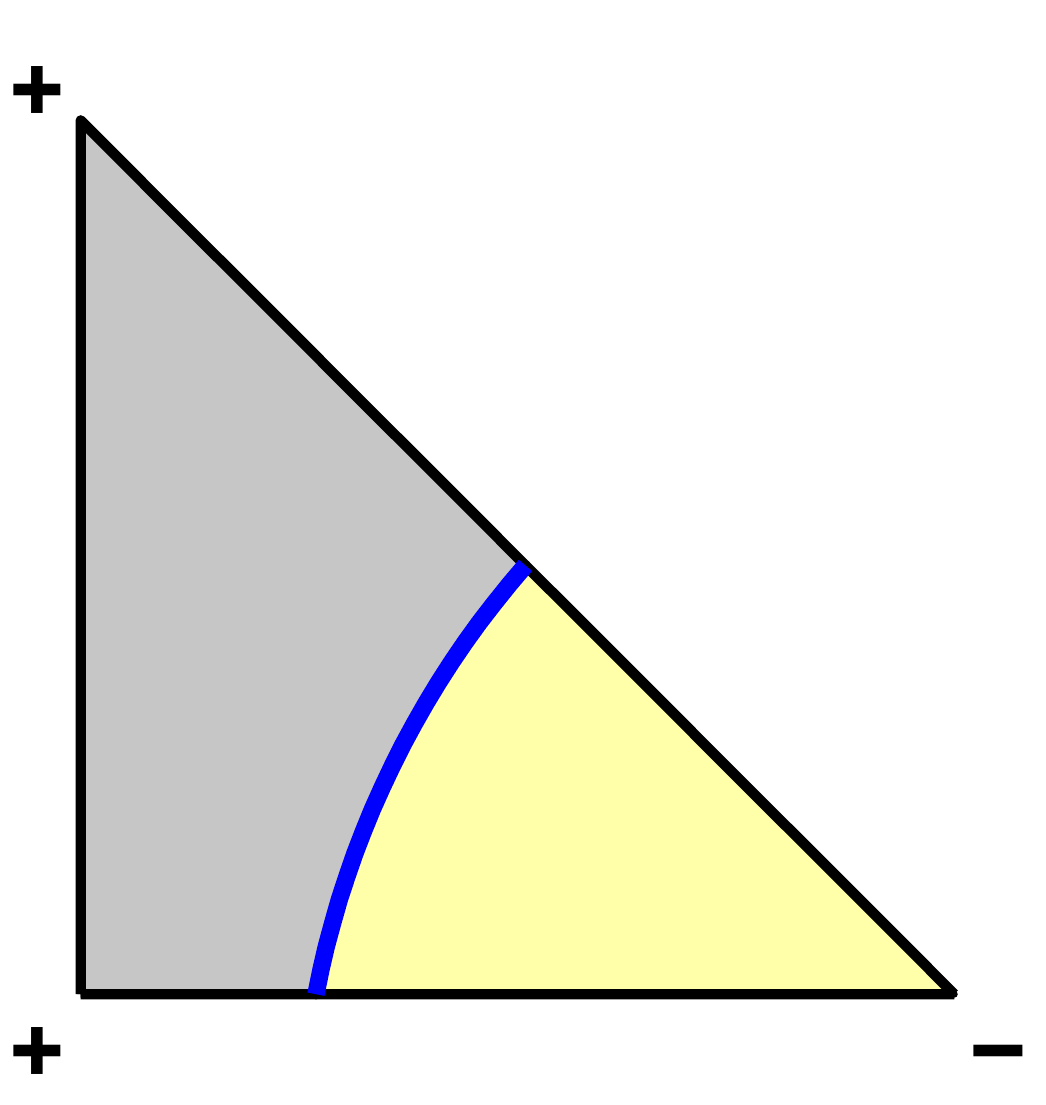}}\qquad\qquad\subfigure[tetrahedron, top.~1]{\includegraphics[height=4cm]{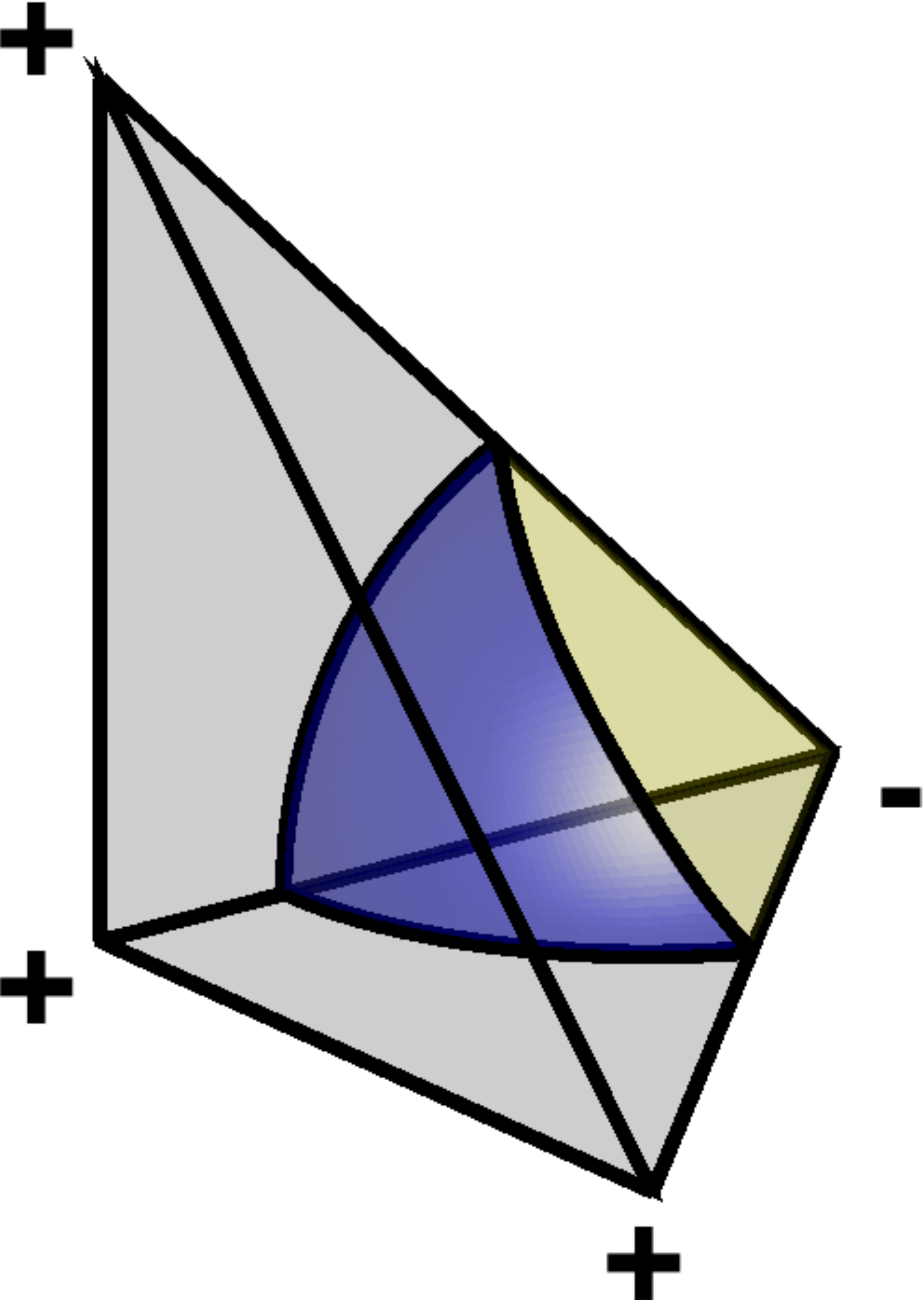}}\qquad\qquad\subfigure[tetrahedron, top.~2]{\includegraphics[height=4cm]{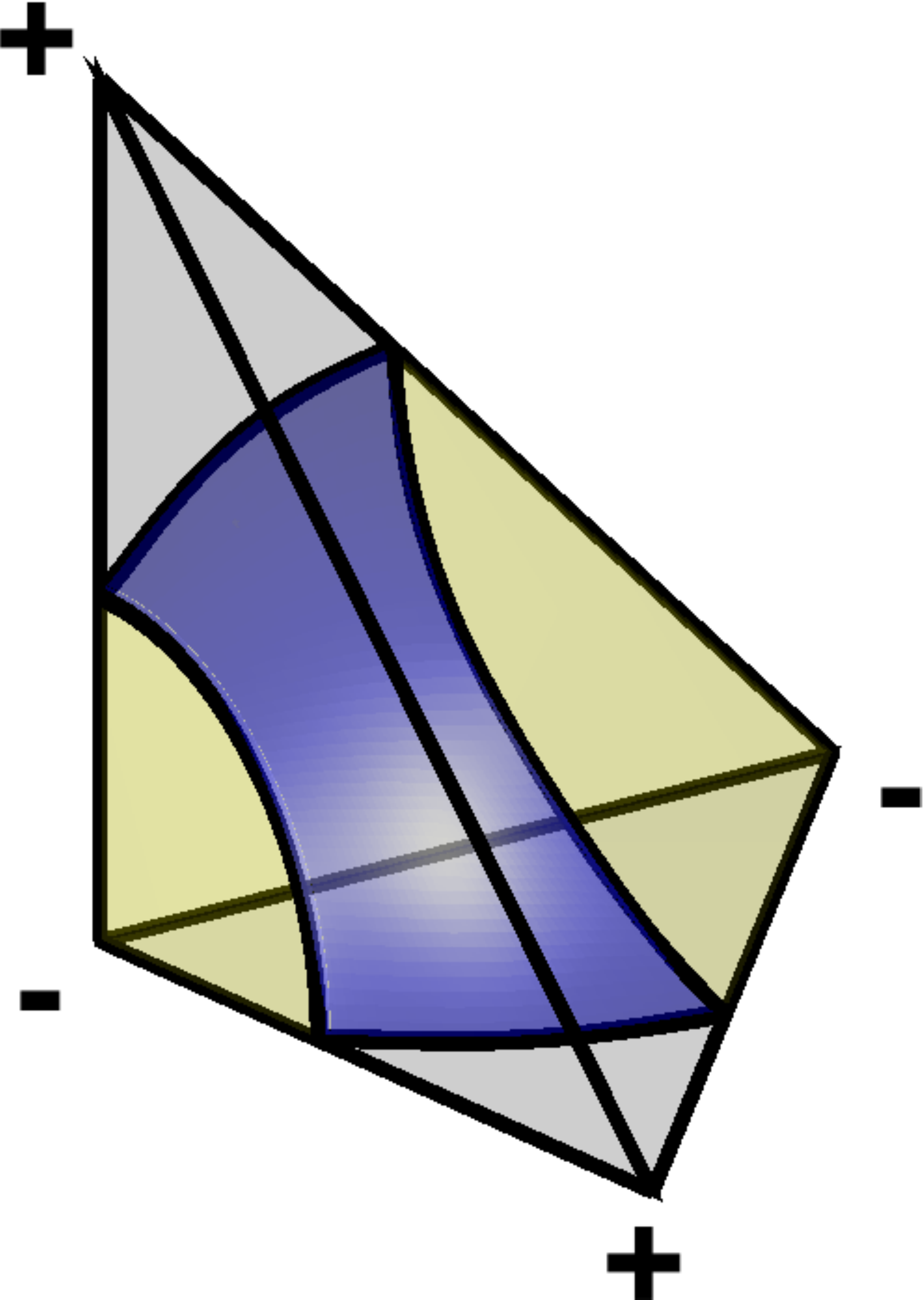}}

\caption{\label{fig:TopologicalCases}Topologically different cut situations
for valid level-set data in triangles and tetrahedra.}
\end{figure}

From now on, it is assumed that valid level-set data is present in
a (refined) reference element. Then, the zero-level set cuts the reference
element in a limited number of topologically different cases, see
Fig.~\ref{fig:TopologicalCases}. Only the signs at the \emph{corner}
nodes of the reference element determine the situation. A cut triangle
is decomposed into one sub-triangle and one sub-quadrilateral; the
latter may be further decomposed into sub-triangles. A cut tetrahedron
is decomposed either into one sub-tetrahedron and one sub-prism (topology
1), see Fig.~\ref{fig:TopologicalCases}(b), or into two sub-prisms
(topology 2), see Fig.~\ref{fig:TopologicalCases}(c). Obviously,
the sub-prisms may be further decomposed into tetrahedra. There is
no need to avoid hanging nodes in the decomposition into sub-elements
for integration purposes.

\subsection{Root search and interface elements\label{sub:InterfaceElements}}

Depending on the topological situation, different algorithms for the
identification of the element nodes of the interface elements are
employed. Each node is a root of the level-set function $\phi^{h}\left(\vek r\right)$
which is obtained by solving a non-linear problem. The root finding
is characterized by (i) the start values of the iteration, (ii) the
search directions in the reference element, and (iii) the iteration
method. The discussion is along the lines of \cite{Fries_2015a}.

\subsubsection{Start values and search directions in 2D\label{sub:StartValuesIn2d}}

\begin{figure}
\centering

\subfigure[11]{\includegraphics[width=4cm]{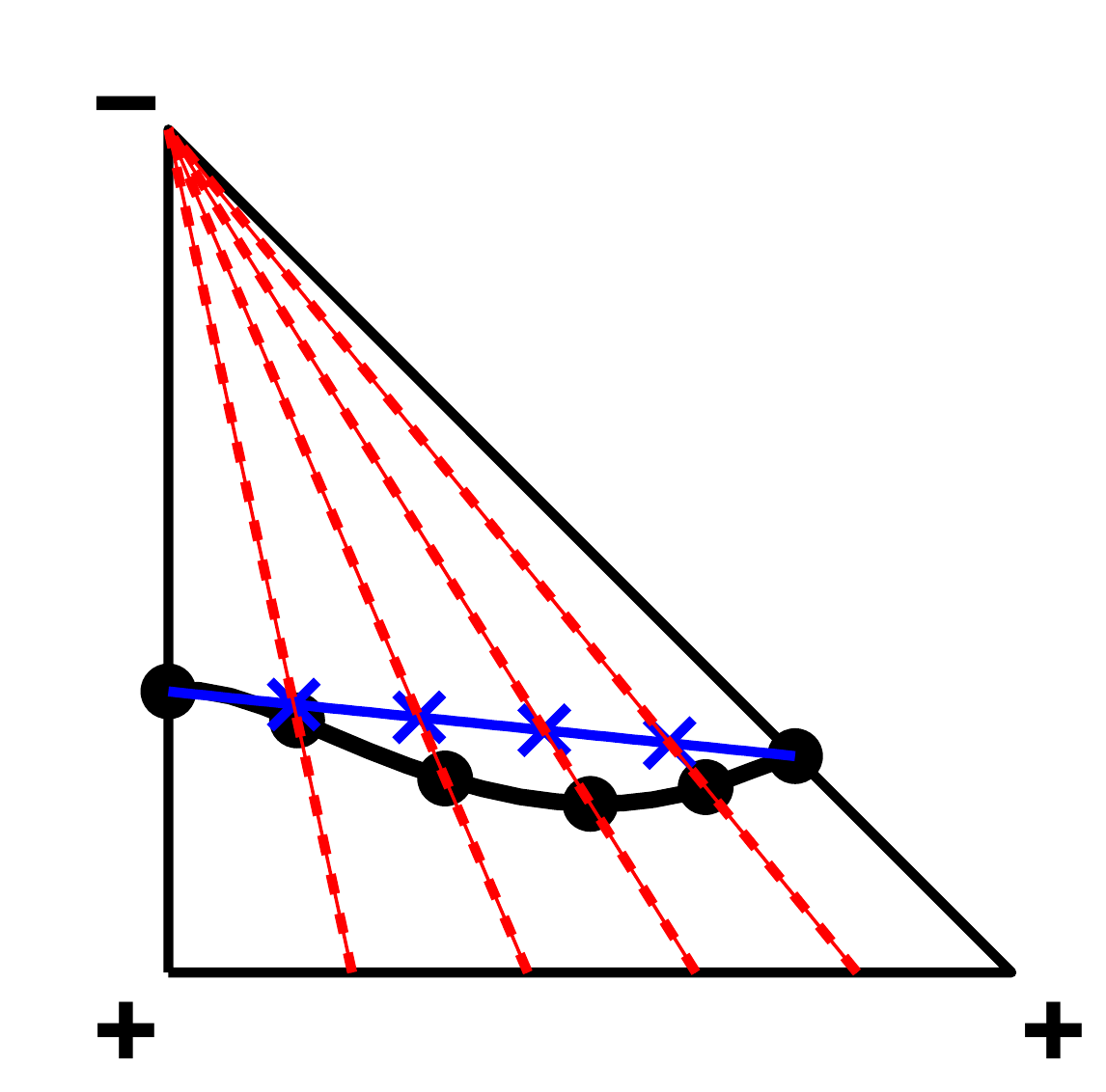}}\subfigure[12]{\includegraphics[width=4cm]{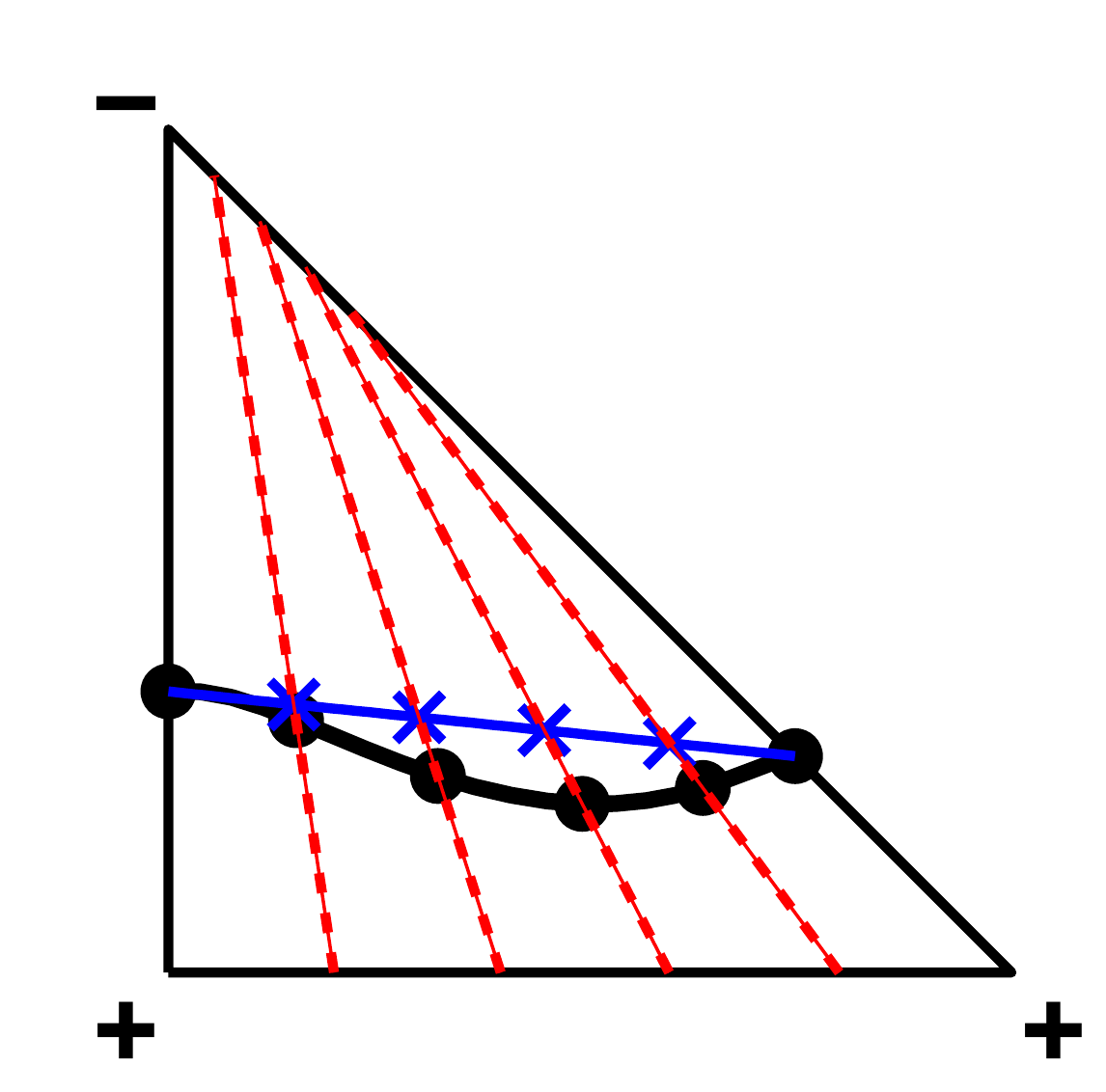}}\subfigure[13]{\includegraphics[width=4cm]{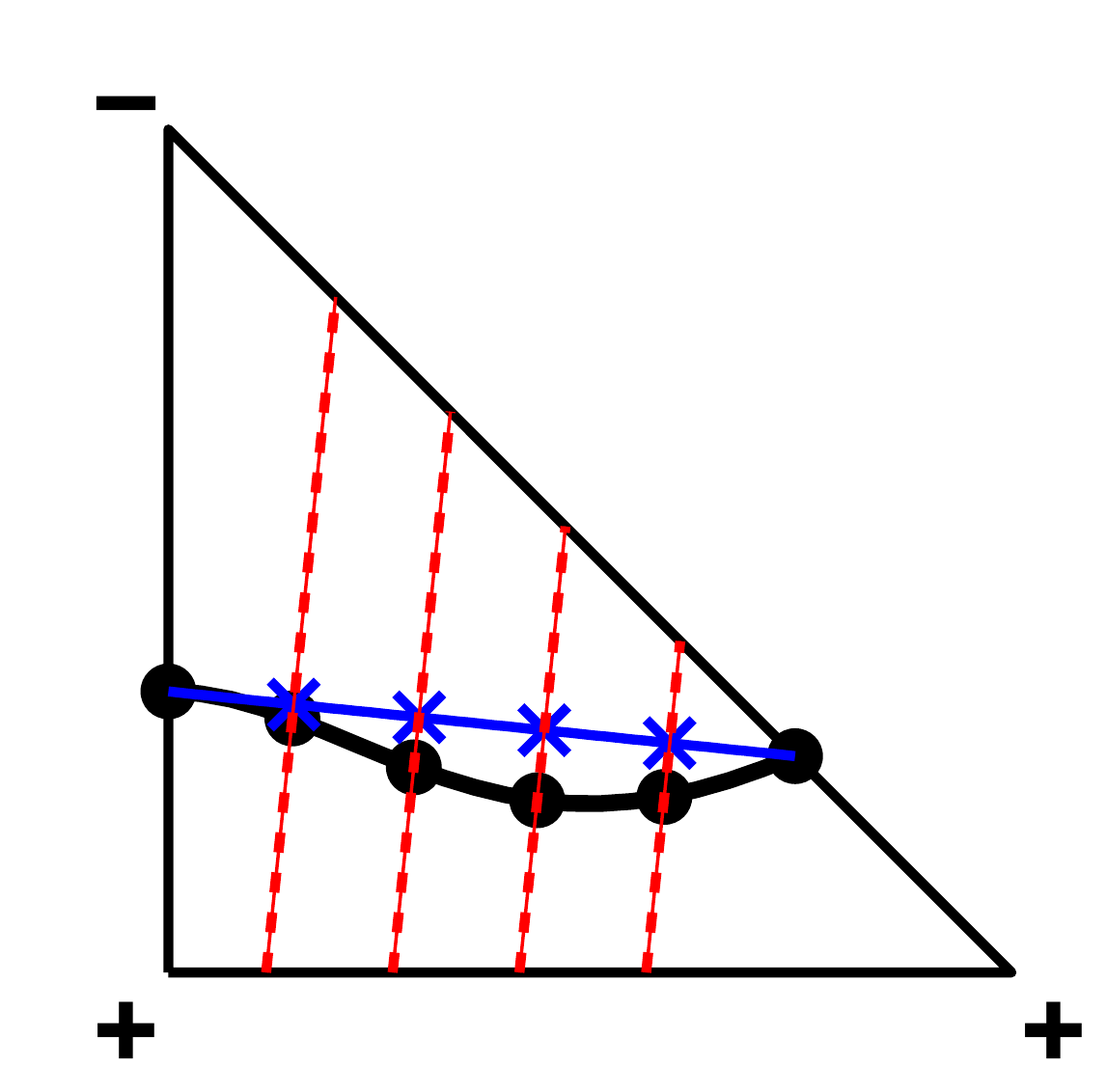}}\subfigure[14]{\includegraphics[width=4cm]{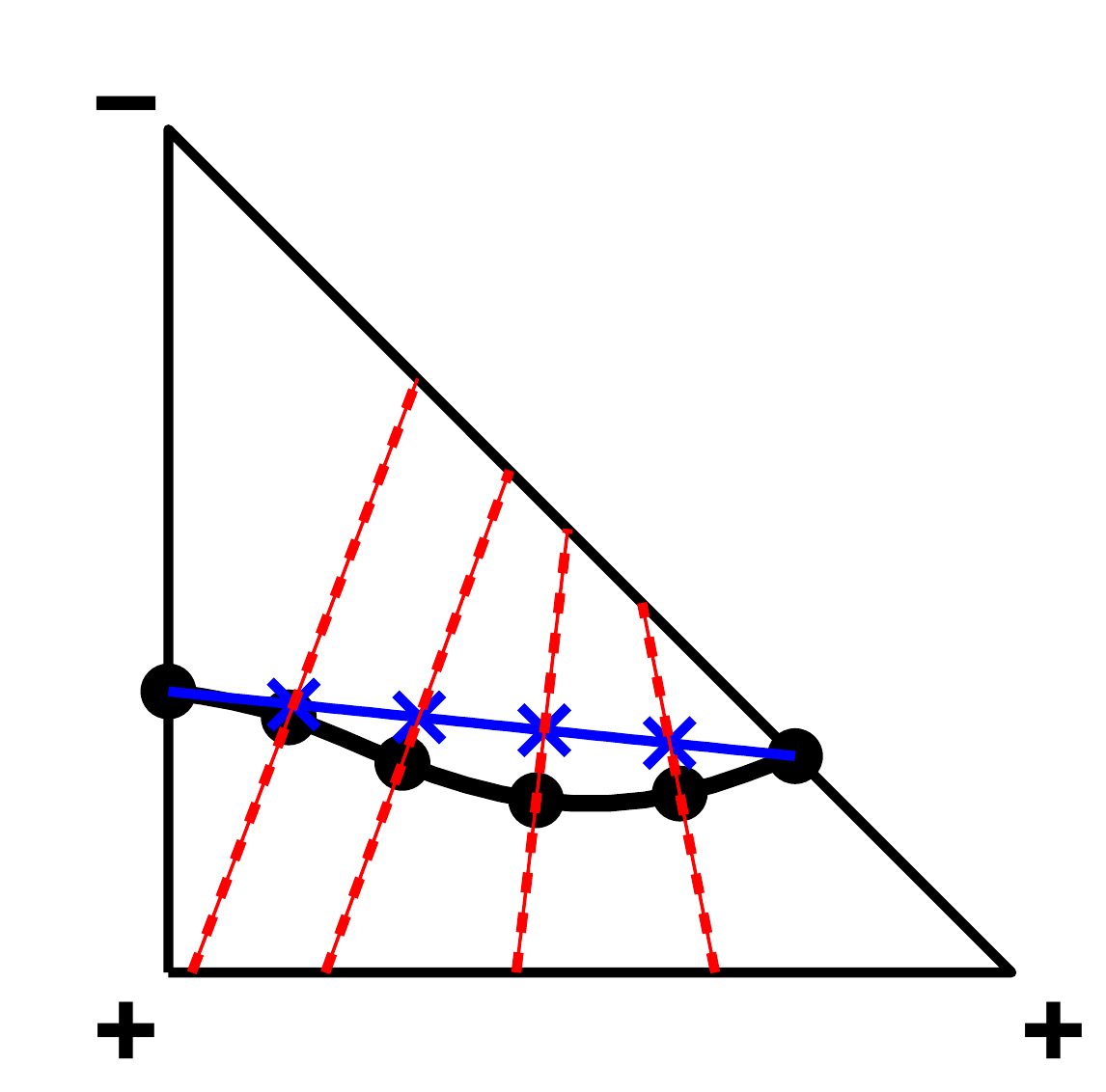}}

\subfigure[21]{\includegraphics[width=4cm]{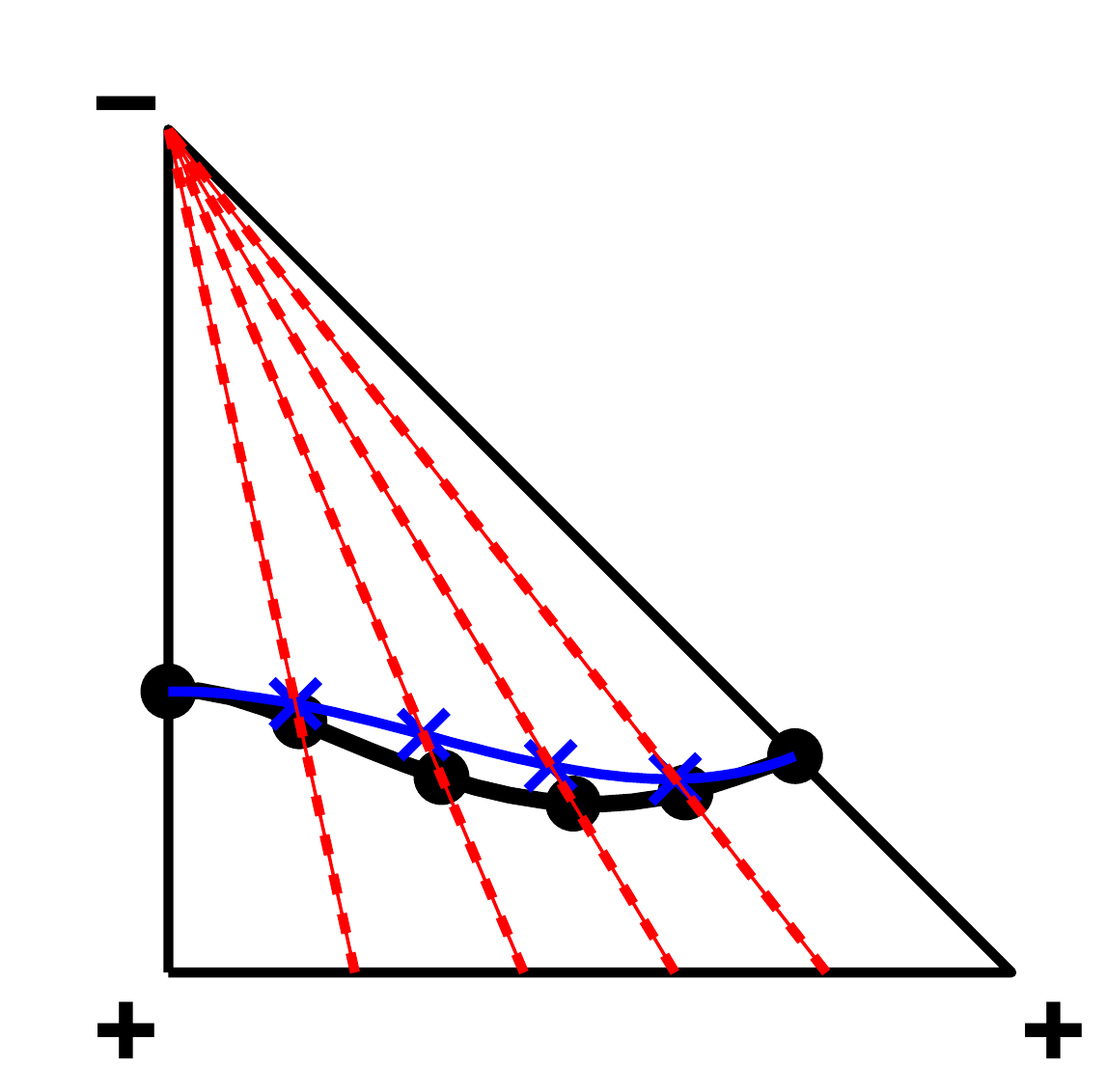}}\subfigure[22]{\includegraphics[width=4cm]{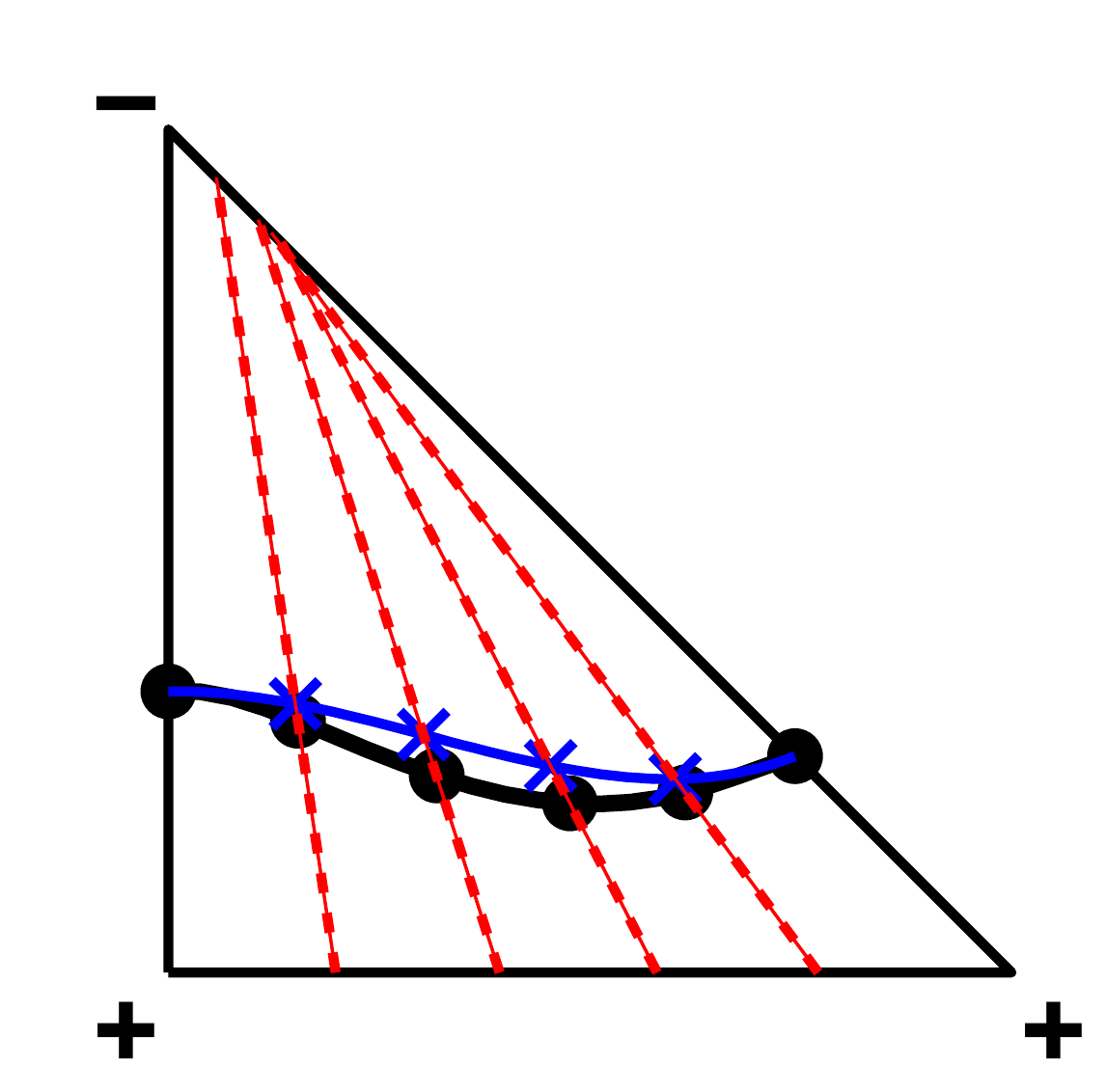}}\subfigure[23]{\includegraphics[width=4cm]{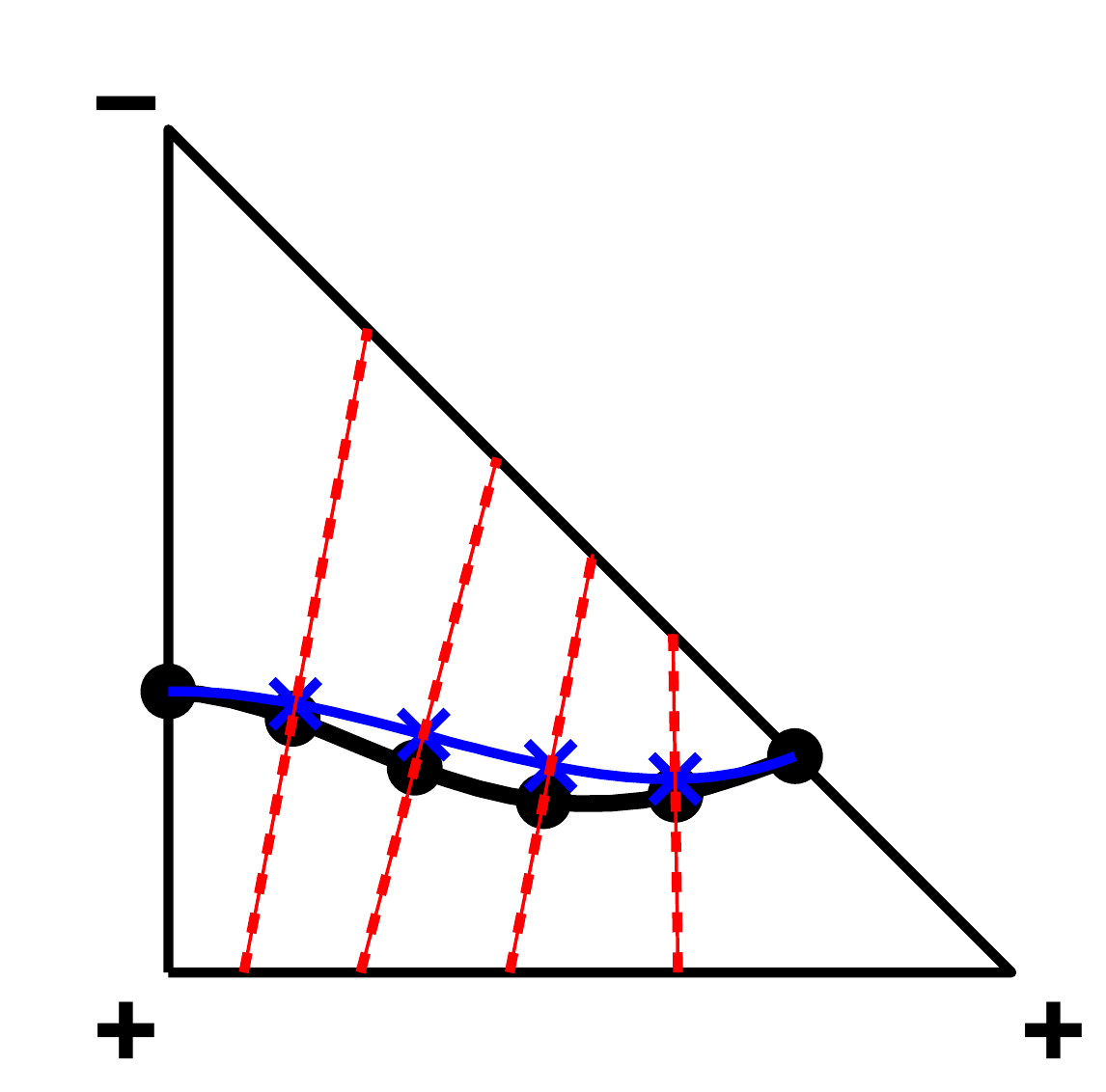}}\subfigure[24]{\includegraphics[width=4cm]{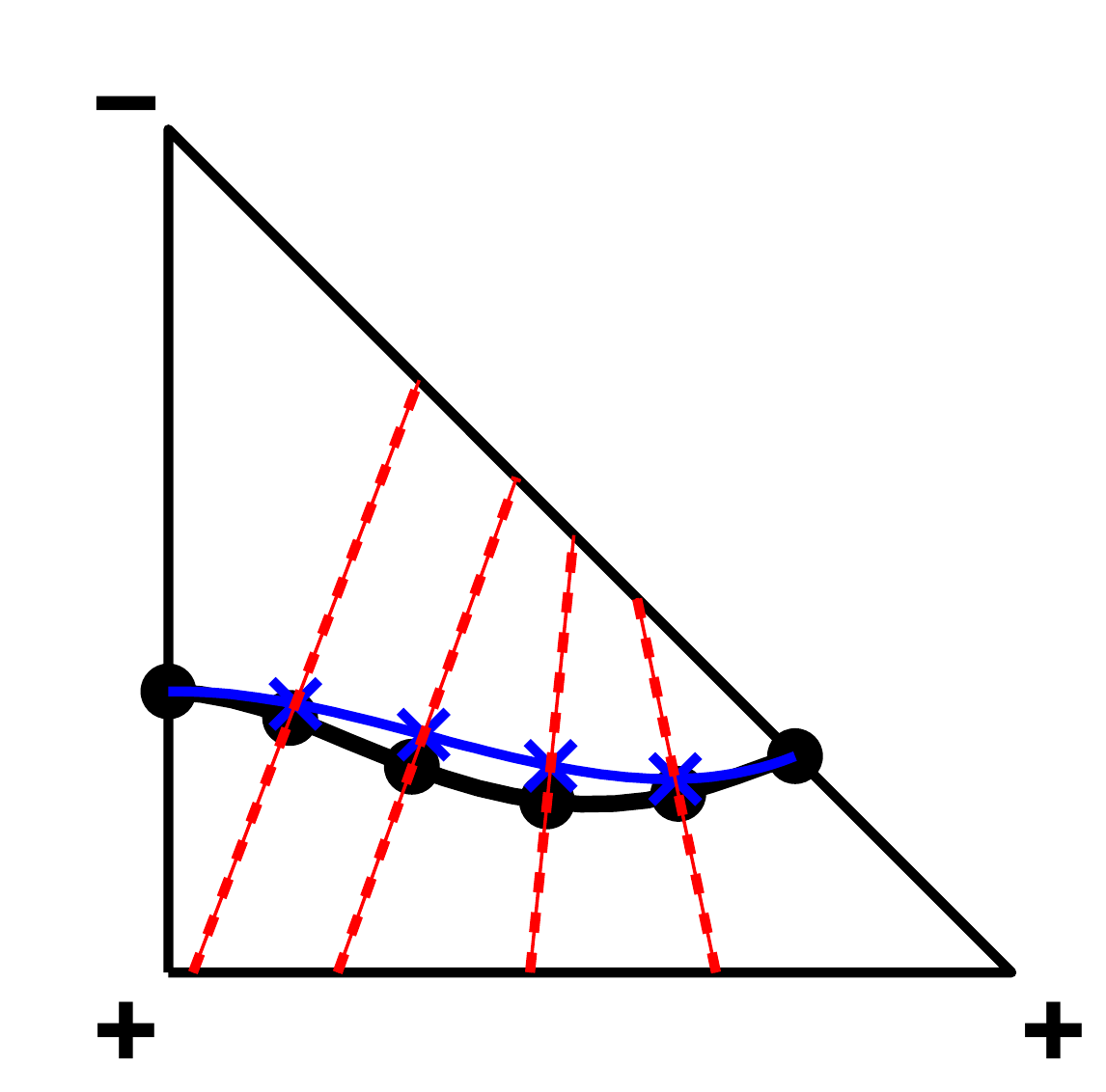}}

\caption{\label{fig:VisSearchDirTri}Starting points (blue), search paths (red)
and reconstructed interface elements (black) in triangular reference
elements. (a) to (d) are based on a linear and (e) to (h) on a Hermite
reconstruction. Four different search directions are investigated
for each case. Each variant is identified by a two-digit number.}
\end{figure}

The first step is to identify the intersection of the zero-level set
with the element edges. In triangular reference elements, these are
two points. One may then define a linear interpolation of these two
points or a Hermite interpolation which takes into account also the
direction of the zero-level set at the intersection points. In order
to obtain the start guesses, an interface element of the desired order
is mapped onto the linear or Hermite intermediate reconstruction.
These are the blue lines and crosses in Fig.~\ref{fig:VisSearchDirTri}
(a) to (d) for the linear and (e) to (h) for the Hermite case. 

For the search directions, four different variants are studied: 
\begin{enumerate}
\item Towards the node on the other side than the other two, see Fig.~\ref{fig:VisSearchDirTri}(a)
and (e).
\item In direction of the interpolated edge directions at the intersection
points, see Fig.~\ref{fig:VisSearchDirTri}(b) and (f).
\item In normal direction to the linear or Hermite reconstruction, see Fig.~\ref{fig:VisSearchDirTri}(c)
and (g).
\item In direction of the gradient of the level-set function, i.e.~$\nabla\phi^{h}\left(\vek r\right)$,
see Fig.~\ref{fig:VisSearchDirTri}(d) and (h).
\end{enumerate}
For later reference in the numerical studies, the 2D variants are
characterized by two-digit numbers: The first digit is either $1$
for a linear or $2$ for a Hermite reconstruction. The second digit
refers to the four variants of the search directions listed above.

\subsubsection{Start values and search directions in 3D\label{sub:StartValuesIn3d}}

In three dimensions, the element nodes of the higher-order interface
element (triangle for topology 1, quadrilateral for topology 2) are
separated into outer and inner nodes. The outer (i.e., edge) nodes
are enforced to remain on the faces of the reference tetrahedron.
In fact, this is achieved by treating one face after the other and
performing a reconstruction in reference triangles as described above.
Thereafter, these nodes are mapped to the corresponding face and define
the edge nodes of the interface element.

For the inner nodes, the start values and search directions depend
on an intermediate reconstruction of the zero-isosurface based on
the mapping defined in \cite{Solin_2003a}. This mapping defines a
surface implied by the three or four curved edges (higher-order line
elements) of the triangular or quadrilateral zero-isosurface within
the cut tetrahedron, see Fig.~\ref{fig:VisSolinReconTetr}. The definition
of the mapping is outlined in the appendix \ref{sub:MappingTriQuad}
of this work. Next, a higher-order reference interface element is
mapped onto the intermediate reconstruction to obtain the start values,
see the blue crosses in Fig.~\ref{fig:VisSolinReconTetr}. The search
directions are then, as above, either the corresponding normal vectors
or the gradients of the level-set function.

It is thus seen that also in 3D, the 2D-variants play a major role
as they fully define the outer nodes and are also important to define
start values for the inner nodes through the intermediate reconstructions.
It is noted that, in contrast to the definition of the start values
in 3D as proposed in \cite{Fries_2015a}, optimal convergence rates
are achieved with the procedure proposed herein.

\begin{figure}
\centering

\subfigure[Topology 1]{\includegraphics[width=5cm]{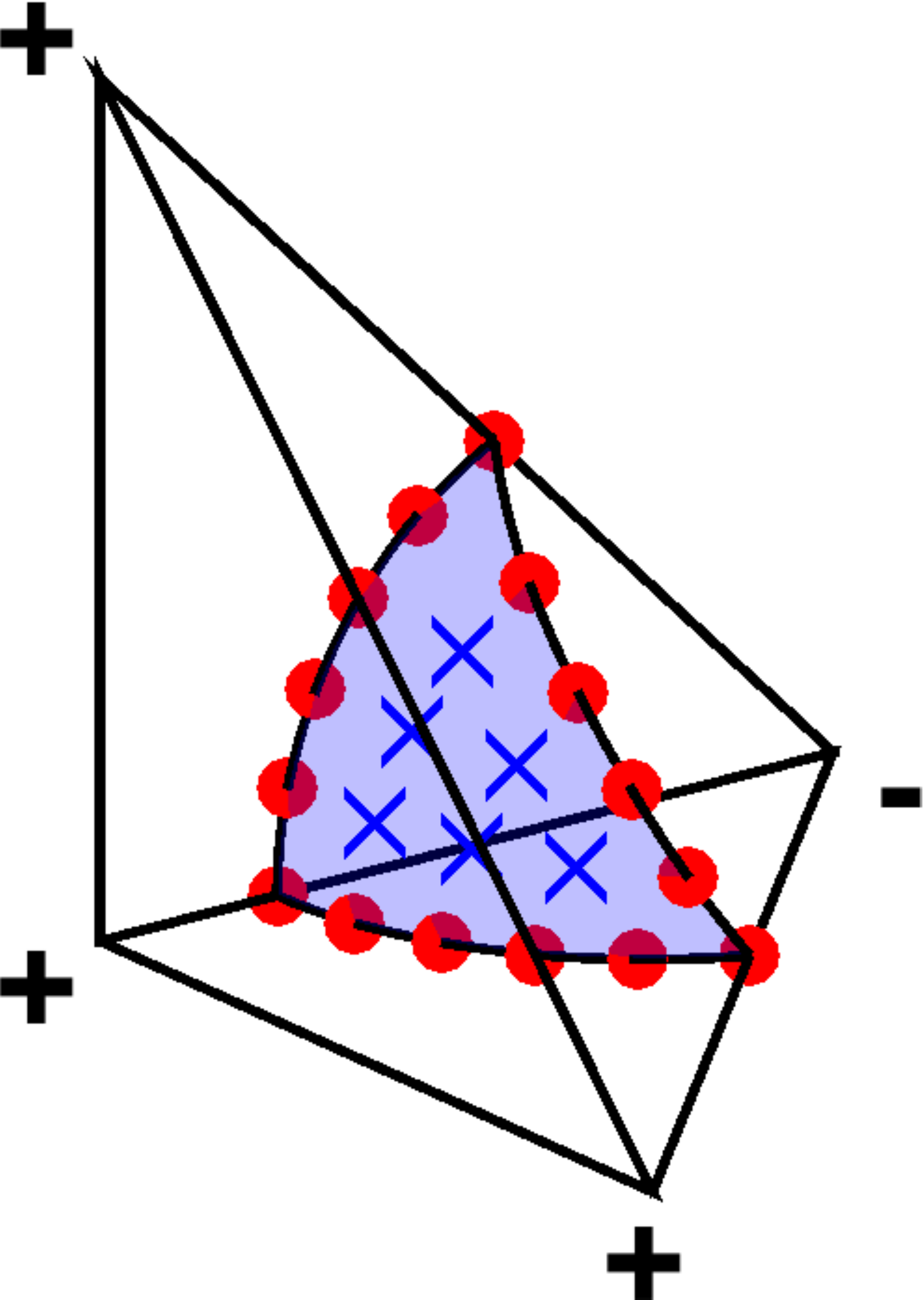}}\qquad\subfigure[Topology 2]{\includegraphics[width=5cm]{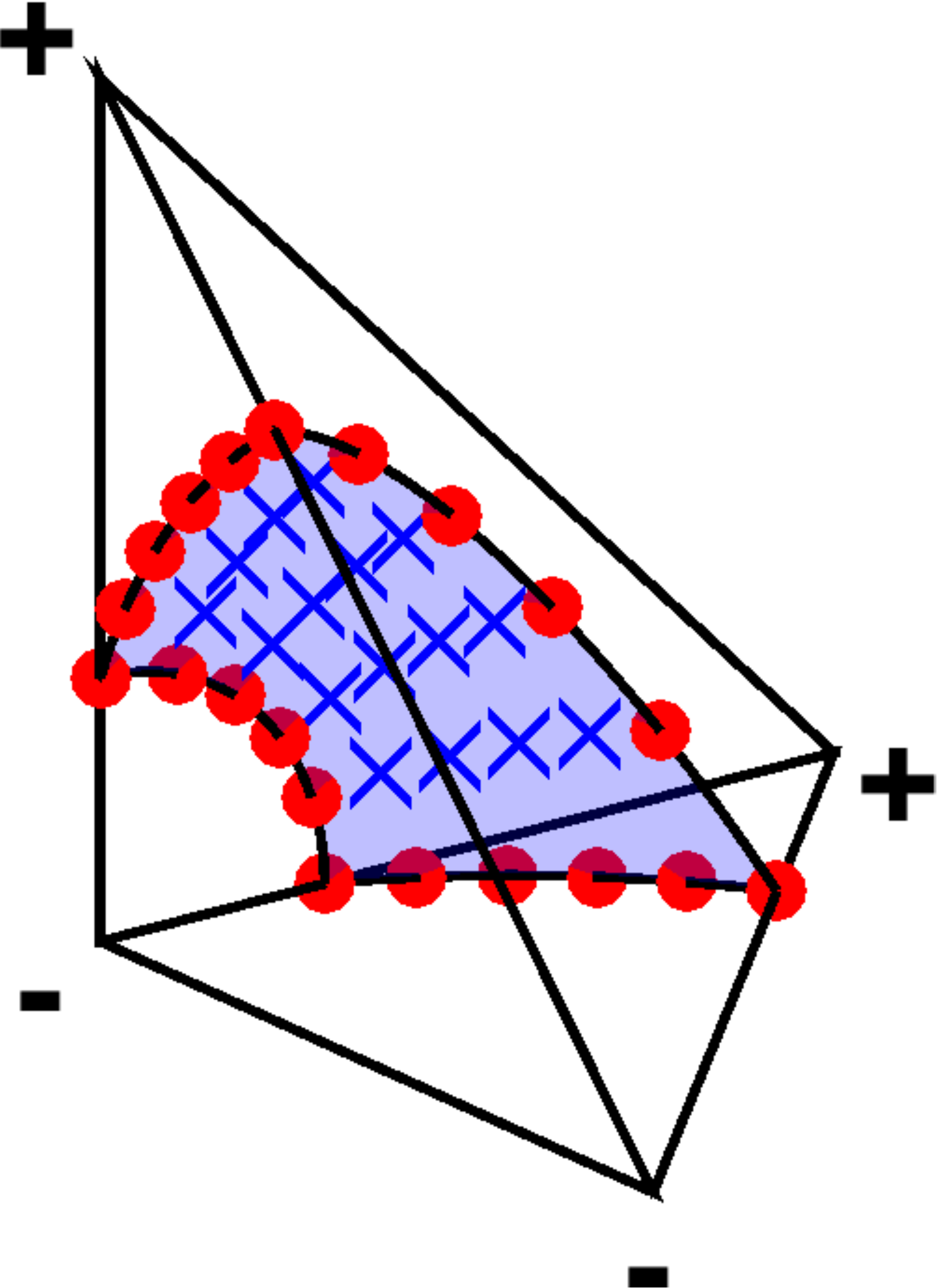}}

\caption{\label{fig:VisSolinReconTetr}Reconstruction of a $5$th-order interface
element in a tetrahedron; this is a triangle for topology 1 and a
quadrilateral for topology 2. The outer (edge) nodes (red dots) on
the faces are reconstructed based on 2D triangles. The start values
for the inner nodes (blue crosses) are obtained by a $5$th-order
reference element mapped to the surface implied by the outer nodes,
see the appendix \ref{sub:MappingTriQuad}.}
\end{figure}

\subsubsection{Iterative procedure\label{sub:IterativeProcedure}}

\begin{figure}
\centering

\subfigure[$\vek N=$normal vector]{\includegraphics[width=3.7cm]{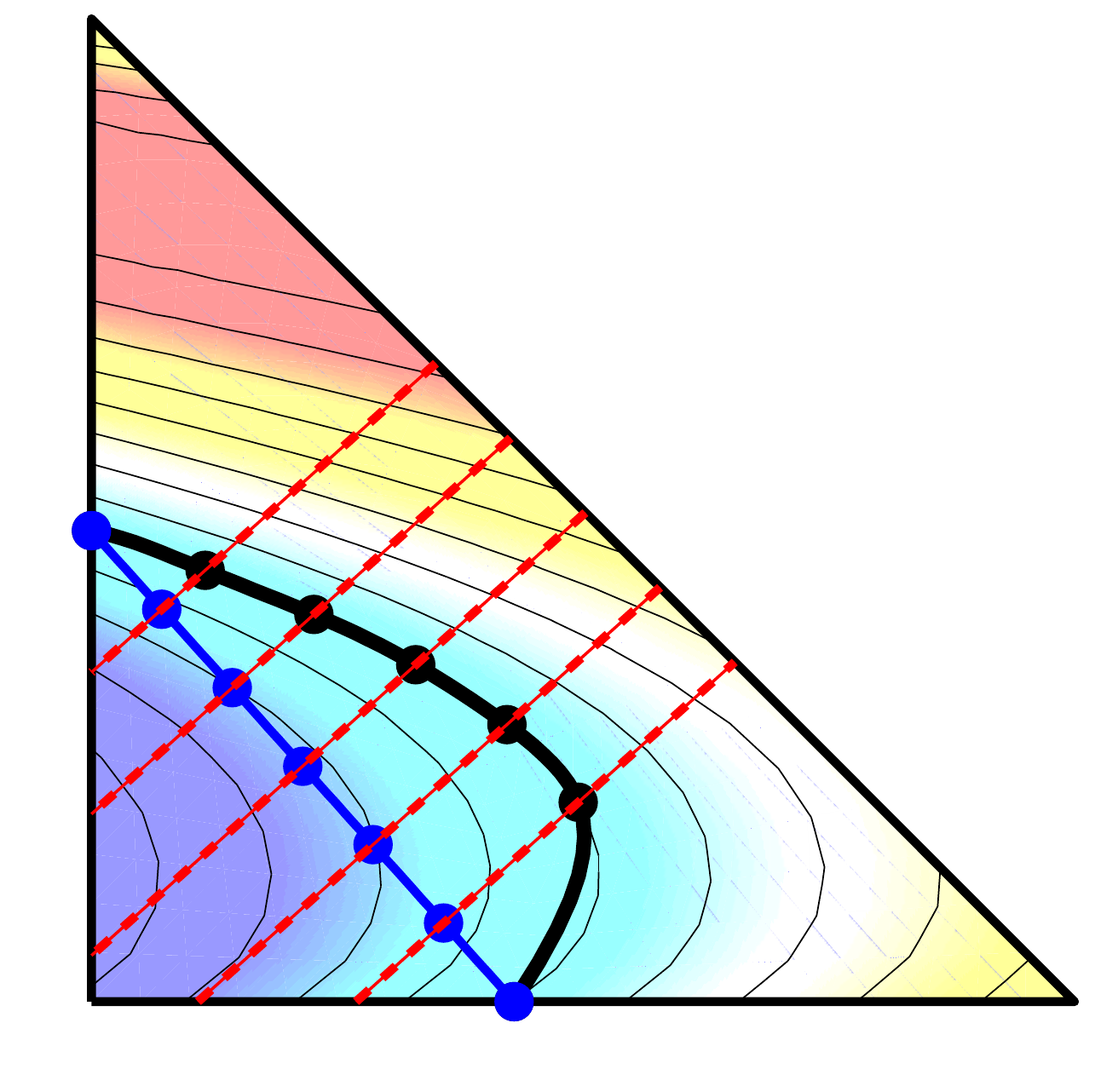}}\qquad\subfigure[$\vek N=\nabla\phi\left(\vek r^{0}\right)$]{\includegraphics[width=3.7cm]{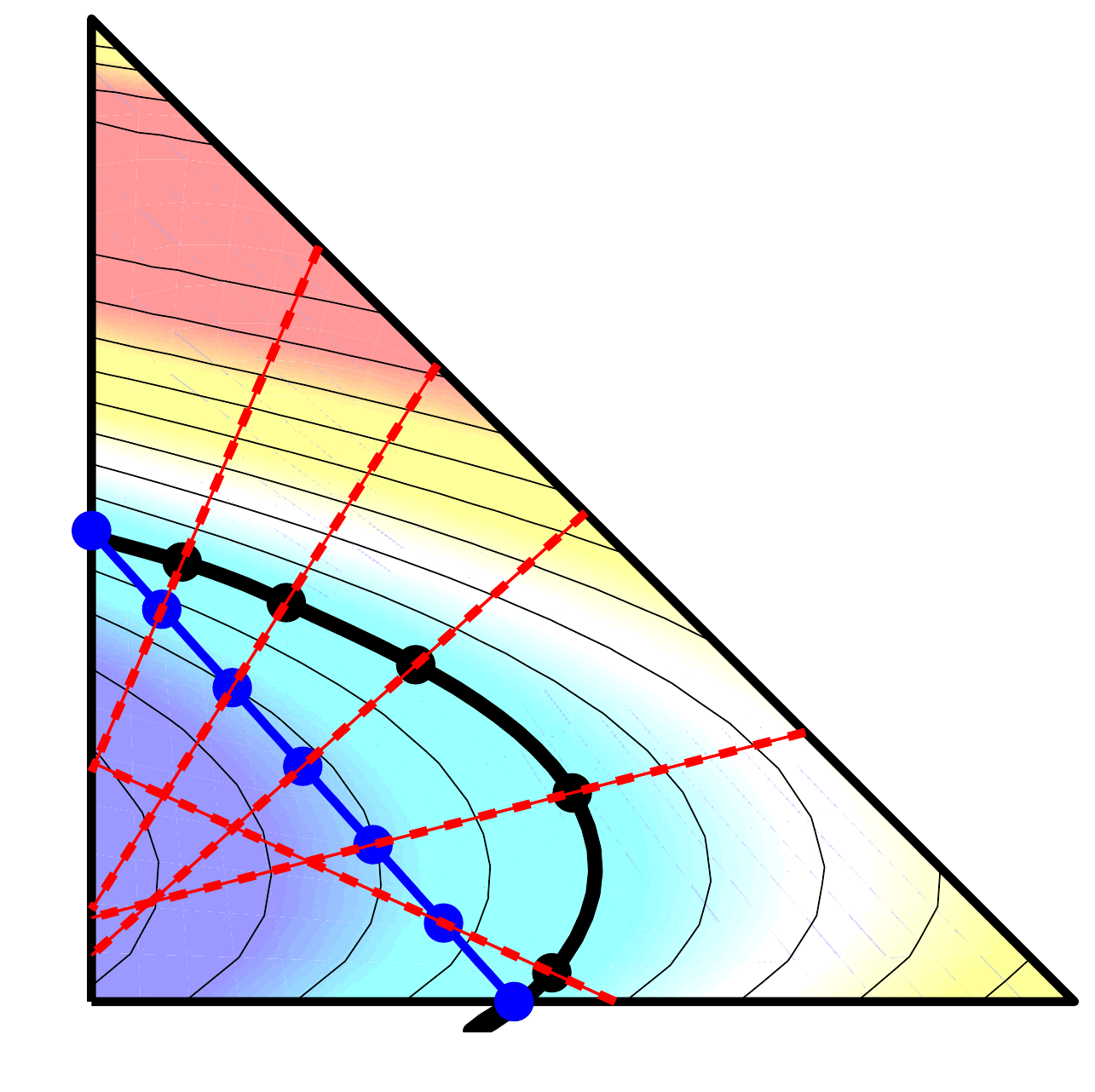}}\quad\subfigure[$\vek N=\nabla\phi\left(\vek r^{i}\right)$]{\includegraphics[width=3.7cm]{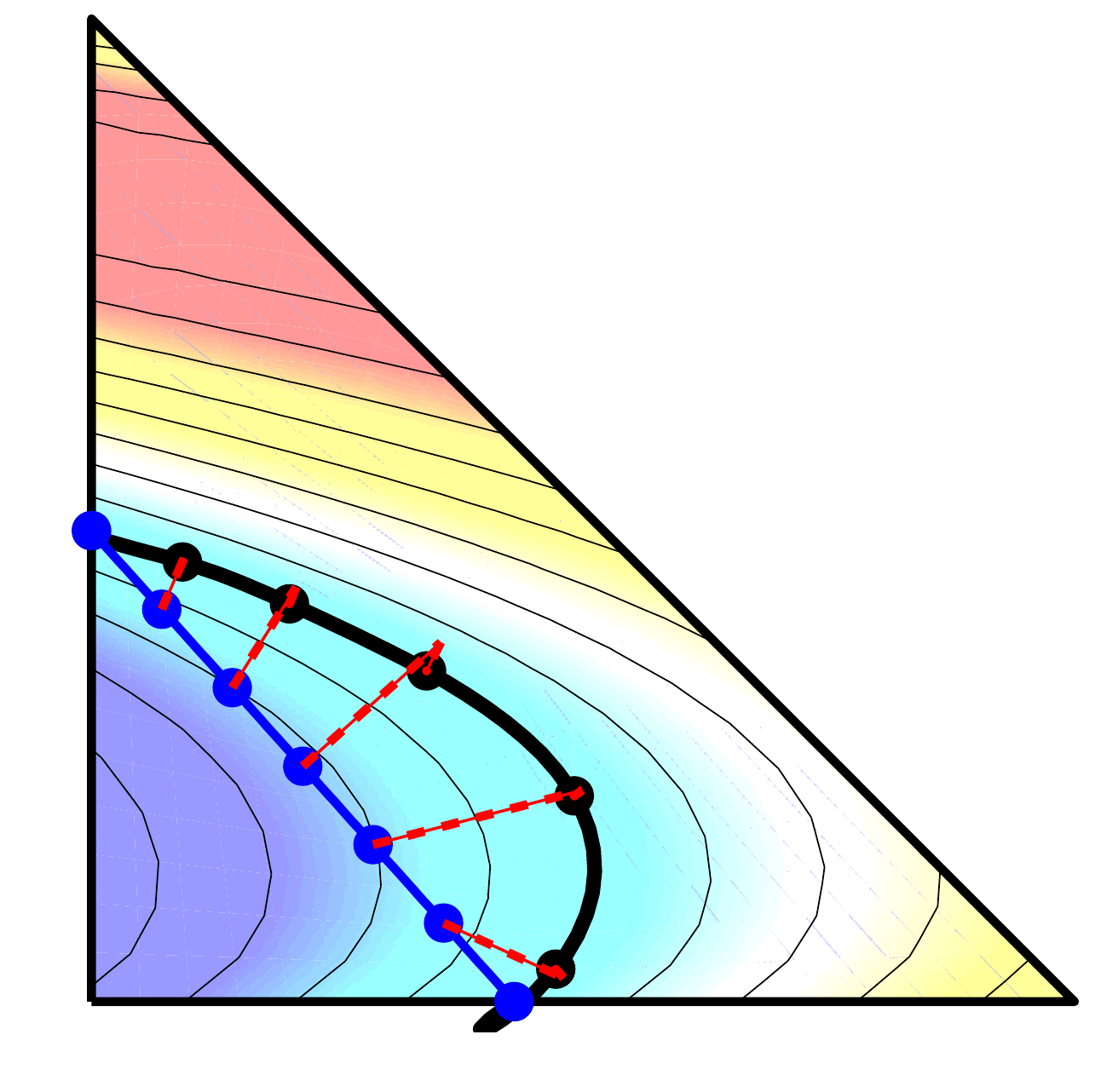}}

\caption{\label{fig:VisSearchDirQuad}(a) and (b) show the iterative method
along fixed search paths, either normal to the linear reconstruction
or in the direction of $\nabla\phi\left(\vek r^{0}\right)$, (c) is
a free search in direction of $\nabla\phi\left(\vek r^{i}\right)$.}
\end{figure}

An iterative procedure is required to identify positions on the zero-level
set from the starting points. One approach is to move strictly along
the straight search paths as defined above, see Fig.~\ref{fig:VisSearchDirQuad}(a)
and (b). The other approach rather uses the gradient of $\phi^{h}\left(\vek r\right)$
at each intermediate position during the iteration, see Fig.~\ref{fig:VisSearchDirQuad}(c).
All approaches yield quadratic convergence rates of the iterative
procedure but, of course, to (slightly) different positions of the
element nodes on the zero-level set.

Mathematically, the algorithm is described as follows. The task is
to find the root of $\phi^{h}\left(\vek r\right)$ in the reference
element, i.e.~some position on the zero-level set. The starting point
of the iterative procedure is labeled $\vek r^{0}$. The Newton-Raphson-type
algorithm for all approaches considered here is based on the following
iteration:
\begin{equation}
\vek r^{i+1}=\vek r^{i}-\frac{\phi^{h}\left(\vek r^{i}\right)}{\nabla\phi\left(\vek r^{i}\right)\cdot\vek N}\cdot\vek N.\label{eq:RootSearch}
\end{equation}

$\vek N$ is the search direction. In the case of a straight search
path, this is e.g.~the normal vector of the intermediate reconstruction
or the gradient at each starting point, hence $\vek N=\nabla\phi\left(\vek r^{0}\right)$.
Otherwise, it may also be the current gradient $\vek N=\nabla\phi\left(\vek r^{i}\right)$.

\subsection{Numerical results for reconstructions in 2D and 3D\label{XX_InterfResults}}

In order to distinguish the quality of the different variants to recover
reconstructions, some numerical results are presented resulting from
an integration and interpolation on the zero-level sets. This is also
justified by the fact that many properties of the second step, i.e.~the
decomposition, as outlined in the next section, are immediately inherited
from the reconstruction. If the reconstruction were not optimal, the
resulting integration points in cut elements would also be sub-optimal.

It is emphasized that in all studies presented herein, the respective
orders of the background elements and the reconstructed surface elements
coincide. Later, in Section \ref{sec:Decomposition}, the same holds
for the order of the resulting sub-elements after the decomposition.
As confirmed in \cite{Fries_2015a}, using different orders would
not be useful as the overall convergence rates would be determined
by the lowest order.

\subsubsection{Numerical results for zero-isolines in 2D\label{XXX_InterfResults2d}}

\begin{figure}
\centering

\subfigure[Setup]{\includegraphics[width=5cm]{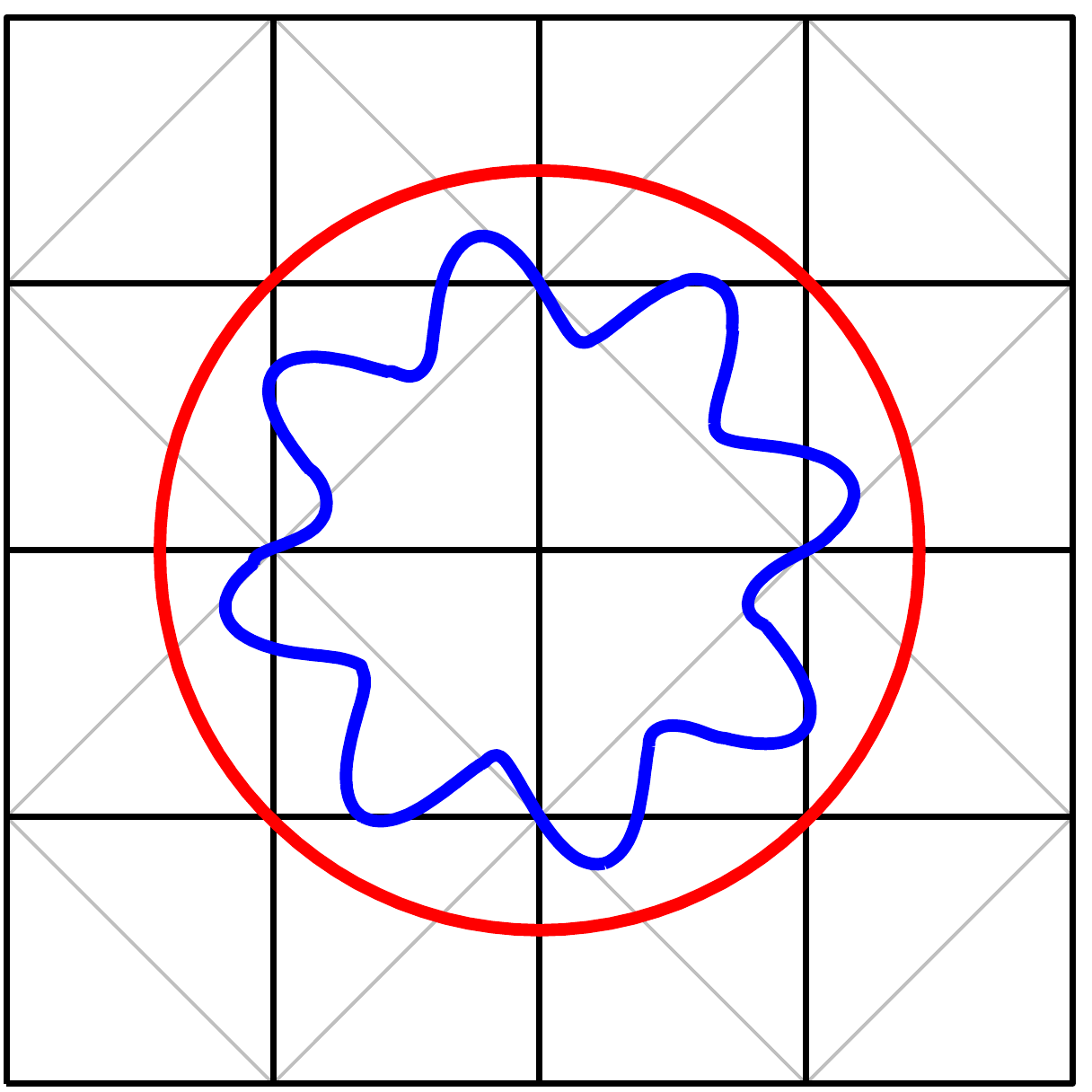}}

\caption{\label{fig:VisRes2dDomain}The background mesh in 2D and the zero-level
sets of $\phi_{1}$ and $\phi_{2}$.}
\end{figure}

We start with the two-dimensional situation. Following \cite{Fries_2016a},
two different level-set functions are considered:
\begin{eqnarray}
\phi_{1}\left(\vek x\right) & = & \sqrt{x^{2}+y^{2}}-r,\label{eq:ExactLevelSet2dA}\\
\phi_{2}\left(\vek x\right) & = & \sqrt{x^{2}+y^{2}}-R\left(\theta\right),\label{eq:ExactLevelSet2dB}
\end{eqnarray}
with $r=0.7123$, $R\left(\theta\right)=0.5+0.1\cdot\sin\left(8\theta\right)$
and $\theta\left(\vek x\right)=\mathrm{atan}\left(y/x\right)$. See
Fig.~\ref{fig:VisRes2dDomain}(a) for a visualization of the zero-level
sets of $\phi_{1}$ and $\phi_{2}$. The zero-level set of $\phi_{1}$
is a circle with radius $r$ and is frequently used in the literature,
e.g.~\cite{Fries_2015a}, and $\phi_{2}$ is similar to \cite{Li_1998a}.
The function $f\left(\vek x\right)$ is defined as
\begin{equation}
f\left(\vek x\right)=\nicefrac{1}{2}\cdot x+\nicefrac{1}{4}\cdot y+x^{2}+2\cdot y^{3}.\label{eq:ExactFctForInt2D}
\end{equation}
and is later integrated on the zero-level sets.

In the convergence studies, $\left\{ 6,10,20,30,50,70,100,150,200,300\right\} $
elements are employed per dimension. We use standard Gauss rules for
the integration on the zero-level set $\Gamma_{0}$ with a rather
high order of $11$ which is kept constant independently of the order
of background and interface elements. Results in 2D are studied in
different error norms which are related to integrating and interpolating
a function $f\left(\vek x\right)$ on $\Gamma_{0}$ \cite{Fries_2016a}:
\begin{eqnarray}
\varepsilon_{1}^{\Gamma} & = & \Big|\Big(\sum_{i}w_{i}\Big)-I_{1}^{\Gamma}\Big|/\left|I_{1}^{\Gamma}\right|\quad\mathrm{with}\quad I_{1}^{\Gamma}=\int_{\Gamma_{0}}1\,\mathrm{d}\Gamma,\label{eq:Error2dInterf1}\\
\varepsilon_{\phi}^{\Gamma} & = & \Big(\sum_{i}w_{i}\cdot\phi\left(\vek x_{i}\right)\Big),\\
\varepsilon_{f}^{\Gamma} & = & \Big|\Big(\sum_{i}w_{i}\cdot f\left(\vek x_{i}\right)\Big)-I_{f}^{\Gamma}\Big|/\left|I_{f}^{\Gamma}\right|\quad\mathrm{with}\quad I_{f}^{\Gamma}=\int_{\Gamma_{0}}f\left(\vek x\right)\,\mathrm{d}\Gamma,\label{eq:Error2dInterf2}\\
\varepsilon_{f_{1D}^{h}}^{\Gamma} & = & \Big|\Big(\sum_{i}w_{i}\cdot f_{1D}^{h}\left(\vek x_{i}\right)\Big)-I_{f}^{\Gamma}\Big|/\left|I_{f}^{\Gamma}\right|,\label{eq:Error2dInterf3}\\
\varepsilon_{f_{2D}^{h}}^{\Gamma} & = & \Big|\Big(\sum_{i}w_{i}\cdot f_{2D}^{h}\left(\vek x_{i}\right)\Big)-I_{f}^{\Gamma}\Big|/\left|I_{f}^{\Gamma}\right|,\label{eq:Error2dInterf4}
\end{eqnarray}
where $\vek x_{i}$ are the integration points in the interface elements
with integration weights $w_{i}$. It is noted that, in contrast to
the \emph{nodes} of the interface elements, the integration points
at $\vek x_{i}$ are only approximately on the zero-level set of $\phi^{h}\left(\vek x\right)$. 

The first error, $\varepsilon_{1}^{\Gamma}$ in Eq.~(\ref{eq:Error2dInterf1}),
is evaluated by summing up integration weights giving the length of
the reconstructed interface which is compared with the exact length
(for $\phi_{1}$, this is the circumference of a circle). The second
error, $\varepsilon_{\phi}^{\Gamma}$, integrates the level-set function
on the zero-level set, which, in the ideal case, would be zero. Nevertheless,
because the nodes of the interface elements are on the zero-level
set of $\phi^{h}$ rather than $\phi$ and, furthermore, the integration
points are only approximately on the zero-level set of $\phi^{h}$,
$\phi\left(\vek x_{i}\right)$ is not exactly zero. The other error
norms integrate either the exact or the interpolated function $f\left(\vek x\right)$
as given in Eq.~(\ref{eq:ExactFctForInt2D}). For $\varepsilon_{f}^{\Gamma}$,
the function $f\left(\vek x\right)$ is evaluated exactly at the integration
points $\vek x_{i}$. Instead, $\varepsilon_{f_{1D}^{h}}^{\Gamma}$
evaluates the function at the element nodes of the \emph{interface}
elements and interpolates at the integration points using the shape
functions of the interface element. Equivalently, $\varepsilon_{f_{2D}^{h}}^{\Gamma}$
uses the element nodes of the \emph{background} elements to evaluate
the function and the corresponding shape functions to interpolate.

\begin{figure}
\centering

\subfigure[]{\includegraphics[width=6cm]{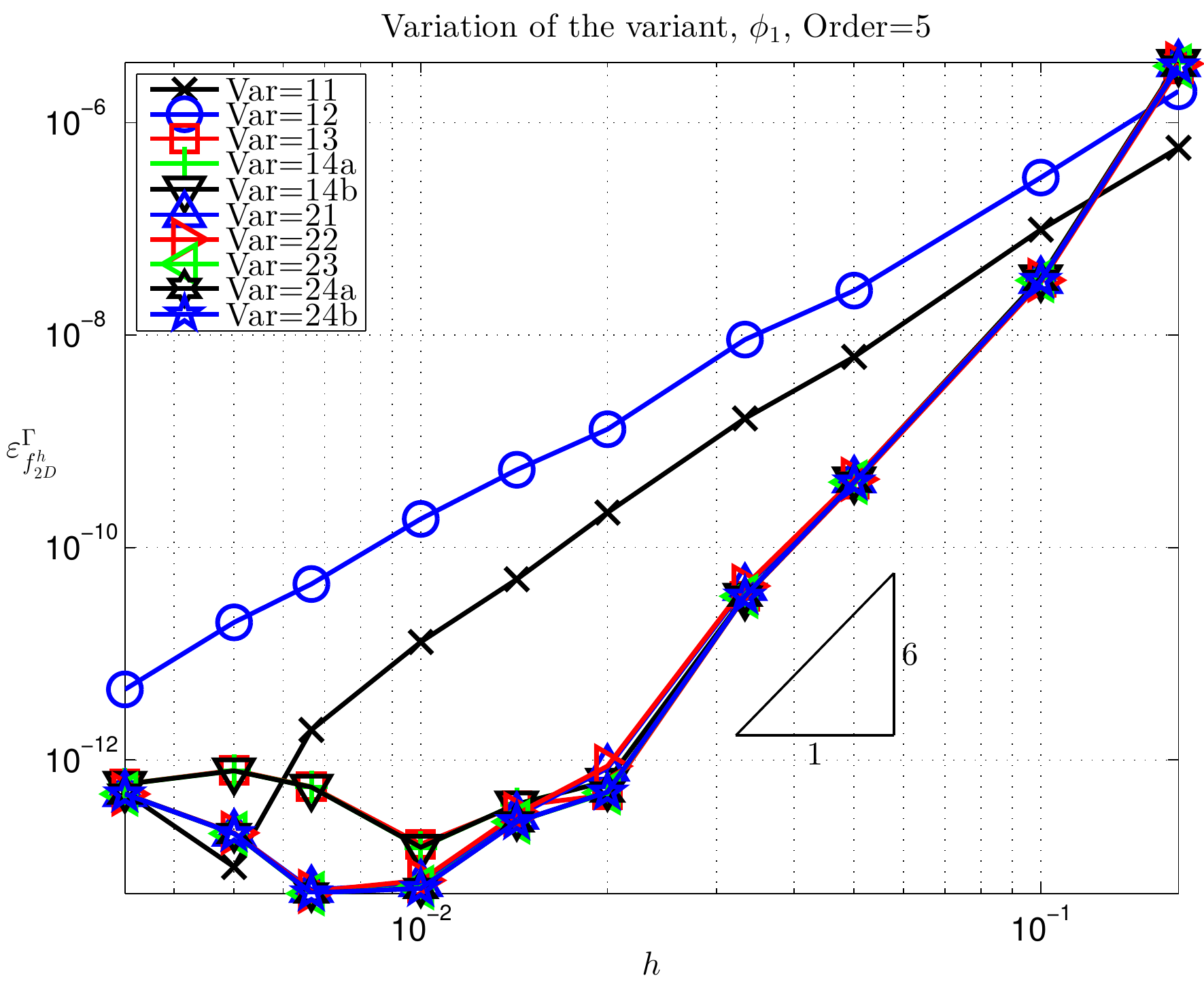}}\quad\subfigure[]{\includegraphics[width=6cm]{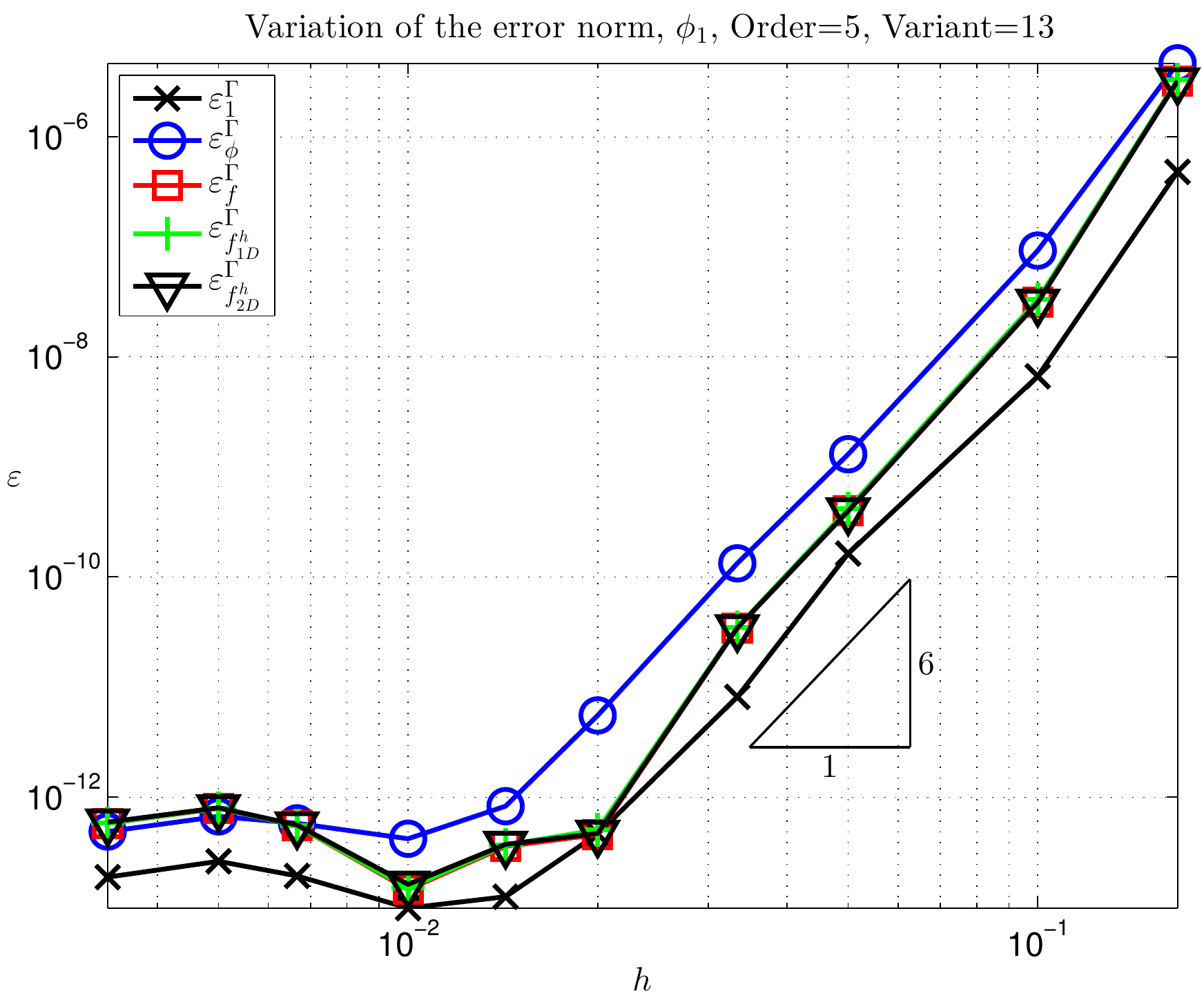}}\\\subfigure[]{\includegraphics[width=6cm]{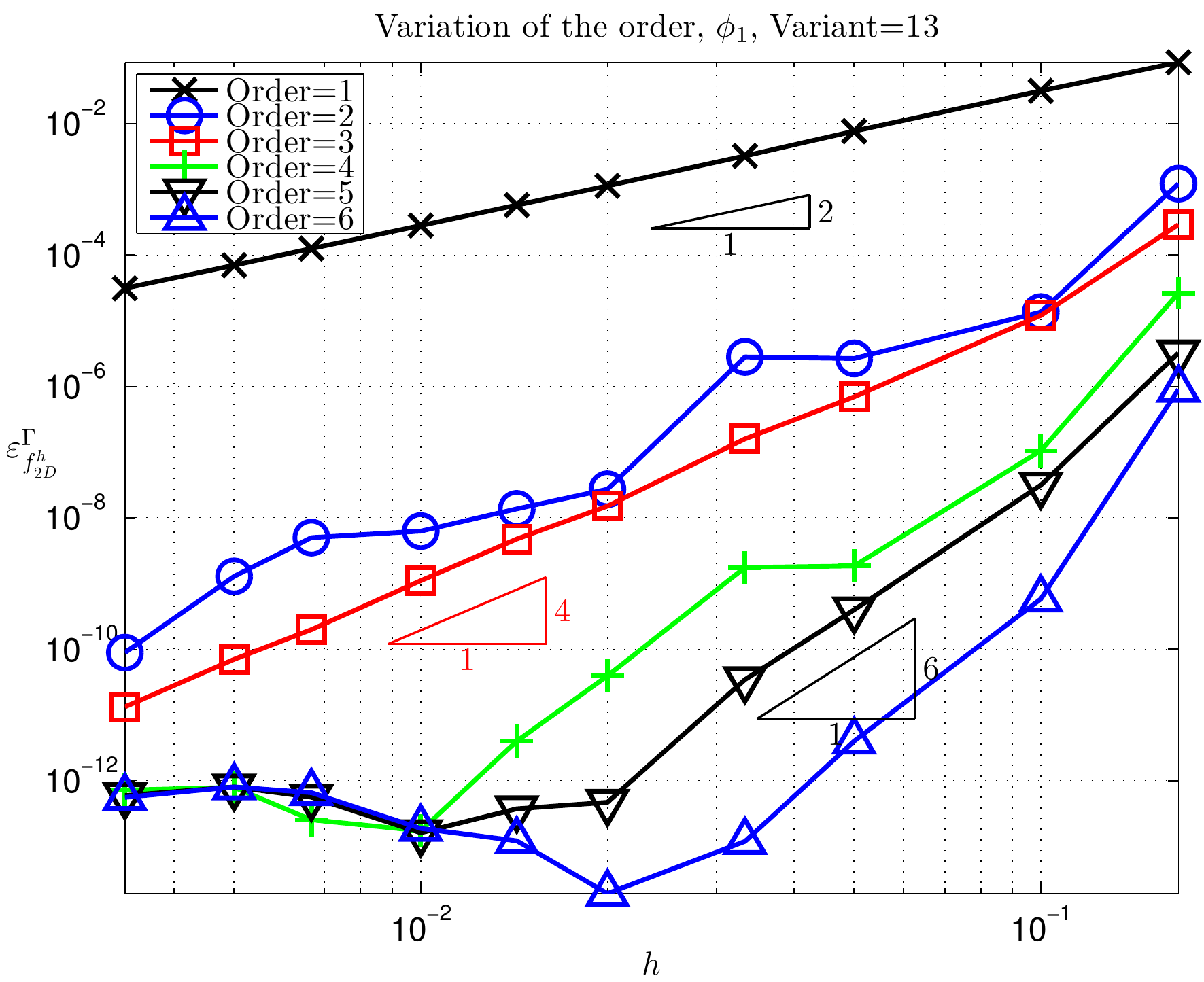}}\quad\subfigure[]{\includegraphics[width=6cm]{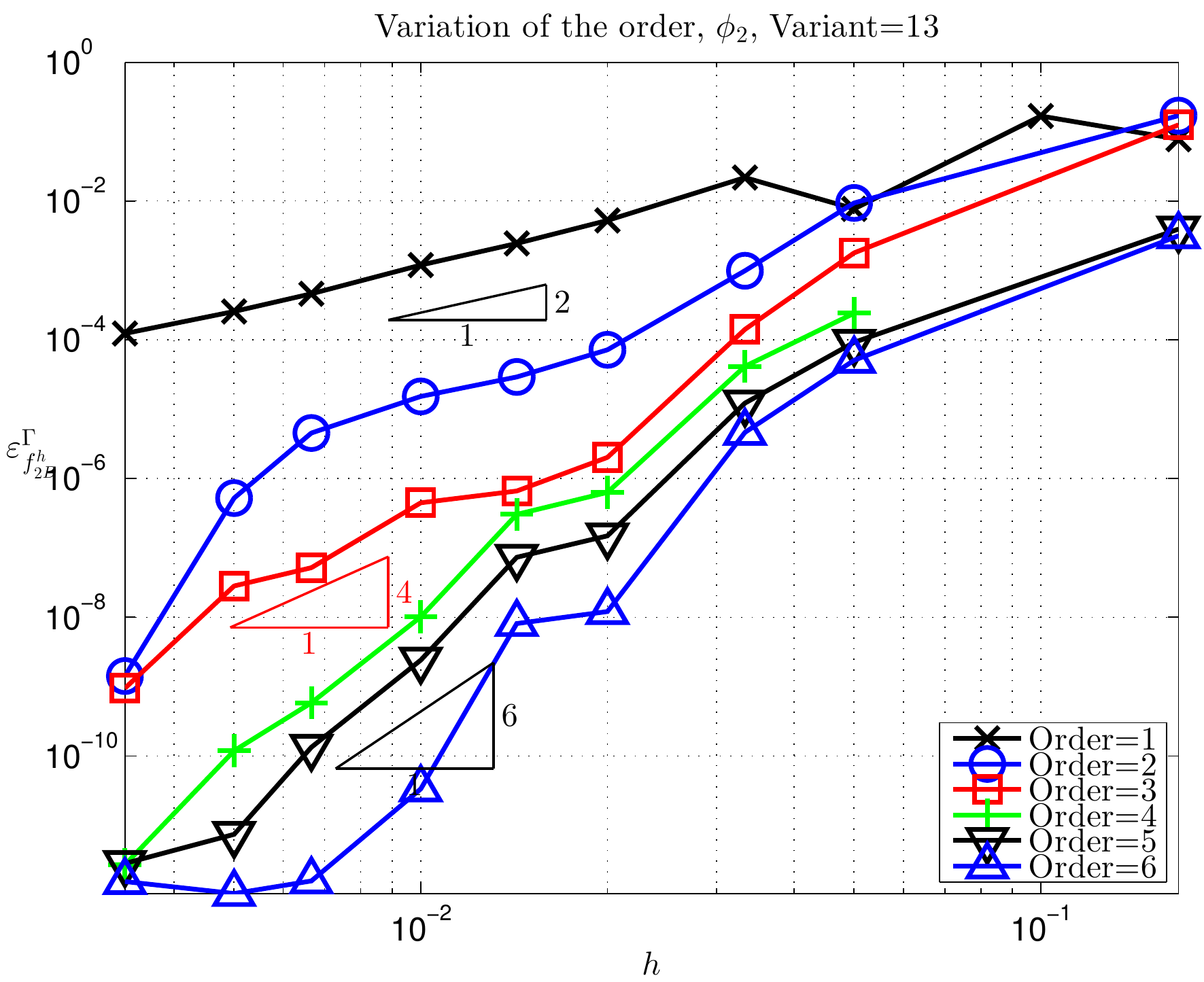}}

\caption{\label{fig:VisRes2dInterf}Convergence results in 2D for integrating
and interpolating on zero-isolines.}
\end{figure}

Results are shown in Fig.~\ref{fig:VisRes2dInterf}. In Fig.~\ref{fig:VisRes2dInterf}(a),
it is seen that most of the variants for the search algorithm (start
values and search directions) perform similar and with optimal convergence.
Variants $14$ and $24$, where the search directions are based on
the gradient of the level-set function are split into two sub-versions
14a/14b and 24a/24b depending on whether the search path is constant,
$\nabla\phi^{h}\left(\vek r^{0}\right)$, or changes during the iteration,
$\nabla\phi^{h}\left(\vek r^{i}\right)$. A few variants actually
underperform, in particular variant $11$ and $12$. In these variants,
the search directions depend on the topological situation only, i.e.~the
search directions are towards the one node on the other side or the
interpolated edge directions but are, otherwise, independent of the
level-set data. It is seen that there is only a negligible difference
between the variants which are based on the linear reconstruction
(first digit in the variant number is $1$) compared to the Hermite
reconstruction (first digit is $2$). Because it is found that this
also applies to the three-dimensional situation, the Hermite reconstruction
is neglected in the following. Of course, Fig.~\ref{fig:VisRes2dInterf}(a)
is just one sample of a large parameter space that was actually studied
where, as an example, $\phi_{1}$ is chosen and an order of $5$ for
the background mesh and the reconstruction.

Fig.~\ref{fig:VisRes2dInterf}(b) shows that the error converges
independently of the error norms in Eqs.~(\ref{eq:Error2dInterf1})
to (\ref{eq:Error2dInterf4}). This was again confirmed for the whole
parameter space that has been studied for the two level-set functions.
Now that it is known which variants perform potentially optimal independently
on the error norm, the full convergence results for different orders
are shown in Fig.~\ref{fig:VisRes2dInterf}(c) and (d) for $\phi_{1}$
and $\phi_{2}$, respectively. The figures look very similar for all
optimal variants. It is confirmed that optimal convergence rates are
achieved. It is noted that only data points are displayed where no
recursive refinement was necessary (because this influences the element
lengths $h$ and complicates the proper placement of the result in
the convergence plots).

\subsubsection{Numerical results for zero-isosurfaces in 3D\label{XXX_InterfResults3d}}

Numerical results are now presented for the integration on zero-isosurfaces
$\Gamma_{0}$ in 3D. The procedure is similar to the studies in two
dimensions. It is recalled that in 3D, inner and outer nodes of each
interface element are treated differently. For the outer nodes, any
of the 2D variants may be used, see Section \ref{sub:StartValuesIn2d}.
Here we restrict ourselves to variant $13$ but all other variants
have been systematically studied as well. For the inner nodes, three
variants are studied, see \ref{sub:StartValuesIn3d}: Variant A and
B use straight search paths, the first based on the normal vector
of the intermediate reconstruction at the start value $\vek r^{0}$
and the second based on the gradient of the level-set function at
$\vek r^{0}$. Variant C is based on the gradient as well, but depending
on the points $\vek r^{i}$ during the iterative procedure.

\begin{figure}
\centering

\subfigure[Zero-level set of $\phi_1$]{\includegraphics[height=5cm]{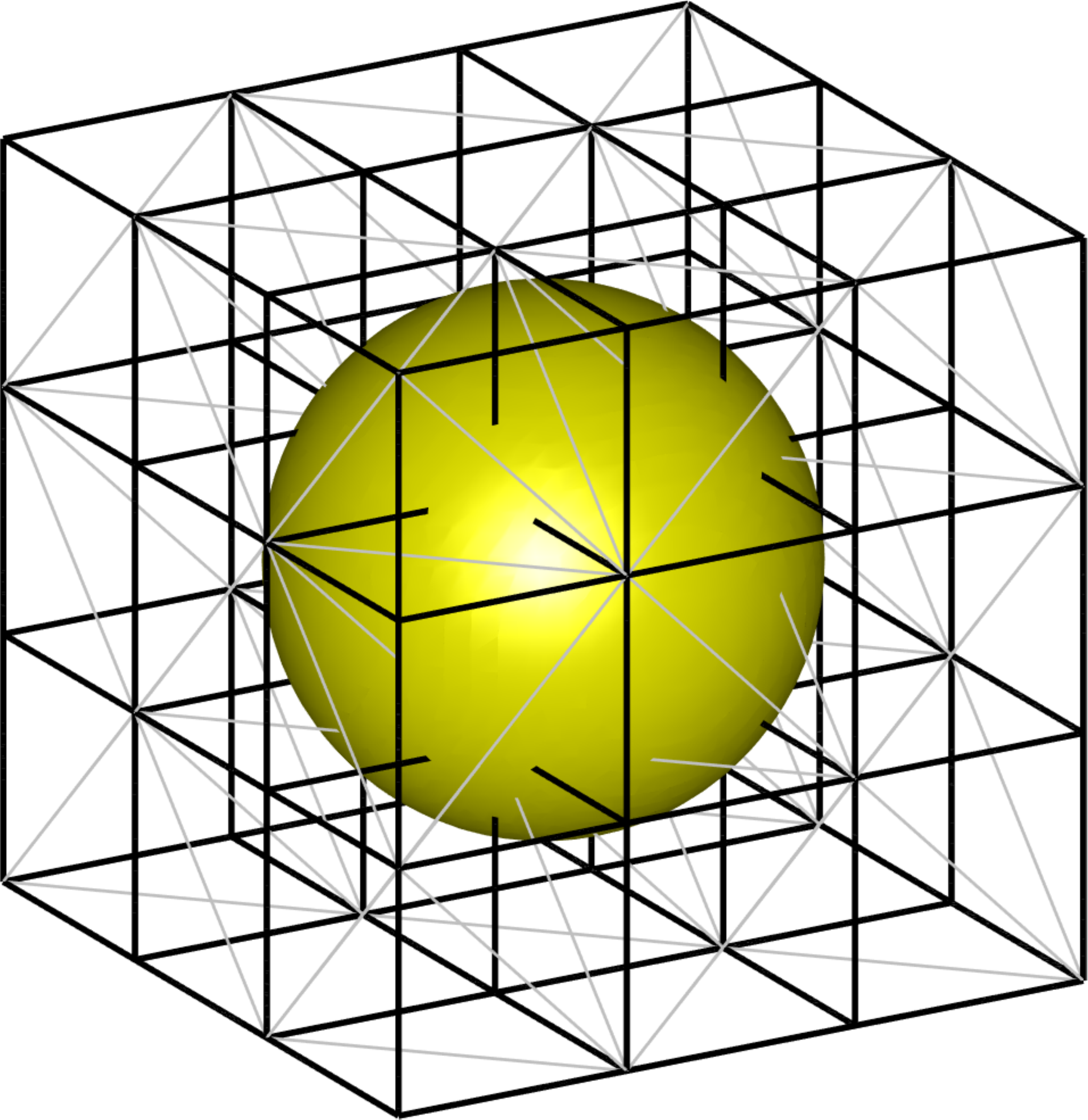}}\qquad\subfigure[Zero-level set of $\phi_2$]{\includegraphics[height=5cm]{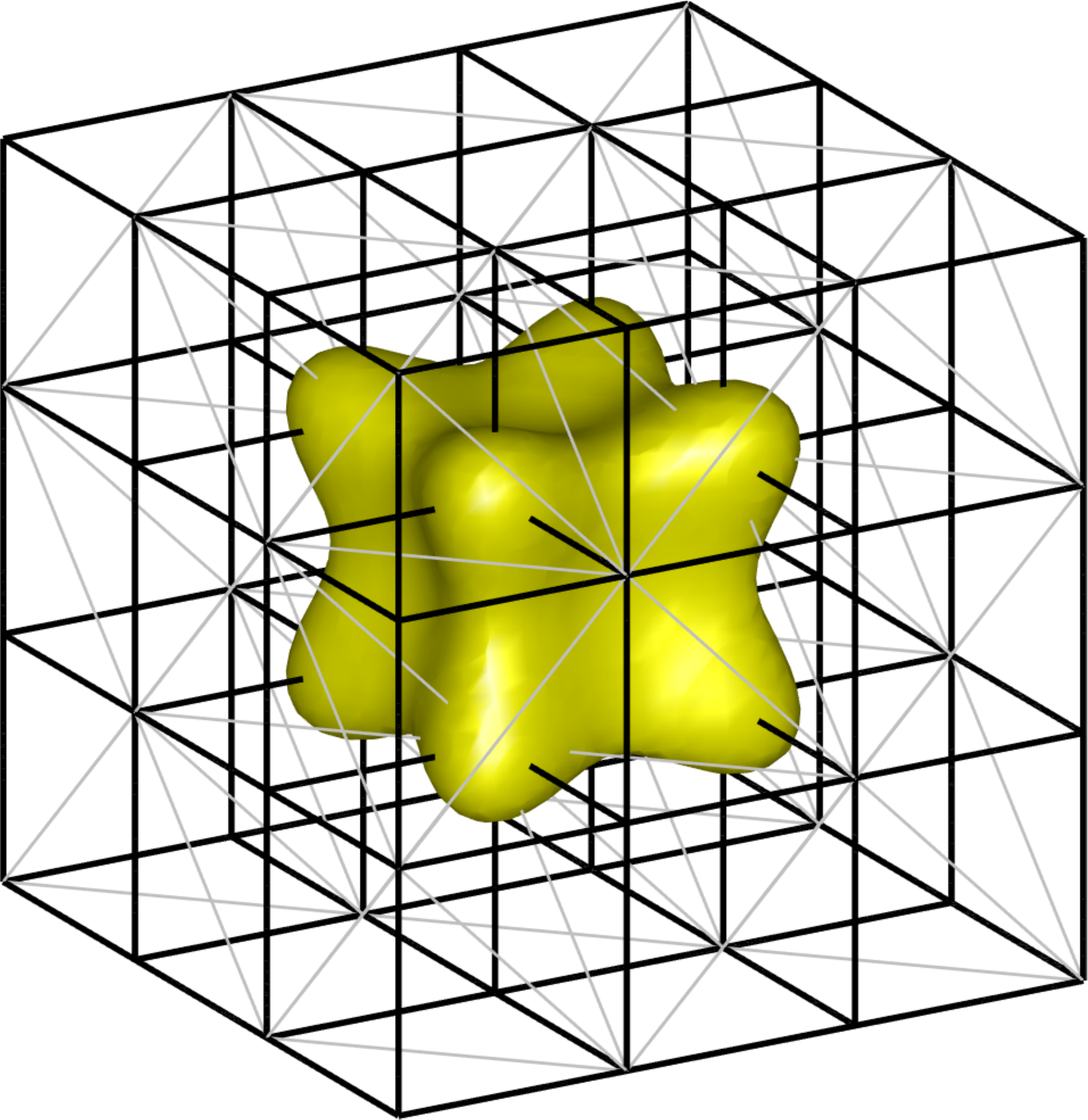}}

\caption{\label{fig:VisRes3dInterf}The background mesh in 3D and the zero-level
sets of (a)$\phi_{1}$ and (b) $\phi_{2}$.}
\end{figure}

The following two level-set functions are considered in 3D:
\begin{eqnarray}
\phi_{1}\left(\vek x\right) & = & \sqrt{x^{2}+y^{2}+z^{2}}-r,\label{eq:ExactLevelSet3dA}\\
\phi_{2}\left(\vek x\right) & = & \phi_{1}\left(\vek x\right)+0.1\cdot\left[\cos\left(2\pi\cdot x\right)+\cos\left(2\pi\cdot y\right)+\cos\left(2\pi\cdot z\right)\right],\label{eq:ExactLevelSet3dB}
\end{eqnarray}
with $r=0.7123$. The zero-isosurfaces of these two level-set functions
are visualized in Fig.~\ref{fig:VisRes3dInterf}. There, also a coarse
tetrahedral background mesh is shown. For the integrand, the function
\begin{equation}
f\left(\vek x\right)=x^{2}+y^{2}+\nicefrac{1}{2}\cdot\cos(z)\label{eq:ExactFctForInt3D}
\end{equation}
is chosen. In the convergence studies, $\left\{ 6,10,14,20,30,50,70,100\right\} $
elements are used per dimension. 

Again, the error is studied in different norms. The three norms in
Eqs.~(\ref{eq:Error2dInterf1}) to (\ref{eq:Error2dInterf2}) are
directly applicable in this 3D study as well, the other two are slightly
modified as
\begin{eqnarray*}
\varepsilon_{f_{2D}^{h}}^{\Gamma} & = & \Big|\Big(\sum_{i}w_{i}\cdot f_{2D}^{h}\left(\vek x_{i}\right)\Big)-I_{f}^{\Gamma}\Big|/\left|I_{f}^{\Gamma}\right|\quad\mathrm{with}\quad I_{f}^{\Gamma}=\int_{\Gamma_{0}}f\left(\vek x\right)\,\mathrm{d}\Gamma,\\
\varepsilon_{f_{3D}^{h}}^{\Gamma} & = & \Big|\Big(\sum_{i}w_{i}\cdot f_{3D}^{h}\left(\vek x_{i}\right)\Big)-I_{f}^{\Gamma}\Big|/\left|I_{f}^{\Gamma}\right|.
\end{eqnarray*}
Therein, $f_{2D}^{h}$ is the integrand interpolated by the shape
functions of the reconstructed interface elements and $f_{3D}^{h}$
is interpolated by the shape functions of the background elements. 

\begin{figure}
\centering

\subfigure[]{\includegraphics[width=6cm]{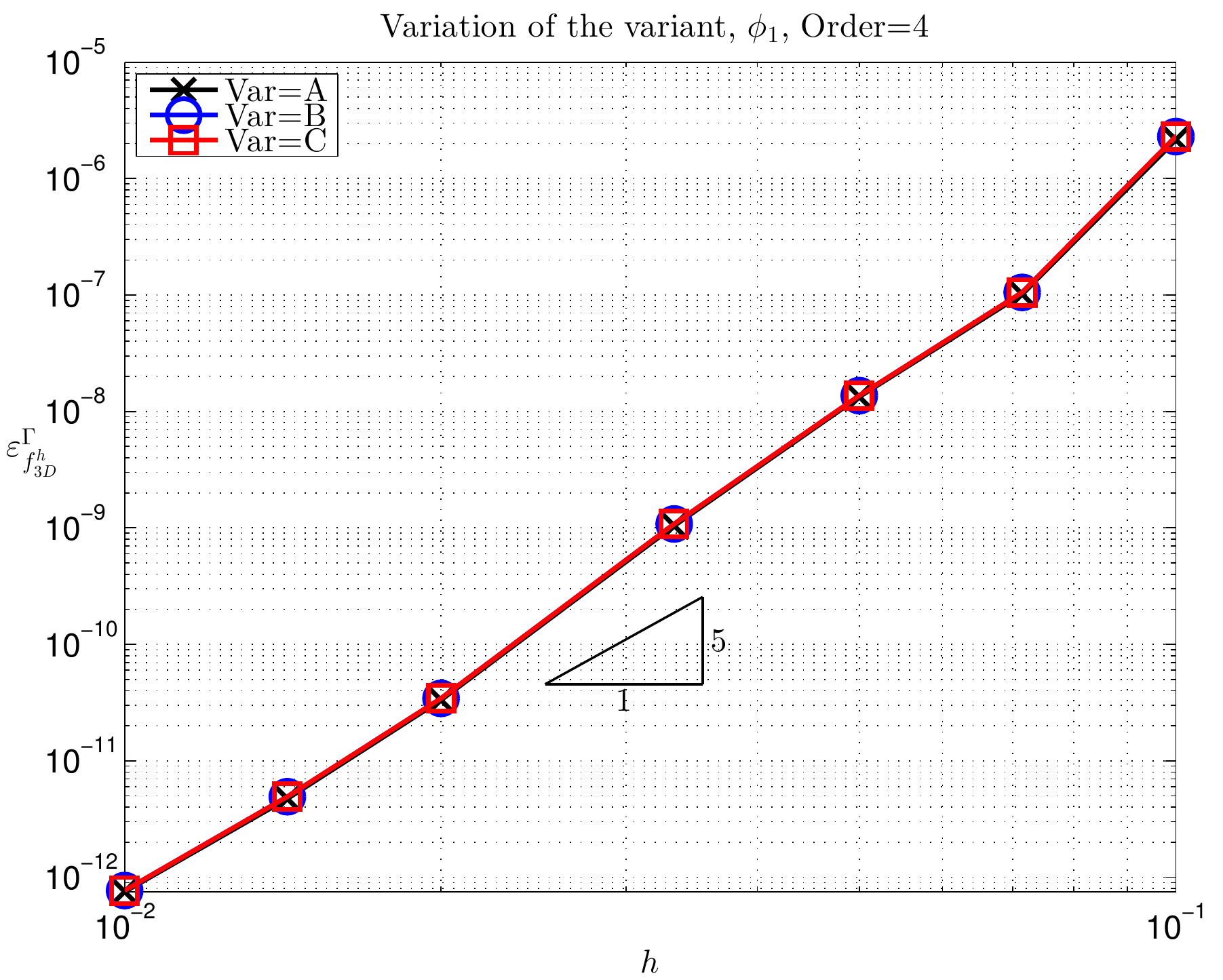}}\quad\subfigure[]{\includegraphics[width=6cm]{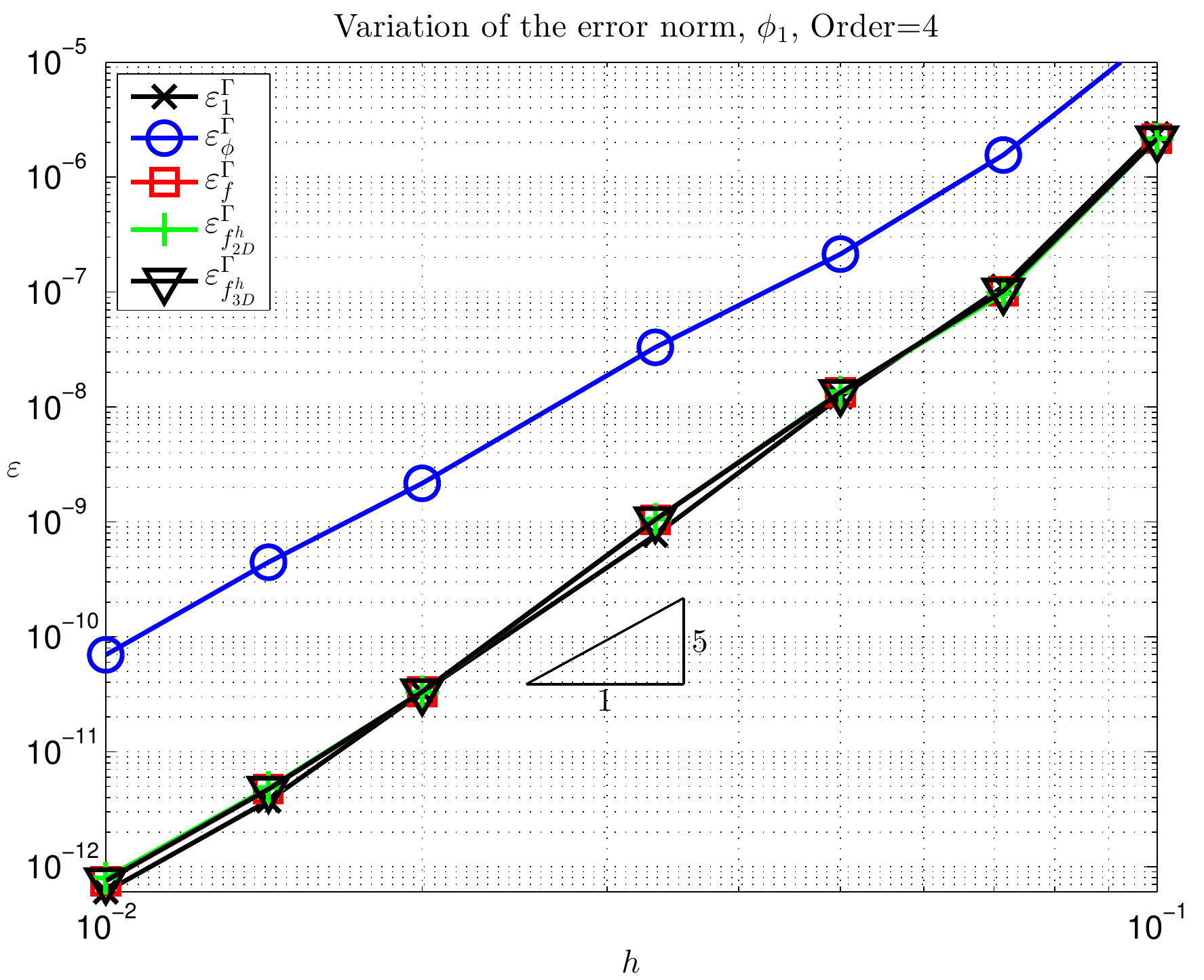}}\\\subfigure[]{\includegraphics[width=6cm]{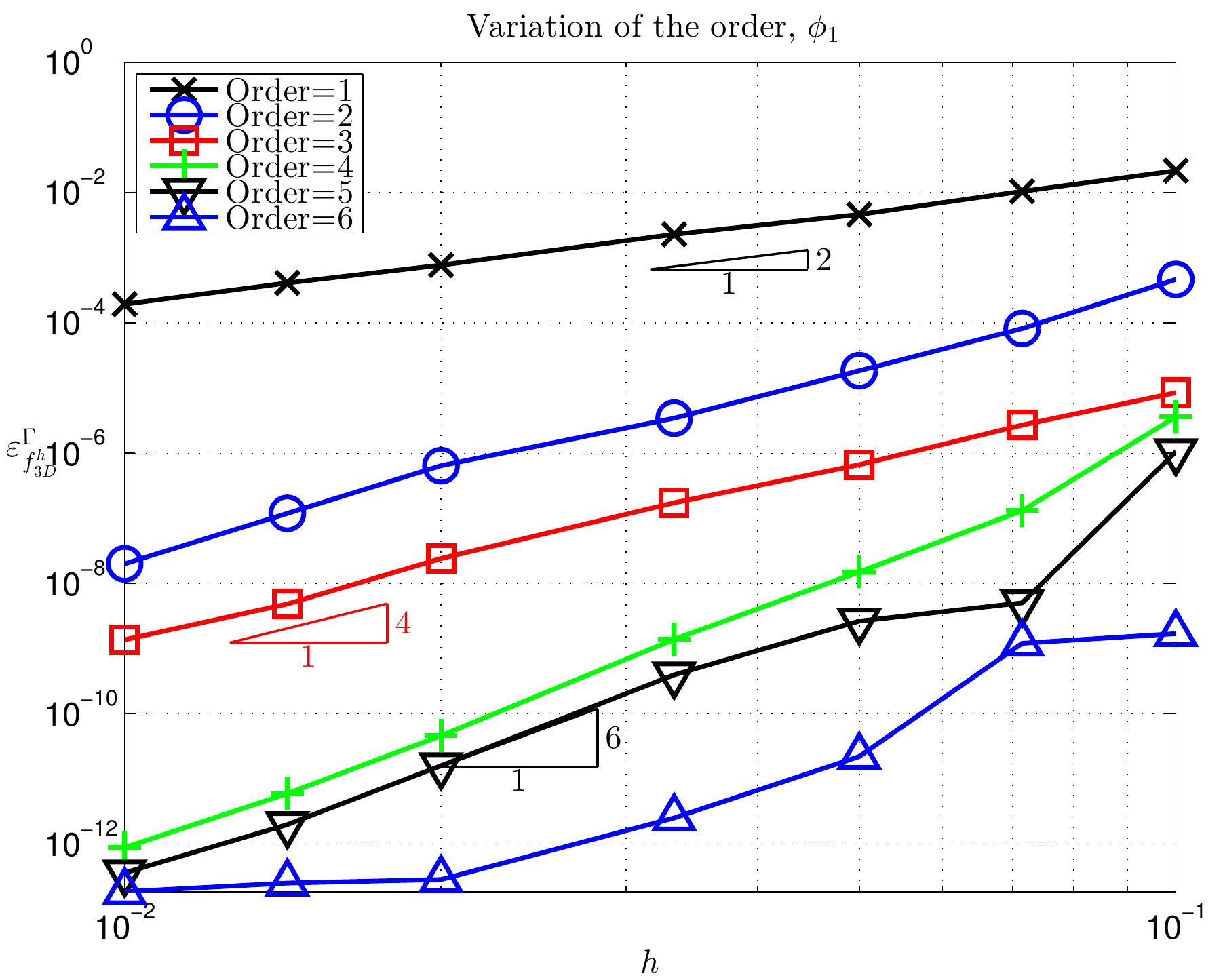}}\quad\subfigure[]{\includegraphics[width=6cm]{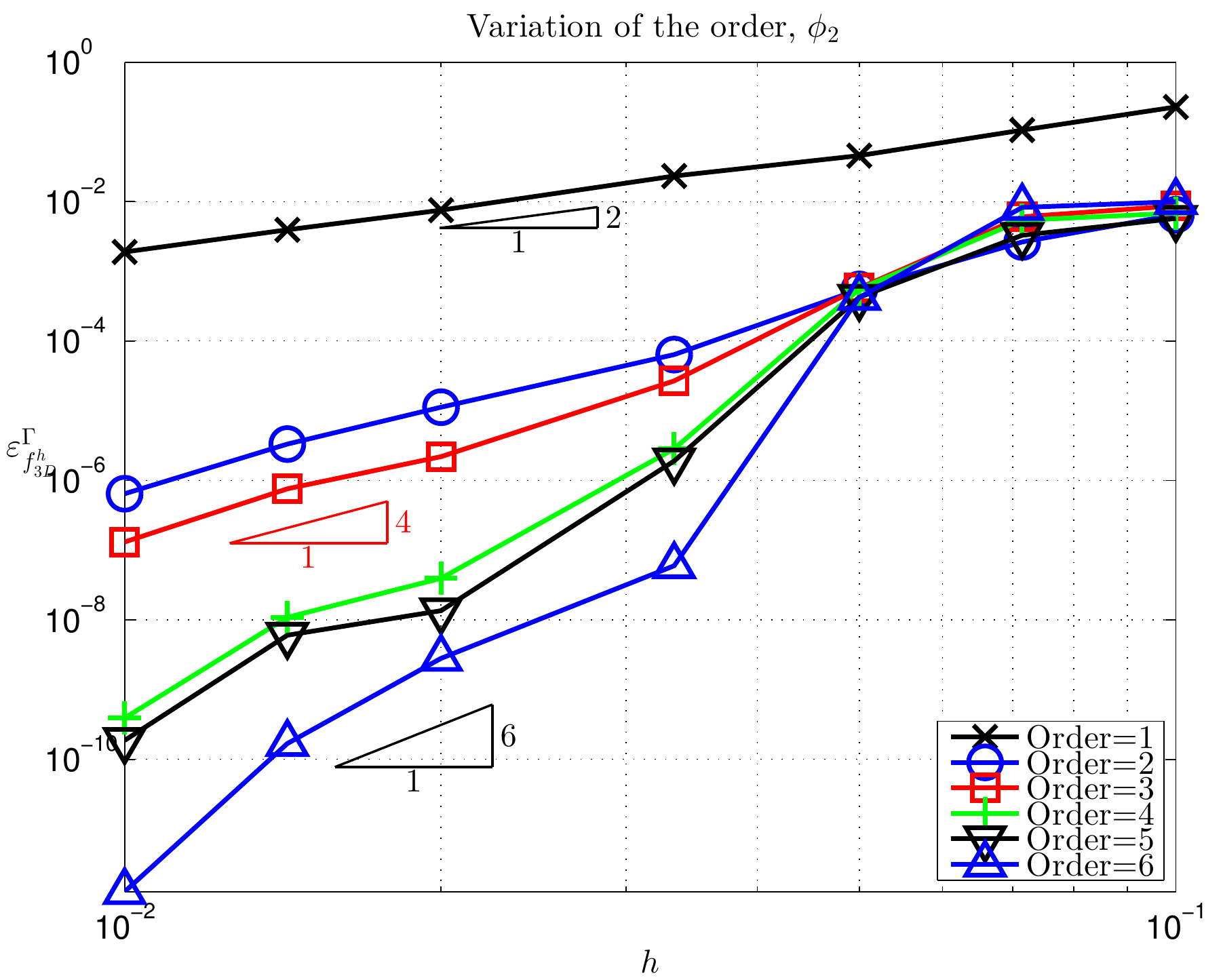}}

\caption{\label{fig:VisRes3dInterfConv}Convergence results in 3D for integrating
and interpolating on zero-isosurfaces.}
\end{figure}

Results are shown in Fig.~\ref{fig:VisRes3dInterfConv}. Fig.~\ref{fig:VisRes3dInterfConv}(a)
shows that all three variants for the inner nodes lead to virtually
identical results. Therefore, we concentrate on the variant A using
the normal vector for the root search. This is consistent to variant
13 for the 2D case, here applied on each of the faces of the tetrahedron.
Fig.~\ref{fig:VisRes3dInterfConv}(b) confirms that the same convergence
rates are expected in all of the five error norms. Finally, Fig.~\ref{fig:VisRes3dInterfConv}(c)
and (d) vary the order and resolution of the background mesh to show
that optimal convergence rates are obtained for integrations on the
zero-isosurfaces of $\phi_{1}$ and $\phi_{2}$.

Hence, optimal convergence rates are found for the reconstruction
regardless of the fact that a number of errors are involved: The reconstruction
is only an approximation of the zero-level set of the interpolated
level-set $\phi^{h}$ which, in fact, is only an approximation of
$\phi$. In some error norms, only the interpolated function $f^{h}$
is integrated rather than $f$. Finally, also the numerical integration
based on Gauss quadrature introduces an integration error. Nevertheless,
it was shown that $f\left(\vek x\right)$ is integrated and interpolated
optimally on zero-level sets implied by $\phi\left(\vek x\right)$.

\section{Decomposition\label{sec:Decomposition}}

Once the interface elements are defined, i.e.~the zero-level set
is meshed with higher-order interface elements, the cut background
elements are to be decomposed into sub-elements depending on the topology
of the cut situation, see Section \ref{sub:Topologies}. This is outlined
in Fig.~\ref{fig:VisOverviewRemesh} for the two- and three-dimensional
situation, respectively. The resulting sub-elements share the property
that they feature \emph{one} higher-order side coinciding with the
reconstructed interface element. All other edges are straight. A key
aspect for the decomposition is the mapping of element nodes into
the special sub-elements on the two sides of the interface.

\subsection{Decomposition in 2D\label{XX_Decomp2d}}

\begin{figure}
\centering

\subfigure[Triangle]{\includegraphics[height=4cm]{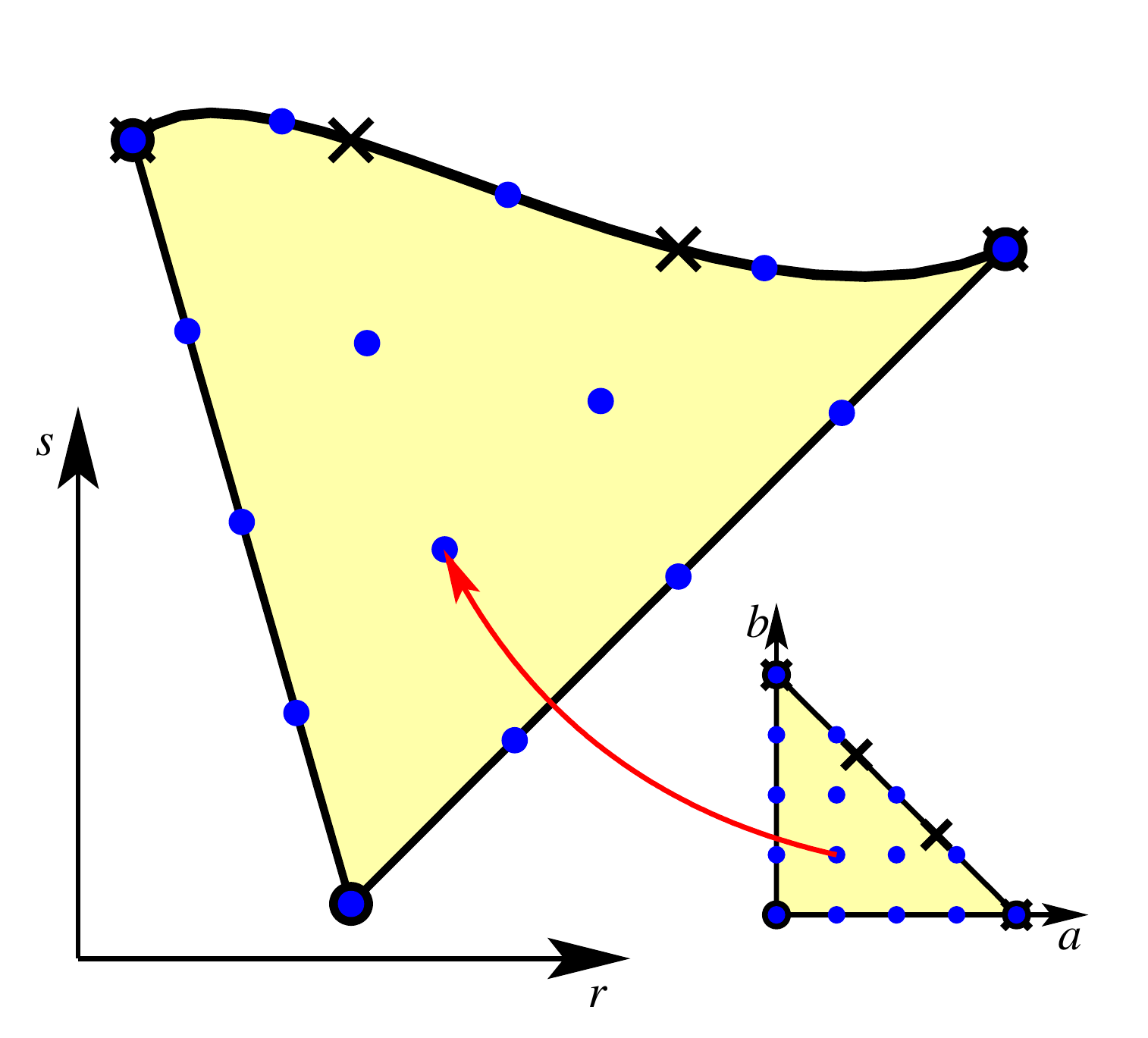}}\quad\subfigure[Quadrilateral]{\includegraphics[height=4cm]{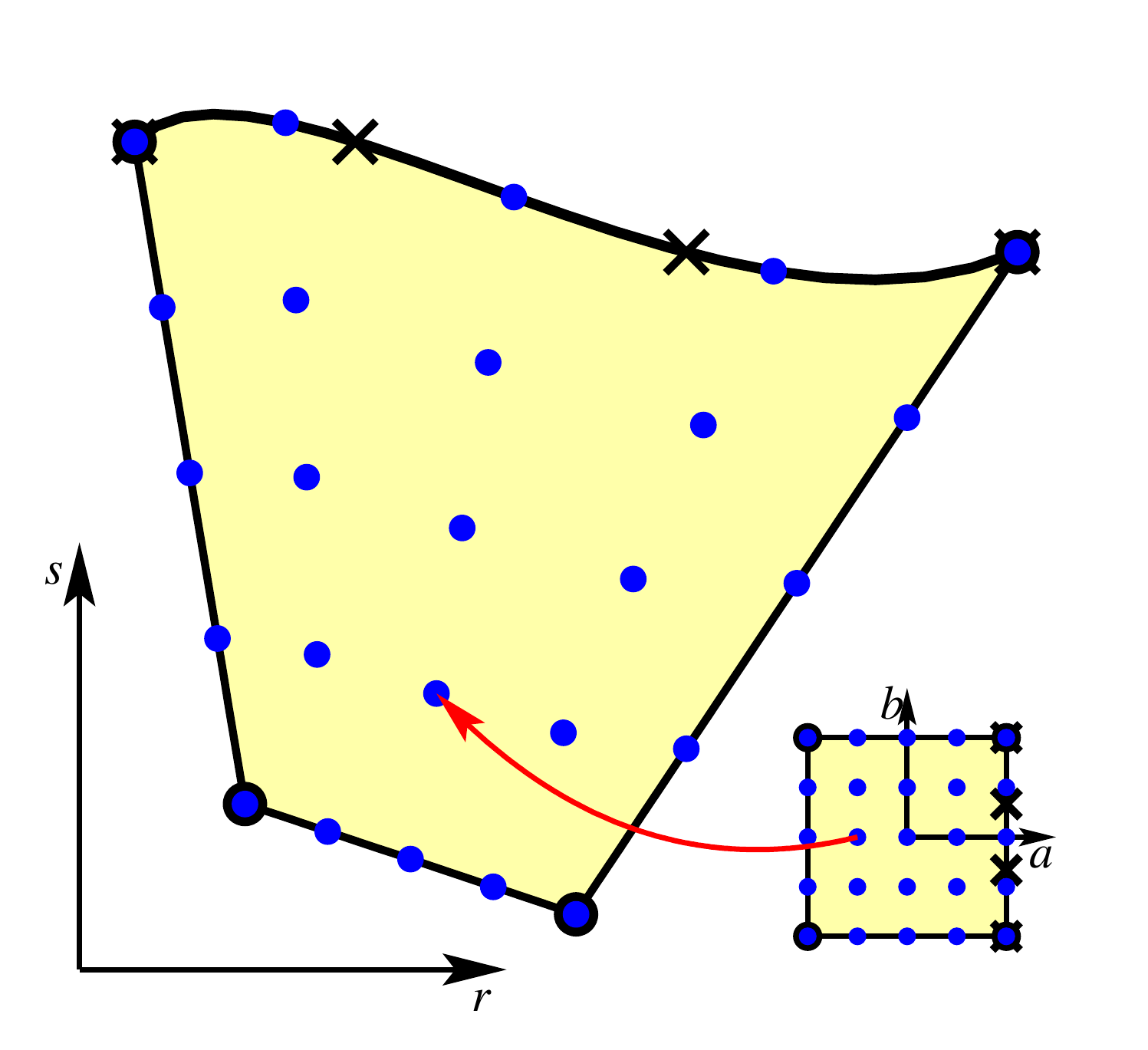}}\quad\subfigure[Different maps]{\includegraphics[height=4cm]{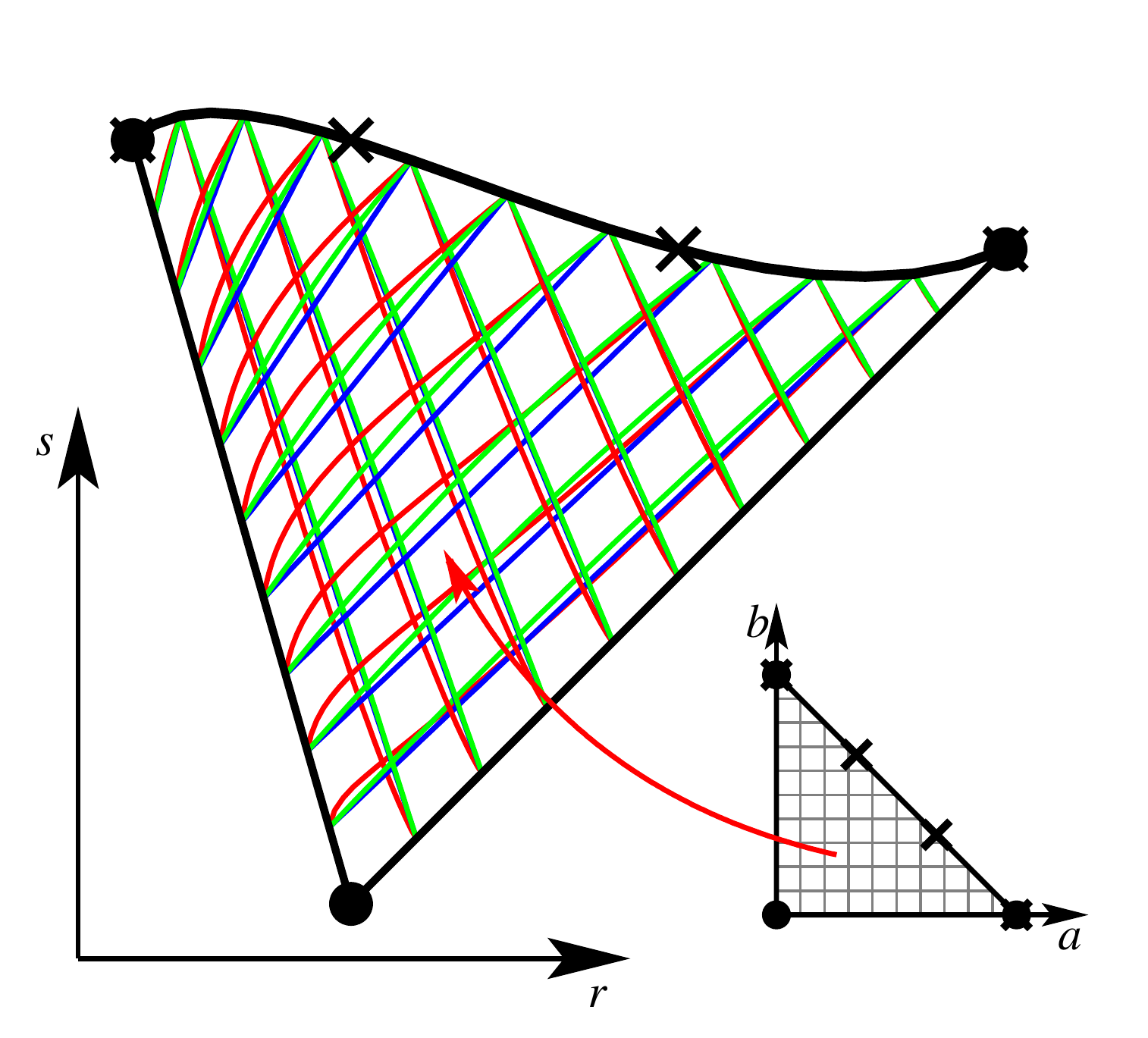}}

\caption{\label{fig:VisMap2dExamples}Examples of resulting higher-order elements
(blue nodes) from the special sub-elements with one higher-order edge
(black crosses), (a) triangular sub-element, (b) quadrilateral sub-element,
(c) some grid-lines of different maps in a triangular element.}
\end{figure}

For the two-dimensional case, a cut reference background element in
the coordinate system $\vek r$ falls into a triangular and a quadrilateral
sub-element with one higher-order side, see Fig.~\ref{fig:VisOverviewRemesh}(a).
In order to define the element nodes of the sub-elements, a mapping
from a higher-order reference triangle or quadrilateral in coordinates
$\vek a$ to each sub-element in coordinates $\vek r$ is sought,
see Fig.~\ref{fig:VisMap2dExamples}(a) and (b) for an example. Such
a mapping is, as before in Section \ref{sub:StartValuesIn3d} in the
context of reconstructions, not unique. However, the properties of
this mapping largely influence the convergence properties of the resulting
sub-elements. 

This is outlined for the case of a triangular sub-element where we
introduce three different mappings $\vek r\left(\vek a\right)$. That
is, the reference element in $\vek a$ is mapped differently to the
coordinate system $\vek r$ as shown in Fig.~\ref{fig:VisMap2dExamples}(c).
Note that the outer contour is the same for the three mappings, however,
different for the inner region. Mapping 1 is defined in \cite{Solin_2003a}
and is outlined in the appendix \ref{sub:MappingTriQuad} for the
case of \emph{three} curved edges in $\mathbb{R}^{3}$. The situation
here is the reduced case where the mapping is to $\mathbb{R}^{2}$
and only one of the edges is curved (the diagonal edge 2 of the reference
triangle). Nevertheless, the same formulas given in the appendix may
be employed, it follows from Eq.~(\ref{eq:EdgeContribution}) that
$\vek r^{\mathrm{edge}\, k}\left(u_{k}\right)\neq\vek0$ only for
$k=2$. 

A second mapping is the original blending function mapping proposed
in \cite{Gordon_1973a,Gordon_1973b,Szabo_2004a} na{\"{i}}vely extended
to triangles. Starting point is the approach in the appendix restricted
to one curved edge. We set $\vek r^{\mathrm{edge}\, k}\left(u_{k}\right)=\vek0$
and $R_{k}=0$ for $k=1,3$. For $k=2$, we define $u_{2}=(b-a)/(a+b)$,
$R_{2}=(a+b)$, and $\vek r^{\mathrm{edge}\,2}\left(u_{2}\right)$
follows from Eq.~(\ref{eq:EdgeContribution}). The third mapping
is the intersection mapping of \cite{Scott_1973a,Gockenbach_2006a}
and is motivated from geometrical considerations. In the reference
element, auxiliary points on the element boundary are computed for
any given inner point $\vek a_{i}$. These points are mapped to the
sub-element in $\vek r$ which is not a problem because they are on
the outer contour. They imply two straight lines and the intersection
of the two lines is the mapped inner point $\vek r_{i}$.

\begin{figure}
\centering

\subfigure[Map 1]{\includegraphics[height=5cm]{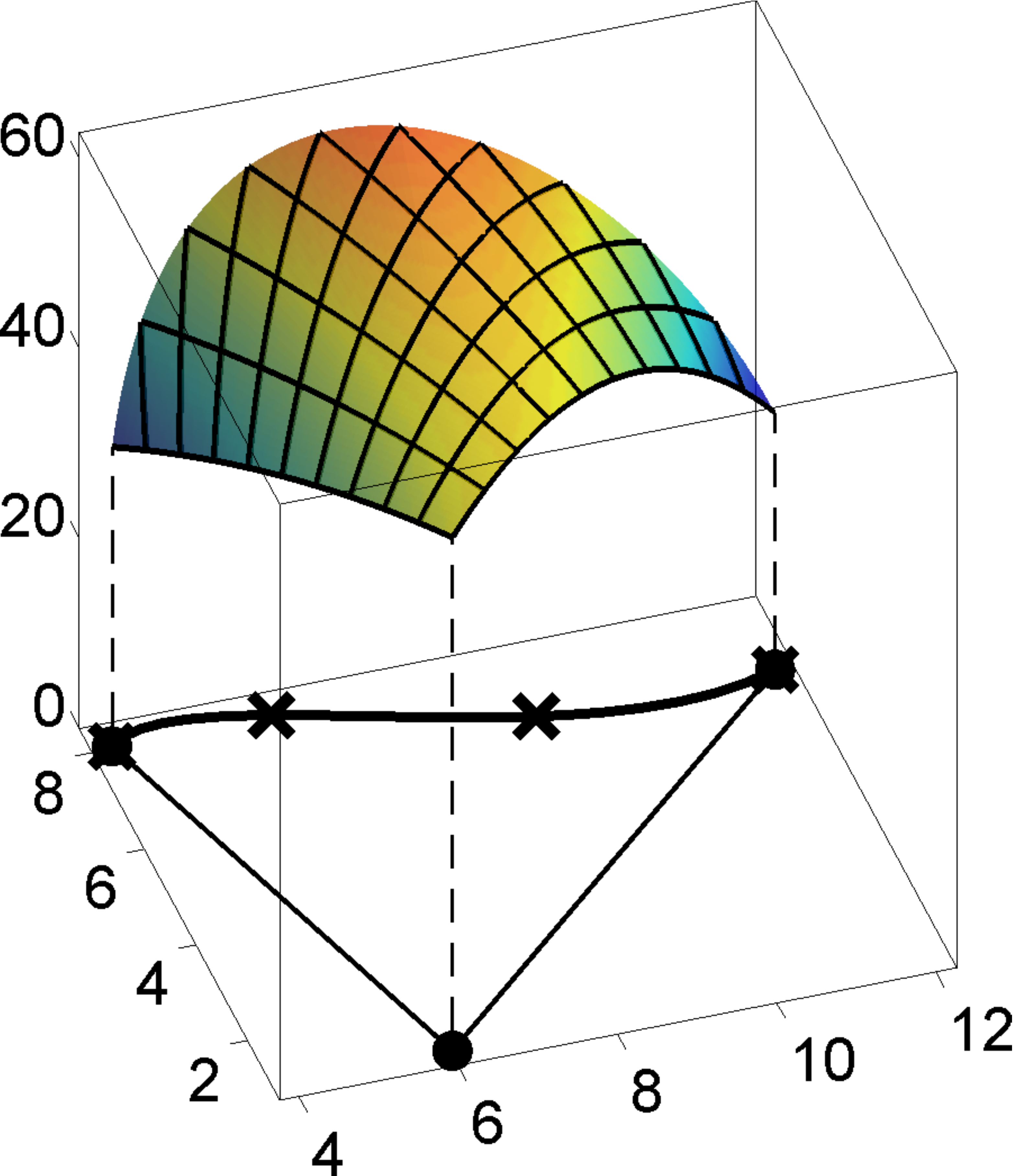}}\qquad\subfigure[Map 2]{\includegraphics[height=5cm]{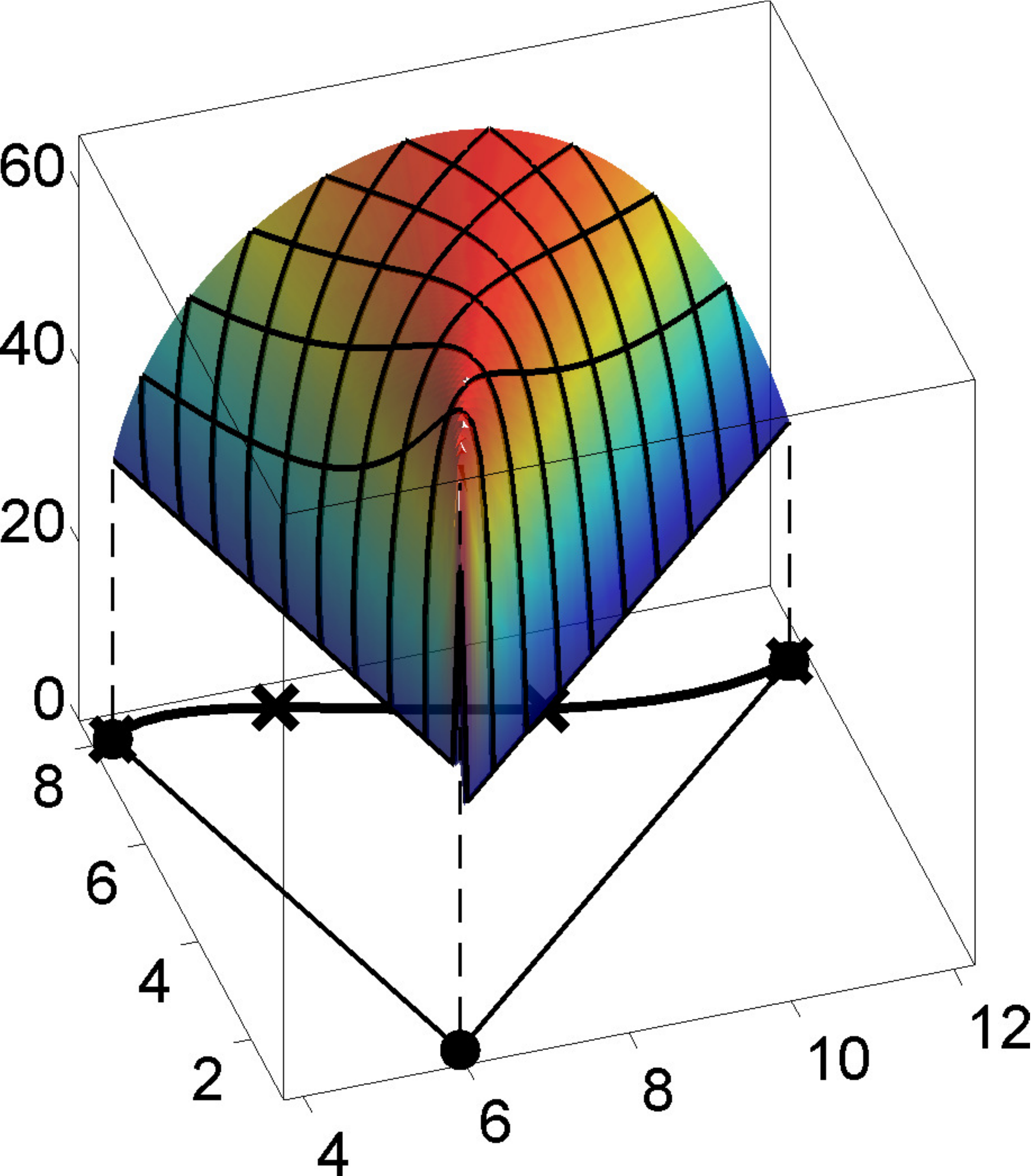}}\qquad\subfigure[Map 3]{\includegraphics[height=5cm]{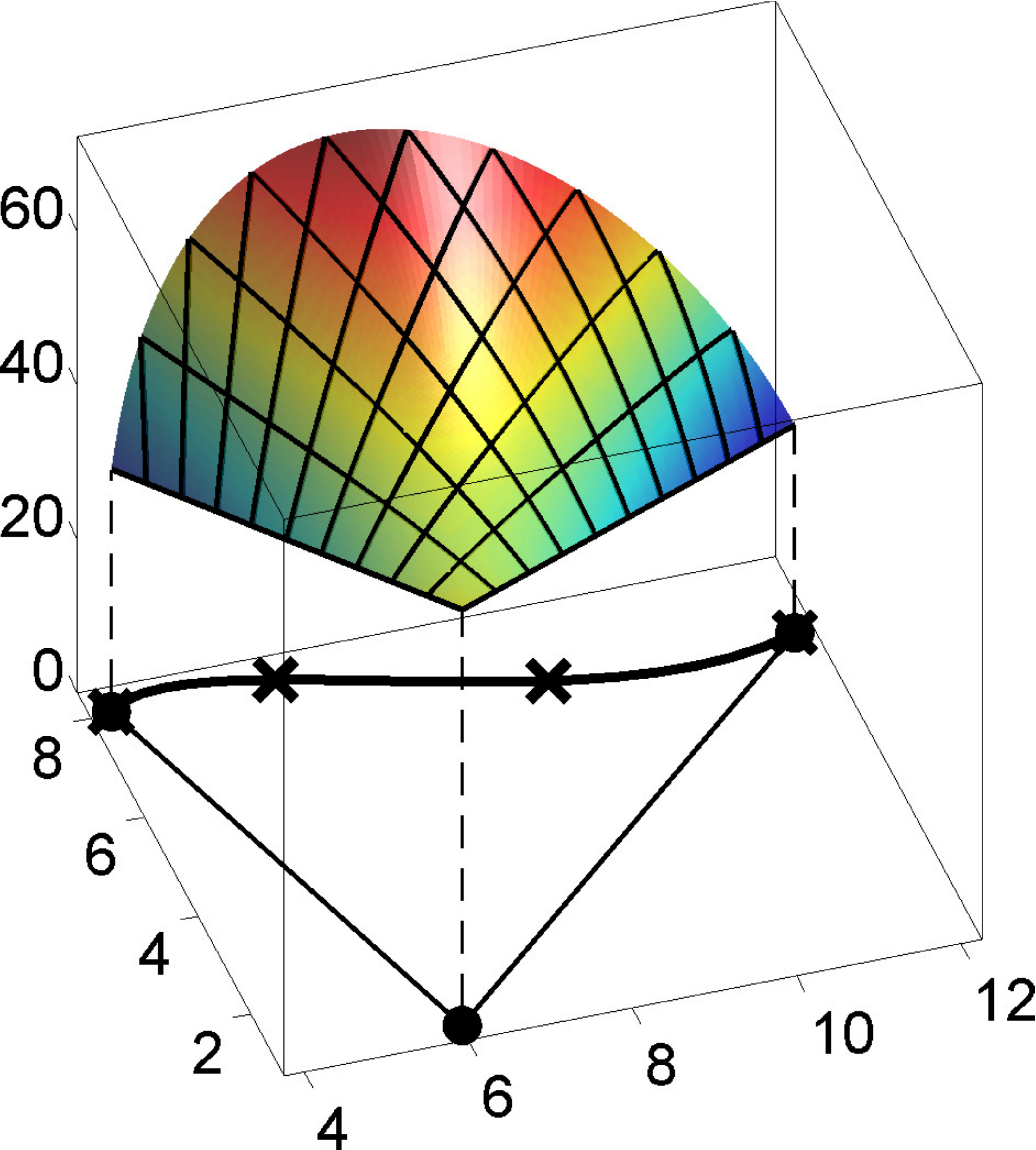}}

\caption{\label{fig:VisMap2dJacobians}Jacobians of the different mappings.}
\end{figure}

Although Fig.~\ref{fig:VisMap2dExamples}(c) indicates that all three
mappings successfully map to the linear triangle with one higher-order
side, it is found that the second map yields sub-optimal accuracy
in the convergence studies. This could be traced back to the fact
that the mapping is not as smooth as needed for higher-order accuracy.
See Fig.~\ref{fig:VisMap2dJacobians} for visualizations of the Jacobians
of the three mappings. It is seen that the second variant does not
feature a smooth Jacobian near the node opposite to the higher-order
side which explains its sub-optimal performance. The intersection
mapping and the map based on \cite{Solin_2003a} and outlined in the
appendix feature smooth Jacobians and perform optimal. However, the
intersection mapping is not easily extended to the three-dimensional
situation (three straight lines do not necessarily intersect in one
point). Therefore, in the following we restrict ourselves only to
the mappings in two and three dimensions as described in \cite{Solin_2003a};
they are general, smooth and lead to optimal element nodes in the
sub-elements as shall be seen below.

The procedure for \emph{quadrilaterals} with one higher order side,
is again an adaption of the more general map outlined in the appendix
\ref{sub:MappingTriQuad}. In this case, the mapping coincides with
the blending function mapping of \cite{Gordon_1973a,Gordon_1973b,Szabo_2004a}.

\subsection{Decomposition in 3D\label{XX_Decomp3d}}

In three dimensions, it was outlined in Section \ref{sub:Topologies}
that a tetrahedron may fall into two different topologies. For topology
1, the tetrahedron was cut into a sub-tetrahedron and a sub-prism
where the curved face is triangular, see also Fig.~\ref{fig:TopologicalCases}(b).
For topology 2, the cut tetrahedron is decomposed into two prisms
where one of the quadrilateral faces is curved in each sub-prism,
see Fig.~\ref{fig:TopologicalCases}(c). The curved faces coincide
with the reconstructed higher-order interface elements as discussed
in Section \ref{sec:Reconstruction}. Each sub-element has straight
edges except for those belonging to the curved face.

\begin{figure}
\centering

\subfigure[Tetrahedron]{\includegraphics[height=4cm]{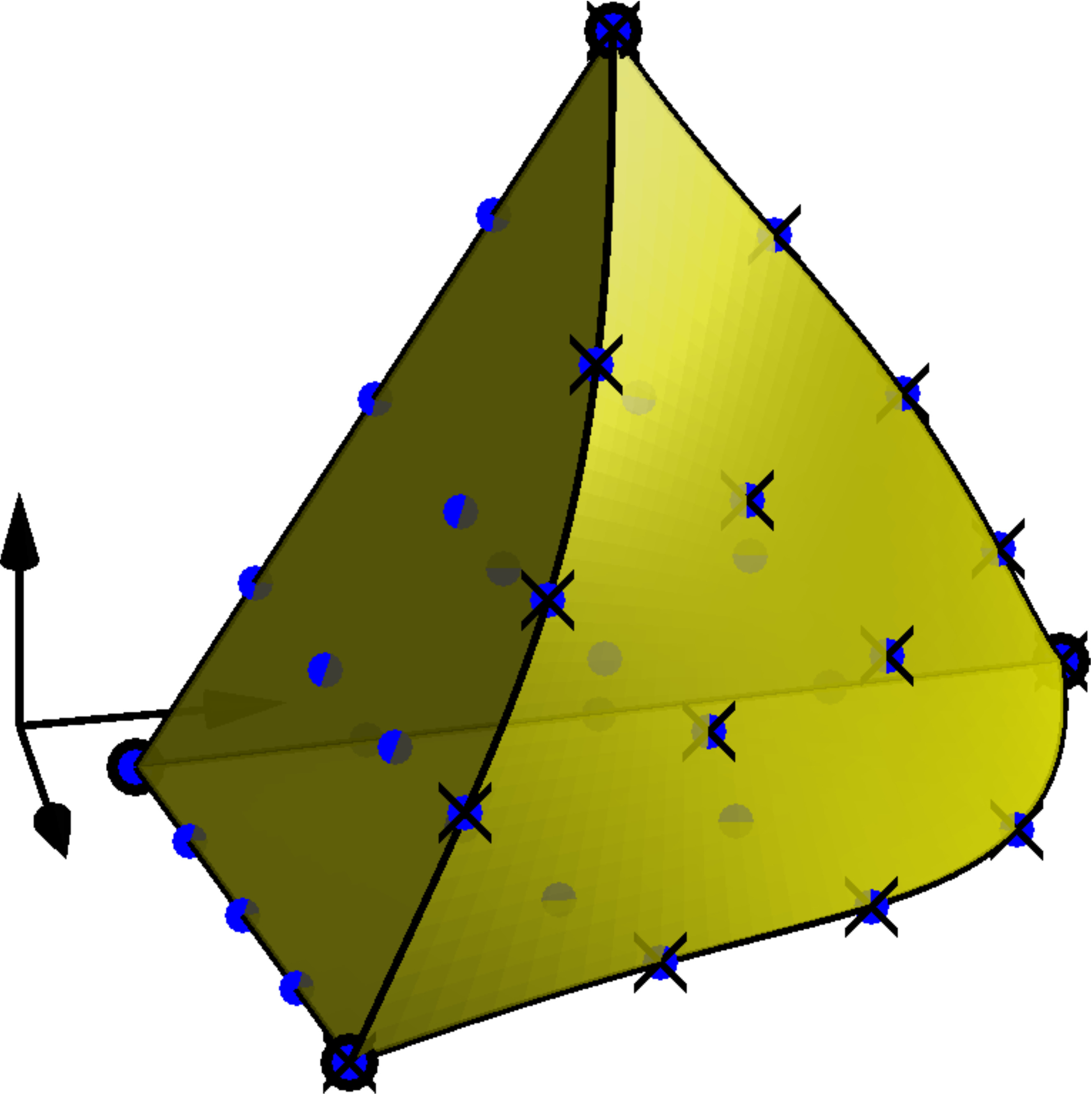}}\quad\subfigure[Prism 1]{\includegraphics[height=4cm]{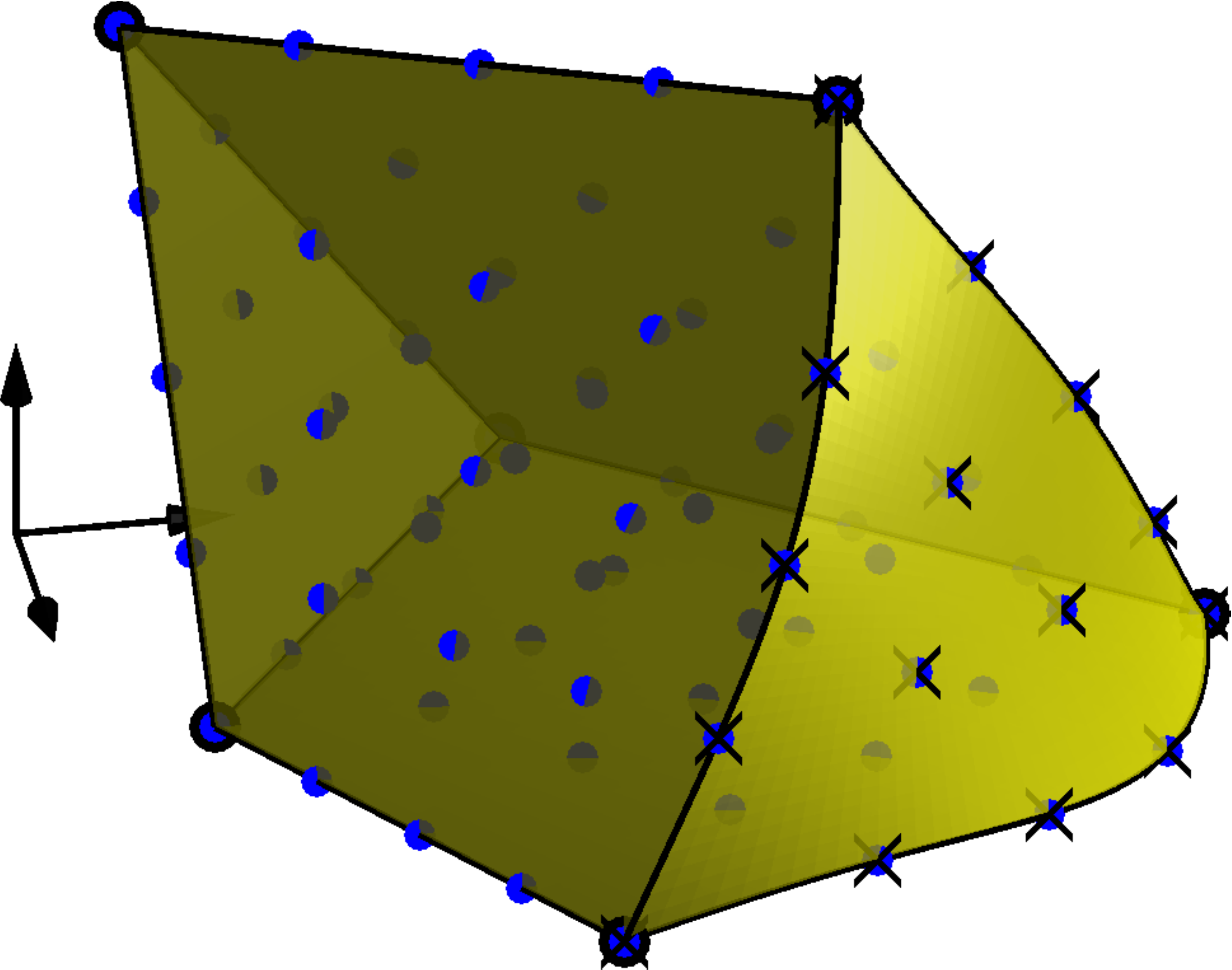}}\quad\subfigure[Prism 2]{\includegraphics[height=4cm]{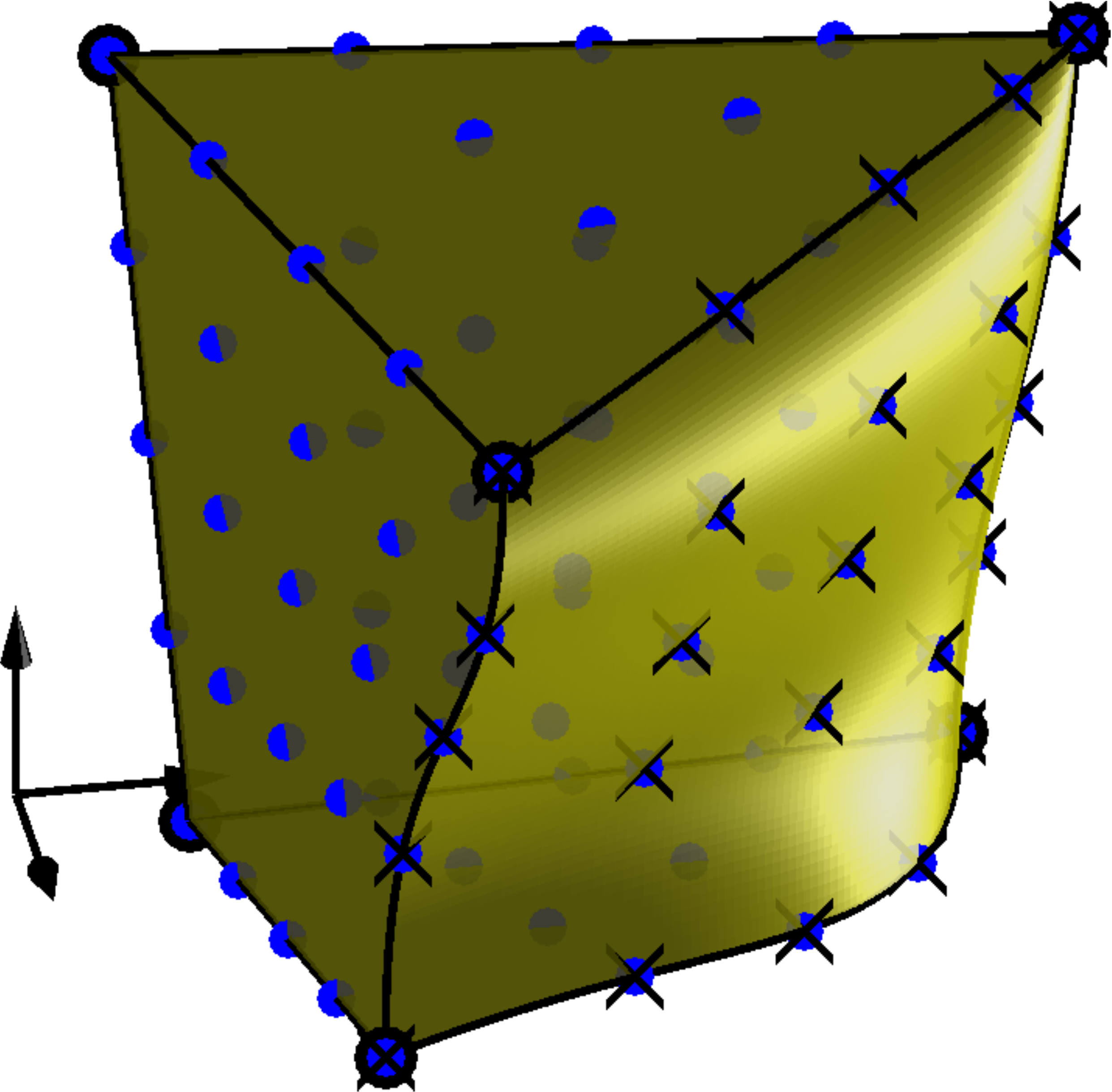}}

\caption{\label{fig:VisMap3dExamples}Examples of resulting higher-order elements
(blue nodes) from the special sub-elements with one higher-order face
(black crosses), (a) tetrahedral sub-element, (b) and (c) prismatic
sub-elements with triangular or quadrilateral higher-order face, respectively.}
\end{figure}

Again, the task is to define maps (of element nodes) from higher-order
reference tetrahedra and prisms in coordinates $\vek a$ to the cut
reference background element in $\vek r$. For some possible situations,
examples are shown in Fig.~\ref{fig:VisMap3dExamples} for all types
of sub-elements. For the sub-elements with one higher-order face,
the corresponding mappings are outlined in appendix \ref{sub:MappingTetra}
for tetrahedra and in appendix \ref{sub:MappingPrisms} for prisms.
The formulas are based on \cite{Solin_2003a} and adapted to the present
situation. The resulting formulas are quite involved wherefore they
are moved to the appendix.

\subsection{Map of integration points}

Once the cut background elements are successfully decomposed into
sub-elements in the reference domain, one may simply map the obtained
element nodes to the physical background element using the isoparametric
concept. It is then simple to map integration points according to
any desired quadrature rule to the sub-elements. Depending on the
particular context for which higher-order accurate integration points
are sought, this may be preferred in the sub-elements of the cut \emph{reference}
or \emph{physical} background element. For the integration studies
performed herein, the integration points are generated in the physical
background elements. However, in the context of FDM and XFEM-related
methods one would generate integration points first in the reference
background elements, then evaluate the corresponding shape functions
and finally map them to the physical background elements using the
isoparametric concept.

\subsection{Numerical results for decompositions in 2D and 3D\label{XX_VolumeResults}}

The numerical results presented here are a direct extension of the
studies of Section \ref{XX_InterfResults} to the integration and
interpolation of implicitly defined areas and volumes.

\subsubsection{Preliminary examples}

\begin{figure}
\centering

\subfigure[]{\includegraphics[width=4cm]{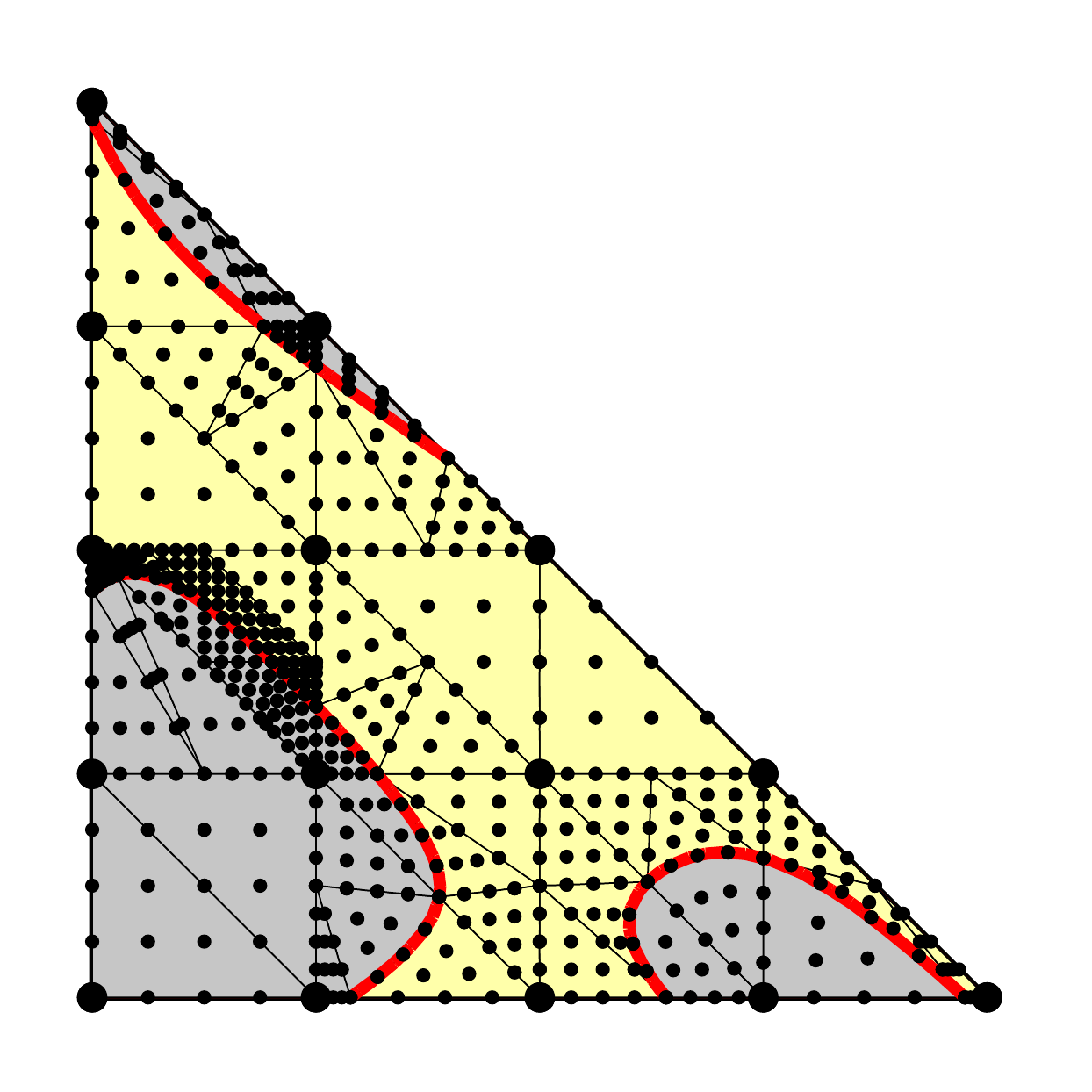}}\subfigure[]{\includegraphics[width=4cm]{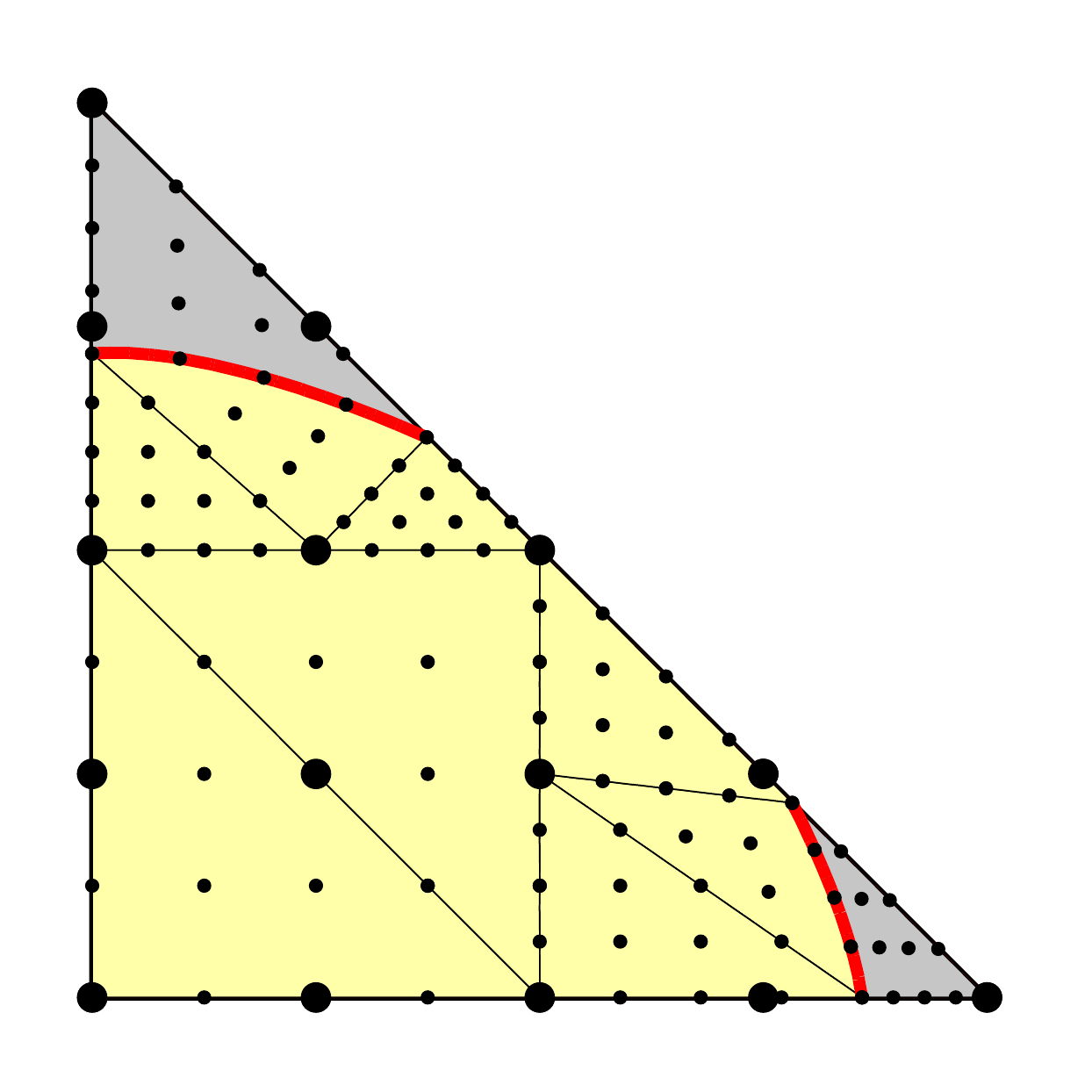}}\subfigure[]{\includegraphics[width=4cm]{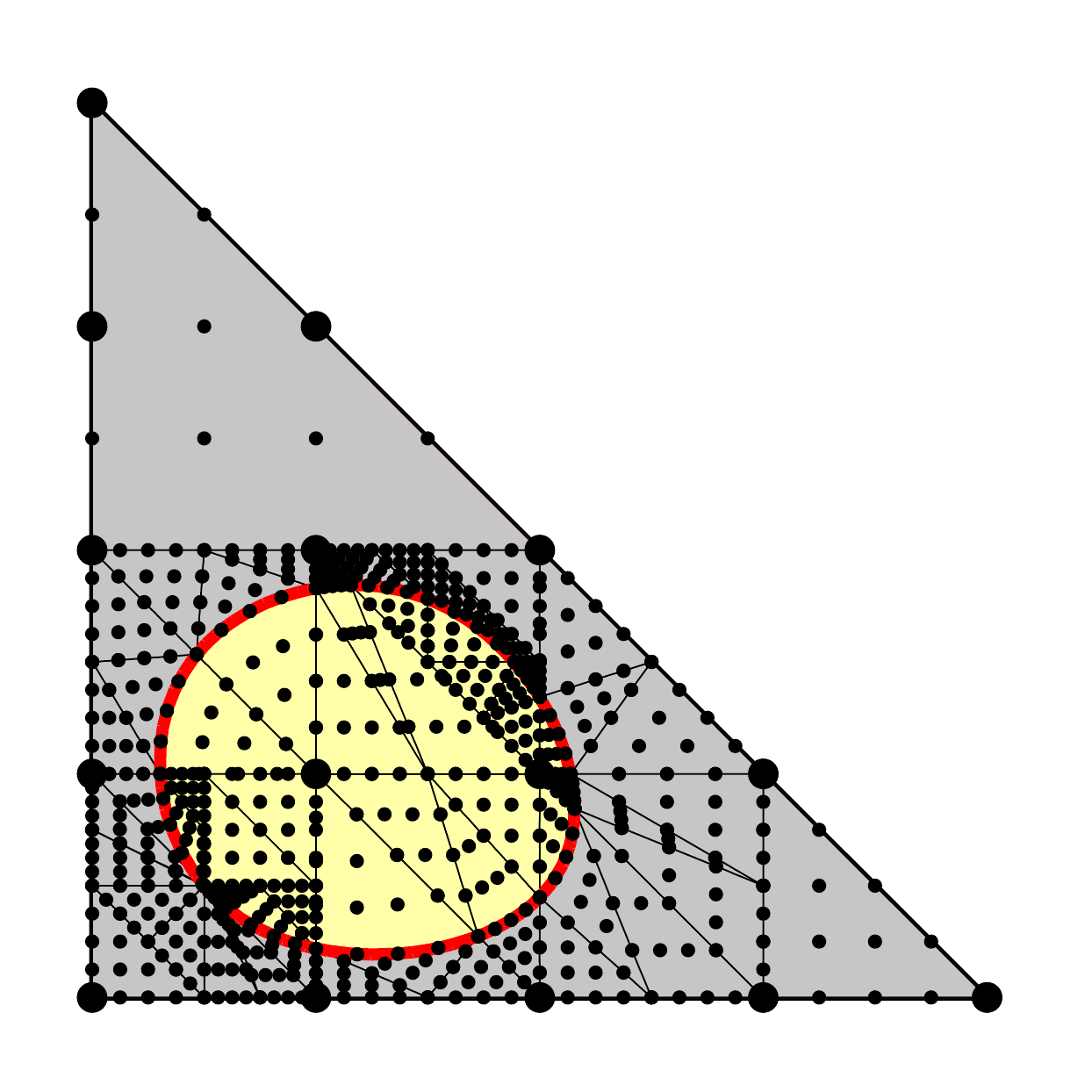}}\subfigure[]{\includegraphics[width=4cm]{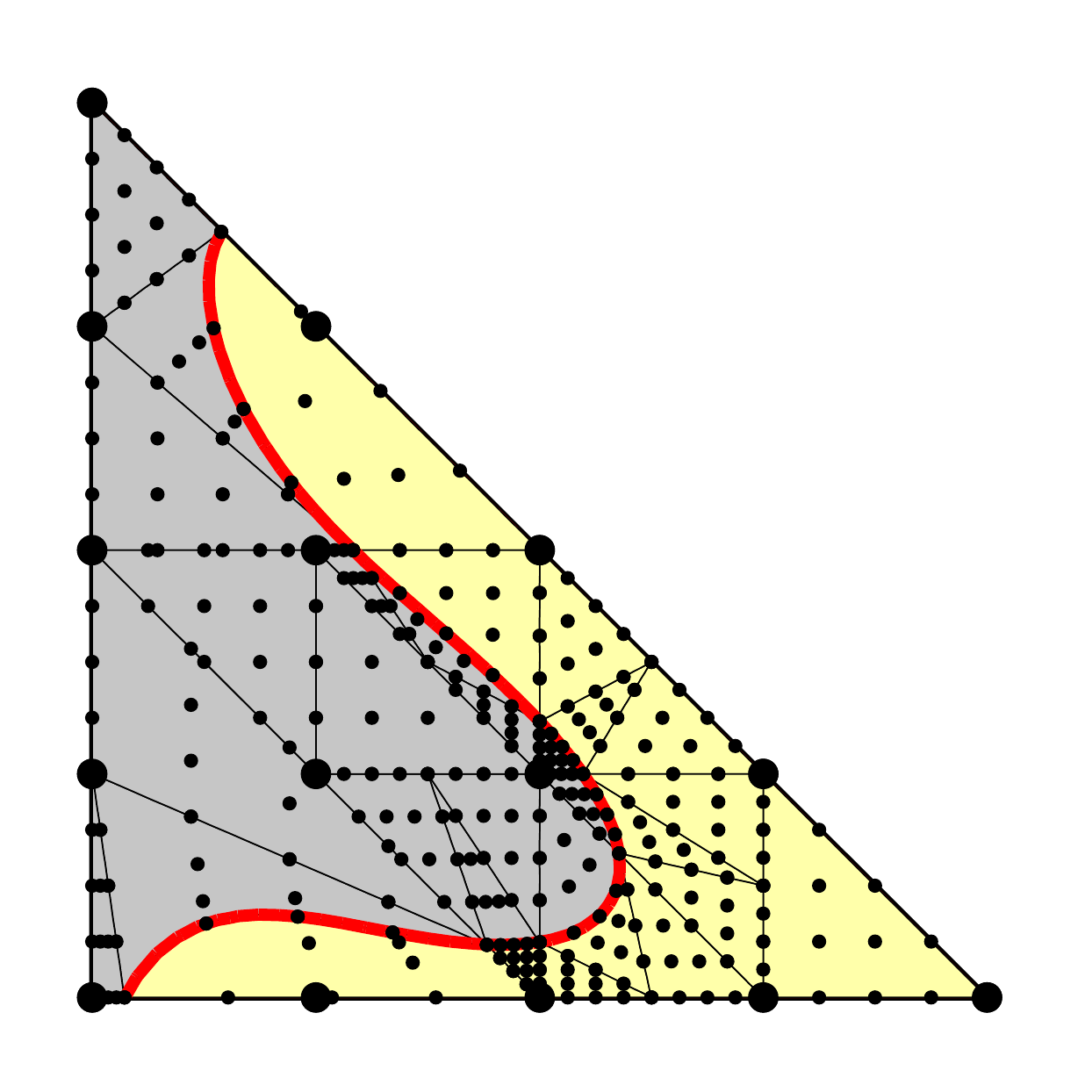}}

\caption{\label{fig:VisDecomp2d}Decomposed reference triangular elements obtained
by recursive refinements.}
\end{figure}

\begin{figure}
\centering

\subfigure[Example 1, inside]{\includegraphics[width=4cm]{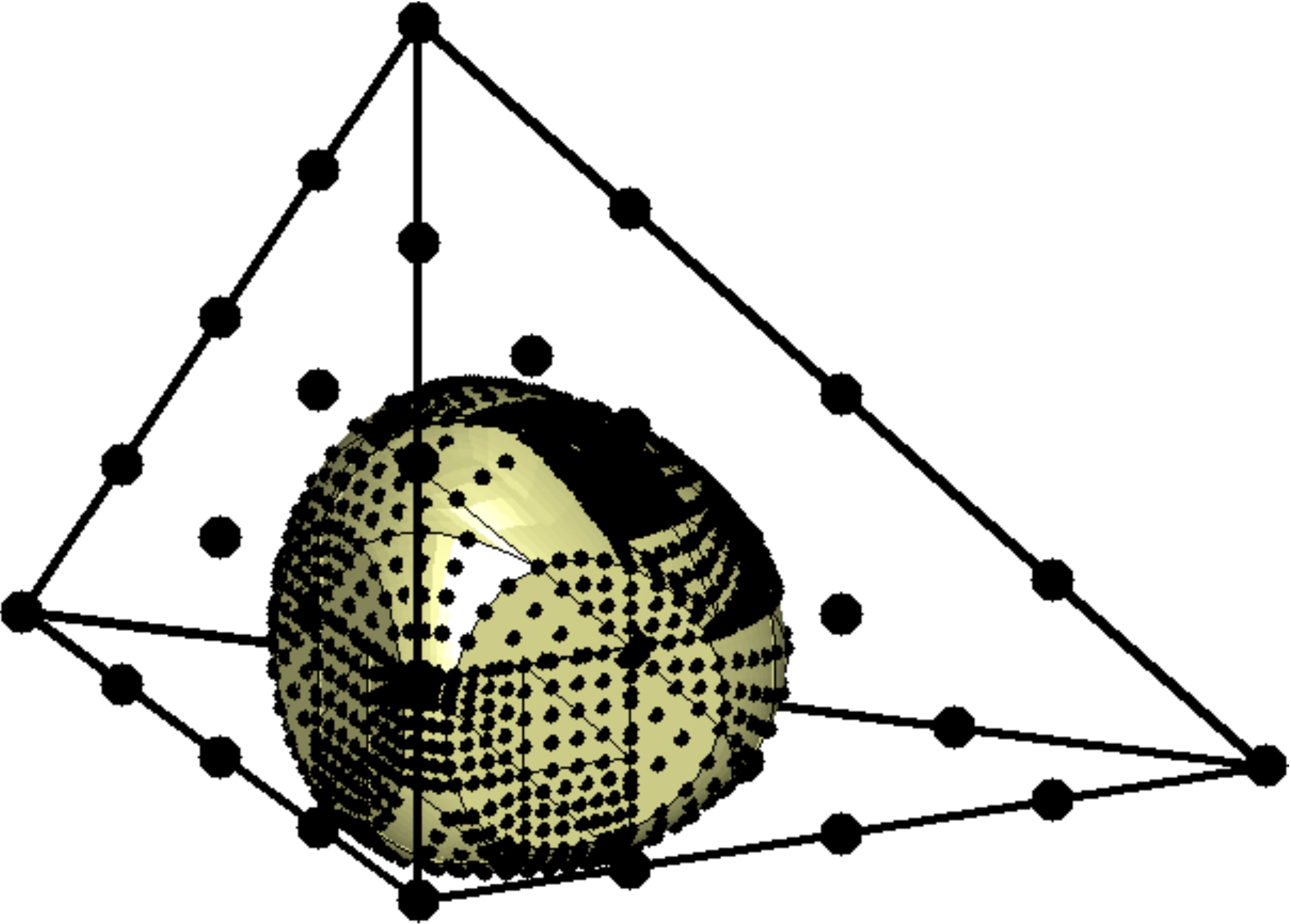}}\subfigure[Example 1, outside]{\includegraphics[width=4cm]{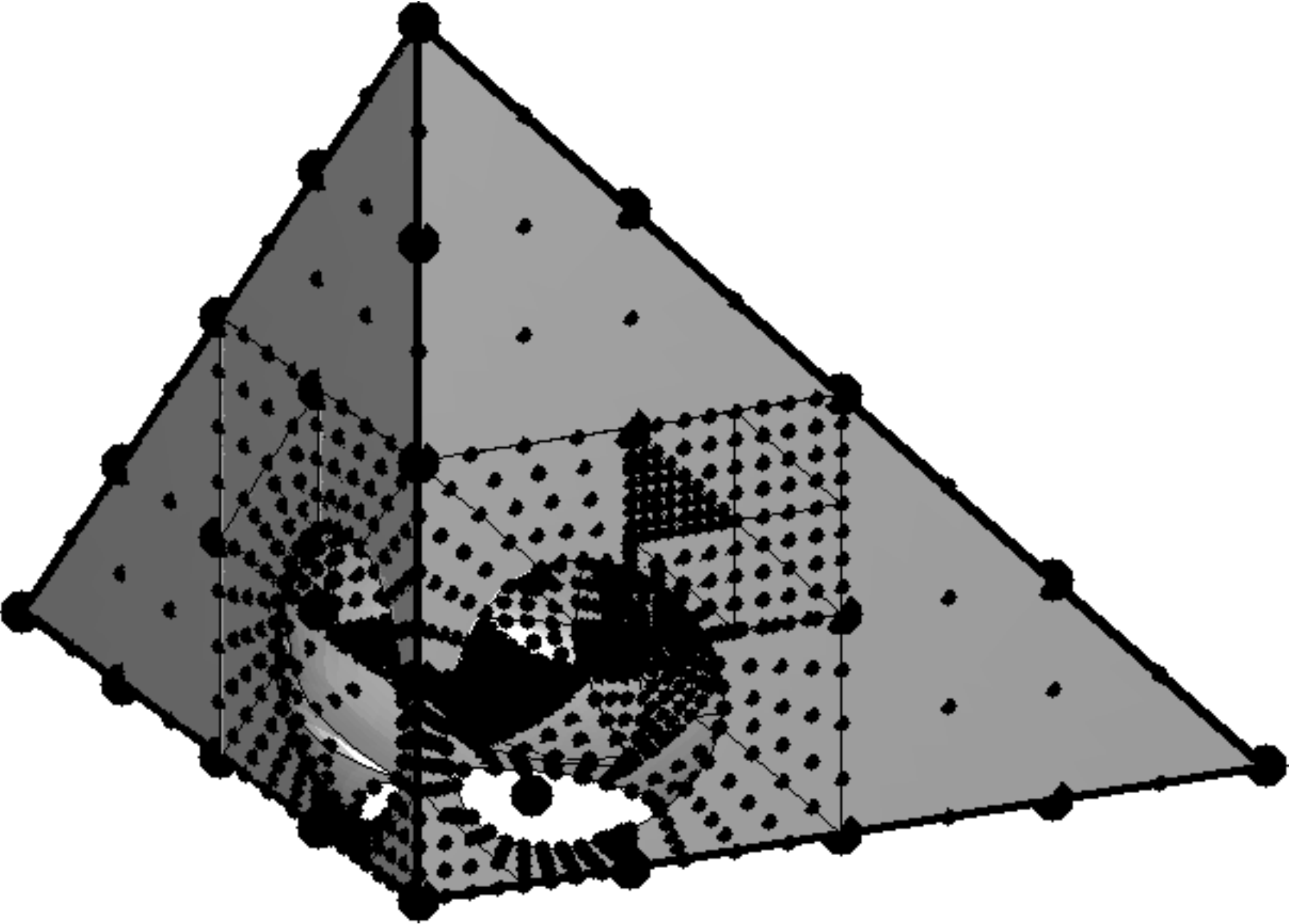}}\subfigure[Example 2, inside]{\includegraphics[width=4cm]{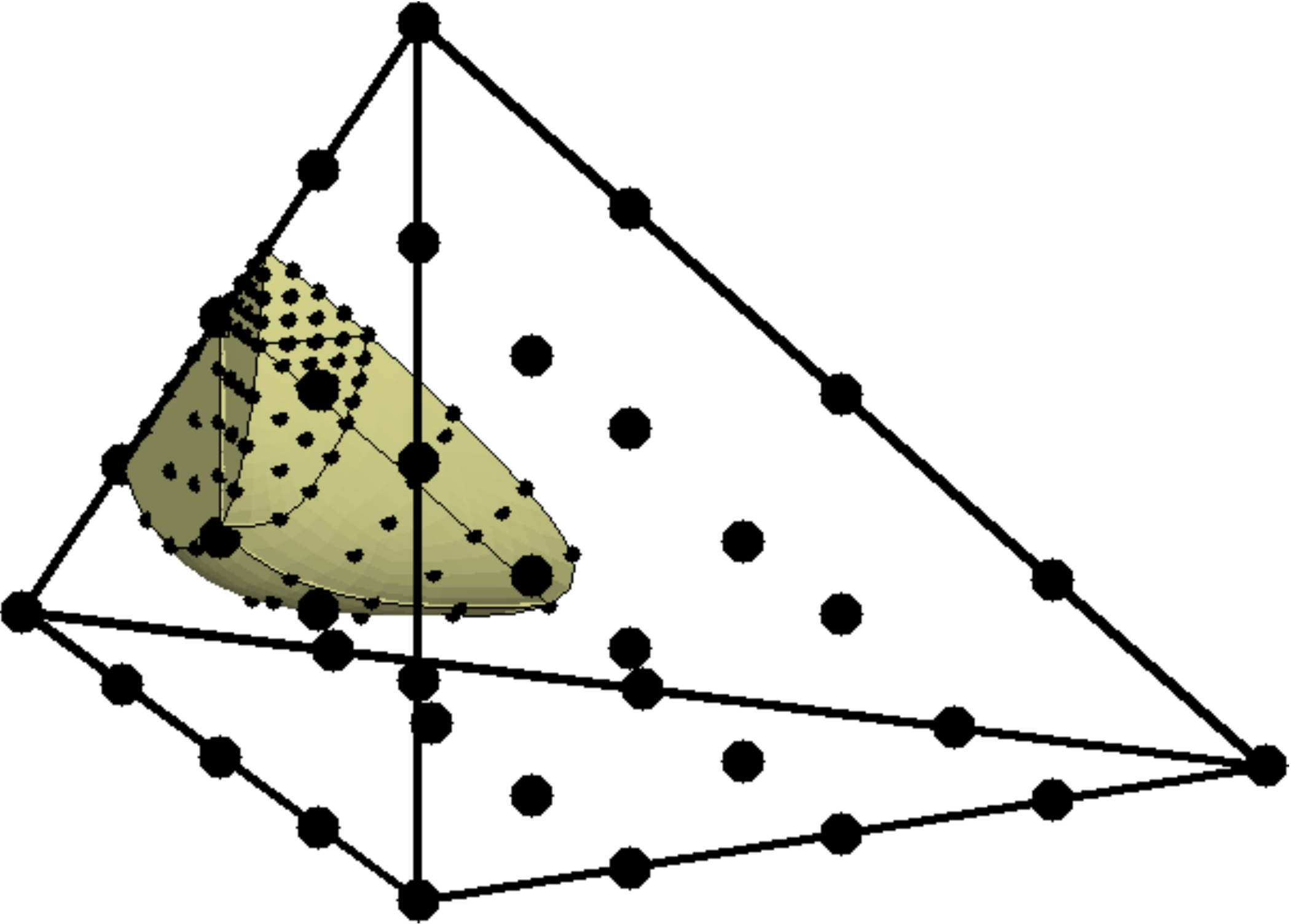}}\subfigure[Example 2, outside]{\includegraphics[width=4cm]{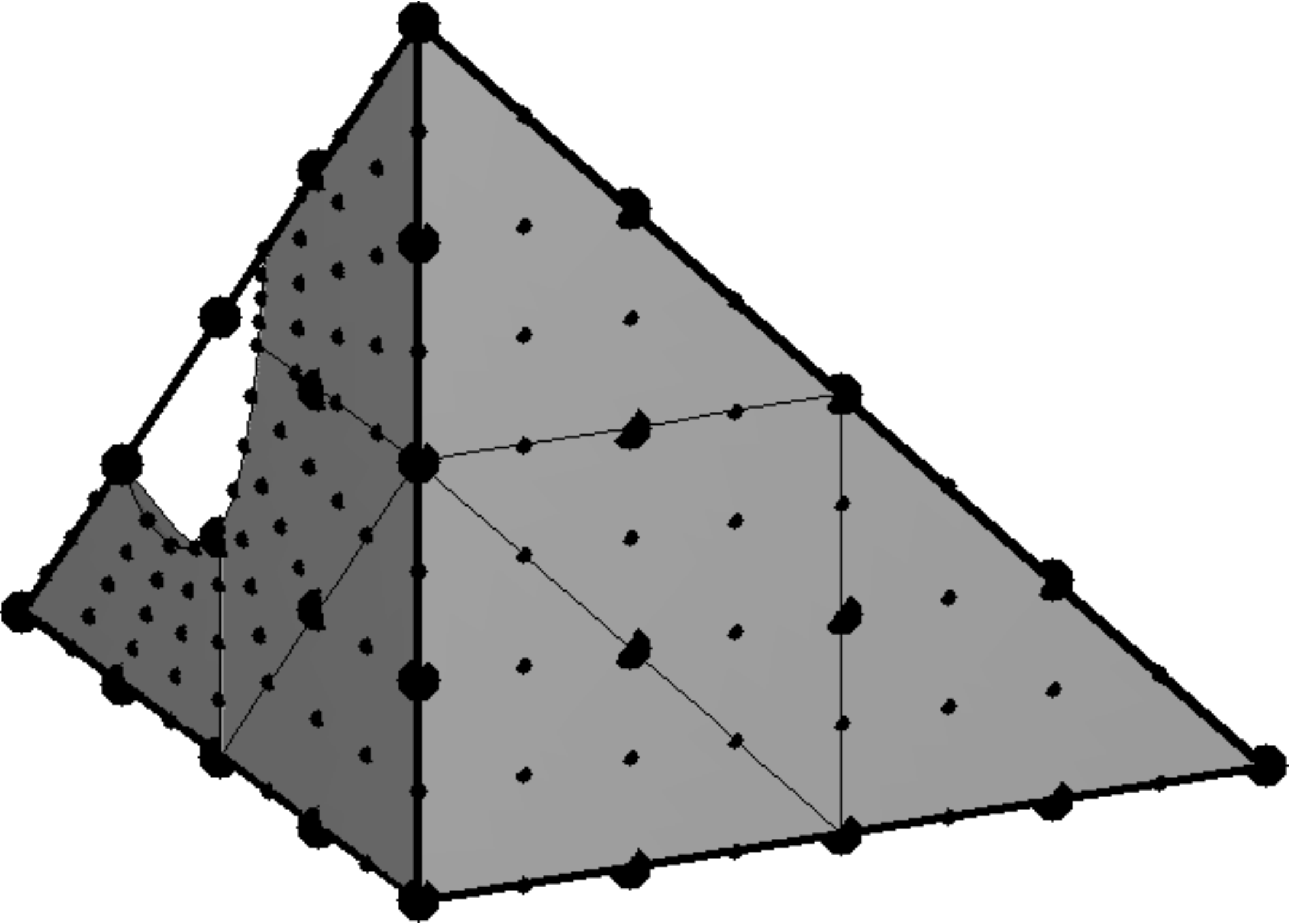}}

\caption{\label{fig:VisDecomp3d}Decomposed reference tetrahedral elements
obtained by recursive refinements.}
\end{figure}

Let us first consider some examples of decompositions. Fig.~\ref{fig:VisDecomp2d}
shows some resulting sub-elements for complex level-set data in 2D.
According to Section \ref{sub:ValididyAndRecursion}, recursive refinements
are necessary to obtain valid level-set data in the refined elements.
Note the relation to Fig.~\ref{fig:Vis2dIsoLinesInElem}. Examples
for decompositions in 3D are shown in Fig.~\ref{fig:VisDecomp3d}
which correspond to zero-level sets visualized in Fig.~\ref{fig:Vis3dIsoSurfInElem}.

\subsubsection{Numerical results for decompositions in 2D\label{XXX_VolumeResults2d}}

The setup is identical to the one described in Section \ref{XXX_InterfResults2d}.
The same background meshes, level-set functions and integrands are
considered here. Three different (relative) error norms are introduced
for the integration in $\Omega^{-}$:
\begin{eqnarray}
\varepsilon_{1}^{\Omega} & = & \Big|\Big(\sum_{i}w_{i}\Big)-I_{1}^{\Omega}\Big|/\left|I_{1}^{\Omega}\right|\quad\mathrm{with}\quad I_{1}^{\Omega}=\int_{\Omega^{-}}1\,\mathrm{d}\Omega,\label{eq:Error2dArea1}\\
\varepsilon_{f}^{\Omega} & = & \Big|\Big(\sum_{i}w_{i}\cdot f\left(\vek x_{i}\right)\Big)-I_{f}^{\Omega}\Big|/\left|I_{f}^{\Omega}\right|\quad\mathrm{with}\quad I_{f}^{\Omega}=\int_{\Omega^{-}}f\left(\vek x\right)\,\mathrm{d}\Omega,\label{eq:Error2dArea2}\\
\varepsilon_{f^{h}}^{\Omega} & = & \Big|\Big(\sum_{i}w_{i}\cdot f^{h}\left(\vek x_{i}\right)\Big)-I_{f}^{\Omega}\Big|/\left|I_{f}^{\Omega}\right|,\label{eq:Error2dArea3}
\end{eqnarray}
where $\vek x_{i}$ are 2D integration points in the special sub-elements
and $w_{i}$ the corresponding weights. 

\begin{figure}
\centering

\subfigure[]{\includegraphics[width=6cm]{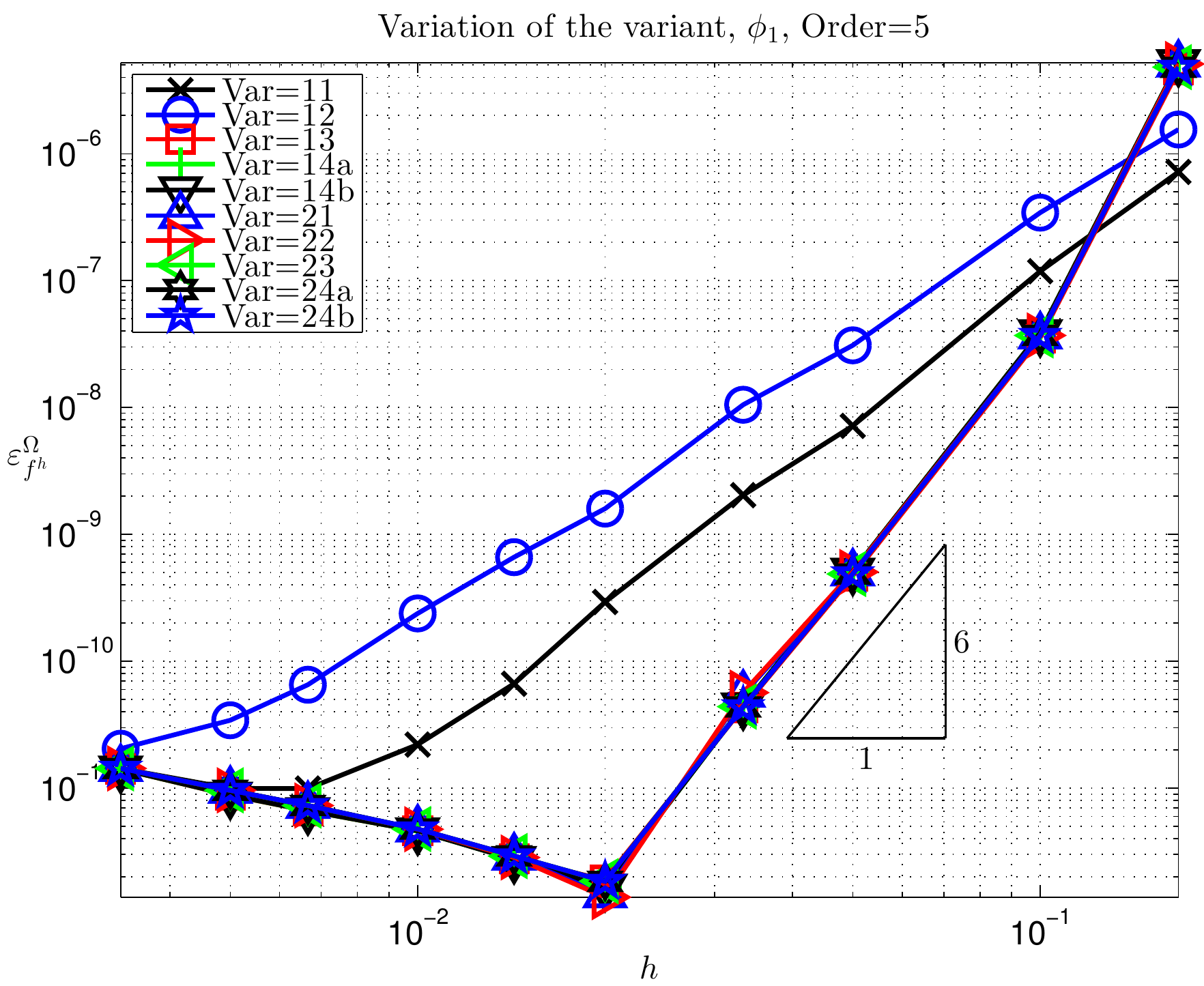}}\quad\subfigure[]{\includegraphics[width=6cm]{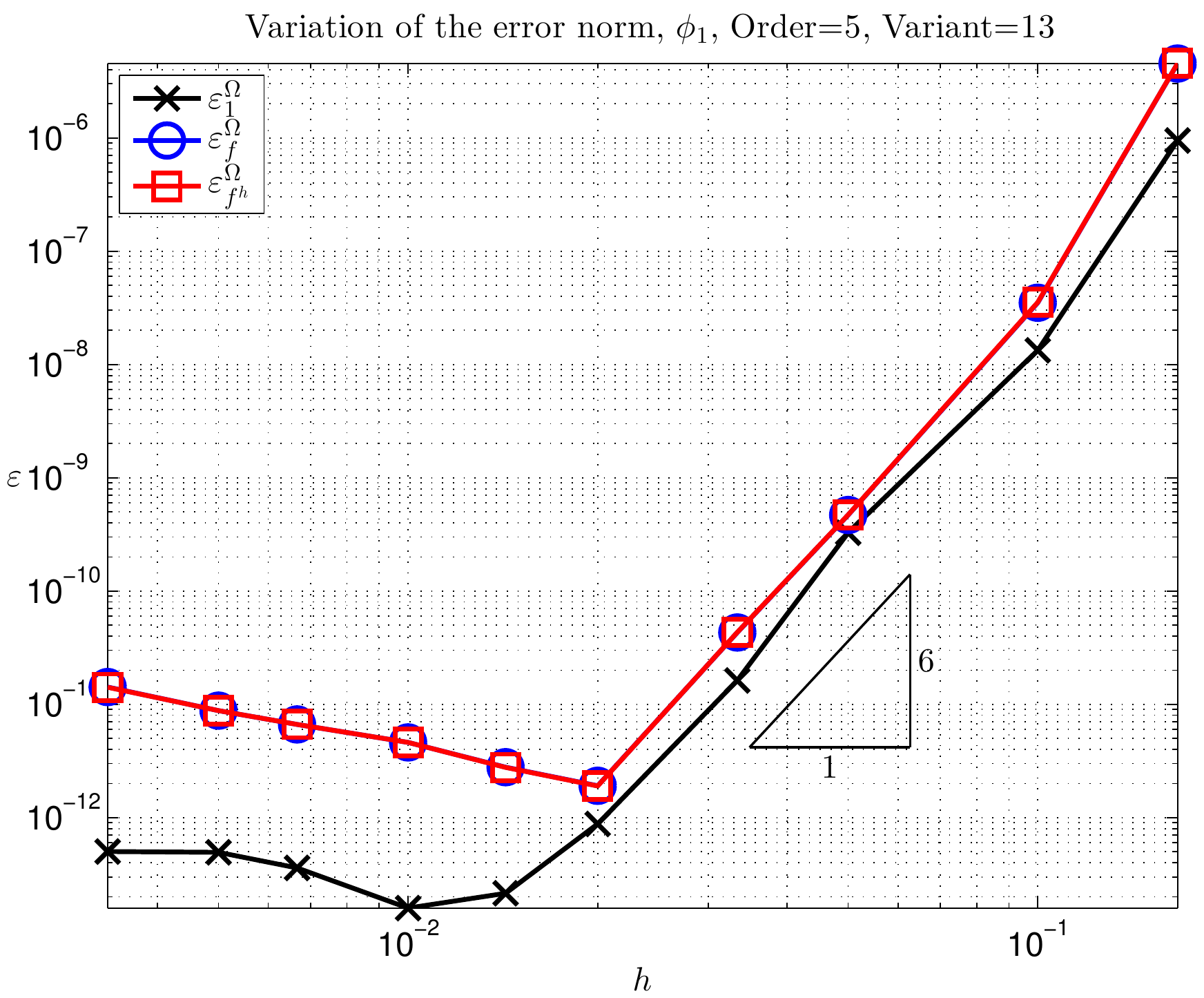}}\\\subfigure[]{\includegraphics[width=6cm]{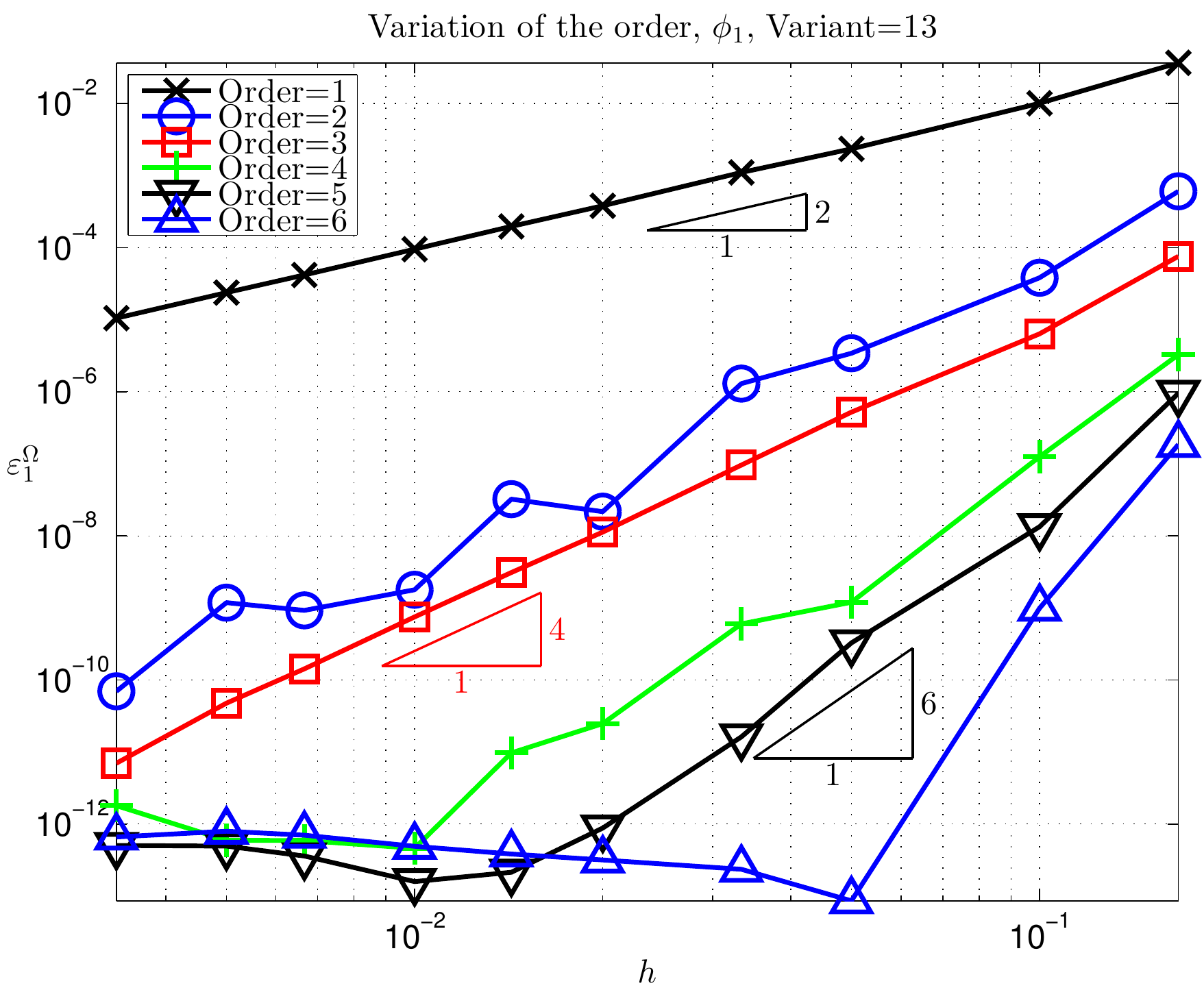}}\quad\subfigure[]{\includegraphics[width=6cm]{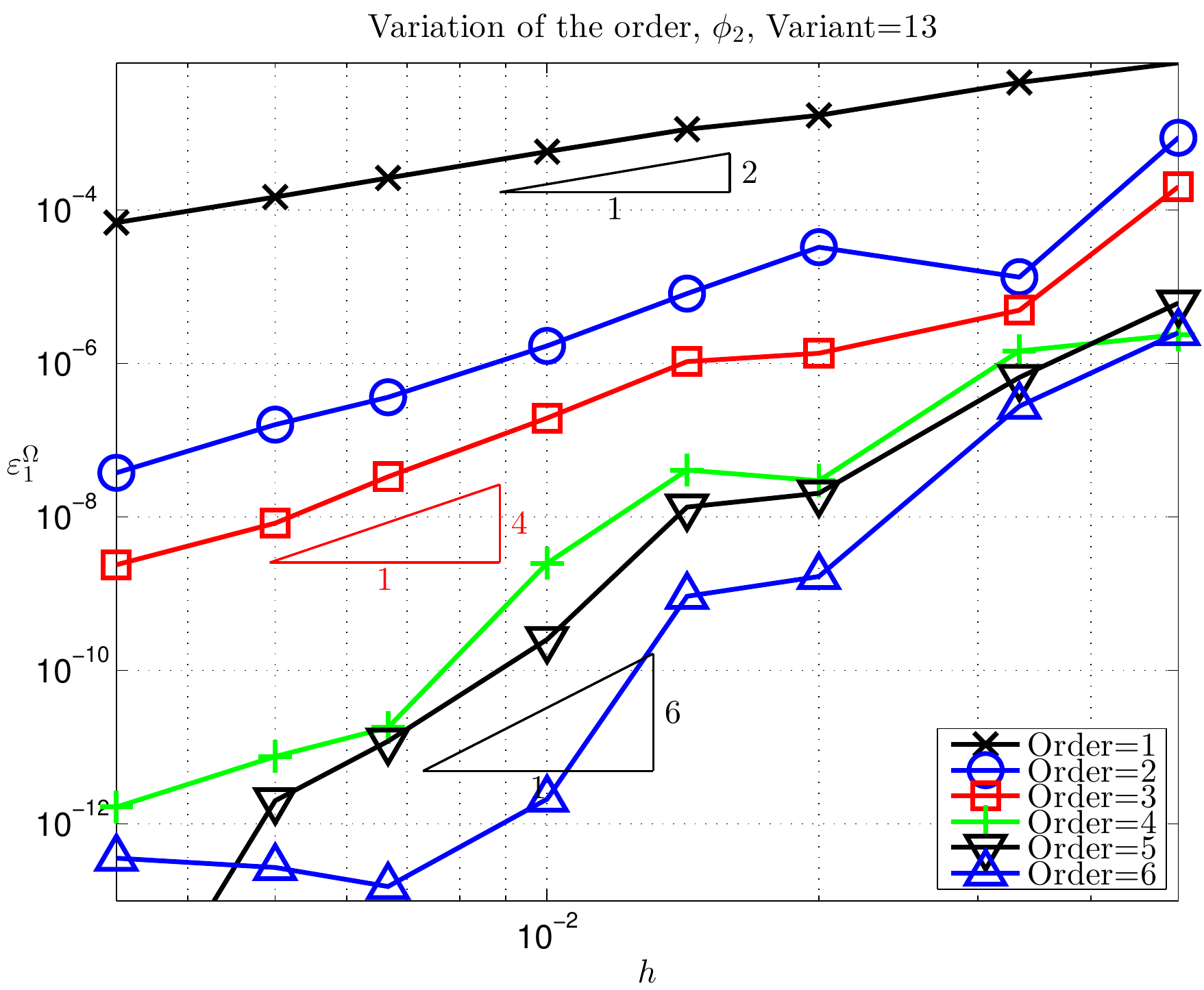}}

\caption{\label{fig:VisRes2dArea}Convergence results in 2D for integrating
and interpolating in implicitly defined areas.}
\end{figure}

Convergence results are shown in Fig.~\ref{fig:VisRes2dArea}. Fig.~\ref{fig:VisRes2dArea}(a)
confirms that most variants of the search directions perform optimally.
A very similar result was shown in Fig.~\ref{fig:VisRes2dInterf}(a)
for the integration on the zero-level sets and the same conclusions
are drawn: Firstly, search directions should be based on the normal
vectors on the intermediate reconstruction or on the gradient of the
level-set function. Secondly, Hermite reconstructions are not superior
to linear ones. Therefore, we suggest variant 13 (based on normal
vectors on linear reconstructions) for its simplicity. Fig.~\ref{fig:VisRes2dInterf}(b)
indicates that all three error norms lead to identical results which
was confirmed for the whole parameter space tested. An exception was
only found for the sub-optimal second variant of the mapping to the
sub-elements described in Section \ref{XX_Decomp2d}: There, it turned
out that results are sub-optimal in $\varepsilon_{f^{h}}^{\Omega}$
(convergence rates are bounded by 4). Optimal results for $\varepsilon_{1}^{\Omega}$
and $\varepsilon_{f}^{\Omega}$ were still possible depending on the
applied quadrature rules (which is also not desirable). Therefore,
as mentioned before, we shall only present results where the mappings
based on \cite{Solin_2003a} are used, see also the appendix. Finally,
the full convergence data for variant 13 for $\phi_{1}$ and $\phi_{2}$
according to Eqs.~(\ref{eq:ExactLevelSet2dA}) and (\ref{eq:ExactLevelSet2dB})
is shown in Fig.~\ref{fig:VisRes2dArea}(c) and (d): Optimal convergence
rates are obtained. Only data points are shown where no recursive
refinements were needed so that there is a clear element length $h$
associated to every result.

\subsubsection{Numerical results for decompositions in 3D\label{XXX_VolumeResults3d}}

\begin{figure}
\centering

\subfigure[]{\includegraphics[width=6cm]{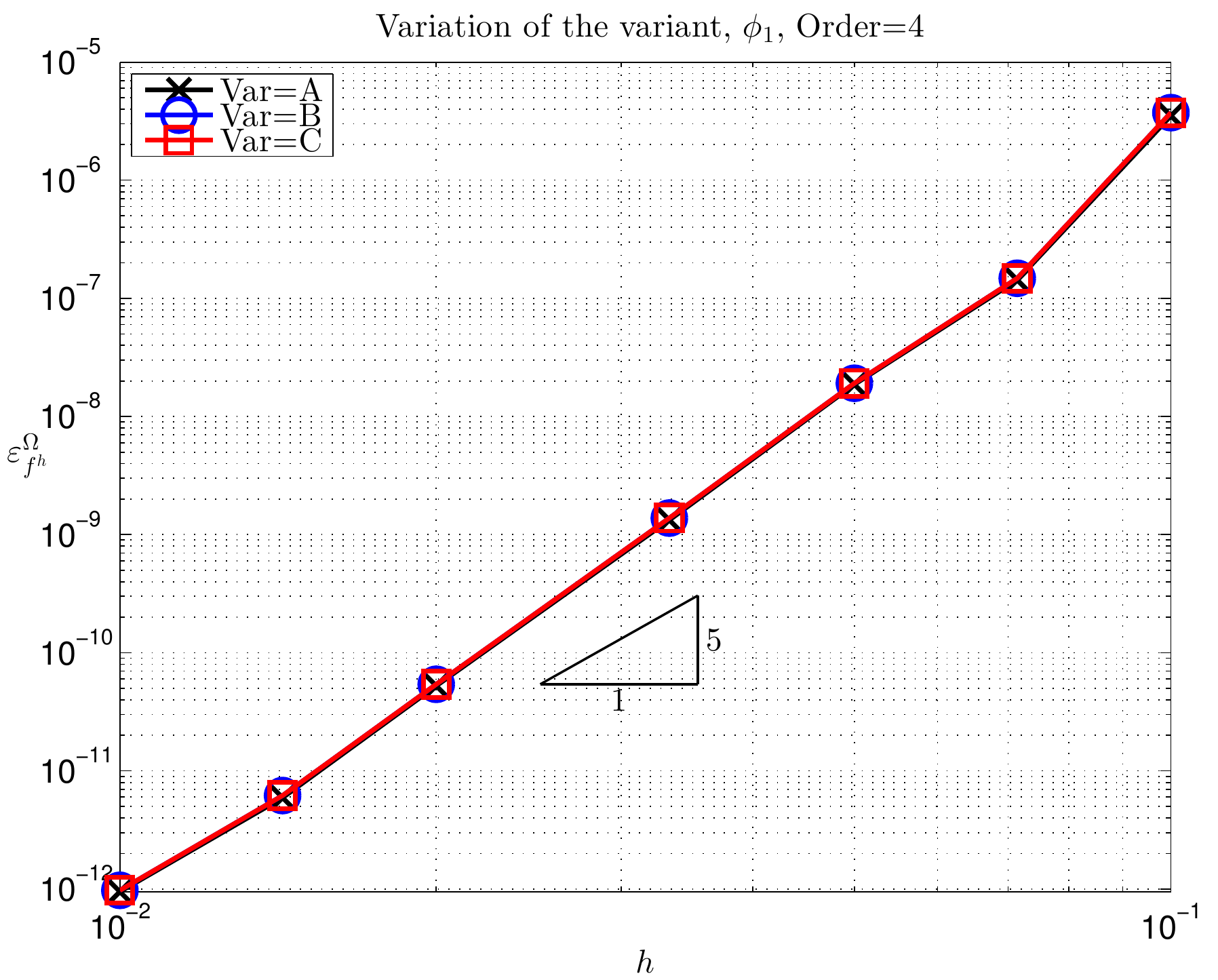}}\quad\subfigure[]{\includegraphics[width=6cm]{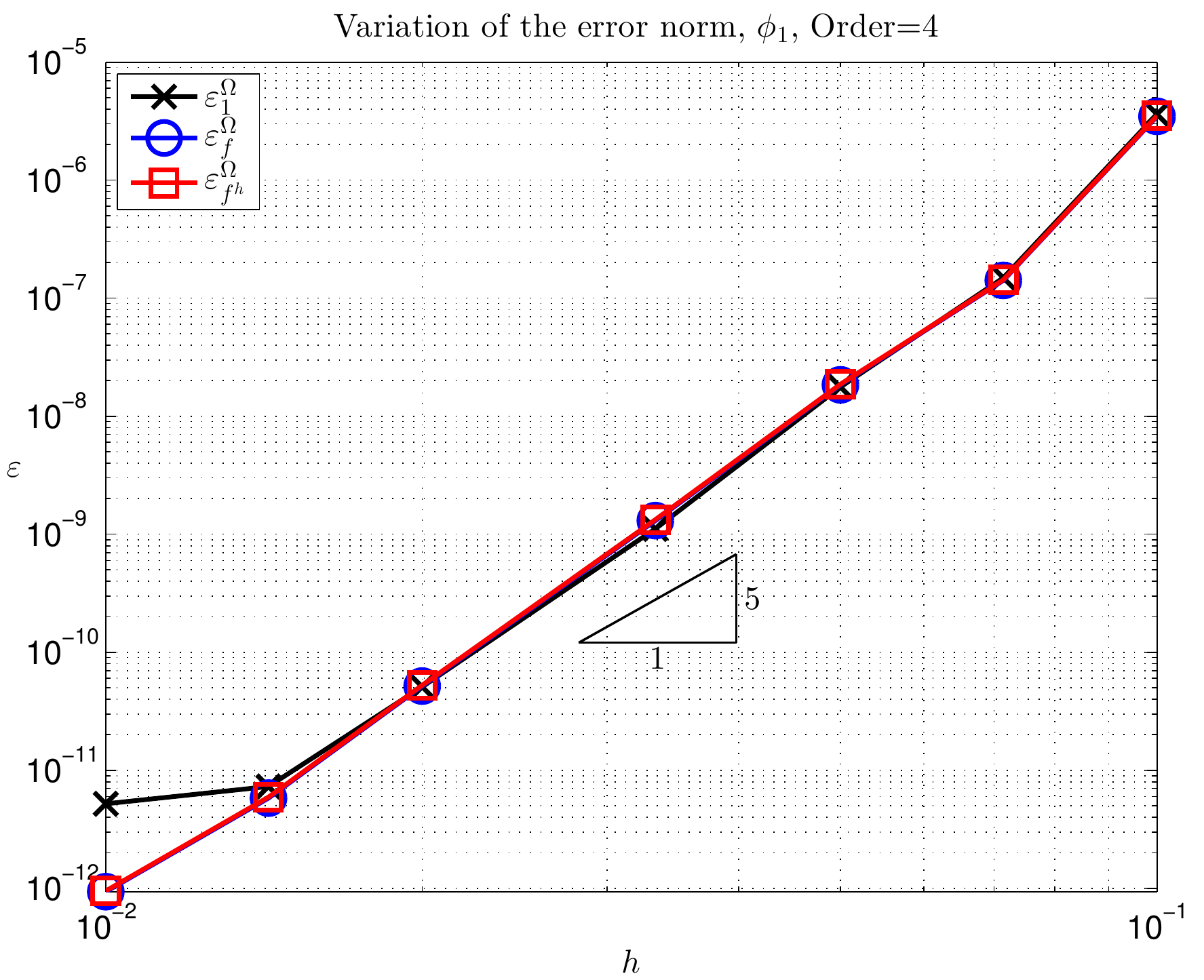}}\\\subfigure[]{\includegraphics[width=6cm]{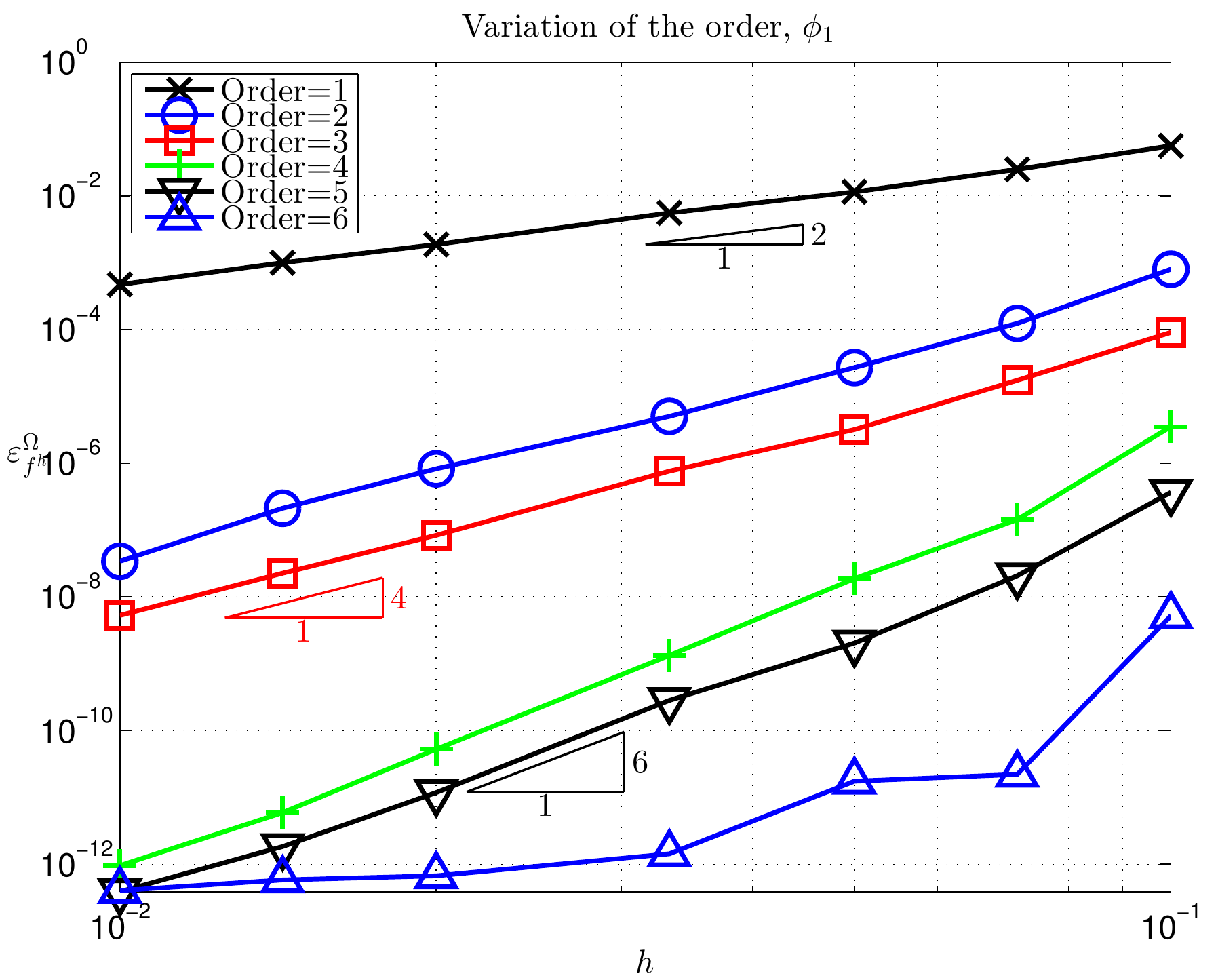}}\quad\subfigure[]{\includegraphics[width=6cm]{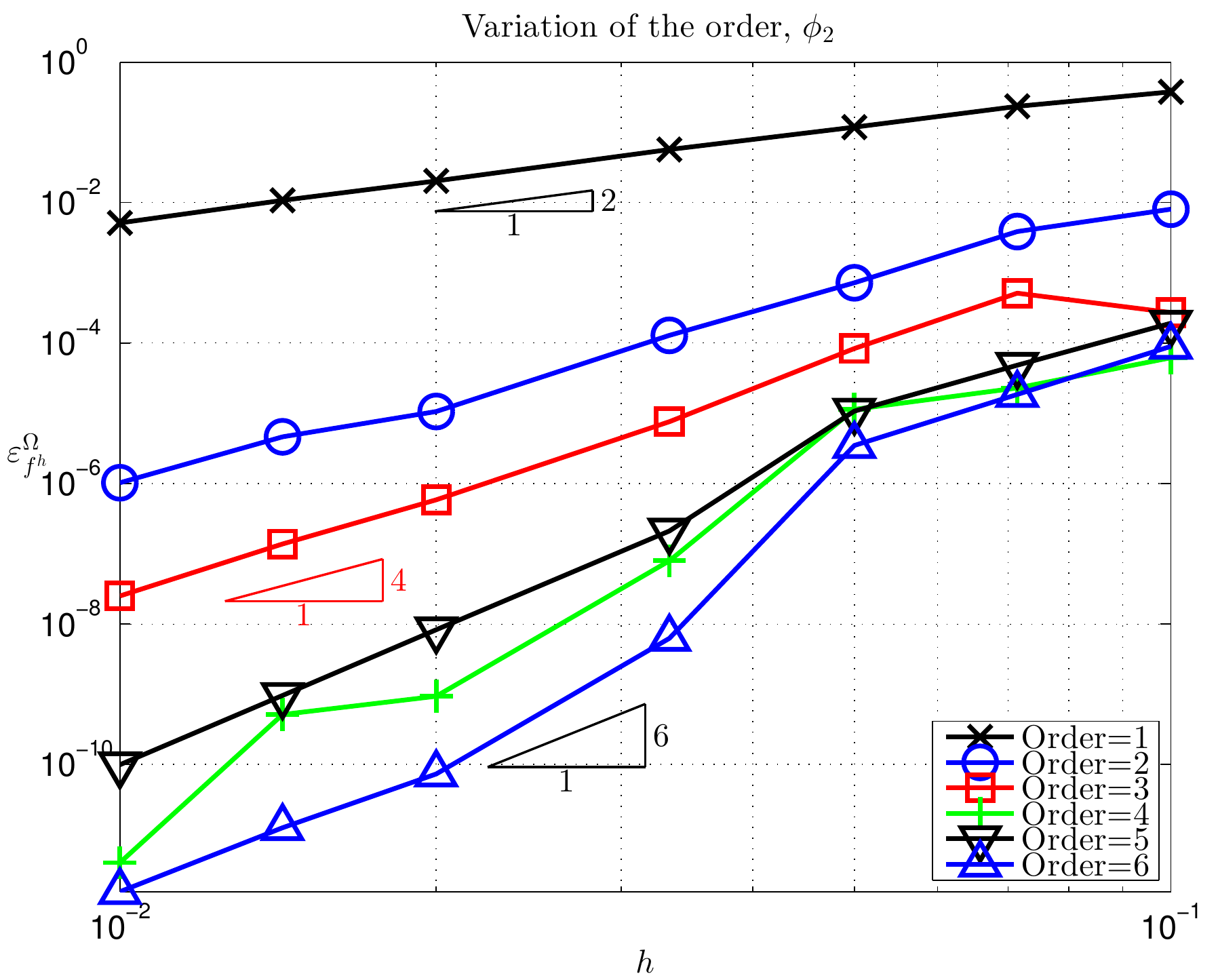}}

\caption{\label{fig:VisRes3dVol}Convergence results in 3D for integrating
and interpolating in implicitly defined volumes.}
\end{figure}

In 3D, the convergence studies follow Section \ref{XXX_InterfResults3d}.
The same definitions of the error norms from Eqs.~(\ref{eq:Error2dArea1})
to (\ref{eq:Error2dArea3}) are used. Results are shown in Fig.~\ref{fig:VisRes3dVol}
where the same form of presentation is employed than in Section \ref{XXX_InterfResults3d}.
Fig.~\ref{fig:VisRes3dVol}(a) uses variant 13 on the faces of the
tetrahedra and variants A, B, C for the inner nodes. The differences
are negligible and we suggest to use variant A because it is based
on the normal vector of the intermediate reconstruction just as variant
13. Of course, all other combinations of search directions for the
inner and outer nodes of the interface elements have been investigated
as well: Those variants for the outer nodes which performed optimal
in Section \ref{XXX_VolumeResults2d} also perform optimal in this
3D context. The fact that all three error norms behave similar is
confirmed in Fig.~\ref{fig:VisRes3dVol}(b). Finally, the full convergence
data (for variant A13) for $\phi_{1}$ and $\phi_{2}$ according to
Eqs.~(\ref{eq:ExactLevelSet3dA}) and (\ref{eq:ExactLevelSet3dB})
is shown in Fig.~\ref{fig:VisRes3dVol}(c) and (d) and, again, optimal
convergence rates are obtained. It is noted that a previous work of
the authors \cite{Fries_2015a} was not able to achieve optimal convergence
rates in 3D because the intermediate reconstructions used there were
not sufficiently smooth.

It is confirmed for all studies that properties of the reconstruction
are inherited by the decomposition. In particular, sub-optimal results
for reconstructions can not be better for decompositions. Of course,
optimal reconstructions may still lead to sub-optimal decompositions
e.g.~if the maps to the sub-elements, as described above, are not
sufficiently smooth.

\section{Corners and edges\label{sec:CornersAndEdges}}

So far, the reconstruction and decomposition of background elements
has been described with respect to \emph{one} level-set function.
It is possible to define a level-set function $\phi\left(\vek x\right)$
whose zero-level set also involves corners and edges. Then, this function
is necessarily $C_{0}$-continuous. However, the meshing procedures
described in the previous sections rely on an \emph{interpolated}
level-set function $\phi^{h}\left(\vek x\right)$. Provided that $\phi\left(\vek x\right)$
is only $C_{0}$-continuous, an accurate interpolation $\phi^{h}\left(\vek x\right)$
poses requirements on the background mesh that we wish to avoid here.
For example, in 2D, for being able to capture a corner, one would
have to impose the constraint that the corner is at least on the edge
of the background mesh. Otherwise, there would be a dramatic loss
in accuracy in the interpolation $\phi^{h}\left(\vek x\right)$ and
a higher-order decomposition would not be possible.

\begin{figure}
\centering

\subfigure[Interface with corners]{\includegraphics[width=4cm]{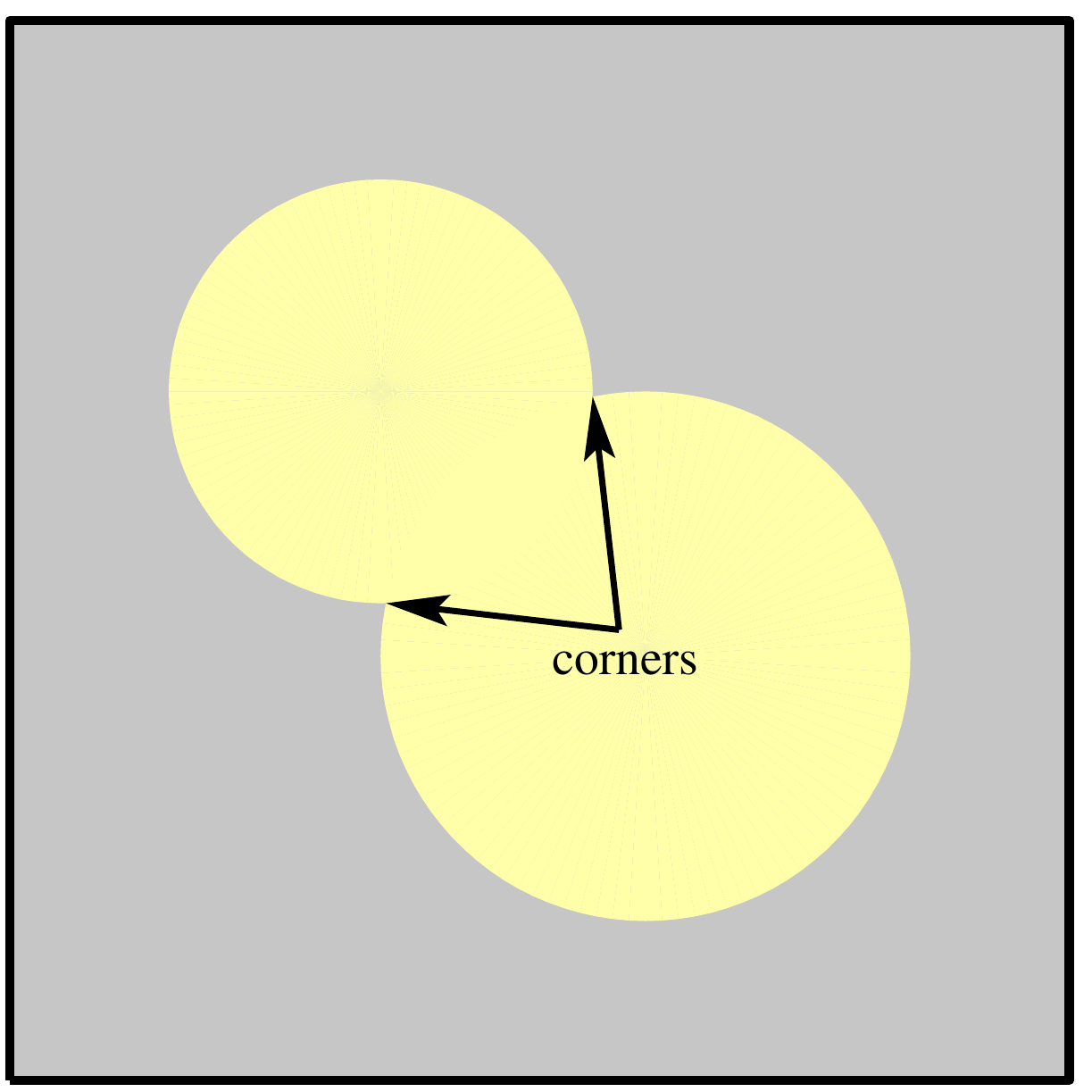}}\;\subfigure[Background mesh]{\includegraphics[width=4cm]{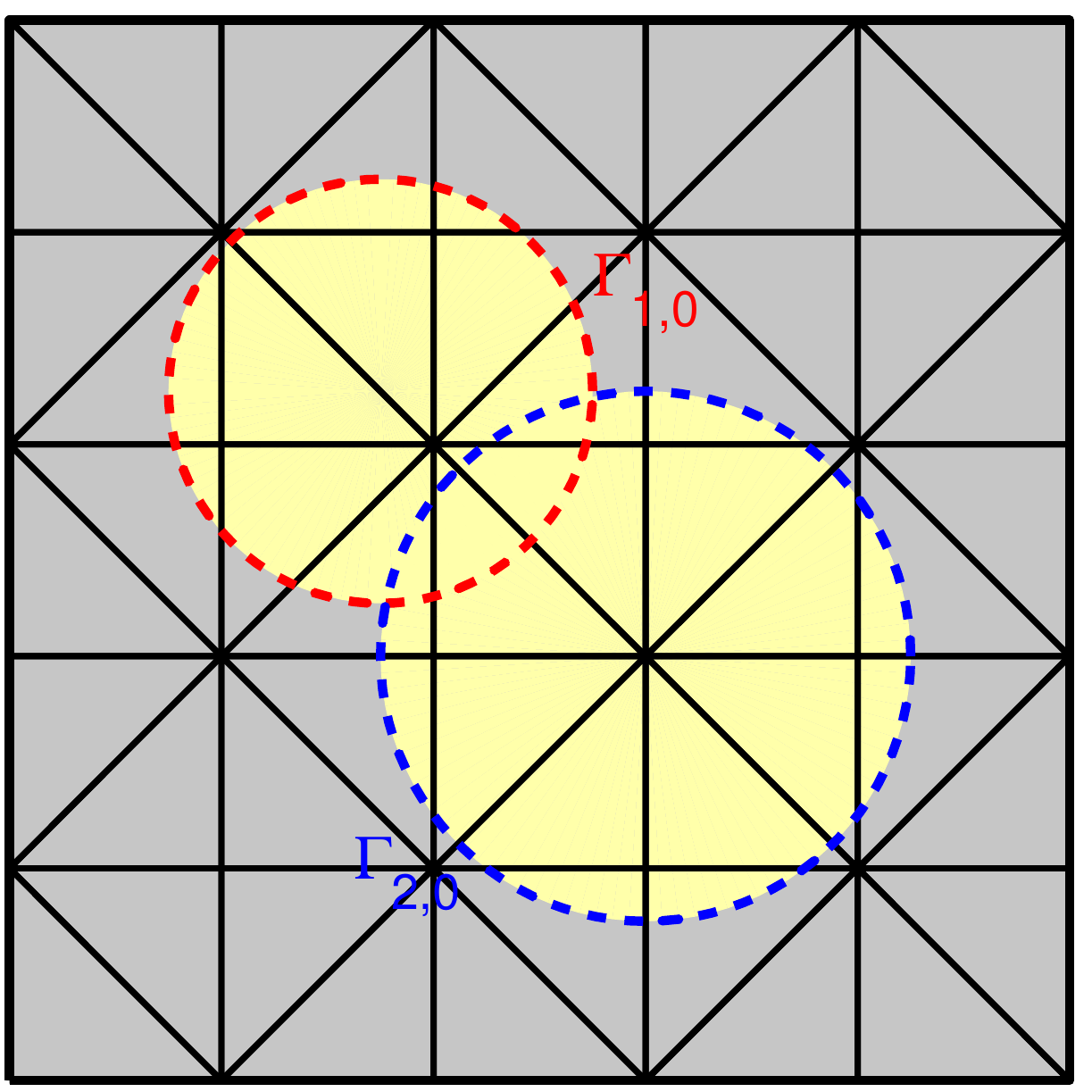}}\;\subfigure[Decomposed background mesh]{\includegraphics[width=4cm]{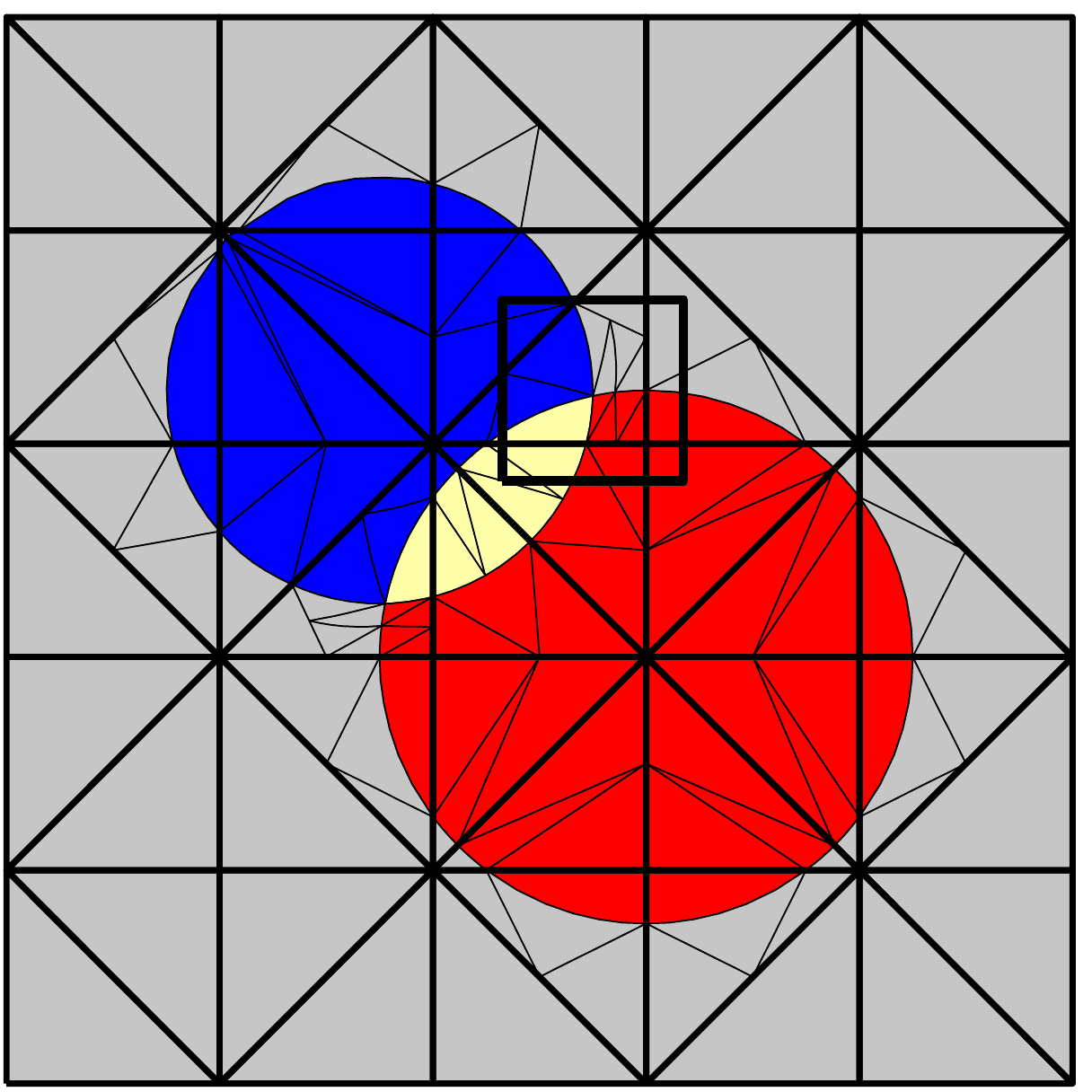}}\;\subfigure[Detail]{\includegraphics[width=4cm]{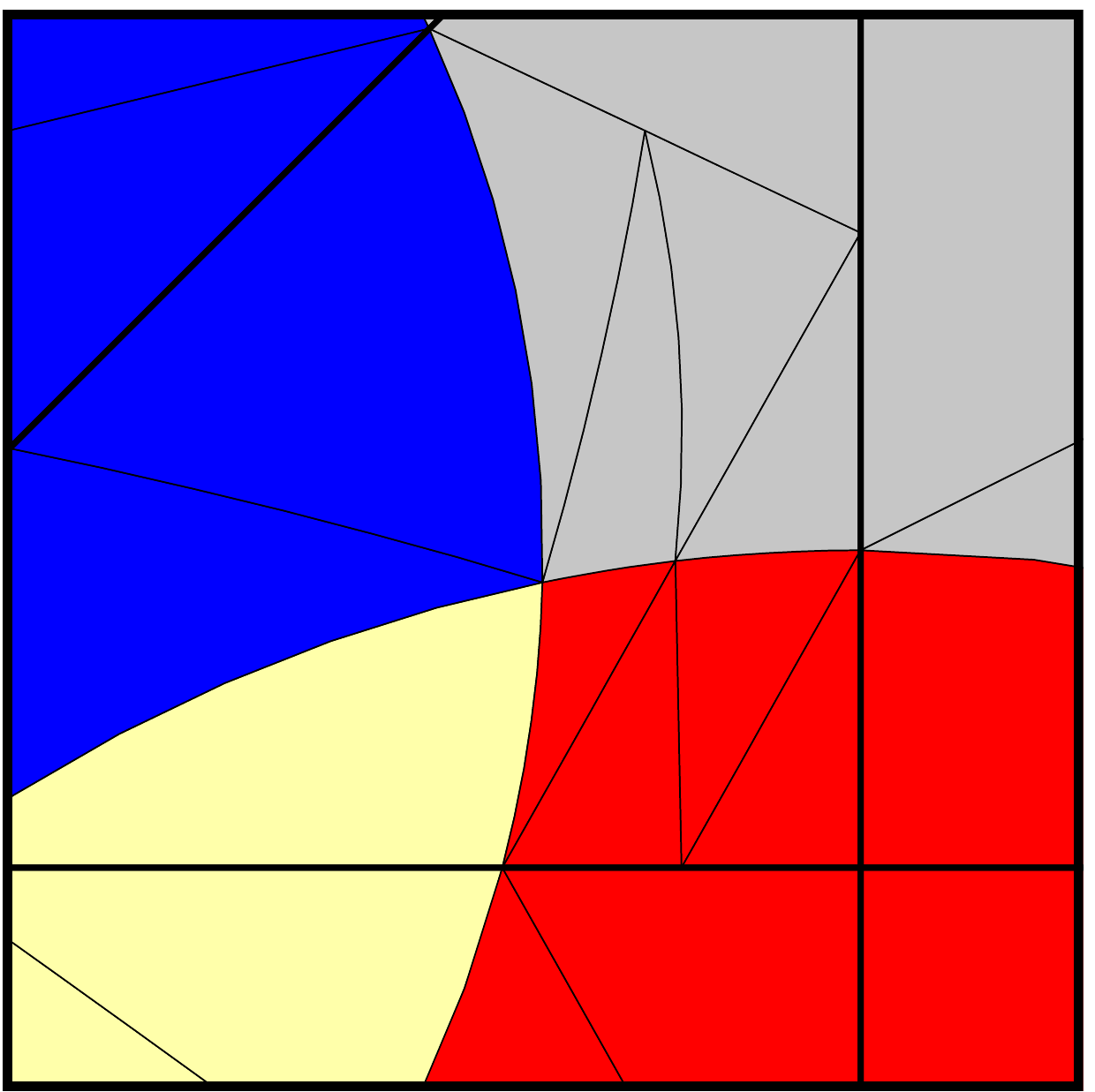}}

\caption{\label{fig:VisDomainCorners2d}In 2D, corners in the interfaces may
be defined by two level-set functions.}
\end{figure}

Therefore, we suggest to use \emph{multiple} level-set functions for
the domain description. Then, several level-set values are present
at the nodes of the background mesh, one for each smooth part of the
interface. In two dimensions, \emph{two} level-set functions define
a corner. See Fig.~\ref{fig:VisDomainCorners2d} for an example.
As seen in Fig.~\ref{fig:VisDomainCorners2d}(a), the domain features
an interface/boundary with kinks, i.e.~corners. Fig.~\ref{fig:VisDomainCorners2d}(b)
shows the two zero-level sets $\Gamma_{1,0}$ and $\Gamma_{2,0}$
of the level-set functions $\phi_{1}$ and $\phi_{2}$, respectively.
Assume the level-set functions are negative inside the circular zero-level
sets. It is seen in Fig.~\ref{fig:VisDomainCorners2d}(c) how the
two level-set functions imply 4 different sub-domains based on their
signs. The gray region in Fig.~\ref{fig:VisDomainCorners2d}(a) is
$\Omega^{\star}=\left\{ \vek x:\phi_{1}\left(\vek x\right)>0\;\mathrm{and}\;\phi_{2}\left(\vek x\right)>0\right\} $.
It is thus seen, how two level-set functions are able to define inner-element
corners in 2D. No additional requirements on the background mesh are
needed to capture the corners accurately. 

\begin{figure}
\centering

\subfigure[3 zero-level sets]{\includegraphics[width=5cm]{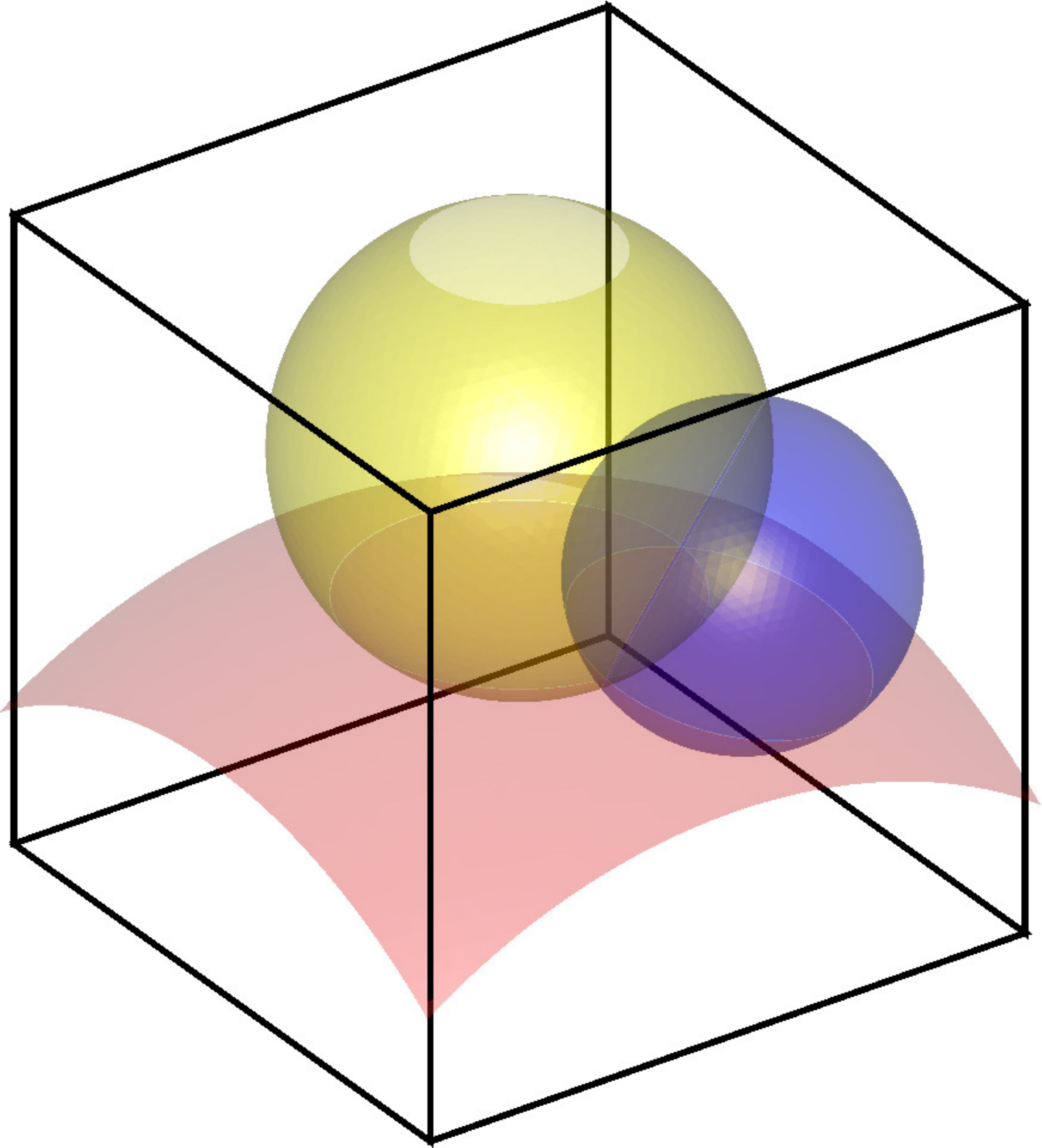}}\quad\subfigure[Example sub-region 1]{\includegraphics[width=5cm]{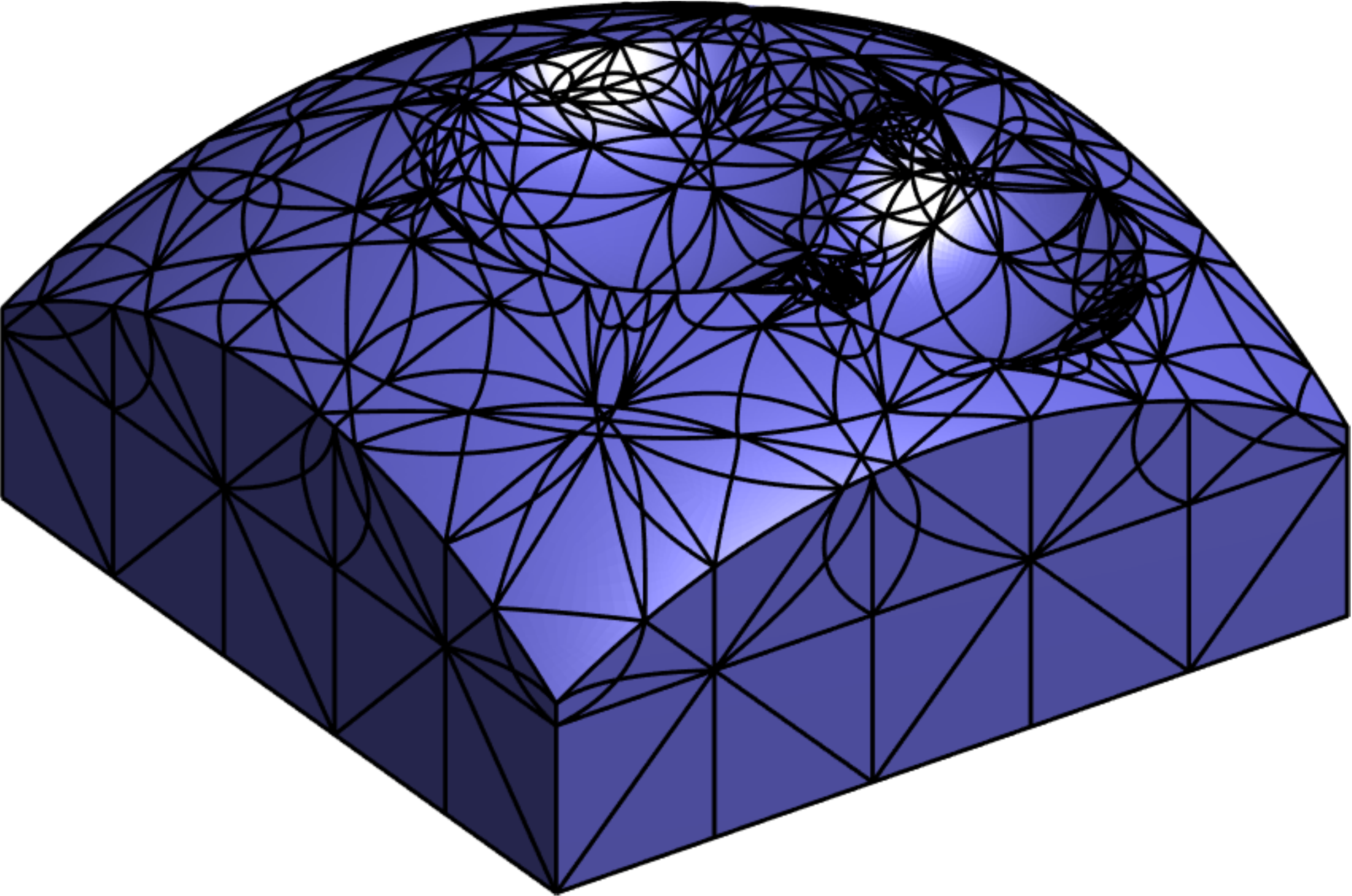}}\quad\subfigure[Example sub-region 2]{\includegraphics[width=5cm]{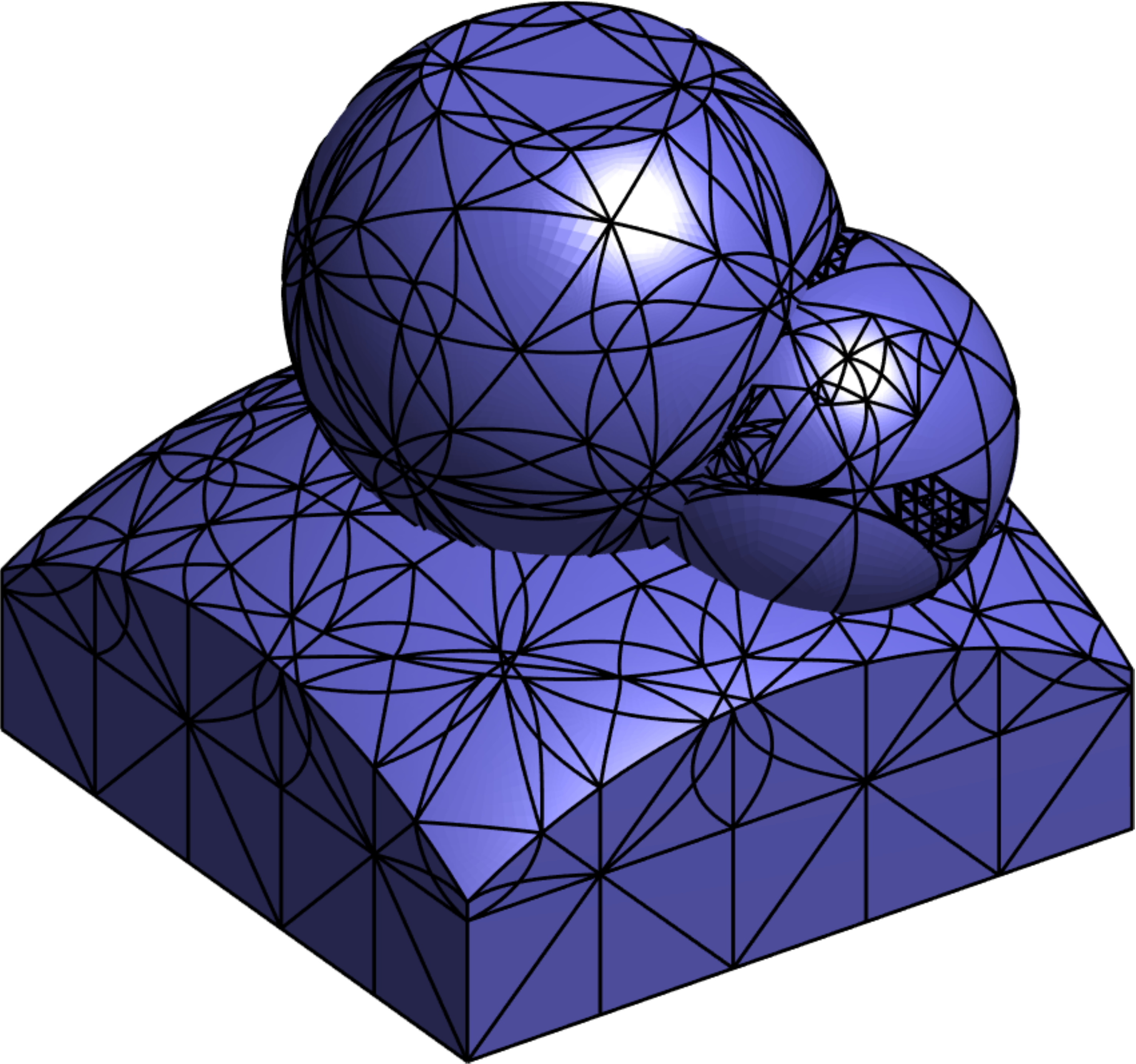}}

\caption{\label{fig:VisDomainCorners3d}Some example domain with implicit corners
and edges implied by three level-set functions.}
\end{figure}

In three dimensions, two level-set functions are able to define an
edge and three functions a corner. See Fig.~\ref{fig:VisDomainCorners3d}
for a graphical representation. Three zero-level sets of $\phi_{1}$,
$\phi_{2}$, and $\phi_{3}$ are shown in Fig.~\ref{fig:VisDomainCorners3d}(a).
It is seen how the intersections of any two level-set functions implicitly
define edges and how three level-set functions imply corners. Of course,
based on the signs of the three level-set functions, one may distinguish
$2^{3}=8$ sub-regions of the background mesh. Depending on the application,
one may now want to compose the domain based on these sub-regions.
Two different examples are shown in Figs.~\ref{fig:VisDomainCorners3d}(b)
and (c); both feature implicit edges and corners.

It is now clear how multiple level-set functions are able to define
edges and corners. It remains to describe how the corresponding reconstructions
and decompositions are realized. The task is to individually mesh
the sub-regions based on the signs of the involved level-set functions.
The resulting meshes are already visualized in Figs.~\ref{fig:VisDomainCorners2d}(c)
and (d) in 2D and Figs.~\ref{fig:VisDomainCorners3d}(b) and (c)
in 3D. Thereafter, it is simple to assemble the domains of interest
(e.g.~for integration, interpolation, or approximation) from these
sub-regions.

\begin{figure}
\centering

\subfigure[Zero level-sets]{\includegraphics[width=4cm]{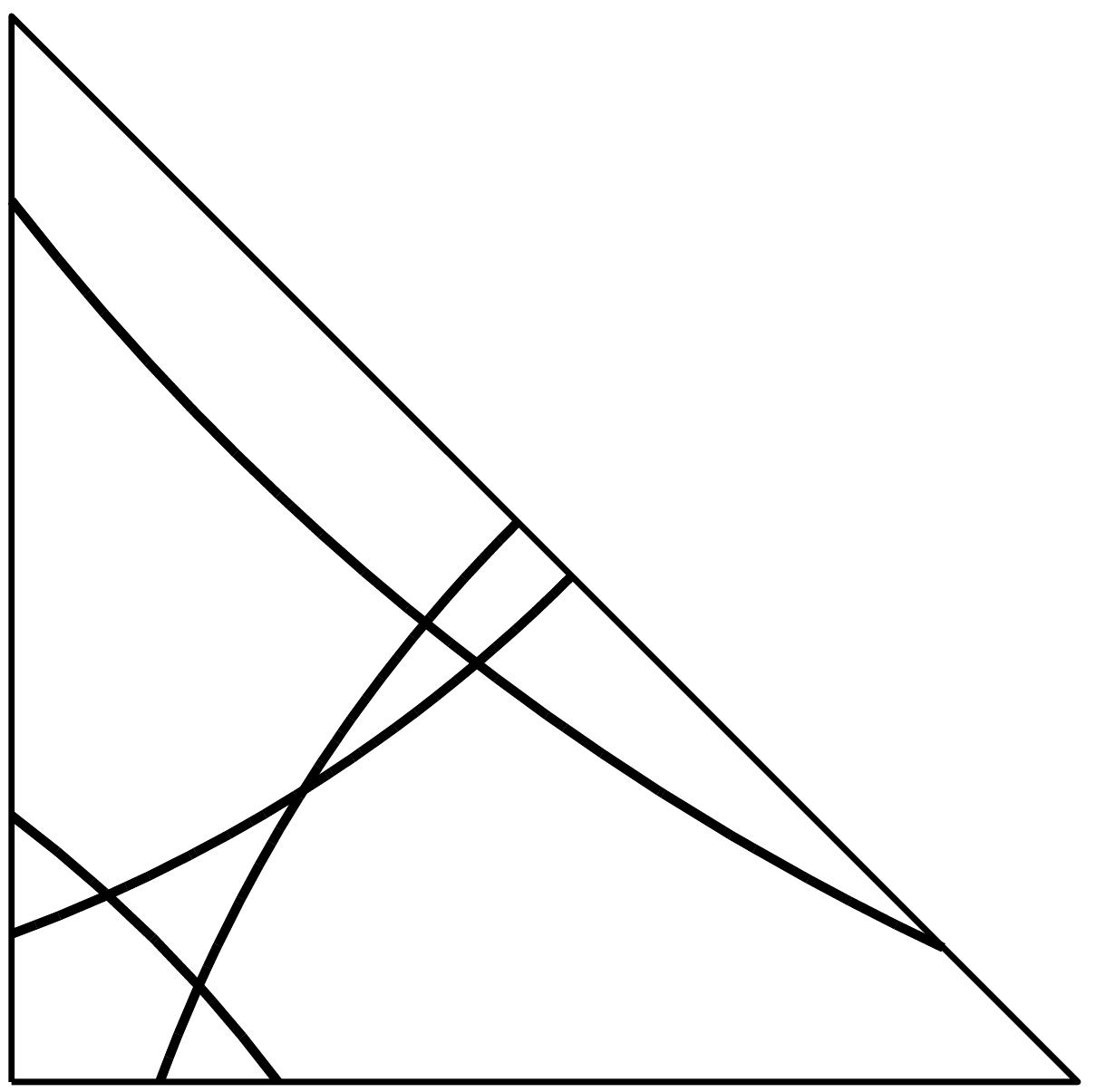}}\quad\subfigure[Decomp. for $\phi_1$]{\includegraphics[width=4cm]{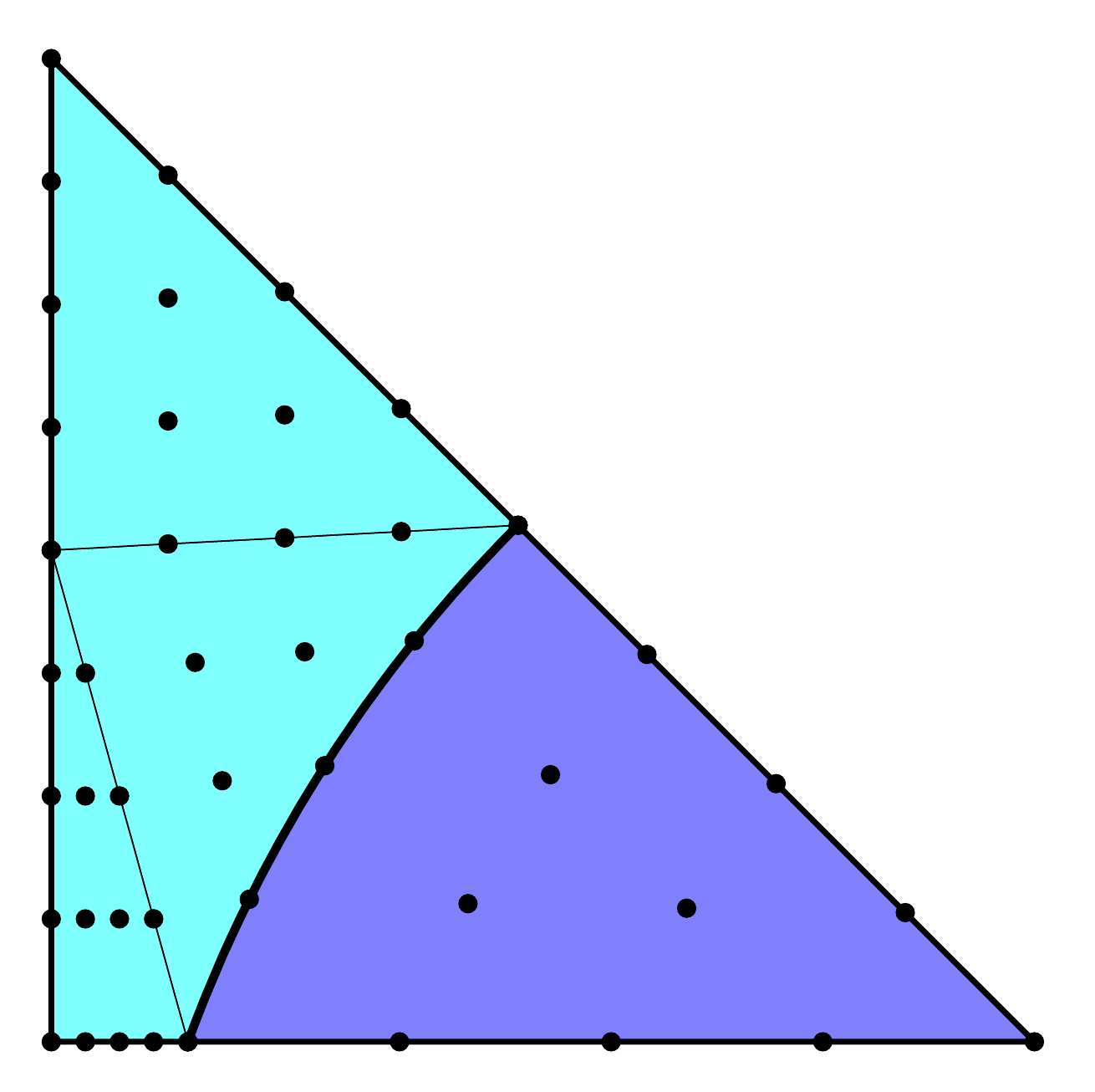}}\quad\subfigure[Decomp. for $\phi_1$ and $\phi_2$]{\includegraphics[width=4cm]{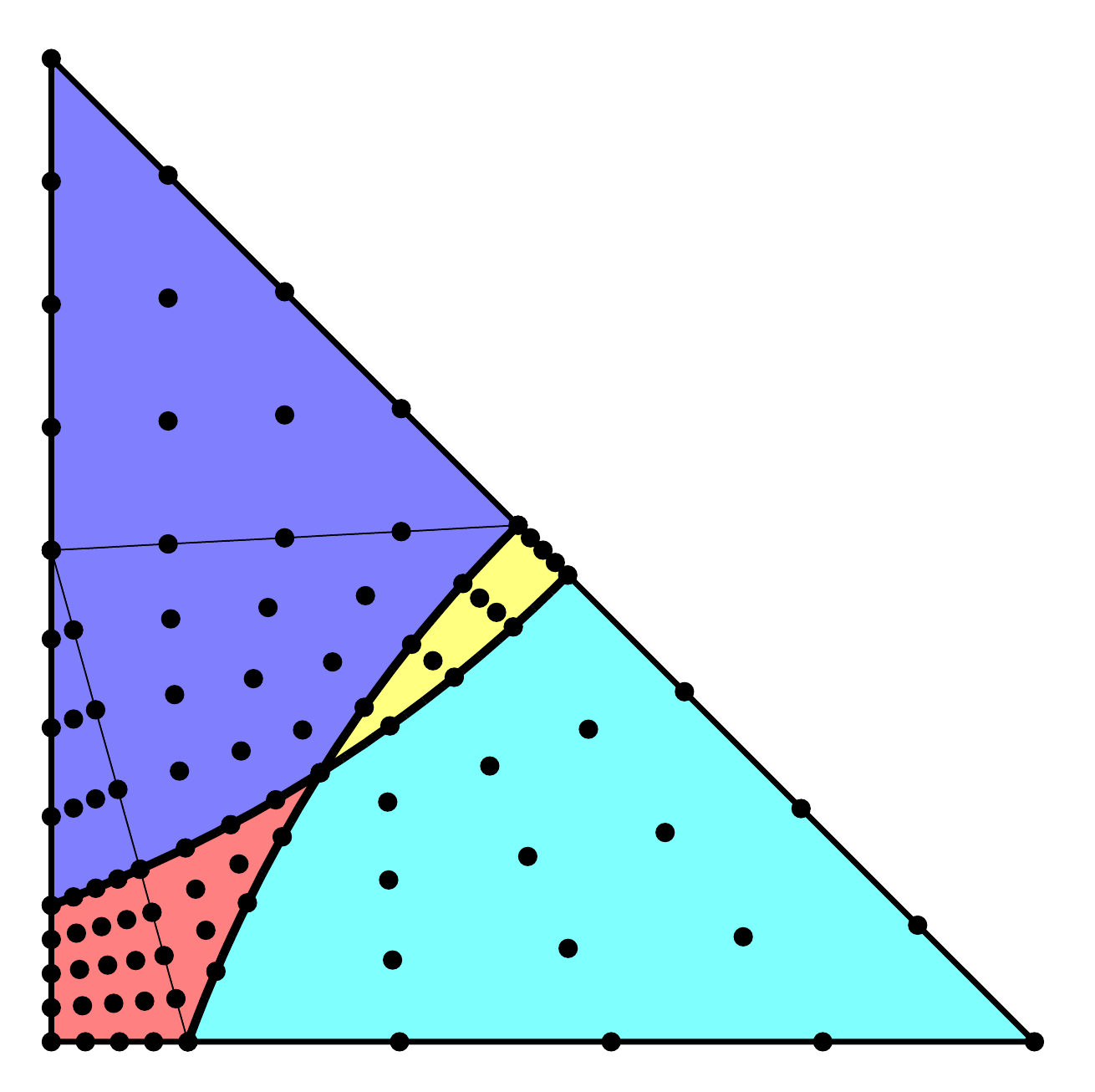}}

\subfigure[Decomp. for $\phi_1$ to $\phi_3$]{\includegraphics[width=4cm]{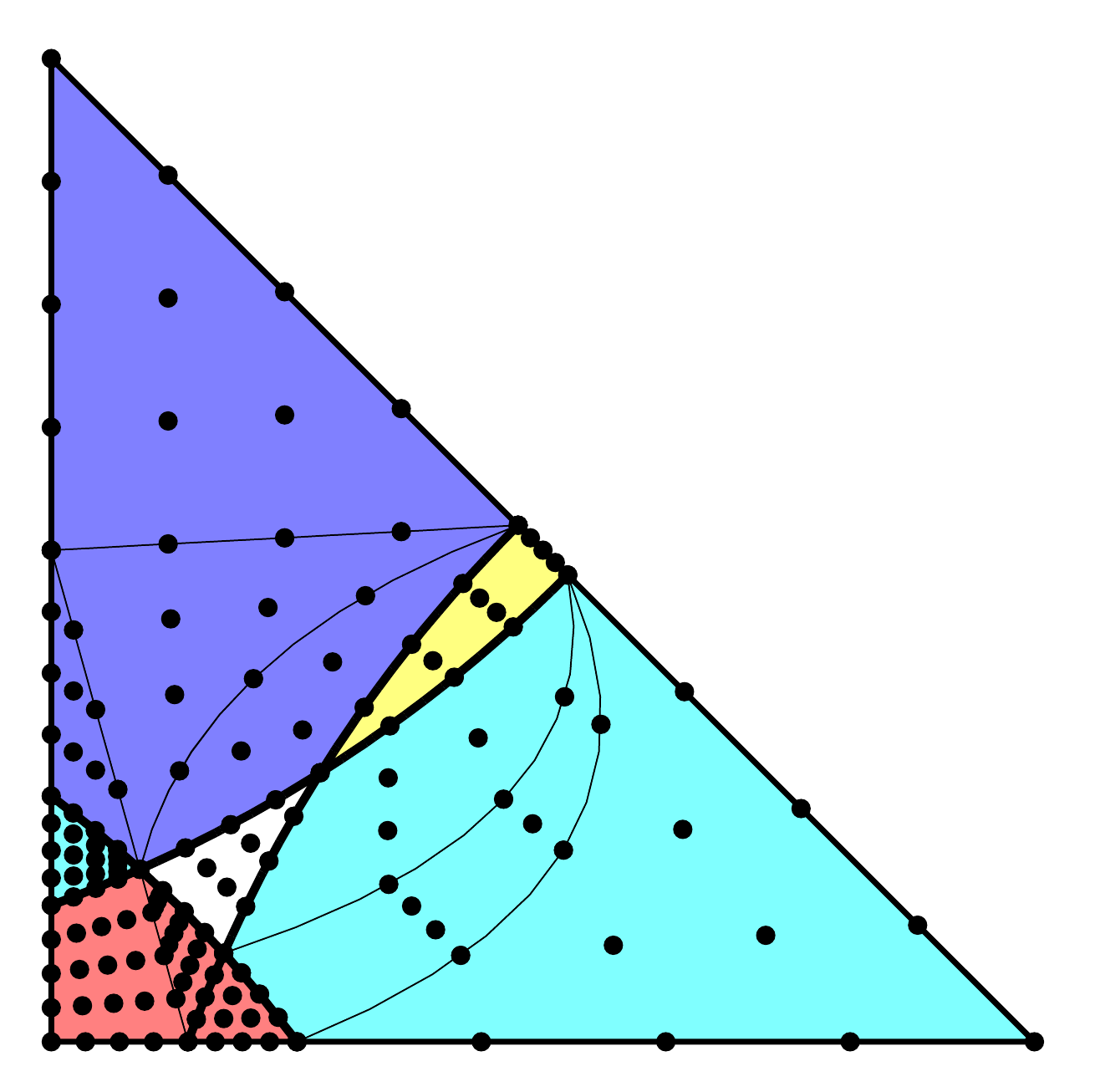}}\quad\subfigure[Decomp. for $\phi_1$ to $\phi_4$]{\includegraphics[width=4cm]{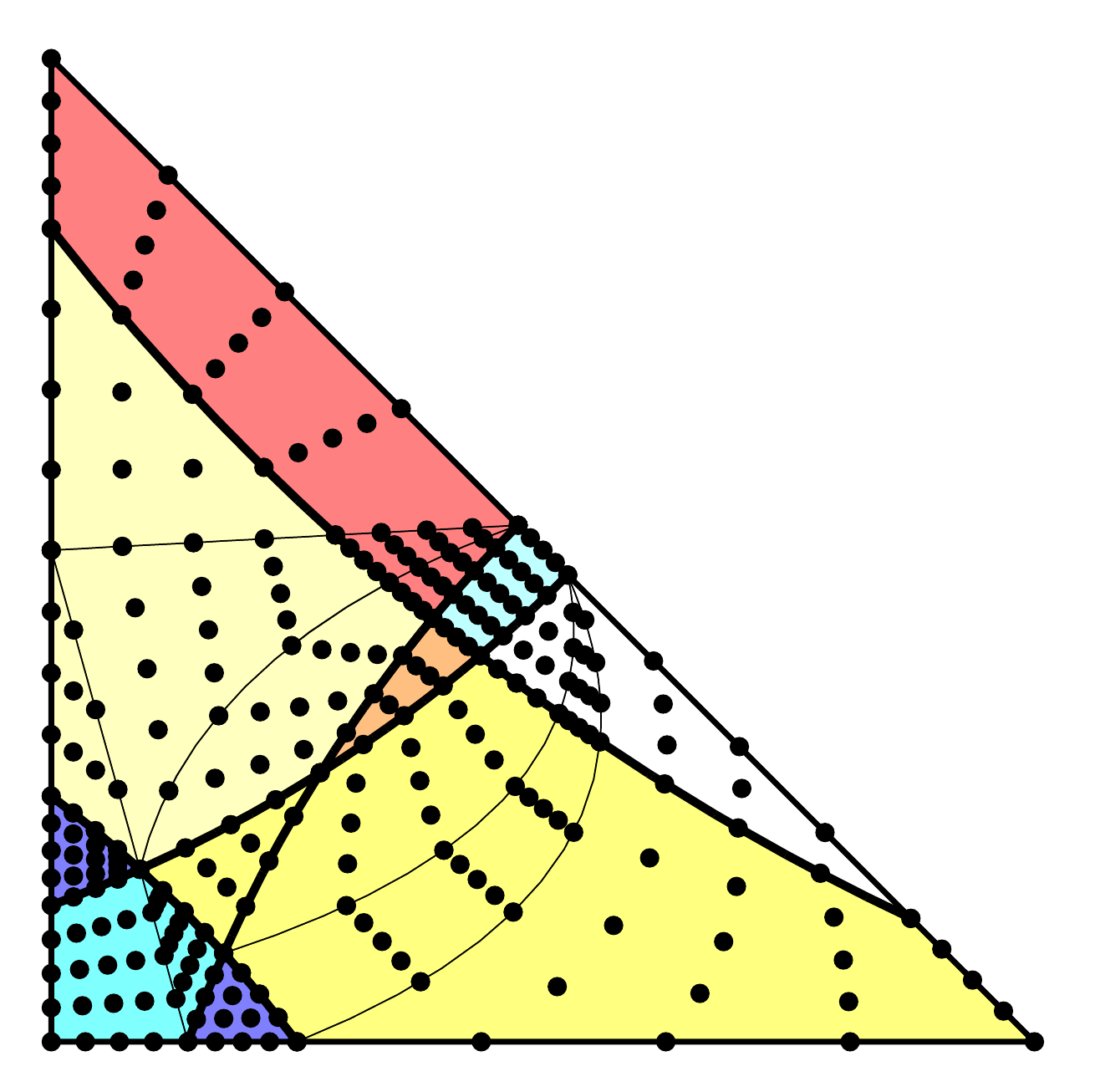}}\quad\subfigure[Sub-domain]{\includegraphics[width=4cm]{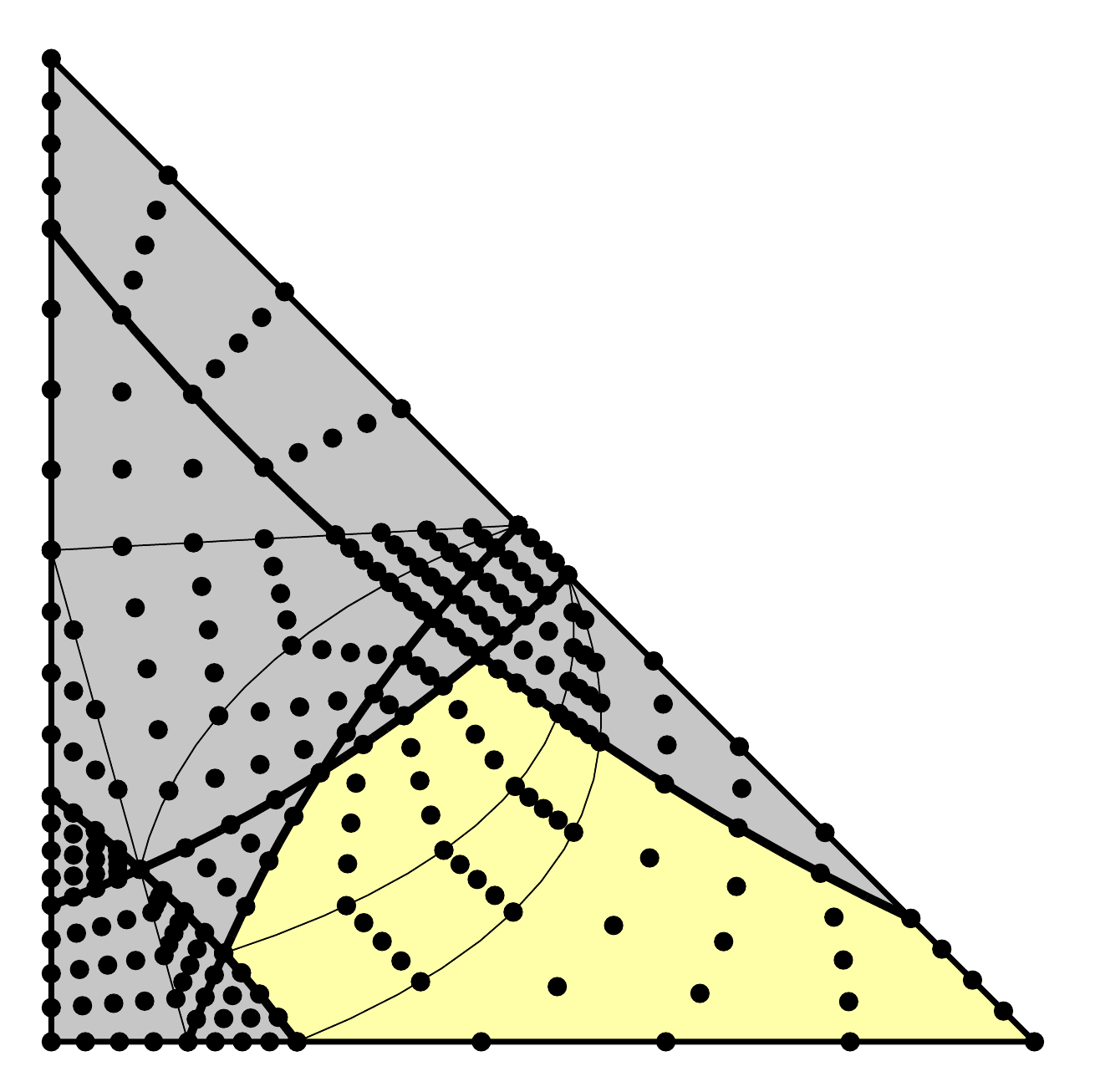}}

\caption{\label{fig:VisMultLevelSet2d}Decomposition of a reference triangular
element with respect to $4$ level-set functions: (a) Shows the zero-level
sets, (b) to (e) show the successive decomposition, (f) shows some
selected sub-domain defined by the signs of the level-sets. This sub-domain
has 3 corners within the reference element.}
\end{figure}

\begin{figure}
\centering

\subfigure[Zero level-sets]{\includegraphics[width=4cm]{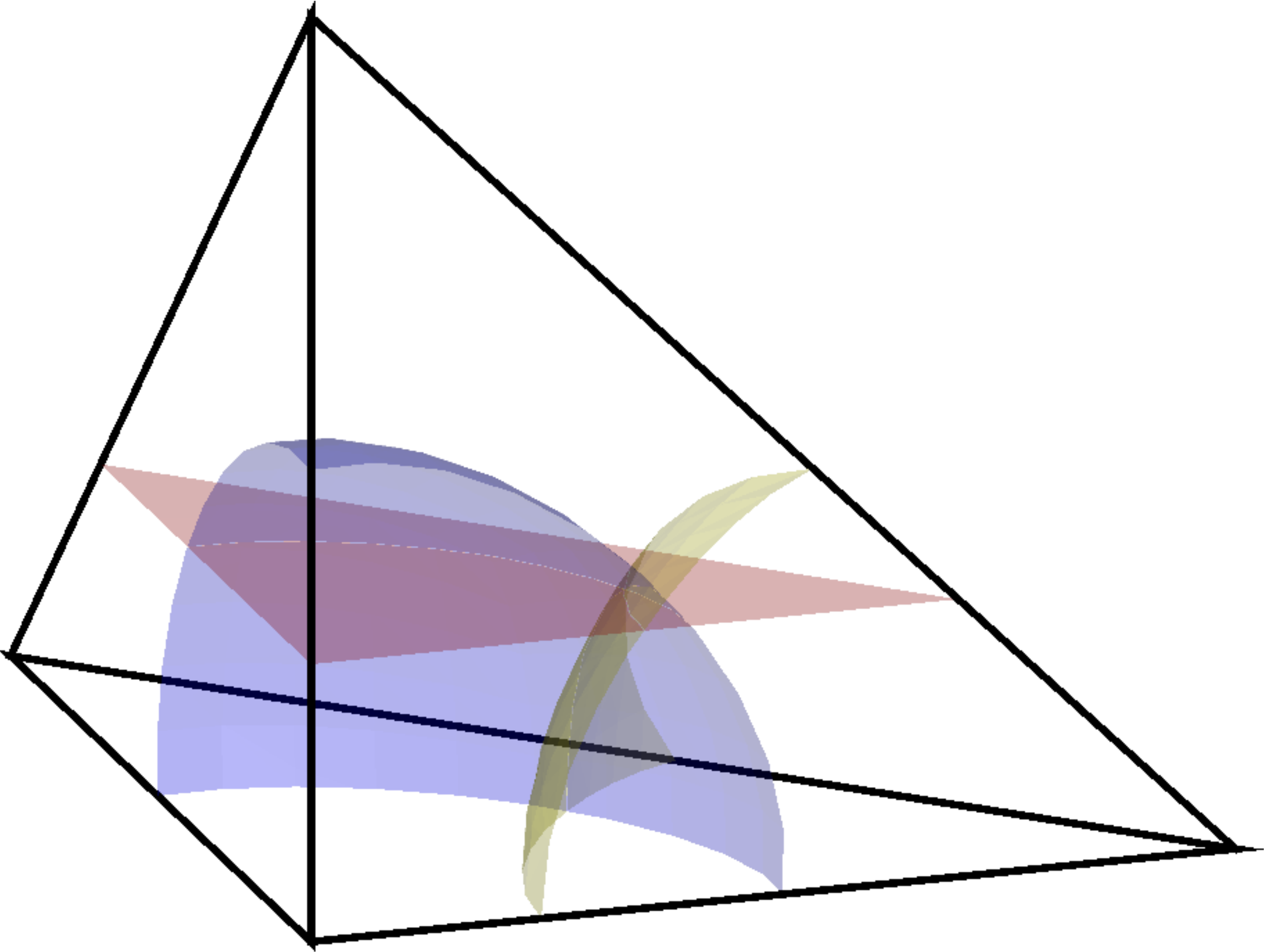}}\quad\subfigure[Decomp. for $\phi_1$]{\includegraphics[width=4cm]{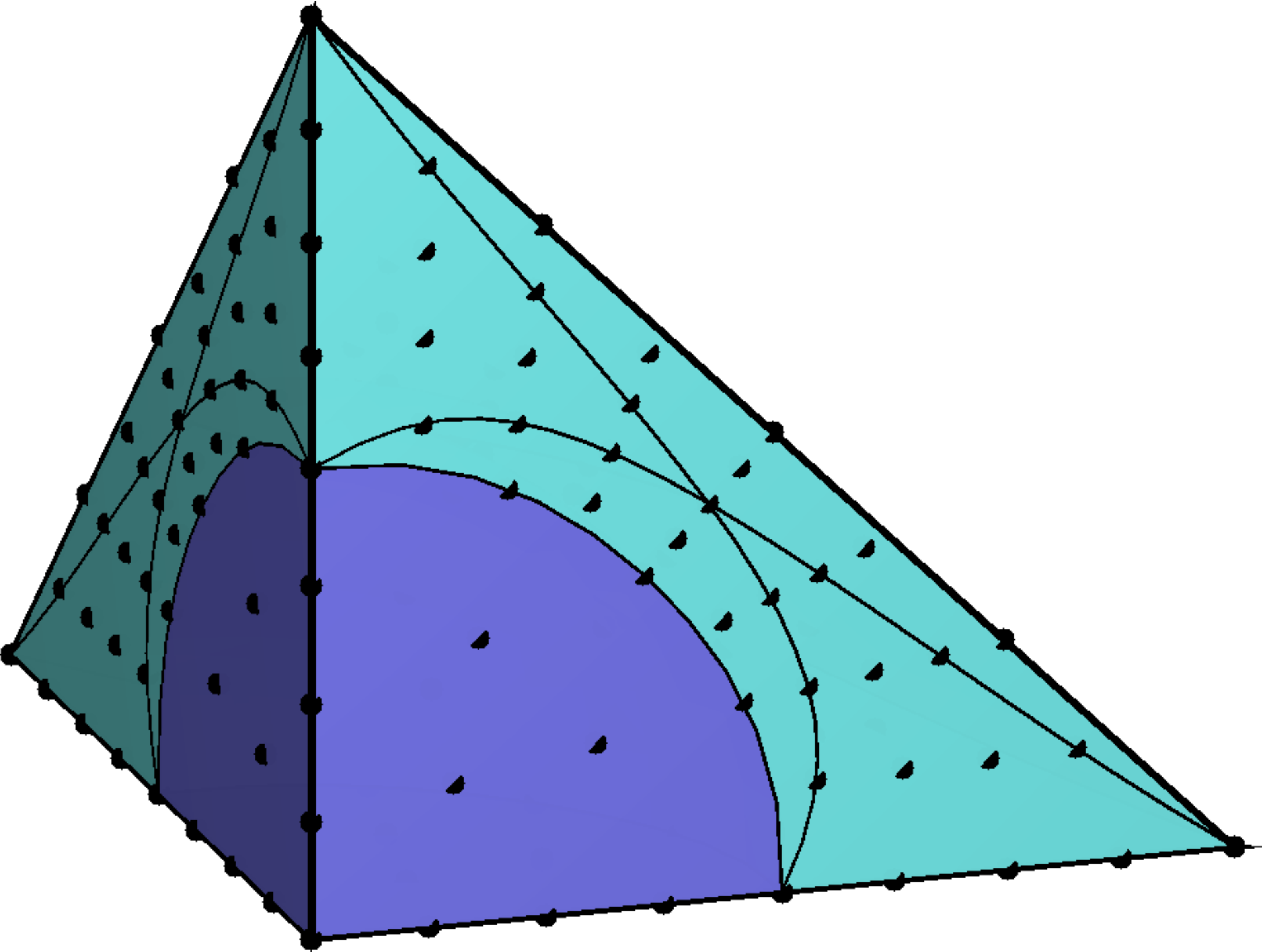}}\quad\subfigure[Decomp. for $\phi_1$ and $\phi_2$]{\includegraphics[width=4cm]{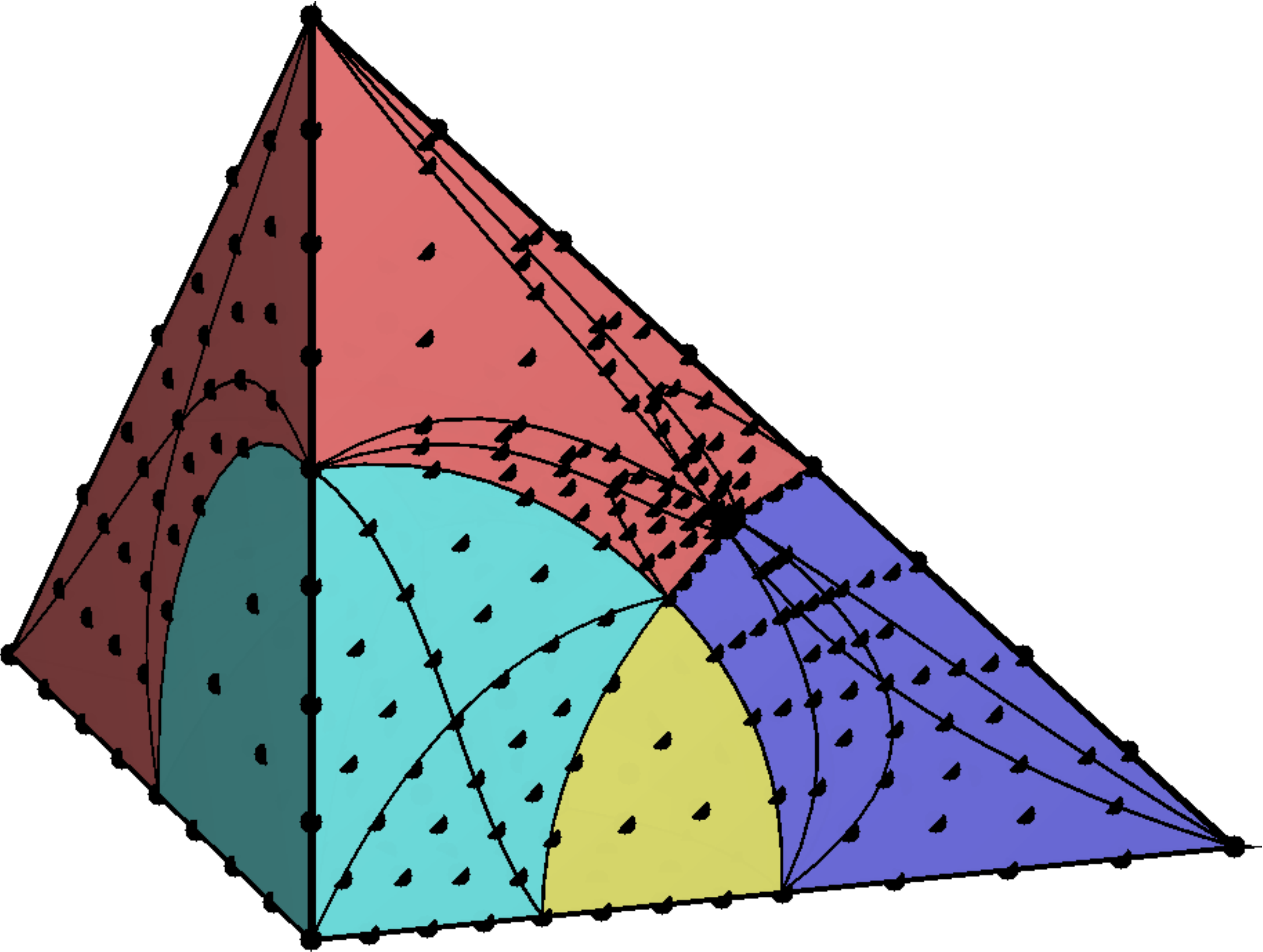}}

\subfigure[Decomp. for $\phi_1$ to $\phi_3$]{\includegraphics[width=4cm]{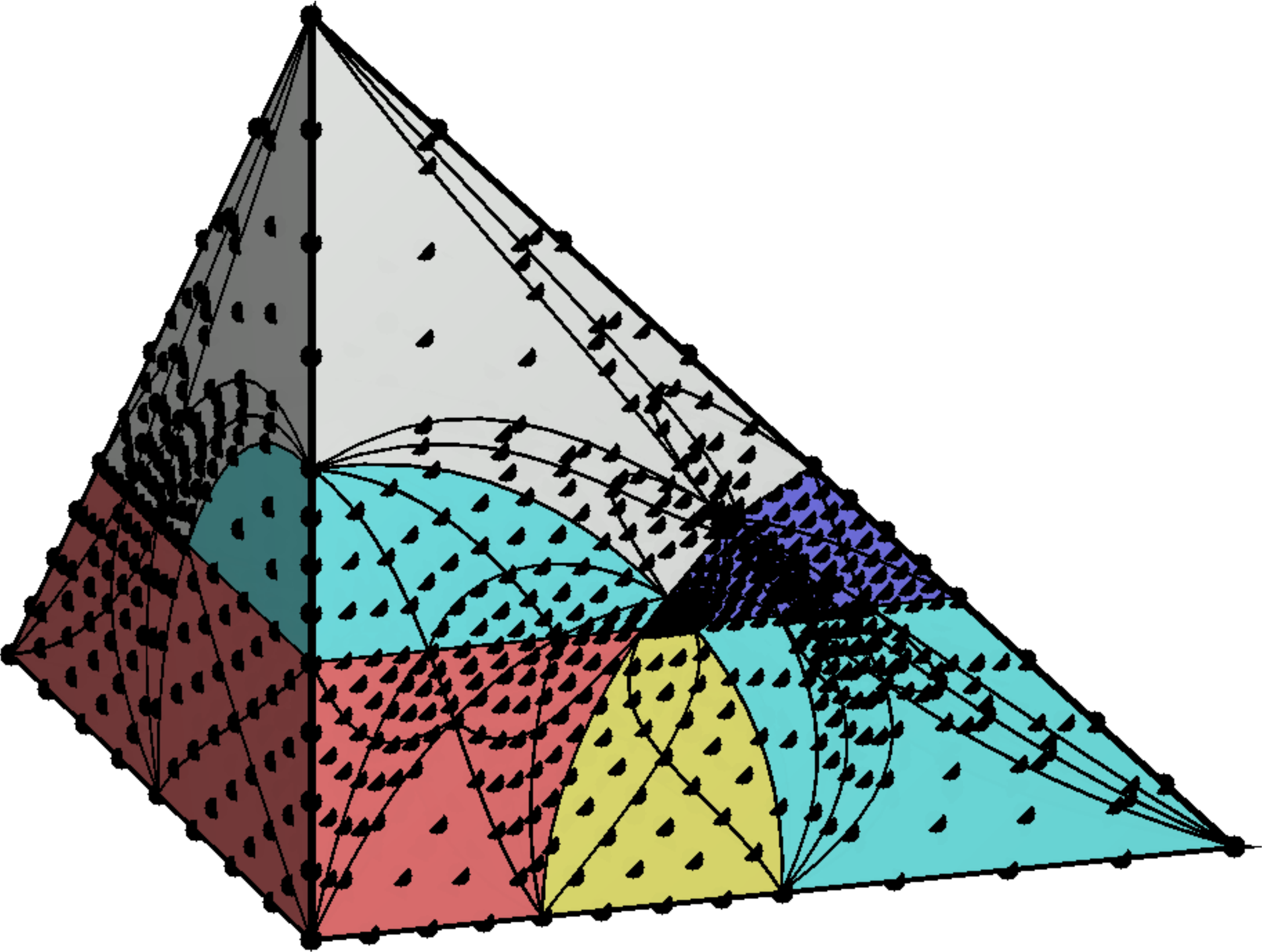}}\quad\subfigure[Selected sub-domain]{\includegraphics[width=4cm]{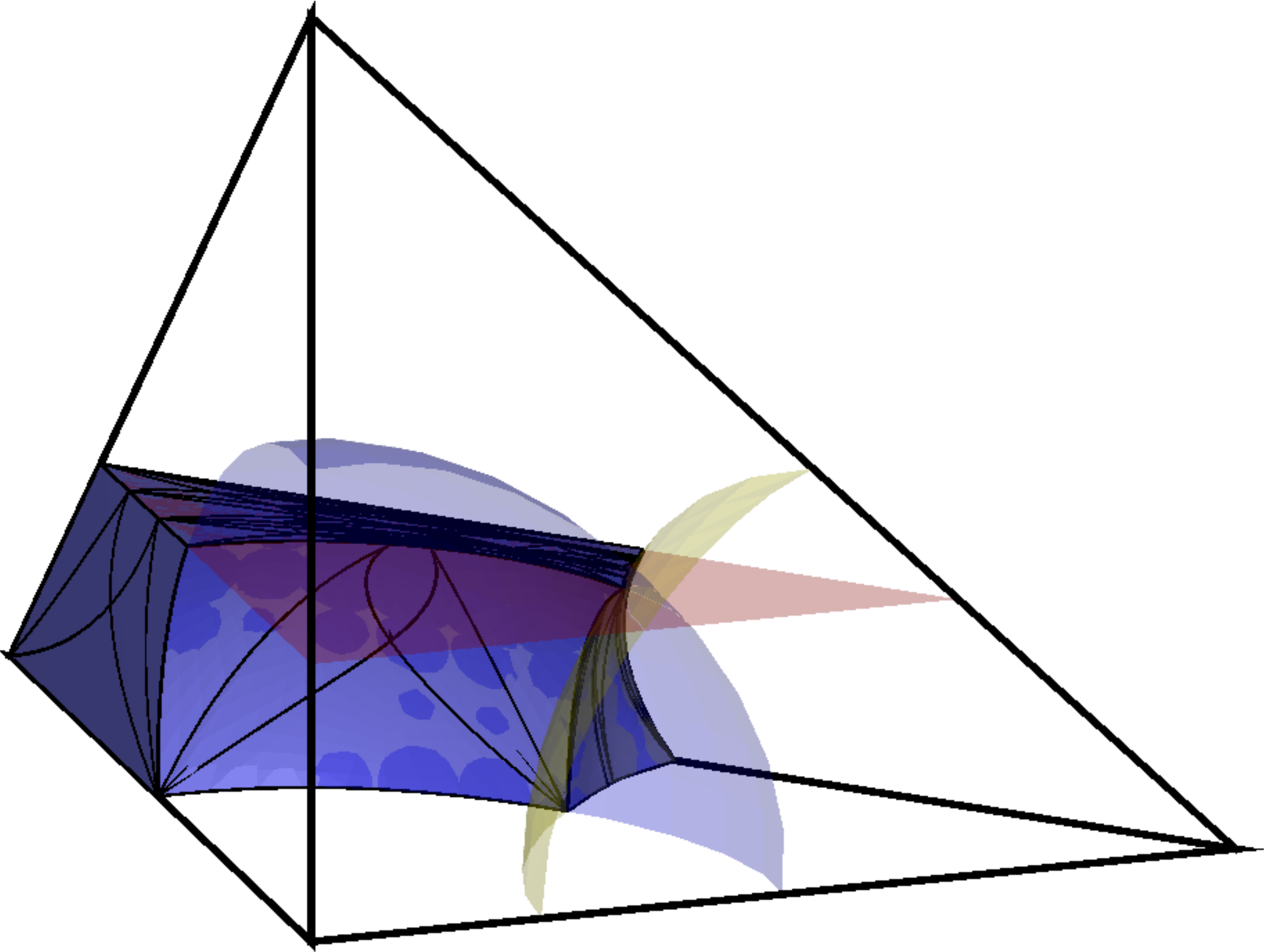}}\quad\subfigure[Selected sub-domain]{\includegraphics[width=4cm]{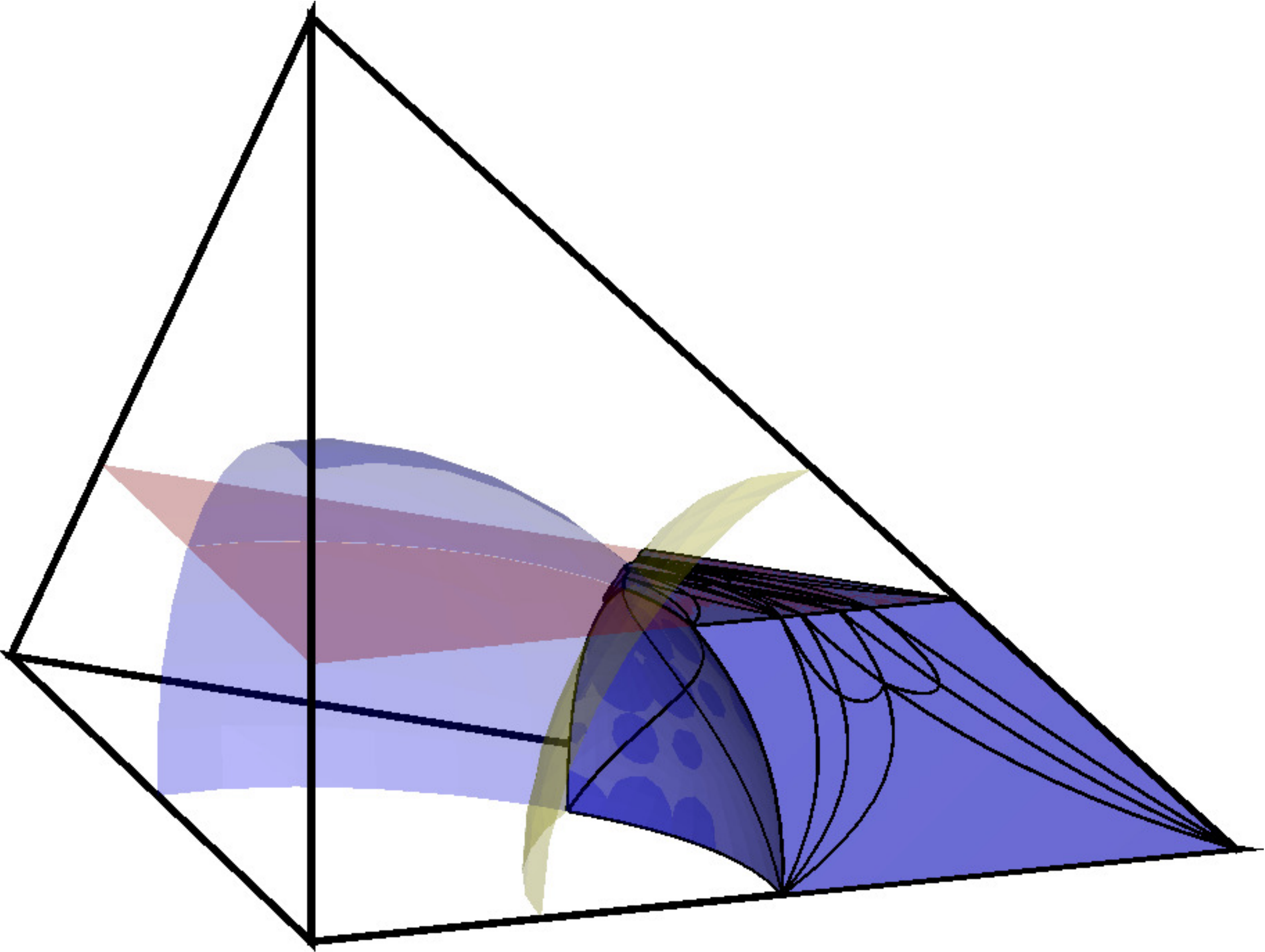}}

\caption{\label{fig:VisMultLevelSet3d}Decomposition of a reference tetrahedral
element with respect to $3$ level-set functions: (a) Shows the zero-level
sets, (b) to (d) show the successive decomposition, (e) and (f) show
some selected sub-domains defined by the signs of the level-sets.
Note the corner inside the reference element.}
\end{figure}

Fortunately, the decomposition of background elements with respect
to several level-set functions is straightforward. The key aspect
is that the decomposition is always carried out in \emph{reference}
background elements. That is, a background element is first decomposed
with respect to the first level-set function $\phi_{1}$. Then, the
remaining level-set functions are interpolated at the new element
nodes of the sub-elements of the cut background element. The decomposition
with respect to the next level-set function is carried out just as
if the already decomposed mesh were the initial mesh. Therefore, it
may be useful to further decompose quadrilateral sub-elements in 2D
into triangles and prismatic sub-elements in 3D into tetrahedra, so
that the newly generated mesh (conforming to the zero-level sets of
previous $\phi_{i}$) only consists of the same element type. Fig.~\ref{fig:VisMultLevelSet2d}
shows an example where a reference element is decomposed with respect
to $4$ level set functions. The resulting decompositions for each
new level-set functions are shown and the different sub-regions that
may be identified by the signs of the level-set functions are color-coded.
Fig.~\ref{fig:VisMultLevelSet2d}(f) shows a sub-domain which features
3 corners within this one element. An example for the situation in
3D is shown in Fig.~\ref{fig:VisMultLevelSet3d}.

Hence, the decomposition with respect to several level-set functions
is a truly straightforward successive application of the meshing procedures
described above. The ability to naturally consider an arbitrary number
of level-set functions, and even the presence of multiple edges and
corners in \emph{one} element seems to be a unique feature of the
proposed method.

\section{Conclusions\label{sec:Conclusions}}

Implicit geometries occur in many applications, in particular in the
context of fictitious domain methods and XFEM-related methods. The
boundaries of a domain of interest or interfaces inside the domain
are then easily defined based on the level-set method. Using a higher-oder
background mesh enables a straightforward path to a higher-order accurate
implicit description. However, during an analysis, integration points
are required, typically, for the integration of the weak form of a
model and accurate quadrature rules in cut background elements become
a critical issue. Herein, a consistent path to obtain higher-order
accurate integration points in cut background elements is described.
There are virtually no requirements on the background mesh, in particular,
there is no need to use Cartesian meshes or elements with straight
edges.

The key idea is to generate higher-order finite elements on the two
sides of an interface which is called decomposition. Before, the zero-level
sets have to be identified and approximated by interface elements
which is called reconstruction. In fact, the reconstruction and decomposition
may be seen as two meshing steps. A number of different variants to
realize the meshing have been studied. It is found that one has to
carefully define the start values and search directions used for the
non-linear detection of the zero-level sets. For the decomposition,
mappings of element nodes into the special sub-elements with one higher-order
side are needed and they have to be sufficiently smooth. It is noted
that depending on the level-set data inside the elements one may have
to recursively refine background elements until valid level-set data
is obtained. Recursive refinements are also needed when the (initial)
reconstruction or decomposition fails e.g.~due to negative Jacobians
of the resulting sub-elements. When the curvatures of the level-set
functions are reasonably adjusted with the resolution of the background
elements, recursive refinements are only needed in a very small number
of elements. The number of generated sub-elements per background element
is then small, thus also the number of additional integration points
in the cut elements.

The proposed meshing procedures will be used in upcoming publications
to approximate boundary values problems without using typical FDM
or XFEM-related methods. In fact, because the generated elements align
with boundares and interfaces, the mesh is used in a classical FEM
context. Major issues are then the quality of the automatically generated
elements which are not necessarily well-shaped. Nevertheless, it shall
be seen that higher-order convergence rates are possible with such
meshes not only in pure integration and interpolation problems as
shown here but also in the approximation of boundary value problems.

The implementation of the proposed meshing procedures is not without
efforts, especially in three dimensions. Therefore, we plan to soon
enable a download of the underlying software on the institute's webpage
at \texttt{www.ifb.tugraz.at}.

\section{Appendix}

Crucial ingredients of the proposed remeshing are the maps of element
nodes (or start values) to the element interiors based on the situation
on the outer contour. For example, for the sub-elements on the two
sides of the zero-level sets, the outer contour is given by linear
elements with one higher-order, typically curved side. For the placements
of start values for the reconstruction in tetrahedra, three or four
higher-order edges of the sought interface element are prescribed.
In the following, the maps needed herein are defined in detail for
triangles, quadrilaterals, tetrahedra and prisms. Because the definitions
of these maps are quite lengthy, we decided to move this to the appendix.
A more general assessment of general transfinite mappings is given
e.g.~in \cite{Solin_2003a}.

\subsection{Mappings for triangles and quadrilaterals\label{sub:MappingTriQuad}}

\begin{figure}
\centering

\subfigure[]{\includegraphics[width=4cm]{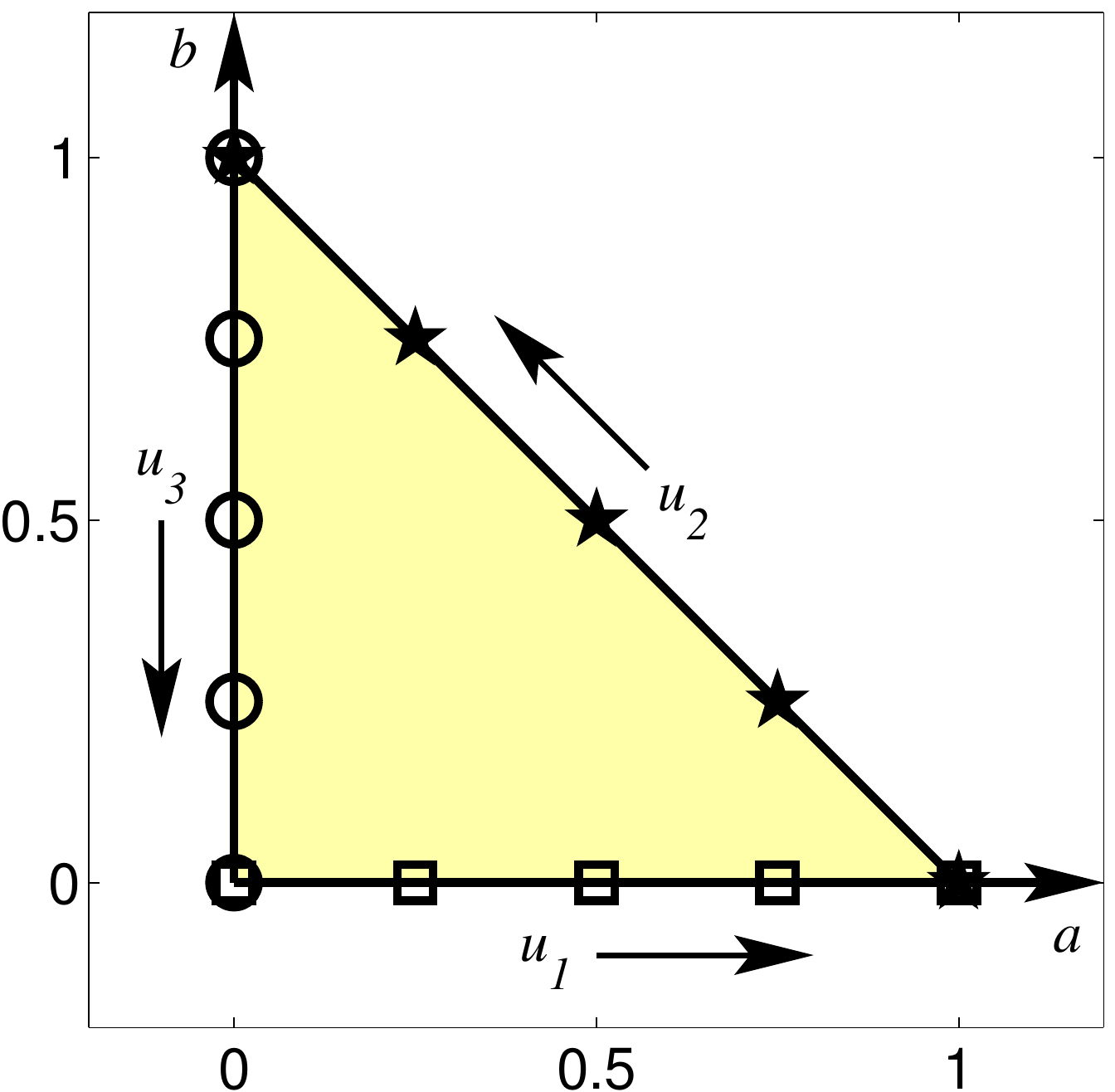}}\subfigure[]{\includegraphics[width=4cm]{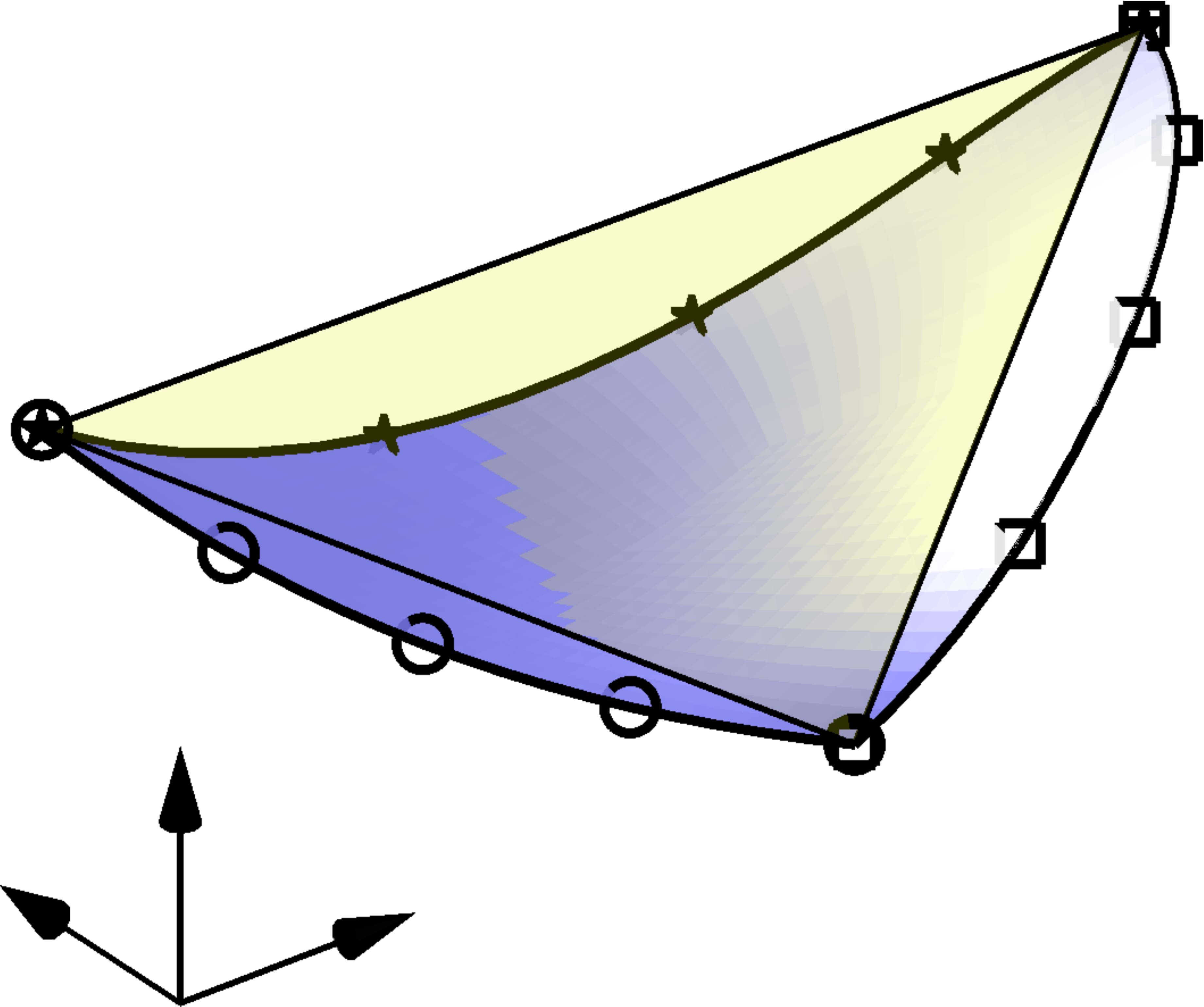}}\subfigure[]{\includegraphics[width=4cm]{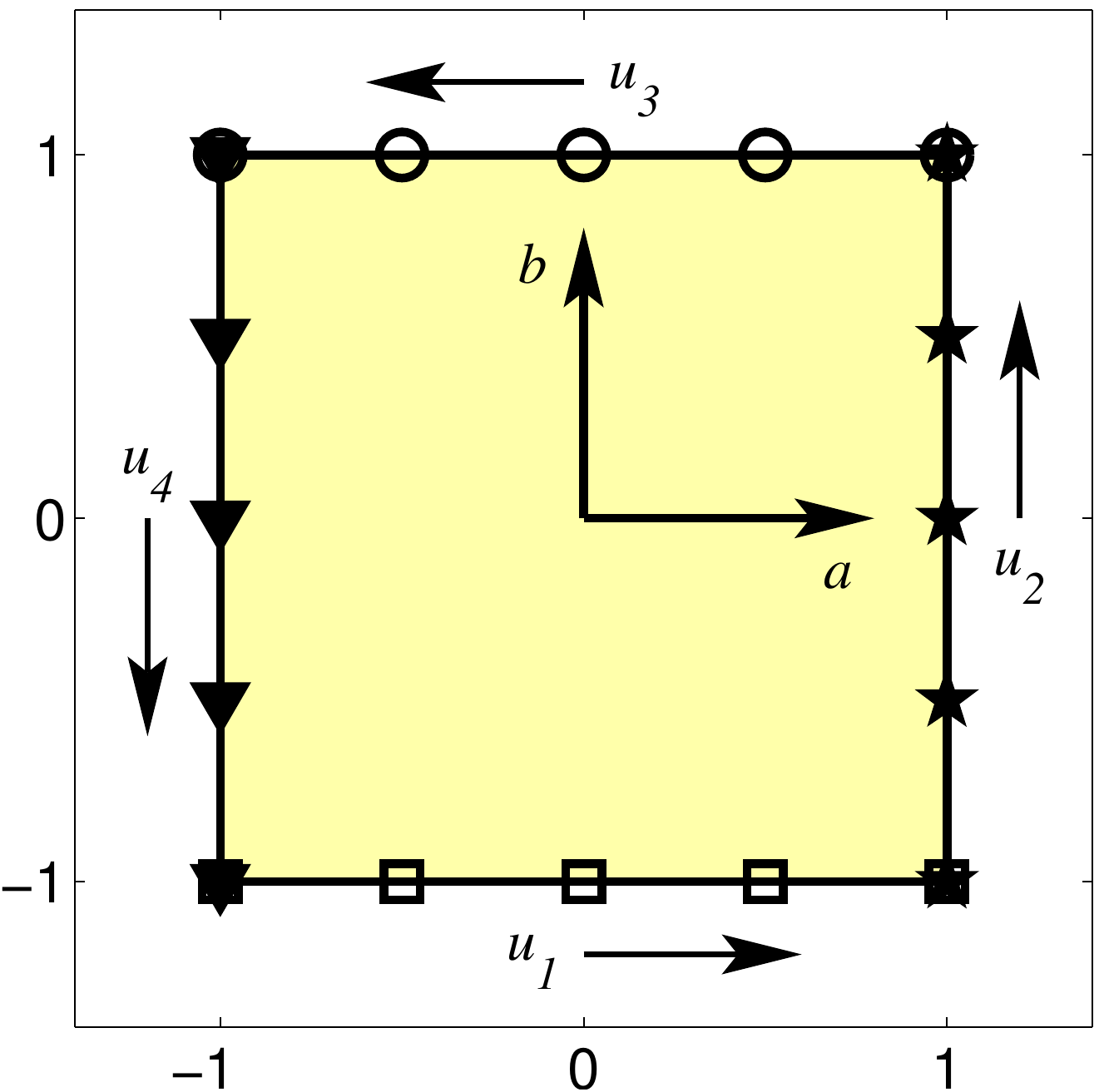}}\subfigure[]{\includegraphics[width=4cm]{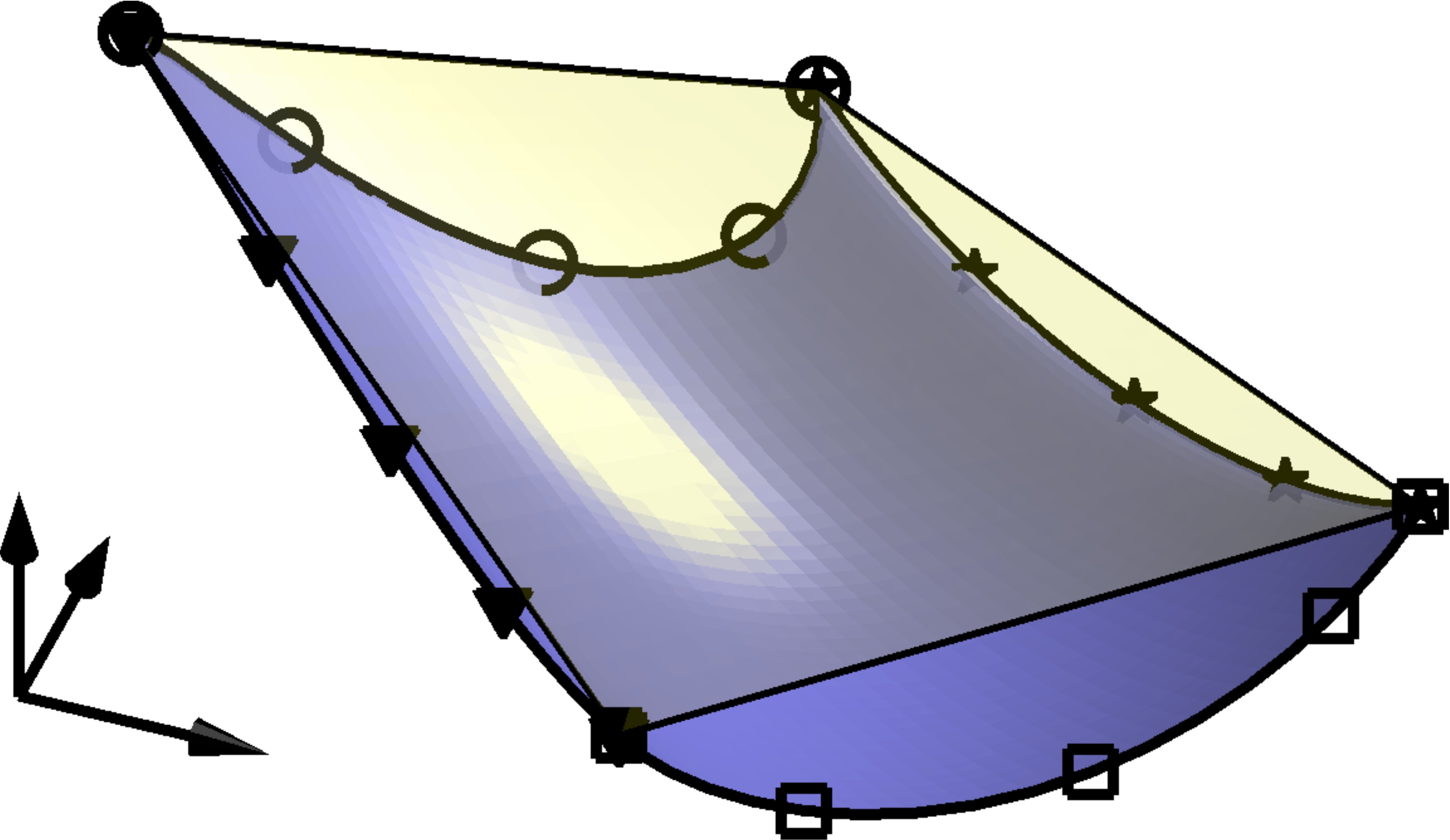}}

\caption{\label{fig:VisRecon3dSolin}Map from a reference triangle or quadrilateral
in $\left(a,b\right)$ to the coordinate system $\left(r,s,t\right)$.
The edges are defined by higher-order line elements and the mapping
of the inner points is sought.}
\end{figure}

Assume that a triangular (or quadrilateral) element is implied by
$n_{e}=3$ (or $n_{e}=4$) curved line elements defining the outer
contour, i.e.~the edges. The line elements are of the same order,
i.e.~the same number of nodes $n_{p}$ and corresponding shape functions
is associated with each edge. The task is to define a smooth mapping
for all points $\left(a,b\right)$ in a reference triangle (or quadrilateral)
to the three-dimensional situation in $\left(r,s,t\right)$. Obviously,
this defines the shape of the triangular (or quadrilateral) surface
in 3D, see Fig.~\ref{fig:VisRecon3dSolin}. Such a mapping is not
unique and the definition outlined below follows the textbook \cite{Solin_2003a}
adapted to the present situation.

The nodal coordinates of the line elements are denoted as $\vek r_{i}^{\mathrm{edge}\, k}\in\mathbb{R}^{3}$
with $i=1,\dots,n_{p}$ and $k=1,\dots,n_{e}$. They must build a
closed contour, i.e.~the start and end coordinates must match properly.
The (bi-)linear shape functions related to the reference triangle
(or quadrilateral) are called $N_{k}^{\star}\left(\vek a\right)$,
more precisely, for triangles
\[
N_{1}^{\star}\left(\vek a\right)=1-a-b,\quad N_{2}^{\star}\left(\vek a\right)=a,\quad N_{3}^{\star}\left(\vek a\right)=b
\]
and for quadrilaterals
\begin{eqnarray*}
N_{1}^{\star}\left(\vek a\right) & = & \nicefrac{1}{4}\left(1-a\right)\left(1-b\right),\; N_{2}^{\star}\left(\vek a\right)=\nicefrac{1}{4}\left(1+a\right)\left(1-b\right),\\
N_{3}^{\star}\left(\vek a\right) & = & \nicefrac{1}{4}\left(1+a\right)\left(1+b\right),\; N_{4}^{\star}\left(\vek a\right)=\nicefrac{1}{4}\left(1-a\right)\left(1+b\right).
\end{eqnarray*}

The higher-order shape functions of the line elements are denoted
as $N_{i}^{\mathrm{HO}}\left(u\right)$ and are defined in a one-dimensional
reference element with coordinates $u\in[-1,1]$. Herein, these are
standard Lagrange shape functions. Furthermore, the linear shape functions
$N_{1}^{\mathrm{lin}}(u)=\nicefrac{1}{2}(1-u)$ and $N_{2}^{\mathrm{lin}}(u)=\nicefrac{1}{2}(1+u)$
are needed. Local coordinates $u_{k}\in\left[-1,1\right]$ along each
of the edges of the triangle or quadrilateral are introduced, see
Fig.~\ref{fig:VisRecon3dSolin}(a) and (c). They have to be related
to the coordinates $\left(a,b\right)$. For triangles, this gives
\[
u_{1}=2a-1,\quad u_{2}=b-a,\quad u_{3}=1-2b
\]
and for quadrilaterals
\[
u_{1}=a,\quad u_{2}=b,\quad u_{3}=-a,\quad u_{4}=-b.
\]
Next, we define
\begin{equation}
\vek r^{\mathrm{edge}\, k}\left(u_{k}\right)=\sum_{i=1}^{n_{p}}N_{i}^{\mathrm{HO}}\left(u_{k}\right)\cdot\vek r_{i}^{\mathrm{edge}\, k}-N_{1}^{\mathrm{lin}}\left(u_{k}\right)\cdot\vek r_{1}^{\mathrm{edge}\, k}-N_{2}^{\mathrm{lin}}\left(u_{k}\right)\cdot\vek r_{n_{p}}^{\mathrm{edge}\, k}\label{eq:EdgeContribution}
\end{equation}
where the first term on the right hand side defines the curve in $\mathbb{R}^{3}$,
and the other two subtract the linear interpolant. Furthermore, a
ramp function is needed which is, for triangles,
\begin{equation}
R_{1}=\frac{N_{1}^{\star}\left(\vek a\right)\cdot N_{2}^{\star}\left(\vek a\right)}{N_{1}^{\mathrm{lin}}(u_{1})\cdot N_{2}^{\mathrm{lin}}(u_{1})},\: R_{2}=\frac{N_{2}^{\star}\left(\vek a\right)\cdot N_{3}^{\star}\left(\vek a\right)}{N_{1}^{\mathrm{lin}}(u_{2})\cdot N_{2}^{\mathrm{lin}}(u_{2})},\: R_{3}=\frac{N_{3}^{\star}\left(\vek a\right)\cdot N_{1}^{\star}\left(\vek a\right)}{N_{1}^{\mathrm{lin}}(u_{3})\cdot N_{2}^{\mathrm{lin}}(u_{3})}\label{eq:RampTri}
\end{equation}
and, for quadrilaterals,
\begin{equation}
R_{1}=N_{1}^{\star}\left(\vek a\right)+N_{2}^{\star}\left(\vek a\right),\, R_{2}=N_{2}^{\star}+N_{3}^{\star},\, R_{3}=N_{3}^{\star}+N_{4}^{\star},\, R_{4}=N_{4}^{\star}+N_{1}^{\star}.\label{eq:RampQuad}
\end{equation}
The overall map $\vek r\left(\vek a\right)$, i.e.~the shape of the
resulting surface element which is implied by the curved line elements,
is then defined as
\[
\vek r\left(\vek a\right)=\sum_{i=1}^{n_{e}}N_{i}^{\star}\left(\vek a\right)\cdot\vek r_{i}^{\star}+\sum_{k=1}^{n_{e}}R_{k}\cdot\vek r^{\mathrm{edge}\, k}\left(u_{k}(\vek a)\right).
\]

The first term on the right hand side is the linear interpolant of
the corner points of the triangle (or quadrilateral) and $\vek r_{i}^{\star}$
is extracted from the start and end points of the curved line elements.
See Fig.~\ref{fig:VisRecon3dSolin} for a graphical representation.

\subsection{Mapping for tetrahedra\label{sub:MappingTetra}}

\begin{figure}
\centering

\subfigure[]{\includegraphics[height=5cm]{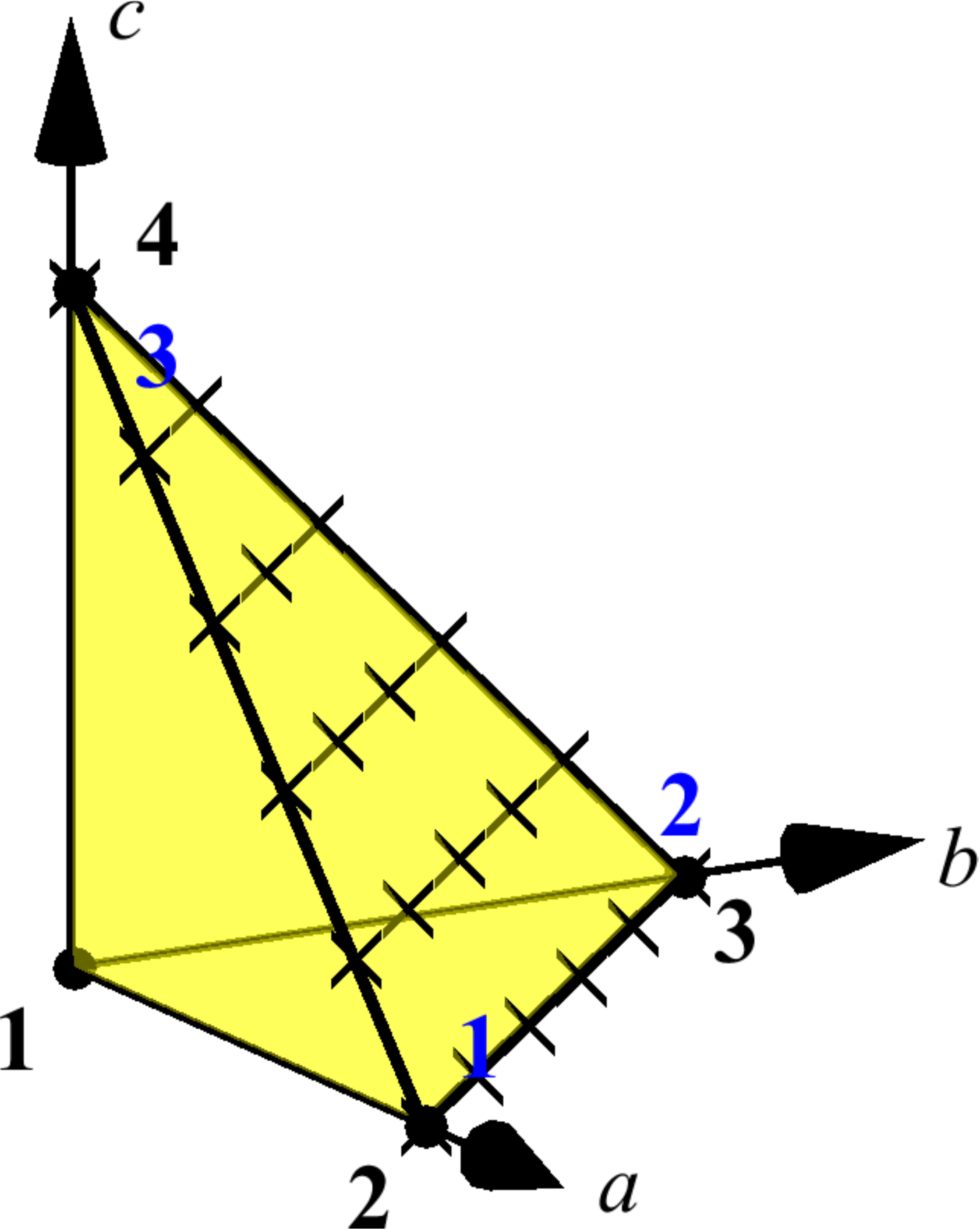}}\quad\subfigure[]{\includegraphics[height=5cm]{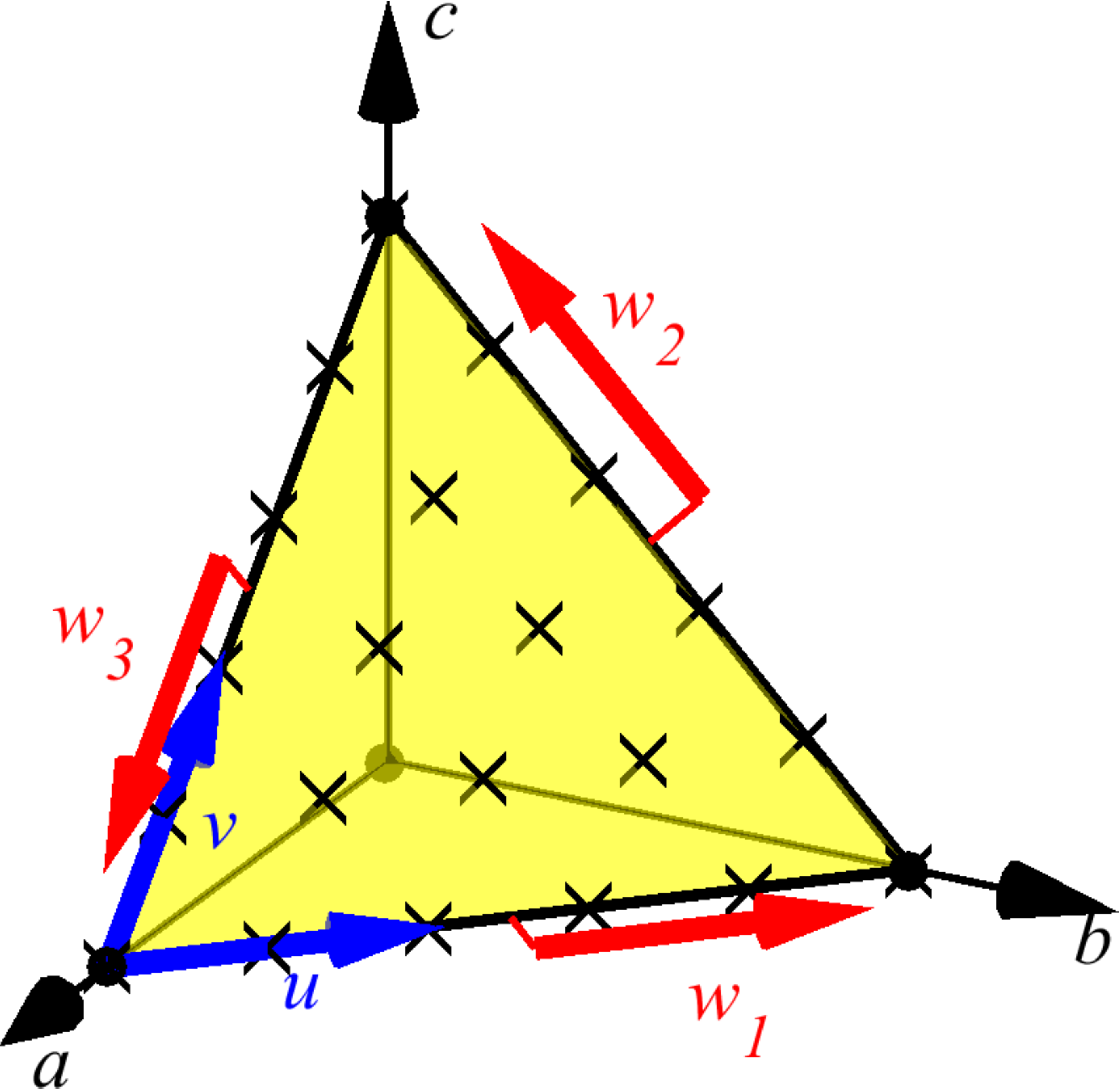}}\quad\subfigure[]{\includegraphics[height=5cm]{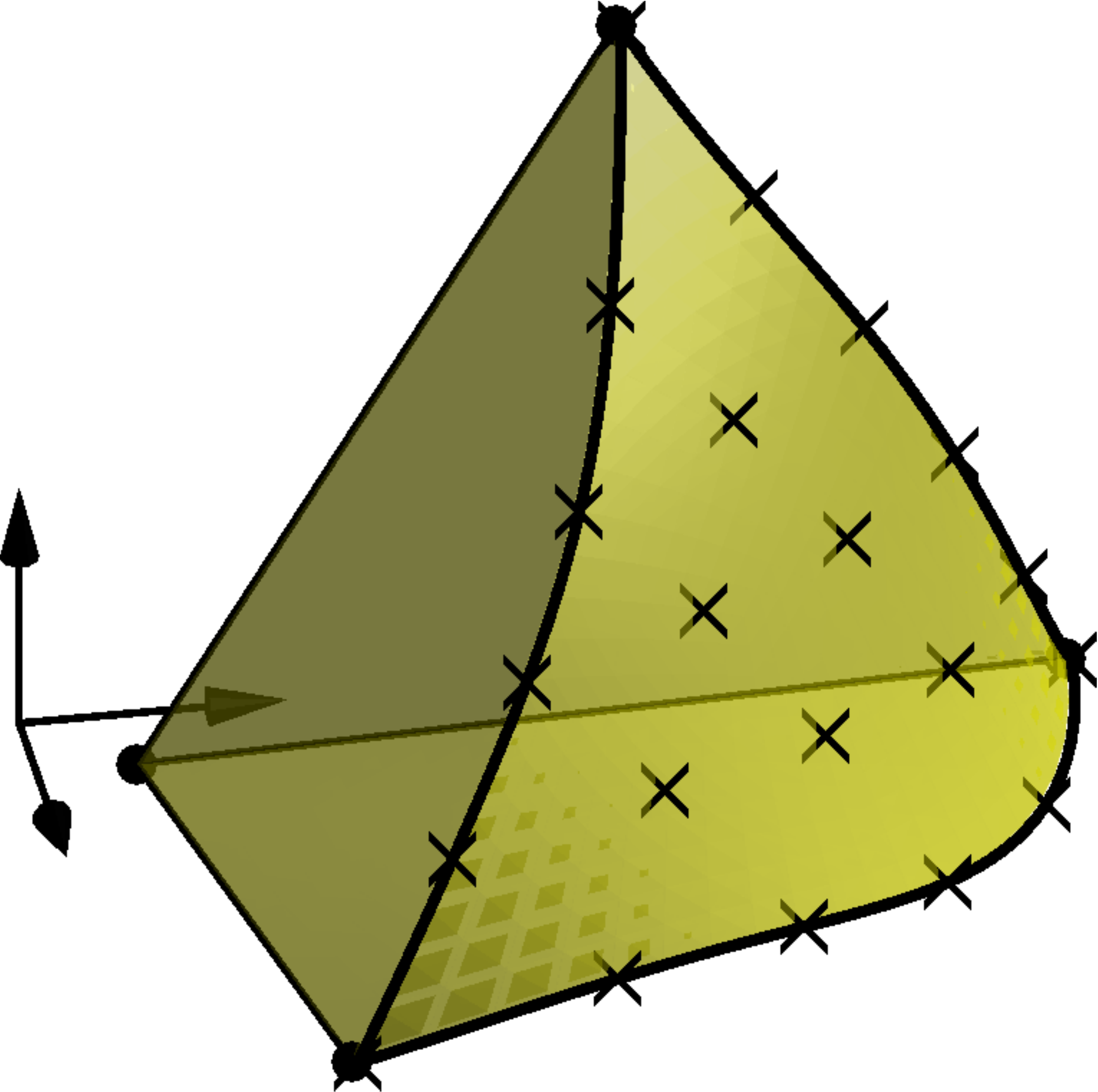}}

\caption{\label{fig:VisRecon3dSolinTetr}(a) Node numbering in the reference
tetrahedron, (b) Local coordinate systems, (c) example of a mapped
tetrahedron in coordinates $\vek r$ with one higher-order side.}
\end{figure}

Let there be a tetrahedron with \emph{one} curved face given by a
higher-order triangular element with $n_{q}$ nodes. The more general
case where all four faces are curved is not needed herein. The curved
face has three curved edges, which are higher-order line elements
with $n_{p}$ nodes each. The other three edges are straight. The
situation is depicted in Fig.~\ref{fig:VisRecon3dSolinTetr}. The
task is to define a smooth mapping for all points $\left(a,b,c\right)$
in a reference tetrahedron to the three-dimensional situation in $\left(r,s,t\right)$.

The nodal coordinates of the triangular face element are denoted as
$\vek r_{i}^{\mathrm{face}}\in\mathbb{R}^{3}$ with $i=1,\dots,n_{q}$.
The three corner nodes of the face are set as the nodes $\vek r_{i}^{\mathrm{tetra}}$
with $i=2,3,4$ of the tetrahedron. The position of the remaining
node of the tetrahedron $\vek r_{1}^{\mathrm{tetra}}$ must be provided
as well. It shall be seen that the resulting mapping may be decomposed
into three parts, which may be associated to vertex, face and edge
contributions, hence,
\begin{equation}
\vek r\left(\vek a\right)=\vek r_{\mathrm{vertex}}\left(\vek a\right)+\vek r_{\mathrm{edge}}\left(\vek a\right)+\vek r_{\mathrm{face}}\left(\vek a\right).\label{eq:MappingTetra}
\end{equation}
The first part is simply the linear mapping of a tetrahdron, i.e.
\begin{equation}
\vek r_{\mathrm{vertex}}\left(\vek a\right)=\sum_{i=1}^{4}N_{i}^{\diamondsuit}\left(\vek a\right)\cdot\vek r_{i}^{\mathrm{tetra}}\label{eq:MappingTetraVertex}
\end{equation}
 with the linear shape functions 
\[
N_{1}^{\diamondsuit}\left(\vek a\right)=1-a-b-c,\quad N_{2}^{\diamondsuit}\left(\vek a\right)=a,\quad N_{3}^{\diamondsuit}\left(\vek a\right)=b,\quad N_{4}^{\diamondsuit}\left(\vek a\right)=c.
\]

For the edge contributions, we proceed similar as in Section \ref{sub:MappingTriQuad}.
That is, local coordinates $w_{k}\in\left[-1,1\right]$ along each
of the three edges belonging to the higher-order face are introduced,
see Fig.~\ref{fig:VisRecon3dSolinTetr}(b). They are related to the
coordinates $\left(a,b,c\right)$ as 
\[
w_{1}=b-a,\quad w_{2}=c-b,\quad w_{3}=a-c.
\]
 It is then simple to define $\vek r^{\mathrm{edge}\, k}\left(w_{k}\right)$
as above in Eq.~(\ref{eq:EdgeContribution}). The ramp function of
Eq.~(\ref{eq:RampTri}) is adapted as
\begin{equation}
R_{1}=\frac{N_{2}^{\diamondsuit}\left(\vek a\right)\cdot N_{3}^{\diamondsuit}\left(\vek a\right)}{N_{1}^{\mathrm{lin}}(w_{1})\cdot N_{2}^{\mathrm{lin}}(w_{1})},\: R_{2}=\frac{N_{3}^{\diamondsuit}\left(\vek a\right)\cdot N_{4}^{\diamondsuit}\left(\vek a\right)}{N_{1}^{\mathrm{lin}}(w_{2})\cdot N_{2}^{\mathrm{lin}}(w_{2})},\: R_{3}=\frac{N_{4}^{\diamondsuit}\left(\vek a\right)\cdot N_{2}^{\diamondsuit}\left(\vek a\right)}{N_{1}^{\mathrm{lin}}(w_{3})\cdot N_{2}^{\mathrm{lin}}(w_{3})}.\label{eq:RampTetrEdge}
\end{equation}
Finally, the overall edge contribution for Eq.~(\ref{eq:MappingTetra})
is defined as
\begin{equation}
\vek r_{\mathrm{edge}}\left(\vek a\right)=\sum_{k=1}^{3}R_{k}\cdot\vek r^{\mathrm{edge}\, k}\left(w_{k}(\vek a)\right).\label{eq:MappingTetraEdge}
\end{equation}

It remains to define the face contribution wherefore we need the local
coordinates $(u,v)$ corresponding to the triangular face element,
see also Fig.~\ref{fig:VisRecon3dSolinTetr}(b), 
\begin{equation}
u\left(\vek a\right)=b+\nicefrac{1}{3}(1-a-b-c),\qquad v\left(\vek a\right)=c+\nicefrac{1}{3}(1-a-b-c).\label{eq:Map u(a), v(a)}
\end{equation}
One needs the higher-order shape functions $N_{i}^{\mathrm{HO}}\left(u,v\right)$
and the linear shape functions $N_{i}^{\mathrm{lin}}\left(u,v\right)$
of the triangular face. Together with the corresponding nodes $\vek r_{i}^{\mathrm{face}}$
and $\vek r_{i}^{\mathrm{tetra}}$ this leads to the map 
\begin{equation}
\vek r^{\mathrm{face234}}\left(\vek u\right)=\sum_{i=1}^{n_{q}}N_{i}^{\mathrm{HO}}\left(\vek u\right)\cdot\vek r_{i}^{\mathrm{face}}-\sum_{i=1}^{3}N_{i}^{\mathrm{lin}}\left(\vek u\right)\cdot\vek r_{i+1}^{\mathrm{tetra}}.\label{eq:FaceContribution}
\end{equation}
We need a bubble function which strictly lives inside the triangular
face. Therefore, one needs to \emph{subtract} the edge contributions
from the face contribution. We have transformed coordinates $\vek a$
to $\vek u$ based on Eq.~(\ref{eq:Map u(a), v(a)}). Next, special
points $\vek a^{\star}$ are generated based on $\vek u$ and the
linear shape functions $N_{i}^{\mathrm{lin}}\left(\vek u\right)$
of the triangle, i.e.
\[
\vek a^{\star}\left(\vek u\left(\vek a\right)\right)=\sum_{i=1}^{3}N_{i}^{\mathrm{lin}}\left(\vek u\right)\cdot\vek a_{i+1}^{\mathrm{tetra}}
\]
where $\vek a_{k}^{\mathrm{tetra}}$, $k=2,3,4$, are simply the coordinates
of the corner nodes of the tetrahedron in the reference configuration
as seen in Fig.~\ref{fig:VisRecon3dSolinTetr}(a). One may then evaluate
the edge contributions (\ref{eq:MappingTetraEdge}) based on $\vek a^{\star}$
obtaining $\vek r_{\mathrm{edge}}\left(\vek a^{\star}\right).$ We
are now ready to define the face contribution to Eq.~(\ref{eq:MappingTetra})
as 
\begin{equation}
\vek r_{\mathrm{face}}\left(\vek a\right)=S\cdot\left(\vek r^{\mathrm{face234}}\left(\vek u\left(\vek a\right)\right)-\vek r_{\mathrm{edge}}\left(\vek a^{\star}\left(\vek u\left(\vek a\right)\right)\right)\right)\label{eq:MappingTetraFace}
\end{equation}
with the ramp function
\[
S=\frac{N_{2}^{\diamondsuit}\left(\vek a\right)\cdot N_{3}^{\diamondsuit}\left(\vek a\right)\cdot N_{4}^{\diamondsuit}\left(\vek a\right)}{N_{1}^{\mathrm{lin}}(\vek u)\cdot N_{2}^{\mathrm{lin}}(\vek u)\cdot N_{3}^{\mathrm{lin}}(\vek u)}.
\]

\subsection{Mapping for prisms\label{sub:MappingPrisms}}

\begin{figure}
\centering

\subfigure[]{\includegraphics[height=5cm]{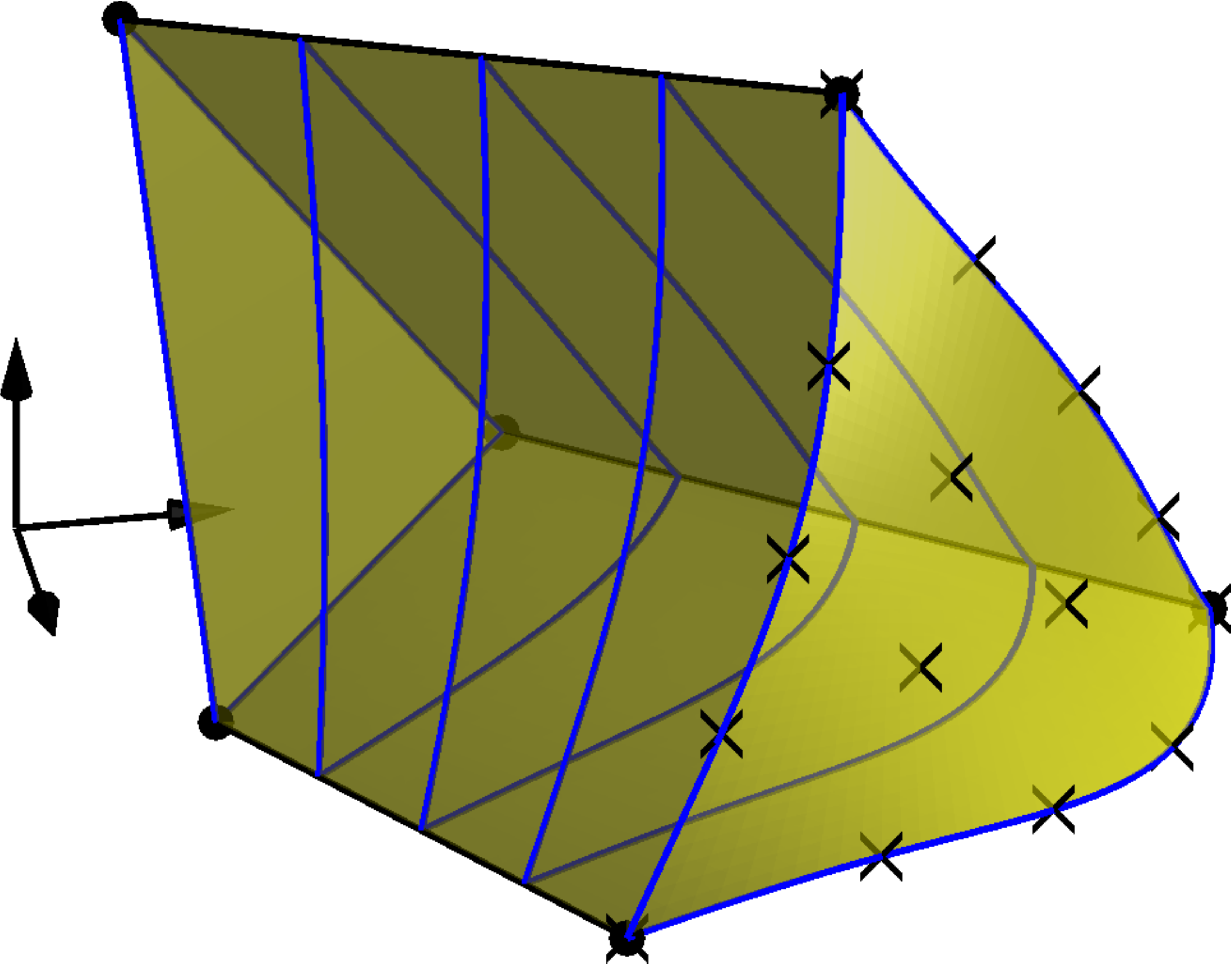}}\qquad\subfigure[]{\includegraphics[height=5cm]{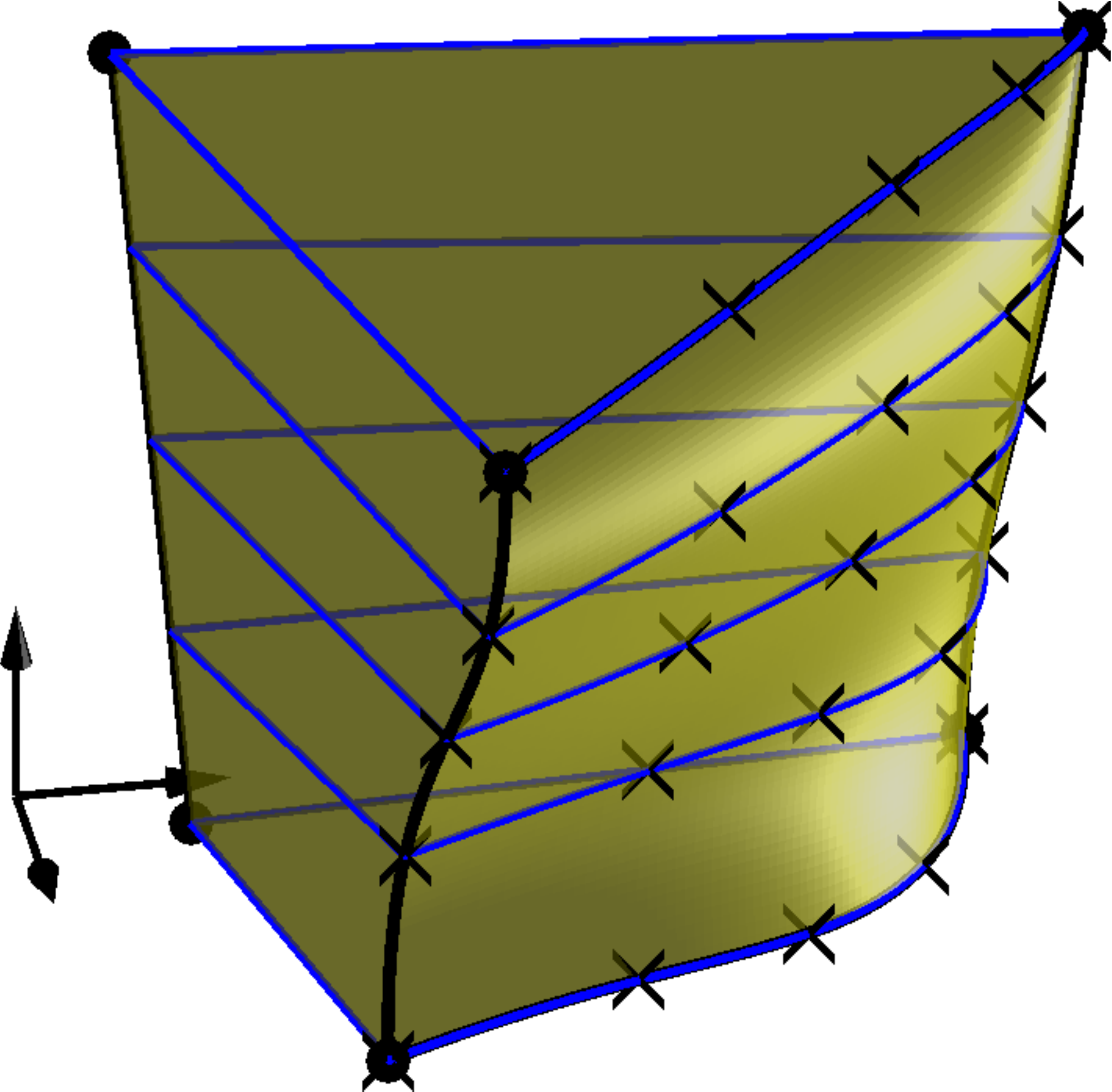}}

\caption{\label{fig:VisRecon3dSolinPrisms}Prisms with one higher-order (a)
triangular or (b) quadrilateral face.}
\end{figure}

The general situation where all faces of the prism are defined by
(curved) higher-order elements is discussed in \cite{Solin_2003a}
but the implementation is rather tedious. Herein, only one of the
faces of the prismatic element may be of higher-order resulting into
two different situations: Case 1 arises from the situation where a
tetrahedral background element is cut into a sub-tetrahedron and a
sub-prism; the higher-order side of the prism is then a triangle,
see Fig.~\ref{fig:VisRecon3dSolinPrisms}(a). Case 2 results when
the tetrahedron is split into two sub-prisms; the higher-order side
of the prism is then a quadrilateral, see Fig.~\ref{fig:VisRecon3dSolinPrisms}(b).
The map for the first case is defined in a straightforward way using
the original idea of \cite{Gordon_1973a,Gordon_1973b}, which is detailed
for the present situation in a previous work of the authors \cite{Fries_2015a}.

In the second case, we wish to avoid the general definition of the
map $\vek r\left(\vek a\right)$ and rather outline the procedure
in an illustrative way: As we are only interested in placing inner
element nodes, it is noted that one may cut the prism into slices
according to the blue lines in Fig.~\ref{fig:VisRecon3dSolinPrisms}(b).
Then, each slice is a triangle with one higher-order side and the
mapping defined in Section \ref{sub:MappingTriQuad} applies. Doing
so in each slice defines all nodes of the sought higher-order prism.

\bibliographystyle{schanz}
\addcontentsline{toc}{section}{\refname}\bibliography{FriesRefs}
 
\end{document}